\documentclass[a4paper,fleqn,11pt]{cas-sc}
\usepackage[square, numbers]{natbib}
\usepackage{algorithm,algpseudocode}
\renewcommand{\algorithmicrequire}{\textbf{Input:}}
\renewcommand{\algorithmicensure}{\textbf{Output:}}
\usepackage{tikz}
\usetikzlibrary{shapes.geometric, arrows}
\usepackage{amsmath}
\numberwithin{equation}{section}
\usepackage{rotating}
\usepackage{lineno}
\usepackage{color,soul}
\usepackage{todonotes}
\usepackage{setspace} 
\usepackage{caption}
\usepackage{subcaption}
\usepackage{arydshln}
\usepackage{schemata}
\newcommand\AB[2]{\schema{\schemabox{#1}}{\schemabox{#2}}}

\def\tsc#1{\csdef{#1}{\textsc{\lowercase{#1}}\xspace}}
\tsc{WGM}
\tsc{QE}
\tsc{EP}
\tsc{PMS}
\tsc{BEC}
\tsc{DE}
\begin{document}
\let\WriteBookmarks\relax
\def\floatpagepagefraction{1}
\def\textpagefraction{.001}
\shorttitle{A parallel branch-and-bound-and-prune algorithm for irregular strip packing with discrete rotations}
\shortauthors{Lastra-Díaz\&Ortuño}

\title [mode = title]{A parallel branch-and-bound-and-prune algorithm for irregular strip packing with discrete rotations}                      
%\tnotemark[1,2]
%\tnotetext[1]{This document is the results of the research project funded by the National Science Foundation.}
%\tnotetext[2]{The second title footnote which is a longer text matter to fill through the whole text width and overflow into   another line in the footnotes area of the first page.}

\author[1]{Juan J. Lastra-Díaz}[orcid=0000-0003-2522-4222]
\cormark[1]
\cortext[cor1]{Corresponding author}
\ead{jlastra@ucm.es}

\author[1]{M. Teresa Ortuño}[orcid=0000-0002-5568-9496]
\ead{mteresa@ucm.es}

\address[1]{Department of Statistics and Operational Research, Institute of Interdisciplinary Mathematics, UCM Research Group HUMLOG, Complutense University of Madrid, Spain}

\begin{abstract}
The irregular strip-packing problem consists of the computation of a non-overlapping placement of a set of polygons onto a rectangular strip of fixed width and the minimal length possible. Recent performance gains of the Mixed-Integer Linear Programming (MILP) solvers have encouraged the proposal of exact optimization models for nesting. The Dotted-Board (DB) MILP model solves the discrete version of the nesting problem by constraining the positions of the polygons to be on a grid of fixed points. However, its number of non-overlapping constraints grows exponentially with the number of dots and types of polygons, which encouraged the proposal of a reformulation called the DB Clique Covering (DB-CC) that sets the current state-of-the-art by significantly reducing the constraints required. However, DB-CC requires a significant preprocessing time to compute edge and vertex clique coverings. Moreover, current knowledge of the stable set polytope suggests that achieving a tighter formulation is unlikely. Thus, our hypothesis is that an ad-hoc exact algorithm requiring no preprocessing might be a better option to solve the DB model than the costly Branch-and-Cut approach. This work proposes an exact branch-and-bound-and-prune algorithm to solve the DB model from the conflict inverse graph based on ad-hoc data structures,  bounding, and forward-checking for pruning the search space. We introduce two 0-1 ILP DB reformulations with discrete rotations and a new lower-bound algorithm as by-products. Our experiments show that DB-PB significantly reduces the resolution time compared to our replication of the DB-CC model. Seventeen open instances are solved up to optimality.
\end{abstract}

%\includegraphics{figs/grabs.pdf}

%\item Two new mixed-integer-programming continuous models without rotations for nesting
%\item A new family of valid cuts, symmetry-breakings, and variable eliminations for nesting
%\item A new family of feasibility cuts among three pieces is introduced for the first time
%\item New models set the state of the art among the family of continuous models for nesting
%\item MIP solver's advances contribute to solving ten open instances, one with 27 pieces

\begin{keywords}
Packing \sep Irregular strip packing \sep Integer programming \sep Branch-and-Bound-and-Prune \sep Forward-checking
\end{keywords}

%\linenumbers
\maketitle
\section{Introduction}

Cutting and packing regular (convex) and irregular (non-convex) polygons onto a rectangular strip with unbounded length is a tedious and usual task in most manufacturing industries based on the cutting of any flat material. For instance, \citet{Milenkovic1991-mo} study the nesting problem for the fashion and apparel industry, whilst \citet{Heistermann1995-qm} and \citet{Whelan1993-hf} do it for leather manufacturing, \citet{Elamvazuthi2009-pu} in furniture, \citet{Han2013-cb} in the glass industry, \citet{Alves2012-ab} in the automotive industry, and \citet{Cheok1991-jy} and \citet{Hamada2019-bx} in shipbuilding. The irregular strip-packing problem is also known as \emph{nesting} or \emph{marker making} in the industry, and it consists in computing a non-overlapping placement of a set of irregular polygons, called \emph{pieces}, onto a fixed-width rectangular strip with unbounded length, called the \emph{board}, whose length is the minimum possible. Another closely related problem, called two-dimensional bin packing \cite{Lodi2002-jz, Iori2021-bi}, is defined as the computation of a non-overlapping placement of a set of polygons onto a larger closed polygon, called the \emph{bin}, to minimize the number of bins required. Regular and irregular strip and bin packing problems, and all their variants concerning the geometry of the pieces or boards, belong to the broader family of Cutting and Packing (C\&P) problems categorized by \citet{Dyckhoff1990-xv} and \citet{Wascher2007-fw}, and extensively reviewed by \citet{Sweeney1992-cb}, \citet{Dowsland1992-su}, \citet{Wang2002-sd}, and \citet{Bennell2013-jd}, among others. 

Research on the irregular strip and bin packing problems can be traced back to the pioneering Linear Programming (LP) models for rectangular bin packing introduced by \citet{Gilmore1965-vd}, and the pioneering heuristic methods for irregular strip packing proposed by \citet{Art1966-vy}, \citet{Adamowicz1976-wv, Adamowicz1976-fp}, and \citet{Albano1980-tq} in the late nineteen sixties and seventies. These early works introduce many of the basic ideas subsequently exploited by all heuristics methods reported in the literature, such as the notion of a feasible non-overlapping region between pieces based on the No-Fit Polygon (NFP) representation \cite{Bennell2008-rc}, and the sequential placement of pieces based on bottom-left heuristics. The boundary and outer region of the no-fit polygon $NFP_{AB}$ set the feasible positions in which polygon $B$ can be placed, without rotating, into a non-overlapping position with respect to polygon $A$. Thus, the no-fit polygon allows one to compute the feasible relative placements between polygons $A$ and $B$ a priori by checking whether their relative position vector $\delta_{AB} = r_B - r_A$ belongs to either the boundary or the outer region of $NFP_{AB}$, where $r_A$ and $r_B$ denote their corresponding reference points. This property has converted the NFP into the most broadly adopted and effective geometric representation for nesting reported in the literature, both for the families of heuristics methods \cite{Bennell2009-cy} and for exact mathematical models \cite{Leao2020-bc}. 

\citet{Fowler1981-vr} show that the irregular strip-packing problem is NP-complete. For this reason, most practical solutions reported in the literature since the pioneering work of \citet{Art1966-vy} are based on sequential placement heuristics to build efficiently feasible solutions that are combined with meta-heuristics for exploring the space of feasible solutions, as shown in most of surveys on nesting \cite{Dowsland1995-qp, Dowsland2002-hx, Hopper2001-rp, Bennell2009-cy, Riff2009-rv}. \citet{Elkeran2013-pe} introduces the current state-of-the-art heuristic method for nesting, called Guided Cuckoo Search (GCS), which defines a two-stage method based on piece clustering and NFP-based bottom-left heuristics \cite{Gomes2002-qj} to build an initial feasible solution that is shrunk by solving an overlap minimization problem using a variant of the cuckoo search meta-heuristics. Although GCS was introduced more than a decade ago, subsequent works have been unable to outperform its results, as shown by \citet[table 3]{Pinheiro2016-bg}, \citet[table 5]{Sato2016-kc}, \citet[table 4]{Cherri2016-zr}, \citet[table 6]{Mundim2017-gs}, \citet[table 2]{Amaro_Junior2017-fi}, and \citet[tables 4-5]{Sato2019-lp}. More recently, other authors have introduced new heuristics for the continuous and discrete versions of the nesting problem, such as the works of \citet{Fang2023-pb}, \citet{Na2023-vn}, \citet{Guo2023-mt}, \citet{Sato2023-cw}, \citet{Queiroz2023-cz}, \citet{Umetani2022-uc}, and \citet{Queiroz2020-oh}. For a recent comparison of heuristics and exact mathematical models for nesting, we refer the readers to the comparative evaluation carried out by \citet{Giovenali2025-hl}.

The recent performance gains achieved by state-of-the-art Mixed-Integer Programming (MIP) solvers \cite{Bixby2002-ry, Achterberg2013-qz, Koch2022-dn} have encouraged the proposal of exact mathematical models for nesting during the last decade, as shown by the literature reviews of \citet{Leao2020-bc} and \citet{Lastra-Diaz2024-il}. The exact mathematical models for nesting can be categorized into mixed-integer linear programming models, Constraint Programming (CP) models, and non-linear programming models, as shown in the categorization of \citet[Fig. 1]{Lastra-Diaz2024-il}. In turn, the family of mathematical programming models for nesting can be categorized into three large families according to the type of decision variables used to represent the position of the pieces as follows. Firstly, the family of continuous models, such as the pioneering MILP model of \citet[\S8]{Li1994-co}, and the subsequent MILP models proposed by \citet[\S5]{Dean2002-wx}, \citet{Fischetti2009-oa}, \citet{Alvarez-Valdes2013-wg}, \citet{Cherri2016-jf}, \citet{Rodrigues2017-dy}, \citet{Lastra-Diaz2024-il}, and \citet{Pantoja-Benavides2024-fv}. Secondly, the family of discrete ILP models for nesting to which this work belongs, whose pioneering model is the Dotted-Board (DB) model introduced by \citet{Toledo2013-oi}, subsequently improved by \citet{Rodrigues2017-zl}, and more recently extended by \citet{Rodrigues_de_Souza_Queiroz2020-tp, Rodrigues_de_Souza_Queiroz2022-jf} to deal with the uncertainty in the demand of pieces by formulating a two-stage stochastic programming model. And thirdly, the semi-continuous MILP model of \citet{Leao2016-xx} whose aim is to overcome the limitations of the grid resolution in the DB model \cite{Toledo2013-oi} by dividing the board into continuous horizontal slices at regular y-axis discrete positions to enlarge the space of feasible positions for the pieces. We refer the reader to the aforementioned works \cite{Leao2020-bc, Lastra-Diaz2024-il} for a comprehensive and updated review of the families of exact methods for nesting. 

The first exact mathematical models for nesting reported in the literature tackled the continuous version of the problem. \citet[\S8]{Li1994-co} introduces the first continuous MIP model for nesting reported in the literature to solve a limitation of the pioneering LP compaction model of \citet{Li1993-uj, Li1995-aj}. Although \citeauthor{Li1994-co}'s model was not experimentally evaluated, it sets the two main features of the family of continuous exact MIP models based on the NFP as follows: (1) the convex decomposition of the outer NFP feasible regions; and (2) the definition of mutually-exclusive binary variables to set the pairwise non-overlapping constraints between pieces defining the feasible regions for their relative placement. For instance, \citet[\S5]{Dean2002-wx} introduces and evaluates for the first time a MIP model for nesting based on a refinement of the \citet{Daniels1994-zn} bin packing model that is essentially identical to the \citeauthor{Li1994-co}'s model. Subsequently, \citet{Fischetti2009-oa} (F\&L) introduce a refinement of \citeauthor{Li1994-co}'s model based on lifting the big-M formulation and a branching-priority algorithm to guide the Branch\&Bound (B\&B) exploration. \citet{Alvarez-Valdes2013-wg} improve the F\&L \cite{Fischetti2009-oa} model by introducing a new MIP model called HS2, which is based on a detailed convex decomposition of the feasible regions into horizontal slices, a lifting for the bound constraints of the continuous variables, six new branching strategies, and the x-axis ordering of identical pieces to remove all symmetric solutions derived from their permutation. \citet{Cherri2016-jf} introduce two continuous MILP models improving the HS2 model \cite{Alvarez-Valdes2013-wg}, together with the first continuous MIP model integrating discrete rotations, which are based on the convex decomposition of the pieces and the definition of the non-overlapping constraints between pieces by using the convex no-fit polygons among their convex parts, together with some families of valid inequalities, and the same x-axis symmetry-breaking for identical pieces proposed by \citet{Alvarez-Valdes2013-wg}. Subsequently, \citet{Rodrigues2017-dy} improve the NFP-CM model by breaking the symmetries of the feasible space for the relative placements between pieces. More recently, \citet{Lastra-Diaz2024-il} introduce two continuous MILP models without rotations called NFP-CM based on Vertical Slices (NFP-CM-VS), which are two reformulations of the NFP-CM model of \citet{Cherri2016-jf} based on a new convex decomposition of the feasible space of relative placements between pieces into vertical slices, together with a large set of logical objects to tighten the formulation derived from the former convex decomposition, such as a new family of valid inequalities, symmetry breakings, variable eliminations, and a family of feasibility cuts among three pieces. The family of NFP-CM-VS MIP models proposed by \citet{Lastra-Diaz2024-il} set the current state-of-the-art in terms of performance among the family of exact continuous mathematical models for irregular strip packing without rotations. Finally, \citet{Pantoja-Benavides2024-fv} introduces a reformulation of the DTM model with rotations of \citet{Cherri2016-jf} that includes a separation constraint and a new set of valid inequalities that allows them to set the state of the art for the exact continuous problem with discrete rotations.

Concerning the families of discrete and semi-continuous MILP models for nesting, \citet{Toledo2013-oi} introduces the pioneering discrete MILP model for nesting, called Dotted-Board (DB), which is based on a set of binary decision variables for constraining the positions of the polygons to be on a grid of fixed points. The DB model has several nice properties, as follows. Firstly, the DB model reduces the definition of non-overlapping constraints among pieces to a set of edge inequalities for pairs of binary variables encoding each infeasible placement between pairs of pieces, which provides one formulation that is independent of both the geometry of the pieces and the board, allowing the representation of irregular pieces with or without holes, as well as any geometry for the board in the same way (e.g. rectangular, leathers with quality zones and defects, etc). Secondly, the DB model preserves the exact geometry of the pieces, unlike other discrete methods which discretize both the pieces and the board, such as the heuristics proposed by \citet{Segenreich1986-bn}, \citet{Oliveira1993-so}, \citet{Babu2001-ax} and \citet{Chehrazad2022-yq} among others. And thirdly, the DB model might be interpreted as a particular case of the more general set packing problem \cite{Balas1972-iz}, whose polytope corresponds to the well-known polytope defined by the ILP formulation of the maximum stable set problem \cite{Giandomenico2013-if, Letchford2020-lj} for an undirected graph $G=(V,E)$, denoted by $STAB(G)$, which can be traced back up to the pioneering work of \citet{Padberg1973-bi}. Thus, this latter property of the DB formulation might potentially allow us to benefit from the large corpus of knowledge and valid inequalities proposed for the Stable Set problem (SSP). However, the DB model has several drawbacks and limitations as follows: (1) its number of non-overlapping constraints grows exponentially with the number of dots and polygon types; (2) it does not consider rotations of the pieces; (3) the quality of the discrete solutions produced by the DB model depends both on the grid resolution, as argued by \citet{Toledo2013-oi} and shown experimentally by \citet{Sato2016-fd}, and the distribution of the points within the board; and finally, (4) it is a mixed-integer programming model, which prevents its resolution with either CP and SAT-based methods or exact implicit enumeration algorithms as done herein. \citet{Cherri2018-rq} bridge the second and third aforementioned drawbacks of the DB model by proposing the use of non-regular grids and one extension of the DB model considering discrete rotations of the pieces. \citet{Rodrigues2017-zl} successfully solve the first drawback detailed above and set the current state of the art of the discrete nesting problem by proposing a Clique Covering MILP reformulation of the DB model that we call DB Clique Covering (DB-CC) model here, which uses edge and vertex clique coverings to reduce significantly the number of constraints and tighten the LP relaxation of the model. However, DB-CC requires a significant preprocessing time to compute the edge and vertex clique coverings, and the lack of a discrete DB formulation without continuous variables is still an open question. Finally, the difficulty of solving exactly the DB model \cite{Toledo2013-oi} has encouraged the proposal of metaheuristics to find near-optimal solutions in reasonable running times, such as the Biased Random-Key Genetic Algorithm (BRKGA) introduced by \citet{Mundim2017-gs} and the Raster Overlap Minimization Algorithm (ROMA) proposed by \citet{Sato2019-lp}. However, no exact alternative algorithm to the traditional Branch-and-Cut approach \cite{Padberg1991-hm, Grotschel1991-qh, Mitchell2002-hg} implemented by all MIP solvers has been proposed yet. 

Concerning the exploration and proposal of new tightening valid inequalities for the DB model \cite{Toledo2013-oi}, the current knowledge about the stable set polytope suggests that achieving a tighter formulation that might be efficiently solved with the standard Branch-and-Cut approach is unlikely, as pointed out by \citet{Giandomenico2013-if} who argue that despite the large set of valid inequalities proposed for the stable set problem, such as clique inequalities \cite{Padberg1973-bi}, odd-hole inequalities \cite{Padberg1973-bi, Nemhauser1992-aw}, web and antiweb inequalities \cite{Trotter1975-tc}, wheel inequalities \cite{Cheng1997-og}, and anti-wheel inequalities \cite{Cheng2002-jl}, ``none of them, to the best of our knowledge, gave any computational outcome" \cite[p.108]{Giandomenico2013-if}. For this reason, this work focuses on proposing exact alternative algorithms to the standard Branch-and-Cut framework to bridge the performance and scalability gaps of the DB model.

The main motivation of this work is to solve efficiently the first drawback of the DB model enumerated above and the preprocessing drawback of the DB-CC model by introducing an exact algorithm called Dotted-Board Parallel Backtracking (DB-PB) to solve the DB model from the conflict inverse graph based on ad-hoc data structures, bounding, and forward-checking for pruning the search space. Likewise, we introduce two 0-1 ILP reformulations with discrete rotations of the current family of DB models \cite{Toledo2013-oi, Rodrigues2017-zl} called the Binary Dotted-Board (BDB) model and the Binary Dotted-Board Clique Covering model (BDB-CC), respectively, to bridge the fourth drawback of the DB model \cite{Toledo2013-oi}. Finally, we introduce a new lower-bound algorithm called (DB-PB-LB) as a by-product. Our main hypothesis is that an ad-hoc exact algorithm requiring no preprocessing might be a better option to solve the DB model than the costly general-purpose Branch-and-Cut or hybrid CP-SAT methods, in addition to avoiding the need for computing edge and vertex clique coverings. A second aim of this work is to carry out a fair and reproducible comparison of our new exact algorithms and 0-1 DB ILP models with the state-of-the-art discrete DB-CC MILP model of \citet{Rodrigues2017-zl} by exactly replicating both their clique covering model and the experiments reported in their article \cite{Rodrigues2017-zl} into a single software and hardware platform based on our C\# implementation of our exact algorithms in the .NET 8 platform and the state-of-the-art Gurobi and Google CPSat solvers. Finally, a third aim is to bridge the lack of reproducibility resources hampering the independent replication and confirmation of previously reported methods and results by providing a detailed reproducibility protocol and dataset as supplementary material to allow the exact replication of all our methods, experiments, and results.

The main research problem tackled by this work is the proposal of efficient exact algorithms for the resolution of the discrete DB MILP model \cite{Toledo2013-oi} for nesting requiring no preprocessing, together with two 0-1 DB ILP reformulations to remove the continuous variables, and a new lower algorithm for discrete nesting. Thus, our main contribution is the introduction of an ad-hoc parallel branch-and-bound-and-prune algorithm for the efficient solution of the DB model, called Dotted-Board Parallel Backtracking (DP-PB) algorithm. A second contribution is the proposal of two O-1 ILP reformulations of the DB model considering discrete rotations, called Binary DB model and Binary DB Clique Covering model, respectively. A third contribution is the proposal of a new lower bound algorithm as a by-product of our new DB-PB algorithm. Finally, a fourth contribution is the exact replication of the DB Clique Covering model and the experiments introduced by \citet{Rodrigues2017-zl}, together with a set of reproducibility resources provided as supplementary material to allow the exact replication of all our methods, experiments, and results.

The rest of the paper is structured as follows. Section \ref{sec:db_models} introduces the formulation and notation of the current family of DB models for nesting to make the comprehension of our new methods and contributions easier. Section \ref{sec:our_binary_DB_model} introduces our two 0-1 ILP reformulations of the former DB  and DB-CC models \cite{Toledo2013-oi, Rodrigues2017-zl}, whilst Section \ref{sec:parallel_backtracking} introduces our exact parallel branch-and-bound-and-prune algorithm to solve efficiently the DB model \cite{Toledo2013-oi} with no preprocessing, and Section \ref{sec:new_lower_bound} introduces our new lower bound algorithm for nesting. Section \ref{sec:evaluation} details our experimental setup and results, whilst section \ref{sec_discussion} introduces our discussion of the results. The final section summarizes our main conclusions and future work. Finally, Appendix A introduces the best solutions obtained for all problem instances in our experiments, whilst Appendix B introduces a detailed reproducibility protocol based on our supplementary dataset \cite{Lastra-Diaz2025-gp} to allow the exact replication of all our models, experiments, and results. Both appendices are provided as supplementary materials.

\section{Preliminary concepts and problem definition}
\label{sec:db_models}

This section introduces the DB model of \citet{Toledo2013-oi} and its reformulation based on clique coverings proposed by \citet{Rodrigues2017-zl} that we called DB-CC in this work, as well as the notation shared by all ILP models and exact algorithms introduced herein, including the decision variables and indexes for the DB model with discrete rotations introduced by \citet{Cherri2018-rq}. The discrete nesting problem with discrete rotations can be formalized as follows. Let $\mathcal{P} = \{P_1,\dots,P_N\}$ be a set of irregular polygons not necessarily different with reference points $r_i \in \mathbb{R}^2$, called \emph{pieces}, and let $\mathcal{T} = \{1,\dots,T\}$ be a set of types of pieces, such that each piece $P_i \in \mathcal{P}$ belongs to some $t \in \mathcal{T}$, and $\Theta = \bigcup\limits_{\forall t \in \mathcal{T}} \Theta_t$ the set of feasible rotations, such that $\Theta_t = \{\theta_t^k \in [0,2\pi], k = 1,\dots, K_t\}$ is the discrete set of feasible rotations for each type of piece $t \in \mathcal{T}$. Let $\mathcal{B} \subset \mathbb{R}^2$ be an open rectangle with fixed width $W$ and length variable $L$ called the \emph{board}, and $\mathcal{D} \subset \mathcal{B}$ a finite set of points defining a discretization of $\mathcal{B}$, as shown in Figure \ref{fig:IFP_dotted-board_model}. Then, the minimization problem (\ref{obj:Abstract}) defines the optimal solution to the discrete irregular strip-packing problem, where $P_i(\theta_i)$ denotes an orientation of the piece $P_i$ respect to its reference point $r_i$ and $P_i(\theta_i) \oplus d_i$ denotes a translation of the piece such that its reference point is positioned at $d_i \in  \mathcal{D}$. Constraints (\ref{ineq:Abstract_non_overlapping}) prevent the overlapping of pieces, whilst Constraints (\ref{ineq:Abstract_innerfit}) force the pieces to be entirely contained in the board.
\begin{align}
\min\limits_{d_i,\theta_i,L} \quad & z = L \label{obj:Abstract} \\
s.t. \quad & int(P_i(\theta_i) \oplus d_i) \cap int(P_j(\theta_j) \oplus d_j) = \emptyset, \; 1 \leq i < j \leq N \label{ineq:Abstract_non_overlapping} \\
& (P_i(\theta_i) \oplus d_i) \subseteq \mathcal{B}, \quad 1 \leq i \leq N \label{ineq:Abstract_innerfit} \\
& d_i \in \mathcal{D}, \; \theta_i \in \Theta, \quad 1 \leq i \leq N \\
& L \in \mathbb{R}^+
\end{align}

\subsection{The Dotted-Board (DB) model}
\label{sec:DB_model}

The Dotted-Board (DB) model \cite{Toledo2013-oi} is defined by the Objective function (\ref{obj:DB}) and the Constraints (\ref{ineq:DB_innerfit}) to (\ref{ineq:DB_length_variable}), together with the decision variables, parameters, and index sets detailed below. The basic DB model is based on a regular discretization of the board; however, it admits any sort of regular or irregular point pattern, as shown by \citet{Cherri2018-rq}. Given any piece of type $t \in \mathcal{T}$, its Inner-Fit Polygon, denoted by $\mathcal{IFP}_t$, defines the set of feasible positions $d \in \mathcal{D}$ in which the piece is completely contained in the board. Concerning the representation of the non-overlapping constraints between pairs of pieces, given two piece instances $P_i, P_j \in \mathcal{P}$ with types $t,u \in \mathcal{T}$ respectively, we call $P_i$ the static piece and $P_j$ the orbiting one, then the non-overlapping relative feasible positions $d,d' \in \mathcal{D}$ between both pieces are defined as the complementary set of their no-Fit polygon, denoted by $\mathcal{NFP}_{tu}^d$. Every grid point $d \in \mathcal{D}$ defines a set of binary decision variables $\delta_t^d$ whose activation in any feasible solution indicates that a piece of type $t \in \mathcal{T}$ is placed in that position. Thus, the number of binary decision variables in the DB model is proportional to $D \times T$, and it grows linearly with the number of points and type of pieces.
\begin{center}
\begin{tabular}{lp{13.5cm}}
\multicolumn{2}{l}{\textbf{Decision variables}}\\
$\delta_t^d$ & Binary variables indicating that one piece of type $t \in \mathcal{T}$ is placed in the dot $d \in \mathcal{D}$.\\
$L$ & Continuous positive decision variable representing the length of the board $\mathcal{B}$. \\
\multicolumn{2}{l}{\textbf{Parameters}} \\
$W$ & Board width \\
$\overline{L}$ & Upper Bound (UB) for the length of the board. \\
$\underline{L}$ & Lower Bound (LB) for the length of the board. \\
$g_x$ & Horizontal x-axis grid resolution. \\
$g_y$ & Vertical y-axis grid resolution. \\
$q_t$ & Number of pieces demanded for each piece type. \\
$T$ & Number of different piece types. \\
$N = \sum\limits_{t = 1}^T q_t$ & Total number of pieces. \\
$x_{t}^m$ & Horizontal distance from the leftmost point of any piece to its reference point $r_i$.\\
$x_{t}^M$ & Horizontal distance from the rightmost point of any piece to its reference point $r_i$. \\
$y_{t}^m$ & Vertical distance from the bottommost point of any piece to its reference point $r_i$.\\
$y_{t}^M$ & Vertical distance from the topmost point of any piece to its reference point $r_i$.\\
$l_t$ & Length of pieces of type $t \in \mathcal{}{T}$. \\
$a_t$ & Area of pieces of type $t \in \mathcal{T}$. \\
$C = \Big\lfloor \frac{\overline{L}}{g_x} \Big\rfloor + 1$ & Number of columns. \\
$R = \Big\lfloor \frac{W}{g_y} \Big\rfloor + 1$ & Number of rows. \\
$D = C \times R$ & Total number of dots.
\end{tabular}
\end{center}

The Objective function (\ref{obj:BDB_obj_function}) sets the minimization of the decision variable $L$ encoding the length of the board. The continuous $L$ decision variable is defined by all current MILP models for nesting to linearize the objective function $z = \min\{\max\{x_i + x_t^M: 1 \leq i \leq N\}\}$, where $x_i$ denotes the position of any piece instances $P_i$ of type $t$ that is placed into the board. Constraints (\ref{ineq:DB_innerfit}) force all pieces to be inside the board. Constraints (\ref{ineq:DB_demands}) set the demand constraints by forcing the placement of all pieces in any feasible solution. Constraints (\ref{ineq:DB_non-overlapping}) encode the non-overlapping constraints among pieces, called edge inequalities in the literature on set packing because they define an edge in the conflict graph $G=(V,E)$, where each vertex $v \in V$ represents a binary decision variable $\delta_t^d$ and there is an edge $e \in E$ for every pair of variables $\delta^{d'}_u, \delta^d_t$ appearing in one Constraint (\ref{ineq:DB_non-overlapping}). Finally, Constraints (\ref{ineq:DB_binary_variables}) and (\ref{ineq:DB_length_variable}) set the domains of the decision variables. The trivial lower bound for the length of the board used by the DB model and most of nesting MILP models is set by the formula (\ref{eq:trivial_lower_bound}) below.
\begin{equation}
\underline{L} = \max\Bigg\{\frac{1}{W}\sum\limits_{t = 1}^T a_t q_t, \max\Big\{\bigcup\limits_{t = 1}^T l_t\}\Big\} \Bigg\} \label{eq:trivial_lower_bound}
\end{equation}
\begin{center}
\begin{tabular}{lp{11cm}}
\multicolumn{2}{l}{\textbf{Index sets}}\\
$\mathcal{P} = \{1,\dots,N\}$ & Indexes for every piece instance $P_i$ to be placed in the board $\mathcal{B}$. \\
$\mathcal{D} = \{0,1,2,\dots,D-1\}$ & Indexes representing every dot $d = c_d \times R + r_d \in \mathcal{D}$ in the board, where $c_d$ and $r_d$ are the column and row indexes corresponding to dot $d \in \mathcal{D}$. \\
$\mathcal{T} = \{1,2,\dots,T\}$ & Indexes denoted by $t \in \mathcal{T}$ representing every different type of piece.\\
$\mathcal{IFP}_{t} = \{d \in \mathcal{D} : P_t\oplus d \subseteq \mathcal{B}\}$ & Indexes for the feasible positions of pieces of type $t \in \mathcal{T}$. \\
\multicolumn{2}{l}{$\mathcal{NFP}_{tu}^d = \{d' \in \mathcal{IFP}_u : int(P_t \oplus d) \cap int(P_u \oplus d') \neq \emptyset\}$} \\
& Set of infeasible relative positions between pieces of type $t$ and $u$ when the former one is placed at point $d \in \mathcal{D}$ .
\end{tabular}
\end{center}
\textbf{Dotted-Board (DB) model} \cite{Toledo2013-oi}
\begin{align}
\min\limits_{\delta^d_t,L} \quad & z = L \label{obj:DB}\\
\text{s.t.} \quad & (c_dg_x + x^M_t) \: \delta^d_t \leq L, \quad \forall d \in \mathcal{IFP}_t, \forall t \in \mathcal{T} \label{ineq:DB_innerfit} \\
& \sum\limits_{\forall d \in \mathcal{IFP}_t} \delta^d_t = q_t, \quad \forall t \in \mathcal{T} \label{ineq:DB_demands} \\
& \delta^{d'}_u + \delta^d_t \leq 1, \quad \forall d' \in  \mathcal{NFP}^d_{tu}, \; \forall (t,u) \in \mathcal{T} \times \mathcal{T}, \; \forall d \in \mathcal{IFP}_t \label{ineq:DB_non-overlapping} \\
& \delta^d_t \in \{0,1\}, \quad \forall d \in \mathcal{IFP}_t, \quad \forall t \in \mathcal{T} \label{ineq:DB_binary_variables} \\
& L \in [\underline{L}, \overline{L}] \subset \mathbb{R}_+  \label{ineq:DB_length_variable} 
\end{align}

\begin{figure}[t!]
\begin{subfigure}[t]{0.48\textwidth}
\centering
\begin{tikzpicture}[scale=1]
% Dibujamos la pieza
\fill[blue!10] (1.75,2) -- (3,1) -- (3.5,1) -- (4,3) -- (3,4) -- (0.5,3.5) -- (0.5,1.5) -- (1.75,2);
% Dibujamos el contorno de la pieza
\draw[thick,blue] (1.75,2) -- (3,1) -- (3.5,1) -- (4,3) -- (3,4) -- (0.5,3.5) -- (0.5,1.5) -- (1.75,2);
% Dibujamos la bbox de la pieza
\draw[dashed,black] (0.5,1) -- (4,1) -- (4,4) -- (0.5,4) -- (0.5,1);
\draw[black] (2.25,2.5) node {$P_i$};
% Dibujamos lo segmentos par marcar los extremosde la bbox
\draw[dashed,black] (0.5,1) -- (0.5,0.5);
\draw[black] (0.75,0.65) node {$x_t^m$};
\draw[dashed,black] (4,1) -- (4,0.5);
\draw[black] (3.75,0.65) node {$x_t^M$};
\draw[dashed,black] (4,1) -- (4.5,1);
\draw[black] (4.25,1.25) node {$y_t^m$};
\draw[dashed,black] (4,4) -- (4.5,4);
\draw[black] (4.25,3.75) node {$y_t^M$};
% Etiquetamos las regiones
\filldraw[black] (0.5,1) circle (0.05);
\draw[black] (0.5,1) circle (0.20);
\draw[black] (1.5,1.35) node {$r_i=(x_i,y_i)$};
\end{tikzpicture}
\caption{Representation of the geometry of piece $P_i$ of type $t \in \mathcal{T}$, where $r_i$ sets its reference point, whilst that $x_t^m$ and $x_t^M$ denotes its minimal and maximal horizontal coordinates regarding $r_i$ and $y_t^m$ and $y_t^M$ denotes its minimal and maximal vertical coordinates. The point $r_i$ is set to the bottom-left point of its bounding box by default.}
\label{fig:piece_vars_DB}
\end{subfigure}
\hfill
\hfill
\begin{subfigure}[t]{0.48\textwidth}
\centering
\begin{tikzpicture}[scale=0.6]
% Dibujamos los dos polígonos
\fill[blue!5] (4,3) -- (6,3) -- (8,6) -- (8,8) -- (2,8) -- (2,6) -- (4,3);
\draw[thick, blue] (4,3) -- (6,3) -- (8,6) -- (8,8) -- (2,8) -- (2,6) -- (4,3);
\draw[thick, blue, ->] (2,6) -- (3,4.5);
\fill[blue!15] (6,5) -- (8,8) -- (4,8) -- (6,5);
\draw[blue] (6,5) -- (8,8) -- (4,8) -- (6,5);
\fill[blue!15] (8.5,4) -- (10.5,4) -- (10.5,6) -- (8.5,6) -- (8.5,4);
\draw[blue] (8.5,4) -- (10.5,4) -- (10.5,6) -- (8.5,6) -- (8.5,4);
\filldraw[black] (6,5) circle (0.1);
\filldraw[black] (8.5,4) circle (0.1);
\draw[black] (6,7) node {$A$};
\draw[black] (9.5,5) node {$B$};
\draw[black] (4,5.5) node {$NFP_{AB}$};
\draw[black] (6,4.7) node {$r_A$};
\draw[black] (8.5,3.7) node {$r_B$};
\end{tikzpicture}
\caption{No-fit polygon $NFP_{AB}$ between the polygons $A$ and $B$ computed in $\mathbb{R}^2$. The no-fit polygons between pair of pieces can be pre-calculated and referenced to the reference point of the static pieces to compute later the NFP dots by intersecting the grid dots with the NFP continuous regions.}
\label{fig:nofit_polygon}
\end{subfigure}
\\
\begin{subfigure}[t]{0.48\textwidth}\centering
\begin{tikzpicture}[scale=0.5]
% Dibujamos el grid
\draw[step=1.0,black!50, very thin] (0,0) grid (10,8);
% Dibujamos las aristas no seleccionadas
\fill[blue!10] (0,0) -- (2,3) -- (4,0) -- (0,0);
\draw[thick, blue] (0,0) -- (2,3) -- (4,0) -- (0,0);
% Punto base
\filldraw [blue] (0,0) circle (3pt);
% Puntos IFP
\filldraw [black] (0,1) circle (2pt);
\filldraw [black] (0,2) circle (2pt);
\filldraw [black] (0,3) circle (2pt);
\filldraw [black] (0,4) circle (2pt);
\filldraw [black] (0,5) circle (2pt);
\filldraw [black!60] (1,0) circle (2pt);
\filldraw [black!60] (1,1) circle (2pt);
\filldraw [black] (1,2) circle (2pt);
\filldraw [black] (1,3) circle (2pt);
\filldraw [black] (1,4) circle (2pt);
\filldraw [black] (1,5) circle (2pt);
\filldraw [black!60] (2,0) circle (2pt);
\filldraw [black!60] (2,1) circle (2pt);
\filldraw [black!60] (2,2) circle (2pt);
\filldraw [black] (2,3) circle (2pt);
\filldraw [black] (2,4) circle (2pt);
\filldraw [black] (2,5) circle (2pt);
\filldraw [black!60] (3,0) circle (2pt);
\filldraw [black!60] (3,1) circle (2pt);
\filldraw [black] (3,2) circle (2pt);
\filldraw [black] (3,3) circle (2pt);
\filldraw [black] (3,4) circle (2pt);
\filldraw [black] (3,5) circle (2pt);
\filldraw [black] (4,0) circle (2pt);
\filldraw [black] (4,1) circle (2pt);
\filldraw [black] (4,2) circle (2pt);
\filldraw [black] (4,3) circle (2pt);
\filldraw [black] (4,4) circle (2pt);
\filldraw [black] (4,5) circle (2pt);
\filldraw [black] (5,0) circle (2pt);
\filldraw [black] (5,1) circle (2pt);
\filldraw [black] (5,2) circle (2pt);
\filldraw [black] (5,3) circle (2pt);
\filldraw [black] (5,4) circle (2pt);
\filldraw [black] (5,5) circle (2pt);
\filldraw [black] (6,0) circle (2pt);
\filldraw [black] (6,1) circle (2pt);
\filldraw [black] (6,2) circle (2pt);
\filldraw [black] (6,3) circle (2pt);
\filldraw [black] (6,4) circle (2pt);
\filldraw [black] (6,5) circle (2pt);
% Texto de filas
\draw (-0.5,0) node {0};
\draw (-0.5,1) node {1};
\draw (-0.5,2) node {2};
\draw (-0.5,3) node {3};
\draw (-0.5,4) node {4};
\draw (-0.5,5) node {5};
\draw (-0.5,6) node {6};
\draw (-0.5,7) node {7};
\draw (-0.5,8) node {8};
% Texto de columnas
\draw (0, -0.5) node {0};
\draw (1, -0.5) node {1};
\draw (2, -0.5) node {2};
\draw (3, -0.5) node {3};
\draw (4, -0.5) node {4};
\draw (5, -0.5) node {5};
\draw (6, -0.5) node {6};
\draw (7, -0.5) node {7};
\draw (8, -0.5) node {8};
\draw (9, -0.5) node {9};
\draw (10, -0.5) node {10};
% Dibujamos las dimensiones del grid
\draw (7.5, 6.5) node {$\mathcal{B}$};
\draw (10.5, 7.5) node {$g_y$};
\draw[dashed] (9,8) -- (9,8.75);
\draw[dashed] (10,8) -- (10,8.75);
\draw (9.5, 8.5) node {$g_x$};
\draw[dashed] (10,8) -- (10.75,8);
\draw[dashed] (10,7) -- (10.75,7);
\end{tikzpicture}
\caption{Regular discretization of the board $\mathcal{B}$ and Inner Fit Polygon (IFP) for a piece of triangle type $\textbf{t}$ regarding the discrete board as defined by the Dotted-Board model \cite{Toledo2013-oi}. The blue point at the origin sets the reference point of the piece, whilst black dots set its feasible positions within the board. The values $g_x$ and $g_y$ denote the horizontal and vertical dimensions of the grid, although any irregular grid might be used \cite{Cherri2018-rq}.}
\label{fig:IFP_dotted-board_model}
\end{subfigure}
\hfill
\begin{subfigure}[t]{0.48\textwidth}
\centering
\begin{tikzpicture}[scale=0.5]
% Dibujamos el grid
\draw[step=1.0,black!50, very thin] (0,0) grid (10,8);
% Dibujamos el triágulo
\fill[blue!10] (0,0) -- (2,3) -- (4,0) -- (0,0);
\draw[thick, blue] (0,0) -- (2,3) -- (4,0) -- (0,0);
% Dibujamos el cuadrado
\fill[blue!10] (6,4) -- (8,4) -- (8,6) -- (6,6) -- (6,4);
\draw[thick, blue] (6,4) -- (8,4) -- (8,6) -- (6,6) -- (6,4);
\filldraw [blue] (6,4) circle (3pt);
% Puntos NFP
\filldraw [blue] (0,0) circle (3pt);
\filldraw [black] (0,1) circle (2pt);
\filldraw [black] (0,2) circle (2pt);
\filldraw [black] (1,0) circle (2pt);
\filldraw [black] (1,1) circle (2pt);
\filldraw [black] (1,2) circle (2pt);
\filldraw [black] (2,0) circle (2pt);
\filldraw [black] (2,1) circle (2pt);
\filldraw [black] (2,2) circle (2pt);
\filldraw [black] (3,0) circle (2pt);
\filldraw [black] (3,1) circle (2pt);
% Texto de filas
\draw (-0.5,0) node {0};
\draw (-0.5,1) node {1};
\draw (-0.5,2) node {2};
\draw (-0.5,3) node {3};
\draw (-0.5,4) node {4};
\draw (-0.5,5) node {5};
\draw (-0.5,6) node {6};
\draw (-0.5,7) node {7};
\draw (-0.5,8) node {8};
% Texto de columnas
\draw (0, -0.5) node {0};
\draw (1, -0.5) node {1};
\draw (2, -0.5) node {2};
\draw (3, -0.5) node {3};
\draw (4, -0.5) node {4};
\draw (5, -0.5) node {5};
\draw (6, -0.5) node {6};
\draw (7, -0.5) node {7};
\draw (8, -0.5) node {8};
\draw (9, -0.5) node {9};
\draw (10, -0.5) node {10};
\end{tikzpicture}
\caption{Black dots define the set of discrete positions belonging to the $NFP^{(0,0)}_{ts}$ defining the infeasible positions for pieces of square type $\textbf{s}$ regarding pieces of triangle type $\textbf{t}$ when a triangle is positioned in the grid point $(0,0)$.}
\label{fig:NFP_dotted-board_model}
\end{subfigure}
\end{figure}

\emph{Computation of the no-fit polygons}. Given two polygons $X,Y \subset \mathbb{R}^2$, their Minkowski sum is defined as $X \oplus Y = \{x + y \::\: x \in X, y \in Y\}$. Using the Minkowski sum, the no-fit polygon between polygons $A$ and $B$ is defined as $NFP_{AB} = A \oplus -B$, where $A$ is the static polygon and $B$ the orbiting one. Currently, there are three main families of methods for the computation of the NFP as follows: (1) the family of orbiting methods, such as those introduced by \citet{Mahadevan1984-to}, \citet{Ghosh1991-eh}, \citet{Dean2006-wp}, \citet{Burke2007-ta}, \citet{Burke2010-hg}, and the method restricted to convex pieces of \citet{Cuninghamegreen1989-si}, whose core idea is to slide the orbiting polygon around the static polygon;  (2) the family of methods based on Minkowski sums, such as those methods proposed by \citet{Bennell2001-wu}, \citet{Bennell2008-ha}, and \citet{Milenkovic2010-pd}, whose core idea is to exploit the definition of the NFP as a Minkowski sum; and finally, (3) the family of decomposition methods, whose core idea is to decompose the problem in three steps as follows: (a) decomposition of the input polygons into convex or star-shaped sub-parts, such as done by \citet{Watson1999-kf}, \citet{Licari2011-as} and \citet[\S4]{Li1995-aj}, respectively; (b) computation of the NFP between convex or star-shaped polygons using any specialized algorithm, such as the orbiting methods of \citet{Cuninghamegreen1989-si} or Minkowski sums \cite[p.299]{De_Berg1997-bi} for convex polygons; and (c) the recombination of the convex or star-shaped NFP parts to obtain the overall NFP, as detailed by \citet[\S4]{Agarwal2002-ip}. To build the non-overlapping constraints (\ref{ineq:DB_non-overlapping}) above, we use the same approach introduced in our previous work \cite{Lastra-Diaz2024-il} in which we decompose all pieces into convex parts using the Greene's convex decomposition algorithm \cite{Greene1983-fl} implemented by the CGAL .NET wrapper \cite{Scrawk2022-on}, and then, we compute all pairwise convex no-fit polygons between convex parts of two different pieces by using the Minkowski sum for convex polygons \cite[Algorithm 1]{Lastra-Diaz2024-il}. However, the readers might use any known NFP algorithm for this aim.

\paragraph{The DB model with discrete rotations.} For the sake of completeness, we introduce below the binary decision variables $\delta_t^{d\theta}$, parameters, and index sets required by the formulation of the DB model with discrete rotations \cite{Cherri2018-rq}, denoted by DB-Rot here, which will be used in our 0-1 DB reformulations. The DB-Rot model \cite{Cherri2018-rq} is easily derived from the basic DB model \cite{Toledo2013-oi} above by substituting the binary variables $\delta_t^d$ by its extended version $\delta_t^{d\theta}$ indicating that a piece of type $t \in \mathcal{T}$ is placed in point $d \in \mathcal{D}$ with orientation $\theta \in \Theta_t$. Finally, it is needed to define the index sets for the discrete rotations $\Theta$, inner-fit polygons $\mathcal{IFP}_t^\theta$, and no-fit polygons $\mathcal{NFP}_{tu}^{d\theta}$, as detailed in eq. (\ref{index:theta_sets}) to (\ref{index:nofit_with_rotations}).

\begin{align} 
\Theta &= \bigcup\limits_{\forall t \in \mathcal{T}} \Theta_t, \; \Theta_t = \{\theta_t^k \in [0,2\pi], k = 1,\dots, K_t\}  \label{index:theta_sets}\\
\mathcal{IFP}_t^\Theta &= \{(d,\theta) \in \mathcal{D} \times \Theta_t: P_t(\theta) \oplus d \subseteq \mathcal{B}\} \label{index:inner-fit_with_rotations} \\
\mathcal{NFP}_{tu}^{d\theta} &= \{(d',\theta') \in \mathcal{IFP}_u^\Theta : int(P_t(\theta) \oplus d) \cap int(P_u(\theta') \oplus d') \neq \emptyset\} \label{index:nofit_with_rotations}
\end{align}
\begin{tabular}{lp{15cm}}
\multicolumn{2}{l}{\textbf{Extended decision variables, parameters, and index sets for the DB model with rotations \cite{Cherri2018-rq}}} \\
$\delta_t^{d\theta}$ & Binary variables indicating that one piece of type $t \in \mathcal{T}$ with orientation $\theta \in \Theta_t$ is placed in $d \in \mathcal{D}$. \\
$x_{t\theta}^M$ & Distance from the rightmost point of one piece of type $t \in T$ with orientation $\theta \in \Theta_t$ to its reference point, as shown in Figure \ref{fig:piece_vars_DB}.
\end{tabular}

\subsection{The Dotted-Board Clique Covering (DB-CC) model}
\label{sec:DB-CC_model}

As mentioned above, the main drawback of the DB model \cite{Toledo2013-oi} to be efficiently solved within the standard Branch-and-Cut framework \cite{Padberg1991-hm, Grotschel1991-qh} is that the number of non-overlapping constraints defined by Constraints (\ref{ineq:DB_non-overlapping}) grows exponentially with the number of dots and types of piece. \citet{Rodrigues2017-zl} solve this latter scalability drawback, and set the current state of the art among the family of discrete mathematical models for nesting, by proposing a new polynomial-size clique-based reformulation of the DB model called Dotted-Board Clique Covering (DB-CC) model here, which is defined by the Objective function (\ref{obj:DB-CC}) and Constraints (\ref{ineq:DB-CC_innerfit}) to (\ref{ineq:DB-CC_length_variable}) detailed below.

The representation of constraints among binary variables as conflict graphs, and the computation of clique covering cuts is a well-known technique in integer programming, as discussed by \citet{Atamturk2000-yv} and efficiently applied to the discrete nesting problem by \citet{Rodrigues2017-zl}. The conflict graph $G=(V,E)$ of any DB problem instance is defined as follows, the vertex set $V$ contains one vertex for each binary decision variable $\delta_t^d$ and one edge $(\delta_t^d,\delta_u^{d'}) \in E$ if the DB model contains one non-overlapping Constraint (\ref{ineq:DB_non-overlapping}) for these pair of variables, as detailed above. On the other hand, the subset of variables $\delta_t^d$ whose associated length bound $L_t^d$ as defined in equation (\ref{eq:DB-CC_lengthUB_perVariable}) is greater than the lower bound for the board's length $\underline{L}$ induces the sub-graph $G_{\underline{L}^+} = (V_{\underline{L}^+}, E_{\underline{L}^+}) \subset G$ that allows defining a vertex clique covering $\mathcal{V}_{\underline{L}^+}$ to partition $V_{\underline{L}^+}$ into a collection of independent vertex sets in $G_{\underline{L}^+}$. One edge clique covering $\mathcal{E}$ of the graph $G$ can be computed with any minimal Edge Clique Covering algorithm, such as those introduced by \citet{Rodrigues2021-ev} and \citet{Conte2020-mo}, whilst one vertex clique covering $\mathcal{V}_{\underline{L}^+}$ of the sub-graph $G_{\underline{L}^+}$ can be computed with any graph coloring algorithm, such as the Recursive Largest First (RLF) heuristics proposed by \citet{Leighton1979-cl} and implemented by \citet{Rodrigues2017-zl} to build the Constraints (\ref{ineq:DB-CC_innerfit}).
\bigskip\\
\textbf{Dotted-Board Clique Covering (DB-CC) model} \; \cite{Rodrigues2017-zl} 
\begin{align}
\min\limits_{\delta^d_t,L}  \quad & z = L \label{obj:DB-CC}\\
\text{s.t.} \quad & \underline{L} + \sum\limits_{\forall (d,t) \in K} (c_dg_x + x_t^M - \underline{L}) \: \delta_t^d \leq L, \quad \forall K \in \mathcal{V}_{\underline{L}^+} \label{ineq:DB-CC_innerfit} \\
& \sum\limits_{\forall d \in \mathcal{IFP}_t} \delta^d_t = q_t, \quad \forall t \in \mathcal{T} \label{ineq:DB-CC_demands} \\
& \sum\limits_{\forall (d,t) \in K} \delta^d_t \leq 1, \quad \forall K \in \mathcal{E} \label{ineq:DB-CC_non-overlapping} \\
& \delta^d_t \in \{0,1\}, \quad \forall d \in \mathcal{IFP}_t, \quad \forall t \in \mathcal{T} \label{ineq:DB-CC_binary_variables} \\
& L \in [\underline{L}, \overline{L}] \subset \mathbb{R}  \label{ineq:DB-CC_length_variable} 
\end{align}

The DB-CC model is based on the definition of the conflict graph $G=(V,E)$ and its corresponding sub-graph $G_{\underline{L}^+}$ of variables with length bounds greater than $\underline{L}$ in expressions (\ref{graph:DB_conflict_graph}) to (\ref{subgraph:DB-CC_edge_set}) and two smart transformations of the DB model detailed below, which contribute to tightening the LP relaxation of the DB-CC model in comparison with the basic DB model \cite{Toledo2013-oi}.
\begin{align} 
L_t^d &= c_dg_x + x_t^M, \quad \forall d \in \mathcal{D} \label{eq:DB-CC_lengthUB_perVariable} \\
G &= (V,E) \label{graph:DB_conflict_graph} \\
V &= \{\delta_t^d, \: \forall d \in \mathcal{IFP}_t, \forall t \in \mathcal{T} \} \label{graph:DB_vertex_set} \\
E &= \{(\delta_t^d,\delta_u^{d'}) \in V \times V : d' \in  \mathcal{NFP}^d_{tu} \wedge (t,u) \in \mathcal{T} \times \mathcal{T} \wedge d \in \mathcal{IFP}_t \}  \label{graph:DB_edge_set}  \\
G_{\underline{L}^+} &=(V_{\underline{L}^+}, E_{{\underline{L}^+}})   \label{subgraph:DB-CC_conflict}  \\
V_{\underline{L}^+} &= \{\delta_t^d \in V : L_t^d > \underline{L} \} \label{subgraph:DB-CC_vertex_set} \\
E_{\underline{L}^+} &= \{\delta_t^d,\delta_u^{d'} \in V_{\underline{L}^+} : \exists (\delta_t^d,\delta_u^{d'}) \in E \} \label{subgraph:DB-CC_edge_set}
\end{align}

\emph{Transformation of the DB non-overlapping constraints}. The non-overlapping Constraints (\ref{ineq:DB_non-overlapping}) of the DB model are substituted by their edge-clique covering reformulation (\ref{ineq:DB-CC_non-overlapping}) derived from the conflict graph $G=(V,E)$ induced by the edge inequalities (\ref{ineq:DB_non-overlapping}) above.

\emph{Transformation of the DB inner-fit constraints}. The inner-fit Constraints (\ref{ineq:DB_innerfit}) of the DB model are substituted by their vertex-clique covering reformulation (\ref{ineq:DB-CC_innerfit}).  This transformation is based on the fact that the inner-fit Constraints (\ref{ineq:DB_innerfit}) of the DB model are always met by all variables for which $L_t^d \leq \underline{L}$. Thus, these latter constraints can be removed from the model and the remaining ones factorized by using one vertex clique covering, as defined by the Constraints (\ref{ineq:DB-CC_innerfit}) above. The variables belonging to the same vertex clique $K \in \mathcal{V}_{\underline{L}^+}$ can be activated at the same time because they belong to an independent set of $G$. Thus, they can be summed in the same inequality. as done in (\ref{ineq:DB-CC_innerfit}).
 
\section{The two new 0-1 ILP reformulations for nesting}
\label{sec:our_binary_DB_model}

This section introduces two 0-1 Integer Linear Programming reformulations of the family of DB models \cite{Toledo2013-oi, Rodrigues2017-zl, Cherri2018-rq} with discrete rotations based on defining a set of binary variables to encode all feasible values for the length upper bounds of the problem. We call these two 0-1 DB ILP reformulations the Binary Dotted-Board (0-1 DB) model and the Binary Dotted-Board Clique Covering model (0-1 DB-CC), respectively. All coefficients of the constraints matrix and objective function of the new 0-1 ILP models take integer values in the set $\{1, 0, -1, N\}$. These two 0-1 ILP formulations have been devised to make the family of DB models amenable to being solved either by state-of-the-art CP solvers, such as Google CP-Sat \cite{Perron2023-nd}, or ad-hoc exact algorithms as the one proposed in section \ref{sec:parallel_backtracking}. Finally, our 0-1 ILP models open the possibility to explore SAT-based reformulations that might be solved with state-of-the-art incremental SAT solvers, such as CaDiCal \cite{Biere2024-oy}, Glucose \cite{Audemard2018-ty} and MiniSAT \cite{Een2004-zw, Sorensson2005-tr} among others. Next, we introduce the extended sets of decision variables, parameters, and index sets used in our formulations.

\textbf{Extended decision variables, parameters, and index sets}
\begin{center}
\begin{tabular}{lp{10cm}}
$z_m$ & Binary variables indicating that the length upper bound $\mathcal{L}_m > \underline{L}$ is enabled in the solution, and thus, solutions with length $L \leq \mathcal{L}_m$ are feasible.\\
$L_t^{d\theta} = c_dg_x + x_{t\theta}^M$ & Length bound induced by every feasible orientation $\theta \in \Theta_t$ of a piece of type $t \in \mathcal{T}$ when it is placed in dot $d \in \mathcal{D}$. \\
$\mathcal{L}^+ = \{\mathcal{L}_m \in (\underline{L}, \overline{L}], 1 \leq m \leq M\}$ & Set of all feasible length upper bounds as defined in (\ref{eq:BDB-feasible_length_bounds}).
\end{tabular}
\end{center}

We use the same notation introduced in Section \ref{sec:db_models} to formalize our new 0-1 DB ILP models. However, we introduce a new set of binary decision variables, denoted by $z_m$, to encode the ordered set of length upper bounds $\mathcal{L}^+$, and the index set $M$ for indexing the variables $z_m$ and their corresponding length upper bounds $\mathcal{L}_m$. To discretize the feasible length values of the problem, the statement (\ref{eq:BDB-feasible_length_bounds}) defines the set of all feasible length upper bounds $\mathcal{L}^+$, which is converted into a total ordered set $(\mathcal{L}^+, \leq_{\mathbb{R}})$ by using the standard order relation of the real numbers. All length upper bounds are sorted in ascending order, such that $\underline{L} < \mathcal{L}_1 < \mathcal{L}_2 < \dots < \mathcal{L}_M$ and one $z_m$ binary decision variable is defined to enable all feasible solutions whose associated length upper bound is equal or greater than $\mathcal{L}_m$.
\begin{align}
\mathcal{L}^+ &:= \Bigg\{\mathcal{L} \in \bigcup\limits_{\forall t \in \mathcal{T}} \Big(\bigcup\limits_{\forall (d,\theta) \in \mathcal{IFP}_t^\Theta} L_t^{d\theta} \in \mathbb{R} \Big) : \mathcal{L} > \underline{L}\Bigg\}, \quad M = |\mathcal{L}^+| \label{eq:BDB-feasible_length_bounds}
\end{align}

\subsection{The Binary Dotted-Board model}

The Binary Dotted-Board (0-1 DB) model is defined by the objective function (\ref{obj:BDB_obj_function}) and the Constraints (\ref{ineq:BDB_1}) to (\ref{vars:BDB_z_variables}), together with the decision variables, parameters, and index sets detailed above. Constraints (\ref{ineq:BDB_1}) enable the length upper bounds represented by the binary decision variable $z_m$ according to the order relationship for the upper bound set $\mathcal{L}^+$. Constraints (\ref{ineq:BDB_2}) set the demand constraints by forcing the placement of all pieces in any feasible solution. Constraints (\ref{ineq:BDB_3}) encode the non-overlapping constraints associated with any infeasible relative positioning and orientation between two pieces. Constraints (\ref{ineq:BDB_4}) enable the feasible placements of the pieces corresponding to each feasible length upper bound $\mathcal{L}_m$. Finally, Constraints (\ref{vars:BDB_delta_variables}) and (\ref{vars:BDB_z_variables}) set the domains of the decision variables.
\bigskip\\
\textbf{Binary DB (0-1 DB) model} 
\begin{align}
\min\limits_{\delta_t^{d\theta}, z_m} \quad & z = \sum\limits_{m = 1}^M z_m \label{obj:BDB_obj_function} \\
s.t.: \quad & z_m \leq z_{m - 1}, \quad 1 < m \leq M \label{ineq:BDB_1} \\
& \sum\limits_{\forall (d,\theta) \in \mathcal{IFP}_t^\Theta} \delta_t^{d\theta} = q_t, \quad \forall t \in \mathcal{T} \label{ineq:BDB_2} \\
& \delta_u^{e\theta'} + \delta_t^{d\theta} \leq 1, \quad \forall (e,\theta') \in \mathcal{NFP}_{tu}^{d\theta}, \forall (d,\theta) \in \mathcal{IFP}_t^\Theta \label{ineq:BDB_3} \\
& \sum\limits_{t = 1}^T \sum\limits_{\forall (d,\theta) \in \mathcal{IFP}_t^\Theta : L_t^{d\theta} \geq \mathcal{L}_m} \delta_t^{d\theta} \leq N z_m, \quad 1 \leq m \leq M \label{ineq:BDB_4} \\
& \delta_t^{d\theta} \in \{0,1\}, \; \forall t \in \mathcal{T}, \forall (d,\theta) \in \mathcal{IFP}_t^\Theta \label{vars:BDB_delta_variables} \\
& z_m \in \{0,1\}, \;  1 \leq m \leq M \label{vars:BDB_z_variables}
\end{align}

Given any feasible solution to the 0-1 DB model above, its length $L$ is defined by the formula (\ref{eq:L_recovering}) below, where $\mathcal{L}_z$ is the length upper bound indexed by the value of the objective function (\ref{obj:BDB_obj_function}).
\begin{equation}
L = \Big\{\begin{array}{cc}
    \underline{L}, & z = 0 \\
    \mathcal{L}_z, & z > 0
\end{array}   \label{eq:L_recovering}
\end{equation}

\subsection{The Binary DB Clique Covering model}
\label{sec:BDB-CC_model}

The Binary Dotted-Board Clique Covering (0-1 DB-CC) model is defined by the Objective function (\ref{obj:BDB-CC_obj_function}) and the Constraints (\ref{ineq:BDB-CC_1}) to (\ref{vars:BDB-CC_z_variables}). Constraints (\ref{ineq:BDB-CC_1}) enable the ordered set $(\mathcal{L}^+, \leq_\mathbb{R})$ of feasible length upper bounds. Constraints (\ref{ineq:BDB-CC_2}) set the demand constraints by forcing the placement of all pieces in any feasible solution. Constraints (\ref{ineq:BDB-CC_3}) encode the non-overlapping constraints associated to the cliques $\mathcal{K} \in \mathcal{E}$ of mutually infeasible relative placements of pieces defined by the conflict graph $G=(V,E)$ of every problem instance, where $\mathcal{E}$ is any edge clique covering of the conflict graph $G$, as explained in Section \ref{sec:DB-CC_model}. Constraints (\ref{ineq:BDB-CC_4}) enable the feasible placements of the pieces corresponding to each feasible length upper bound $\mathcal{L}_m$. Finally, Constraints (\ref{vars:BDB-CC_delta_variables}) and (\ref{vars:BDB-CC_z_variables}) set the domains of the decision variables.
\bigskip\\
\textbf{Binary DB Clique Covering (0-1 DB-CC) model}
\begin{align}
\min\limits_{\delta_t^{d\theta}, z_m} \quad & z = \sum\limits_{m = 1}^M z_m  \label{obj:BDB-CC_obj_function} \\
s.t.: \quad & z_m \leq z_{m - 1}, \quad 1 < m \leq M \label{ineq:BDB-CC_1} \\
& \sum\limits_{\forall (d,\theta) \in \mathcal{IFP}_t^\Theta} \delta_t^{d\theta} = q_t, \quad \forall t \in \mathcal{T} \label{ineq:BDB-CC_2} \\
& \sum\limits_{\forall (t,d,\theta) \in \mathcal{K}} \delta_t^{d\theta} \leq 1, \quad \forall \mathcal{K} \in \mathcal{E} \label{ineq:BDB-CC_3} \\
& \sum\limits_{t = 1}^T \sum\limits_{\forall (d,\theta) \in \mathcal{IFP}_t^\Theta : L_t^{d\theta} \geq \mathcal{L}_m} \delta_t^{d\theta} \leq N z_m,  \label{ineq:BDB-CC_4} \\
& \delta_t^{d\theta} \in \{0,1\}, \; \forall t \in \mathcal{T}, \forall (d,\theta) \in \mathcal{IFP}_t^\Theta \label{vars:BDB-CC_delta_variables} \\
& z_m \in \{0,1\}, \;  1 \leq m \in M\label{vars:BDB-CC_z_variables}
\end{align}

\section{The parallel branch-and-bound-and-prune algorithm for nesting}
\label{sec:parallel_backtracking}

This section introduces our parallel branch-and-bound-and-prune algorithm for the discrete nesting problem, which is based on a parallel backtracking enumeration with bounding and forward-checking to prune the solution space of any Dotted-Board problem instance, together with a set of ad-hoc data structures for their efficient implementation. Our Dotted-Board Parallel Backtracking (DB-PB) algorithm is made up of the following components: (1) the ad-hoc global data structures detailed in Table \ref{tab:DBPB_data_objects}; (2) the initialization procedure detailed in Algorithm \ref{alg:initialization}; (3) the main solving function detailed in Algorithm \ref{alg:DBPBSolve}; (4) the branch-and-bound-and-prune algorithm detailed in Algorithm \ref{alg:BranchAndBoundAndPrune}; (5) an ad-hoc forward-checking \cite{Haralick1980-tl} algorithm detailed in Algorithm \ref{alg:ForwardChecking}; and (6) the thread-safe incumbent updating and dominated variable pruning detailed in Algorithm \ref{alg:UpdateIncumbent}. The efficiency of DB-PB is achieved by an appropriate integration of the features detailed in the paragraphs below.

\begin{table}[h!]
\centering
\caption{Global data structures shared and required by all algorithms introduced herein.}
\label{tab:DBPB_data_objects}
\begin{tabular}{lp{11.5cm}}
\textbf{Global object} & \textbf{Description} \\
\hline
$\pi_{\delta_t^{d\theta}}^V$ & Mapping function from a DB decision variable $\delta_t^{d\theta}$ to its corresponding vertex index $v \in V$ in the conflict inverse graph $\tilde{G}$. \\
$\pi_V^{\delta_t^{d\theta}}$ & Inverse mapping from a vertex index $v \in V$ to its corresponding variable. \\
$\Pi_{\tilde{G}} \in \{0,1\}^{V \times V}$ & BitArray matrix corresponding to the adjacency matrix of inverse graph $\tilde{G}$ of the conflict graph $G = (V,E)$ together with the symmetry-breaking for identical pieces detailed in Algorithm 1, where each vertex $v \in V$ corresponds to a binary decision variable $\delta_t^{d\theta}$, such that $\Pi_{\tilde{G}}^{ij} = 0,$ if $\exists (v_i,v_j) \in  E$, or $\Pi_{\tilde{G}}^{ij} = 1$, otherwise. \\
$\pi_{inst} \in \mathbb{Z}^N$ & Integer array storing the indexes of piece instances denoted by $\pi_{inst} \in \mathbb{Z}^N$. \\
$\pi_{idVar} \in \mathbb{Z}^{V \times 2}$  & Integer matrix storing the index intervals defining the decision variables corresponding to each piece type within the BitArray $\pi_{nonDomVar}$ representing the values for all binary decision variables. \\
$\pi_{type} \in \mathcal{T}^N$ & Integer array encoding the piece type for each piece instance $P_i \in \mathcal{P}$. \\
$\pi_{nonDomVar} \in \{0,1\}^V$ & BitArray encoding the feasibility state of all binary decision variables $\delta_t^{d\theta}$.\\
$\pi_{pieceTypeVarMask} \in \{0,1\}^{\mathcal{T} \times V}$ & BitArray vector encoding the variable masks for each piece type $t \in \mathcal{T}$. \\
$L_t^V \in \mathbb{R}^V$ & Real-value array encoding the length upper bound $L_t^{d\theta}$ of each indexed decision variable $\delta_t^{d\theta}$, such that exists an unique vertex index $v \in V$ for each binary decision variable. \\
$x_{inc} \in V^N$ & Integer array encoding the index set of the binary decision variables selected in the incumbent solution for each piece instance. \\
$L_{inc} \in \mathcal{L}^+ \cup \underline{L}$ & Objective value for the incumbent solution $x_{inc}$. \\
$opt_{proof} \in \{T,F\}$ & Boolean variable indicating that the optimality has been proven.
\end{tabular}
\end{table}

\textbf{BitArray representation}. All non-overlapping constraints of any DB problem instance are encoded in the adjacency matrix $\Pi_{\tilde{G}}$ of the complement or inverse graph $\tilde{G} = (V,\tilde{E})$ of the conflict graph $G=(V,E)$, as detailed in Table \ref{tab:DBPB_data_objects}. $\Pi_{\tilde{G}}$ is represented in memory by an array of BitArray objects with only a memory space complexity of $\mathcal{O}(\frac{V^2}{8})$ bytes, which allows DB-PB implementing a very efficient binary constraint propagation by using the fast bitwise AND operation. Thus, our BitArray representation makes unnecessary any clique cover reformulation of the DB model \cite{Rodrigues2017-zl}, because our algorithms can efficiently solve any DB problem instance directly from the inverse of the conflict graph, denoted by $\tilde{G} = (V,\tilde{E})$. For instance, a DB problem instance with half a million binary decision variables and thousands of millions of constraints might be represented and solved with roughly 30 GB of RAM. However, state-of-the-art general-purpose MIP and hybrid CPSat solvers might find serious difficulties in representing and solving these large problem instances, even their clique cover reformulations. On the other hand, for even larger problem instances, we might extend our data structures and algorithms to use sparse BitArray representations or virtual memory.

\textbf{Symmetry-breaking for identical pieces}. Algorithm \ref{alg:initialization} details the initialization of the main data structures supporting our algorithms, among which the most important one is the adjacency matrix $\Pi_{\tilde{G}}$ of the conflict inverse graph, because each row of this matrix encodes the constraints induced by the selection of any feasible position and orientation (binary variable) of a piece instance in the solution. Every time a piece instance is placed on the board at one decision level of our DB-PB algorithm by selecting a feasible binary variable in the solution, the constraints encoded by its corresponding row are propagated, setting to 0 all conflicting decision variables in subsequent decision levels. To break the symmetries generated by any permutation of the binary decision variables of the same piece type, we disable all binary variables of the same piece type corresponding to the precedent positions of the activated variable, according to the non-decreasing length bound order of the variables of each piece type, as shown in step 19 of Algorithm \ref{alg:initialization}. These latter symmetries appear in our sequential solving and any CP formulation, but not in the family of DB integer programming models \cite{Toledo2013-oi, Rodrigues2017-zl, Cherri2018-rq}.

\begin{algorithm}[h]
\caption{Our initialization algorithm to build the data structures detailed in Table \ref{tab:DBPB_data_objects} that are used by DB-PB to solve any DB problem instance, which includes the symmetry-braking of all solutions induced by any permutation of the binary decision variables $\delta_t^{d\theta}$ for identical pieces during the solving stage.}
\label{alg:initialization}
\begin{algorithmic}[1]
\Require Dotted-Board problem instance
\Ensure Initialization of the main global data structures detailed in Table \ref{tab:DBPB_data_objects}.
\Procedure{Initialization}{}
\State $\pi_{\delta_t^{d\theta}}^V \gets $ sortVariableVector() \Comment{non-decreasing order of the decision variables according to $L_t^{d\theta}$}
\State $\pi_V^{\delta_t^{d\theta}} \gets (\pi_{\delta_t^{d\theta}}^V)^{-1}$  \Comment{gets the inverse mapping}
\State $\pi_{nonDomVar} \gets \textbf{1}^V$ \Comment{enabling all decision variables}
\State $L_t^V \gets \textbf{0}^V$
\For {$\forall v \in V$}
\State $L_t^V(v) \gets L_t^{d\theta} : index(\delta_t^{d\theta} = v)$ \Comment{initializes length bounds of variables}
\EndFor
\State $\Pi_{\tilde{G}} \gets \textbf{1}^{V \times V}$ \Comment{initializes the adjacency (constraint) matrix}
\For {$\forall (v, w) \in V \times V: v \neq w$}    \Comment{fills the adjacency matrix of the inverse conflict graph $\tilde{G}$}
    \If {$\exists (\delta_t^{d\theta}(v), \delta_t^{d\theta}(w)) \in E$} \Comment{$(v,w) \in V \times V$ is an edge of the conflict graph $G=(V,E)$}
        \State $\Pi_{\tilde{G}}(v,w) \gets 0$
        \State $\Pi_{\tilde{G}}(w,v) \gets 0$
    \EndIf
\EndFor
\For {$\forall t \in \mathcal{T}$}  \Comment{symmetry-breaking for identical pieces}
    \For {$i = \pi_{idVar}(t,0)$ to $\pi_{idVar}(t,1) - 1$}
        \For {$j = \pi_{idVar}(t,0)$ to $i$}
            \State $\Pi_{\tilde{G}}(i,j) \gets 0$   \Comment{disables all previous feasible variables, including the current one}
        \EndFor
    \EndFor
\EndFor
\EndProcedure
\end{algorithmic}
\end{algorithm}

\textbf{Piece-driven Binary Constraint Propagation (BCP)}. Algorithm \ref{alg:DBPBSolve} details the solving function of our DB-PB algorithm that solves the problem by using an exact backtracking enumeration with $N$ decision levels, one for each piece instance, which is based on an ad-hoc piece-driven BCP approach. Unlike the general-purpose MIP, CP, and SAT solvers, whose decision level is a single binary decision variable, the overall efficiency of DB-PB is achieved because the BCP induced by each decision level can be efficiently implemented by a single bitwise AND operation between the input feasibility BitArray and the row of $\Pi_{\tilde{G}}$ matrix corresponding to the selected decision variable, as shown in step 10 of Algorithm \ref{alg:BranchAndBoundAndPrune}, which starts with the current set of non-dominated variables represented by the BitArray $\pi_{nonDomVar}$, as shown in step 14 of Algorithm \ref{alg:DBPBSolve} and defined in Table \ref{tab:DBPB_data_objects}. This latter feature allows DB-PB to ensure that every time one decision variable is included in any partial solution, hundreds or thousands of conflicting decision variables in subsequent decision levels can be disabled in microseconds. Likewise, our piece-driven BCP approach allows DB-PB to carry out a complete and very efficient pruning of infeasible partial solutions at the earliest decision level by forward-checking, as detailed in Algorithm \ref{alg:ForwardChecking}.

\textbf{Parallel backtracking}. As mentioned above, Algorithm \ref{alg:DBPBSolve} details the solving function of our DB-PB algorithm which solves the problem by using a parallel backtracking enumeration with $N$ decision levels based on partitioning the binary decision variables corresponding to the piece type of the first piece instance, called the root level, and calling the recursive \emph{BranchAndBoundAndPrune} function detailed in Algorithm \ref{alg:BranchAndBoundAndPrune} from multiple threads, one for each subset of root decision variables. During the initialization detailed in Algorithm \ref{alg:initialization}, all piece types and instances are first ordered in non-increasing area order, and later, the decision variables for each piece type are ordered in non-decreasing length bound order. Thus, the feasible positions of the root piece disable a larger number of feasible positions for the subsequent decision levels.

\begin{algorithm}[h]
\caption{The main function of our exact Dotted-Board Parallel Backtracking (DB-PB) algorithm for discrete nesting. This function randomly partitions all feasible solutions for the first piece among all available threads, which are subsequently extended in parallel to complete a feasible solution by Algorithm \ref{alg:BranchAndBoundAndPrune}.}
\label{alg:DBPBSolve}
\begin{algorithmic}[1]
\Require The global data structures detailed in Table \ref{tab:DBPB_data_objects} and initialized by Algorithm \ref{alg:initialization} \ref{alg:initialization}.
\Ensure The incumbent solution $x_{inc}$, incumbent objective value $L_{inc}$, and optimality proof $opt_{proof}$.
\Function{DBPBSolve}{}
    \State $opt_{proof} \gets F$ \Comment{initialization of optimality proof (global variable)}
    \State $x_{inc} \gets \emptyset$  \Comment{initialization of incumbent (global variable)}
    \State $L_{inc} \gets \infty$ \Comment{initialization of incumbent objective value (global variable)}
    \State mutex $\gets$ new Mutex() \Comment{initialization of the global mutex object}
    \State $p_{root} \in [N] \gets \pi_{inst}(1)$ \Comment{gets the index of the root piece instance $P_i \in \mathcal(P)$ at decision level $1$}
    \State $t_{root} \in \mathcal{C} \gets \pi_{type}(p_{root})$ \Comment{gets the piece type for the root piece instance}
    \State $V_{root} \gets \{v \in V : \pi_{idVar}(t_{root},0) \leq v < \pi_{idVar}(t_{root},1)\}$ \Comment{gets the indexes for all variables of root piece}
    \State $V_{root} \gets randSort(V_{root})$ \Comment{random sorting of the root-level decision variable indexes}
    \While {$opt_{proof} \neq T$}
    \For {\textbf{parallel} $\forall v \in V_{root}$} \Comment{parallel solution space exploration from the root}
        \If {$\pi_{nonDomVar}(v)$} \Comment{prunes the dominated variables}
            \State $x \gets \emptyset$ \Comment{initializes an empty stack to store the variable indexes in the solution}
            \State $x \gets x \cup v$ \Comment{start a new partial solution at the root decision level (push)}
            \State $\phi_1 \gets \pi_{nonDomVar} \:\&\:  \Pi_{\tilde{G}}(v)$ \Comment{propagates all constraints induced by the variable $v$}
            \If {|x| < N}
                \If {ForwardChecking($\phi_1, 1$)} \Comment{pruning by detecting any future dead end}
                    \State BranchAndBoundAndPrune($1, \phi_1, x, L_t^V(v), \pi_{inst}$) \Comment{next decision level}
                \EndIf
            \Else
                \State UpdateIncumbent($x, L_t^V(v)$) \Comment{thread-safe updates incumbent}
            \EndIf
            \State $x \gets x \setminus v$  \Comment{removes the latest tested variable index $v$ (pop)}
            \If {(currentTime > timeLimit) \textbf{or} ($opt_{proof} = T$)}
                \State \textbf{break}
            \EndIf
        \EndIf
    \EndFor
    \EndWhile
    \State \Return \{$x_{inc}, L_{inc}, opt_{proof}$\}
\EndFunction
\end{algorithmic}
\end{algorithm}

\begin{algorithm}[h]
\caption{The Branch-and-Bound-and-Prune algorithm that visits all feasible solutions at each decision level using bounding and forward-checking for pruning the search space at the topmost decision level.}
\label{alg:BranchAndBoundAndPrune}
\begin{algorithmic}[1]
\Require (1) the decision level $\Delta \in [N]$; (2) a BitArray $\phi_{\Delta - 1} \in \{0,1\}^V$ setting the feasibility vector for all decision variables $\delta_t^{d\theta}$ at previous decision level, where $V$ denotes the total number of binary variables; (3) the input partial solution $x \in V^{\Delta - 1}$ represented by a stack containing the integer indexes of the variables selected in the solution; (4) the length upper bound $L_{\Delta - 1}^x \in \mathbb{R}$ set by input partial solution $x$; and (5) the indexes of the piece instances $\pi_{inst} \in \mathbb{Z}^N$ sorted in non-increasing area order.
\Procedure{BranchAndBoundAndPrune}{$\Delta, \phi_{\Delta-1}, x, L_{\Delta - 1}^x, \pi_{inst}$}
\If {$L_{\Delta - 1}^x < L_{inc}$} \Comment{prune current partial solution by comparing to the best objective value $L_{inc}$}
    \State $p_i \in \mathbb{Z} \gets \pi_{inst}(\Delta)$ \Comment{gets the index of the piece instance $P_i \in \mathcal(P)$ at decision level $\Delta$}
    \State $t \in \mathcal{C} \gets \pi_{type}(p_i)$ \Comment{gets the piece type for the piece instance to solve at decision level $\Delta$}
    \Statex \Comment{the for loop below tests all level decision variables in non-decreasing order}
    \Statex \Comment{regarding their length upper bound $L_t^v$, where index $v$ refers a variable $\delta_t^{d\theta}$}
    \For {$v = \pi_{idVar}(t,0)$ \textbf{to} $\pi_{idVar}(t,1) - 1$ \textbf{while} $\pi_{nonDomVar}(v)$ \textbf{and} $(L_{\Delta - 1}^x < L_{inc})$}
        \If {$\phi_{\Delta - 1}(v)$} \Comment{test that the variable $\delta_t^{d\theta}$ with index $v$ is enabled}
            \State $L_\Delta^{x \cup v} = \max\{L_{\Delta-1}^x,  L_t^V(v)\}$ \Comment{gets the length UB for the extended solution $x = x \cup v$}
            \If {$L_\Delta^{x \cup v} < L_{inc}$} \Comment{prune next partial solution by bounding before propagating it}
                \State $x \gets x \cup v$ \Comment{expands the partial solution (push)}
                \State $\phi_\Delta \gets \phi_{\Delta - 1} \:\&\:  \Pi_{\tilde{G}}(v)$ \Comment{propagates all constraints induced the variable $v$}
                \If {$|x| < N$} \Comment{tests if the solution needs to be completed}
                    \If {ForwardChecking($\phi_\Delta, \Delta$)} \Comment{pruning by detecting any future dead end}
                        \State BranchAndBoundAndPrune($\Delta + 1, \phi_\Delta, x, L_\Delta^{x \cup v}, \pi_{inst}$) \Comment{next decision level}
                    \EndIf
                \Else
                    \State UpdateIncumbent($x, L_\Delta^{x \cup v}$) \Comment{thread-safe incumbent updating}
                \EndIf
                \State $x \gets x \setminus v$  \Comment{removes the latest tested variable index $v$ (pop)}
            \EndIf
        \EndIf
    \EndFor
\EndIf
\EndProcedure
\end{algorithmic}
\end{algorithm}

\textbf{Branch-and-Bound-and-Prune}. Algorithm \ref{alg:BranchAndBoundAndPrune} implements the core function of our DB-PB algorithm, which introduces a chronological backtracking enumeration empowered with the following pruning techniques. First, pruning by bounding on the length of the current partial solutions in steps 2, 5 and 8. Second, pruning of dominated variables by escaping from the loop in step 5, because all variables $\delta_t^{d\theta}$ encoding the feasible positions and orientations of each piece type $t \in \mathcal{T}$ are ordered in non-decreasing length-bound order, according to their $L_t^{d\theta}$ bounds. Thus, the loop 5 can be aborted whenever any dominated variable is found. And third, pruning by forward-checking in step 12 to avoid the expansion of any partial solution that would lead to a dead end in any subsequent decision level. As mentioned above, our BitArray representation allows both a very efficient constraint propagation in step 10 and efficient forward-checking in step 12. The branching order of DB-PB is fixed a priori by the length-bound ordering of the variables, which allows saving computational time and escaping from the loop in step 5. Despite other branching heuristics are possible, they might require an extra computational cost that would negatively impact the global performance of our algorithms.

\textbf{Forward-checking}. Algorithm \ref{alg:ForwardChecking} details our ad-hoc forward-checking algorithm based on pruning dead ends by checking in step 7 that the current partial solution has not disabled all decision variables for the subsequent decision levels defined by their corresponding piece types. The input feasibility vector $\phi_\Delta$ encodes the propagation of all conflicts from the root to the current decision level $\Delta$. The current performance of our forward-checking algorithm avoids DB-PB's need to implement costly conflict-driven learning approaches to filter recurrent conflict sets of infeasible variables, as commonly done by most of CP \cite{Dechter1990-ay, Frost1994-wz, Chai2003-wt, Perron2023-nd}, SAT \cite{Marques-Silva2009-ij, Ohrimenko2009-ak, Feydy2009-ky, Stuckey2010-qh}, MIP \cite{Achterberg2007-sy, Achterberg2013-qz, Witzig2017-qn, Witzig2021-is, Nieuwenhuis2023-vu}, and hybrid CP-Sat solvers \cite{Perron2023-nd}.

\begin{algorithm}[t]
\caption{Our forward-checking algorithm to prune the expansion of partial solutions leading to infeasible solutions (dead ends) by checking that there still are feasible solutions for all pieces at the next decision levels.}
\label{alg:ForwardChecking}
\begin{algorithmic}[1]
\Require (1) a BitArray $\phi_\Delta \in \{0,1\}^V$ setting the feasibility vector for all decision variables $\delta_t^{d\theta}$, where $V$ is equal to the total number of decision variables; and (2) $\Delta \in [N]$ denotes the current decision level.
\Ensure Boolean value indicating that the current partial solution is feasible or sets a conflict (dead end).
\Function{ForwardChecking}{$\phi_\Delta, \Delta$}
\State feasible $\gets T$
\State $\mathcal{T}_{next} \subset \mathcal{T} \gets \pi_{nextPieceTypes}(\Delta)$ \Comment{retrieves the set of piece types for next decision levels}
\For {$\forall t \in \mathcal{T}_{next}$}
\State $\phi_{t} \in \{0,1\}^V \gets clone(\pi_{pieceTypeVarMask}(t))$ \Comment{gets the variable mask for the piece type}
\State $\phi_t^\Delta \gets \phi_{t} \:\&\: \phi_\Delta$ \Comment{binary constraint forward-propagation for a child piece type}
\If {$\phi_t^\Delta = \textbf{0}^V$} \Comment{check that there exist feasible solutions for $t \in \mathcal{T}$}
    \State feasible $\gets F$
    \State \textbf{break}
\EndIf
\EndFor
\State \Return feasible
\EndFunction
\end{algorithmic}
\end{algorithm}

\begin{algorithm}[h]
\caption{Thread-safe incumbent update function called from any thread.}
\label{alg:UpdateIncumbent}
\begin{algorithmic}[1]
\Require (1) a new incumbent candidate $x \in V^N$, and (2) its associated objective incumbent value $L_x \in \mathbb{R}$.
\Procedure{UpdateIncumbent}{$x, L_x$}
\State mutex.Wait() \Comment{waits until the incumbent data is released}
\If {$opt_{proof} = F$ \textbf{and} $L_x < L_{inc}$} \Comment{checks that none thread had improved the incumbent}
    \For {$\forall v \in V$}
        \If {$L_t^V(v) >= L_x$}
            \State $\pi_{nonDomVar}(v) \gets F$ \Comment{disables all dominated decision variables}   
        \EndIf
    \EndFor
    \For {$\forall t \in \mathcal{T}$} \Comment{disable all dominated variables in the mask of each piece type}
        \State $\pi_{pieceTypeVarMask}(t) \gets \pi_{pieceTypeVarMask}(t) \:\&\: \pi_{nonDomVar}$  
    \EndFor
    \State $L_{inc} \gets L_x$ \Comment{updates incumbent objective value}    
    \State $x_{inc} \gets x$ \Comment{updates incumbent solution}
    \If {$L_{inc} = \underline{L}$} \Comment{optimality test}
        \State $opt_{proof} \gets T$
    \EndIf
    \If {$opt_{proof} = T$}
        \State $\pi_{nonDomVar}(v) \gets \textbf{0}^V$ \Comment{disables all variables}
    \EndIf
\EndIf
\State mutex.Release() \Comment{releases the mutex}
\EndProcedure
\end{algorithmic}
\end{algorithm}
\clearpage

\begin{algorithm}[h]
\caption{Our DB Parallel Backtracking Lower Bound (DB-PB-LB) algorithm, which searches for a tight lower bound by iteratively solving the problem instance from the trivial lower bound (\ref{eq:trivial_lower_bound}) up to the optimal.}
\label{alg:lower_bound}
\begin{algorithmic}[1]
\Require The global data structures detailed in Table \ref{tab:DBPB_data_objects} and initialized by Algorithm \ref{alg:initialization}. 
\Ensure A tighter lower bound value $\underline{L}$, the optimal solution $x_{inc}$, and optimality proof $opt_{proof}$, if found.
\Function{DBPBSolveLowerBound}{}
    \State $opt_{proof} \gets F$ \Comment{initialization of optimality proof (global variable)}
    \State $x_{inc} \gets \emptyset$  \Comment{initialization of incumbent (global variable)}
    \State $L_{inc} \gets \infty$ \Comment{initialization of incumbent objective value (global variable)}
    \State $\underline{L} \gets \underline{L}_{trivial}$ \Comment{Initialization to the trivial LB value in equation (\ref{eq:trivial_lower_bound})}
    \State mutex $\gets$ new Mutex() \Comment{initialization of the global mutex object}
    \State $p_{root} \in [N] \gets \pi_{inst}(1)$ \Comment{gets the index of the root piece instance $P_i \in \mathcal(P)$ at decision level $1$}
    \State $t_{root} \in \mathcal{C} \gets \pi_{type}(p_{root})$ \Comment{gets the piece type for the root piece instance}
    \While {$opt_{proof} = F$}
        \For {$\forall v \in V$}  \Comment{enables only the non-dominated decision variables}
            \If {$L_t^V(v) \leq \underline{L}$}
                \State $\pi_{nonDomVar}(v) \gets T$
            \Else
                \State $\pi_{nonDomVar}(v) \gets F$
            \EndIf
        \EndFor
        \State $V_{root} \gets \{v \in V : \pi_{idVar}(t_{root},0) \leq v < \pi_{idVar}(t_{root},1)\}$ \Comment{gets all root decision variables}
        \State $V_{root} \gets randSort(V_{root})$ \Comment{random sorting of the root-level decision variable indexes}
        \For {\textbf{parallel} $\forall v \in V_{root}$} \Comment{parallel solution space exploration from the root}
            \If {$\pi_{nonDomVar}(v)$} \Comment{prunes the dominated variables with $L_t^V(v) \geq L_{inc}$}
                \State $x \gets \emptyset$ \Comment{initializes an empty stack to store the variable indexes in the solution}
                \State $x \gets x \cup v$ \Comment{start a new partial solution at the root decision level (push)}
                \State $\phi_1 \gets \pi_{nonDomVar} \:\&\:  \Pi_{\tilde{G}}(v)$ \Comment{propagates all constraints induced by the variable $v$}
                \If {|x| < N}
                    \If {ForwardChecking($\phi_1, 1$)} \Comment{pruning by detecting any future dead end}
                        \State BranchAndBoundAndPrune($1, \phi_1, x, L_t^V(v), \pi_{inst}$) \Comment{next decision level}
                    \EndIf
                \Else
                    \State UpdateIncumbent($x, L_t^V(v)$) \Comment{thread-safe update incumbent}
                \EndIf
                \State $x \gets x \setminus v$  \Comment{removes the latest tested variable index $v$ (pop)}
                \If {(currentTime > timeLimit) \textbf{or} ($opt_{proof} = T$)}
                    \State \textbf{break}
                \EndIf
            \EndIf
        \EndFor
        \If {$opt_{proof} = F$} \Comment{test the infeasibility of current lower bound $\underline{L}$}
            \State $\underline{L} \gets \arg\min\limits_{\lambda}\{\lambda \in (\mathcal{L}^+, \leq_{\mathbb{R}}): \lambda > \underline{L}\}$) \Comment{raises the lower bound to the next feasible value}
        \EndIf
    \EndWhile    
    \State \Return \{$\underline{L}, x_{inc}, opt_{proof}$\}
\EndFunction
\end{algorithmic}
\end{algorithm}
\clearpage

\textbf{Incumbent update}. Algorithm \ref{alg:UpdateIncumbent} details the thread-safe function responsible of updating the incumbent values in steps 12 to 13 and of disabling all dominated variables in step 6, as well as the masks of decision variables for each piece type in step 10. Finally, this latter function checks the optimality of the new incumbent in step 15 and disables all variables whenever the optimal solution has been found, which produces the automatic suspension of all backtracking loops in all threads. We use a mutually-exclusive (mutex) object to update the incumbent values and feasibility state of all decision variables in a thread-safe manner. The \emph{UpdateIncumbent} is the only synchronization point among threads to avoid costly interruptions.

\section{The new lower bound algorithm for nesting}
\label{sec:new_lower_bound}

Algorithm \ref{alg:lower_bound} details our DB-PB-LB lower bound algorithm, which is based on using our base DB-PB algorithm introduced in Section \ref{sec:parallel_backtracking} to solve any problem instance by setting $\overline{L}$ to the trivial lower bound to disable all dominated variables and iteratively increasing the lower bound up to find the optimal solution or reaching the time limit. Algorithm \ref{alg:lower_bound} iteratively fixes the length lower and upper bounds of our DB-PB data structures to each feasible length bound in the ordered set $(\mathcal{L}^+ \cup \underline{L}, \leq_{\mathbb{R}})$, starting with the trivial lower bound $\underline{L}$ in formula (\ref{eq:trivial_lower_bound}). Then, DB-PB-LB disables all dominated variables $v \in V$, such that $L_t^V(v) > \underline{L}$, as done in steps 9 to 13. Next, it searches for any feasible solution in one iteration of the main loop detailed in steps 8 to 37. If any feasible is found, then it is optimal, otherwise DB-PB-LB proves the infeasibility for the current length lower bound, and thus, $\underline{L}$ can be increased to the next feasible length upper bound in $\mathcal{L}_m \in \mathcal{L}^+$ for the next iteration, as done in step 35. In short, DB-PB-LB uses our basic DB-PB algorithm to solve one problem instance from the trivial lower bound (\ref{eq:trivial_lower_bound}) up to the optimal.
\section{Evaluation}
\label{sec:evaluation}

The goals of the experiments in this section are as follows: (1) to evaluate the performance of our new exact DB-PB algorithm; (2) to test our main hypothesis, claiming that an ad-hoc exact algorithm requiring no preprocessing might be a better option to solve the DB
model than the costly general-purpose Branch-and-Cut or hybrid CP-SAT methods by comparing our new DB-PB algorithm with our replication of the DB-CC model \cite{Rodrigues2017-zl} and our new reformulation 0-1 DB-CC with the state-of-the-art Gurobi and Google CPSat solvers; (3) to evaluate the performance of our new 0-1 ILP reformulation of the family of DB models \cite{Toledo2013-oi, Rodrigues2017-zl}; (4) to carry out a fair comparison of the performance of our new exact methods with current state-of-the-art exact methods for discrete nesting by implementing them into the same hardware and software platform; (5) to replicate the DB-CC model \cite{Rodrigues2017-zl} from scratch; (6) to develop a reproducible benchmark of state-of-the-art exact methods for the discrete besting problem based on our software implementation of all models and exact algorithms evaluated herein into the same software library, which is provided as supplementary material (see Appendix \ref{sec:appendix_B}); (7) the independent replication and confirmation of previous findings and results reported in \cite{Rodrigues2017-zl}; and finally, (8) to elucidate the current state of the art of the problem in a sound and reproducible way.

\begin{table}[h!]
\centering
\caption{Implementation details for all mathematical programming models and algorithms evaluated herein. All models are built by ordering the pieces by non-increasing area to build the models.}
\label{tab:models_evaluated}
\begin{tabular}{lp{2cm}p{11cm}}
Method & Reference & Implementation details \\
\hline
\\
DB-CC & \makecell[l]{Rodrigues and \\ Toledo \cite{Rodrigues2017-zl}} & Exact replication of the Dotted-Board Clique Covering model \cite{Cherri2016-jf} as detailed in Section \ref{sec:DB-CC_model}, which is solved by Gurobi 11.0.3 with default parameters. \\
0-1 DB-CC & this work & Our reformulation of the DB-CC model \cite{Rodrigues2017-zl} called 0-1 DB-CC, which is introduced in Section \ref{sec:BDB-CC_model} and solved by Google CPSat 9.11.4210  with default parameters.\\
DB-PB & this work & Exact Dotted-Board Parallel Backtracking (DP-PB) algorithm introduced in Section \ref{sec:parallel_backtracking}, which is implemented in C\# language and .NET 8 platform. \\
DB-CC-LB & \makecell[l]{Rodrigues and \\ Toledo \cite{Rodrigues2017-zl}} & Lower bound algorithm proposed by \citet[Algorithm 1]{Rodrigues2017-zl} as a by-product of their DB-CC model whose results are copied for comparison herein.\\
DB-PB-LB & this work & New lower-bound algorithm for the discrete nesting problem introduced in Section \ref{sec:new_lower_bound}, which is implemented in C\# language and .NET 8 platform.
\end{tabular}
\end{table}

\subsection{Experimental setup}
\label{sec:experimental_setup}

We exactly replicate the same experiments carried out by \citet[Table 1]{Rodrigues2017-dy} to evaluate all mathematical programming models and exact algorithms detailed in Table \ref{tab:models_evaluated} in the same set of eighty-one instances proposed by the former authors. Table \ref{tab:models_evaluated} lists the state-of-the-art integer programming models and exact algorithms for the discrete nesting problem evaluated herein and the details of our software implementation. We use the same iterative LB values for the length of the problem instances that were generated by the lower-bound algorithm of \citet[Algorithm 1]{Rodrigues2017-zl} and used in their experiments, as reported in \cite[Table 2]{Rodrigues2017-dy}, whilst the length upper bounds defined by the length of the board for each problem instance in Table \ref{tab:features_dataset} were also obtained by \citet{Rodrigues2017-dy} with a heuristics.

All our experiments are based on our software implementation of all mathematical programming models and exact algorithms evaluated herein in the same software and hardware Windows 10 Pro platform. Thus, all our results can be directly compared in a fair and unbiased way. Our software is coded in C\# programming language and .NET 8 platform on Windows 10. However, it is written in portable C\# .NET 8 platform, which might be used in Windows or UBUNTU platforms. The DB-CC problem instances of \citet{Rodrigues2017-zl} are solved with Gurobi 11.0.3, whilst the problem instances of our new 0-1 DB-CC model are solved with Google CPSat 9.11.4210. Google CPSat \cite{Perron2023-nd} is a state-of-the-art lazy-clause generation hybrid CP solver \cite{KrupkeUnknown-fe} based on combining four solving strategies as follows: a traditional Branch-and-Cut MIP solver, a Constraint Programming (CP) solver on top of an incremental SAT engine, and a Large Neighbourhood Search (LNS) matheuristics based on RINS \cite{Danna2005-vt}. There are many works reporting competitive or outperforming results of CP formulations regarding their equivalent state-of-the-art MILP formulations in scheduling problems \cite{Wang2015-ls, Lunardi2020-yd, Naderi2023-ok}, resource availability \cite{Kreter2018-jg}, and guillotine rectangular cutting \cite{Polyakovskiy2023-bu}. For instance, \citet{Polyakovskiy2023-bu} introduce a CP formulation for a guillotine rectangular cutting problem that outperforms state-of-the-art MILP formulations solved with Gurobi when it is solved with Google CPSat. For these reasons, we argue that the comparison of our new ad-hoc DB-PB exact algorithm with the state-of-the-art MILP models for discrete nesting solved with Gurobi and Google CPSat is a sound way of testing our main hypothesis. 

Our source code and a pre-compiled version of our software are publicly available in our reproducibility dataset \cite{Lastra-Diaz2025-gp}. We have replicated the state-of-the-art DB-CC model \cite{Rodrigues2017-zl} from scratch by integrating all required steps to evaluate all problem instances as follows: (1) preprocessing to decompose the pieces into convex parts based on the CGAL implementation of the Greene's algorithm \cite{Greene1983-fl} to compute subsequently the no-fit polygons between the convex parts of pieces; (2) computation of the vertex clique covering of every problem instance using our own software implementation of the RLF algorithm \cite{Leighton1979-cl}, as proposed in \cite{Rodrigues2017-zl}; (3) computation of the edge clique covering of every problem instance using the ECC8 software library of \citet{Conte2020-mo} instead of the MAX1\_EK algorithm of \citet{Rodrigues2021-ev} proposed in \cite{Rodrigues2017-zl}, because no software implementation of MAX1\_EK is available; (4) in-memory building of the DB-CC and 0-1 DB-CC models by calling the C\# API of Gurobi 11.0.3 and Google CPSat 9.11.4210, respectively; and (5) resolution of the optimization models by calling the proper solver's functions with its default parameters. All our experiments were implemented on a Windows 10 Pro laptop computer with an Intel Core i7-10510U CPU@1.8GHz and 48 GB of RAM, which is a modest performance CPU with 4 cores (up to 8 parallel threads) that is 19\% slower\footnote{\url{https://cpu.userbenchmark.com/Compare/Intel-Core-i7-2600-vs-Intel-Core-i7-10510U/620vsm891469}} than the Intel Core i7-2600 CPU@3.4 GHz used by \citet{Rodrigues2017-zl} in their experiments. Thus, any performance improvement in our experiments replicating the DB-CC \cite{Rodrigues2017-zl} model regarding the results reported by \citet[Table 4]{Rodrigues2017-zl} can be only attributed to the technological advances in Gurobi 11.0.3.

\subsection{Reproducing our benchmarks}

All our experiments were generated by running the \emph{DBLNSNestingApp} program distributed with our DBLNSNesting V1R1 software library \cite{Lastra-Diaz2025-gp} with four reproducible benchmark files in XML file format that include all problem instances and methods evaluated herein. All problem instances are defined in ESICUP XML file format and provided with our raw output data in our reproducibility dataset \cite{Lastra-Diaz2025-gp}. The evaluation of every exact method in Table \ref{tab:models_evaluated} generates a Scalable Vector Graphics (SVG) file and another LaTex TikZ picture text file providing an image of the solution, and a raw output file in (*.csv) file format reporting the following data for each problem instance: (1) name of the problem instance; (2) number of pieces; (3) nesting efficiency; (4-5) lower and upper bounds of the solution; (6) solution gap; (7) number of binary variables; (8) number of B\&B, CP, or backtracking nodes; (9) running time in seconds; (10) the number of constraints; (11) time in seconds to get the best solution (only for DB-PB method); (12) backtracking node where the best solution was found; and (13) building time of the model in seconds. All the data reported here is automatically generated by running the \emph{benchmark\_report} R-language script file on the collection of raw output files. Finally, all our methods, experiments, and results can be exactly reproduced by following the instructions detailed in Appendix B, which also explains how to evaluate the exact methods models detailed in Table \ref{tab:models_evaluated} in other problem instances provided in the ESICUP XML file format. The aim of our reproducibility resources is not only to allow our results to be easily reproduced, but also to evaluate unseen instances provided in the ESICUP XML file format. Likewise, to allow the direct evaluation of all problem instances with the command-line tools provided by Gurobi, Google CPSat, or other MIP solvers without using our software, our reproducibility dataset \cite{Lastra-Diaz2025-gp} provides the DB-CC \cite{Rodrigues2017-zl} and 0-1 DB-CC models generated for all problem instances in the LP and MPS file formats by Gurobi, and all the corresponding CPSat models in *.txt file format generated for the 0-1 DB-CC problem instances.

\subsection{Evaluation metrics}

To evaluate the performance of the exact mathematical programming models and algorithms detailed in table \ref{tab:models_evaluated}, we compare the Upper Bound (UB), Lower Bound (LB), terminating gap in \%, as well as the building time of the model, resolution time, and total time in seconds obtained by each exact method in the evaluation of all problem instances. To provide a fair and unbiased comparison of the performance of all exact methods evaluated here, we adopt the same approach proposed by \citet{Cherri2016-jf} and \citet{Lastra-Diaz2024-il}, which uses performance profiles \cite{Dolan2002-wd} based on the ratios between the computation time of each method and the best computation time obtained by any method as a performance metric. Let $\Phi$ and $M$ be the sets of problem instances and exact methods evaluated in this study, respectively. Then the ratio of computation times $t_{\phi,m}$  for each exact method $m \in M$ is defined by $r_{\phi,m}$ in formula (\ref{eq:ratio}), considering that the computation time is infinite whenever the methods are unable to solve the problem instance up to optimality. We assume an arbitrary parameter $r_M \geq r_{\phi,m} \: \forall (\phi,m) \in \Phi \times M$, such that $r_{\phi,m} = r_M$ if and only if the method $m$ is unable to solve the problem $\phi$. \citet{Dolan2002-wd} define a \emph{performance profile} of a solver or optimization model as the Cumulative Distribution Function (CDF) of the computation time ratios $r_{\phi,m}$, denoted by $\rho_m(\tau)$, and defined in formula (\ref{eq:performance_profile}) below. The performance profiles defined by the $\rho_m(\tau)$ function set a well-founded and broadly adopted metric to compare exact optimization methods by avoiding any bias derived from a particular set of problem instances and dealing with those cases in which the methods are unable to solve the problem up to optimality.
\begin{linenomath*}
\begin{align}
r_{\phi,m} &= \frac{t_{\phi,m}}{\min \{t_{\phi,i}: i \in M \}}, \quad \forall (\phi, m) \in \Phi \times M \label{eq:ratio} \\
\rho_m(\tau) &= \frac{1}{n_\phi} \{\phi \in \Phi : r_{\phi,m} \leq \tau\}, \quad \forall (m,\tau) \in M \times [1, r_M] \label{eq:performance_profile}
\end{align}
\end{linenomath*}

\subsection{Results obtained}

\begin{figure}[h]
\centering
% SOLUTION for problem instance = threep3w9
% Scaling factor for coordinates = 1
\begin{tikzpicture}[scale = 0.3]
% Draw the convex parts of piece2
\fill[black!15](9,6) -- (9,3) -- (12,3) -- (12,6) -- (9,6);
% Draw the boundary of piece2
\draw[black](9,6) -- (9,3) -- (12,3) -- (12,6) -- (9,6);
% Draw the convex parts of piece2
\fill[black!15](9,3) -- (9,0) -- (12,0) -- (12,3) -- (9,3);
% Draw the boundary of piece2
\draw[black](9,3) -- (9,0) -- (12,0) -- (12,3) -- (9,3);
% Draw the convex parts of piece2
\fill[black!15](8,9) -- (8,6) -- (11,6) -- (11,9) -- (8,9);
% Draw the boundary of piece2
\draw[black](8,9) -- (8,6) -- (11,6) -- (11,9) -- (8,9);
% Draw the convex parts of piece1
\fill[black!15](3,7) -- (5,5) -- (7,7) -- (5,9) -- (3,7);
% Draw the boundary of piece1
\draw[black](3,7) -- (5,5) -- (7,7) -- (5,9) -- (3,7);
% Draw the convex parts of piece1
\fill[black!15](2,3) -- (4,1) -- (6,3) -- (4,5) -- (2,3);
% Draw the boundary of piece1
\draw[black](2,3) -- (4,1) -- (6,3) -- (4,5) -- (2,3);
% Draw the convex parts of piece1
\fill[black!15](0,6) -- (2,4) -- (4,6) -- (2,8) -- (0,6);
% Draw the boundary of piece1
\draw[black](0,6) -- (2,4) -- (4,6) -- (2,8) -- (0,6);
% Draw the convex parts of piece3
\fill[black!15](5,4) -- (9,4) -- (7,7) -- (5,4);
% Draw the boundary of piece3
\draw[black](5,4) -- (9,4) -- (7,7) -- (5,4);
% Draw the convex parts of piece3
\fill[black!15](4,0) -- (8,0) -- (6,3) -- (4,0);
% Draw the boundary of piece3
\draw[black](4,0) -- (8,0) -- (6,3) -- (4,0);
% Draw the convex parts of piece3
\fill[black!15](0,0) -- (4,0) -- (2,3) -- (0,0);
% Draw the boundary of piece3
\draw[black](0,0) -- (4,0) -- (2,3) -- (0,0);
% Draw the bounding box of the solution
\draw[black](0,0) -- (12,0) -- (12,9)-- (0,9) -- (0,0);
\draw (6,-1) node {threep3w9};
\draw (6,-2) node {$L^* = 12$};
\draw (6,-3.5) node {$\Theta_1=\{0\},\Theta_2=\{0\},\Theta_3=\{0\}$};
\end{tikzpicture}
\hspace{1cm}
% SOLUTION for problem instance = threep3w9
% Scaling factor for coordinates = 1
\begin{tikzpicture}[scale = 0.3]
% Draw the convex parts of piece2
\fill[black!15](8,4) -- (8,1) -- (11,1) -- (11,4) -- (8,4);
% Draw the boundary of piece2
\draw[black](8,4) -- (8,1) -- (11,1) -- (11,4) -- (8,4);
% Draw the convex parts of piece2
\fill[black!15](5,3) -- (5,0) -- (8,0) -- (8,3) -- (5,3);
% Draw the boundary of piece2
\draw[black](5,3) -- (5,0) -- (8,0) -- (8,3) -- (5,3);
% Draw the convex parts of piece2
\fill[black!15](4,6) -- (4,3) -- (7,3) -- (7,6) -- (4,6);
% Draw the boundary of piece2
\draw[black](4,6) -- (4,3) -- (7,3) -- (7,6) -- (4,6);
% Draw the convex parts of piece1
\fill[black!15](7,6) -- (9,4) -- (11,6) -- (9,8) -- (7,6);
% Draw the boundary of piece1
\draw[black](7,6) -- (9,4) -- (11,6) -- (9,8) -- (7,6);
% Draw the convex parts of piece1
\fill[black!15](1,2) -- (3,0) -- (5,2) -- (3,4) -- (1,2);
% Draw the boundary of piece1
\draw[black](1,2) -- (3,0) -- (5,2) -- (3,4) -- (1,2);
% Draw the convex parts of piece1
\fill[black!15](0,5) -- (2,3) -- (4,5) -- (2,7) -- (0,5);
% Draw the boundary of piece1
\draw[black](0,5) -- (2,3) -- (4,5) -- (2,7) -- (0,5);
% Draw the convex parts of piece3
\fill[black!15](9,9) -- (5,9) -- (7,6) -- (9,9);
% Draw the boundary of piece3
\draw[black](9,9) -- (5,9) -- (7,6) -- (9,9);
% Draw the convex parts of piece3
\fill[black!15](3,6) -- (7,6) -- (5,9) -- (3,6);
% Draw the boundary of piece3
\draw[black](3,6) -- (7,6) -- (5,9) -- (3,6);
% Draw the convex parts of piece3
\fill[black!15](5,9) -- (1,9) -- (3,6) -- (5,9);
% Draw the boundary of piece3
\draw[black](5,9) -- (1,9) -- (3,6) -- (5,9);
% Draw the bounding box of the solution
\draw[black](0,0) -- (11,0) -- (11,9)-- (0,9) -- (0,0);
\draw (5.5,-1) node {threep3w9};
\draw (5.5,-2) node {$L^* = 11$};
\draw (5.5,-3.5) node {$\Theta_1=\{0,45^{\circ}\},\Theta_2=\{0\},\Theta_3=\{0,180^{\circ}\}$};
\end{tikzpicture}
\caption{Two solutions obtained by DB-PB for the same problem instance without rotations and with them.}
\label{fig:rotation_example}
\end{figure}
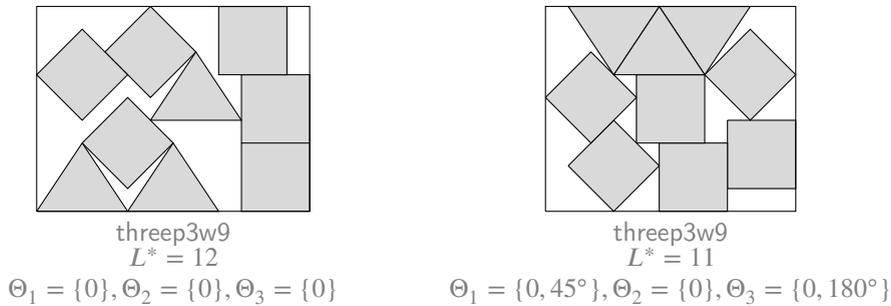

Table \ref{tab:features_dataset} details the features of the eighty-one problem instances used for comparing all exact methods evaluated herein, which exactly correspond to the set of problem instances proposed by \citet[Table 1]{Rodrigues2017-zl} for the evaluation of their DB-CC model. Table \ref{tab:lower_bound_results} compares the results of the iterative LB algorithm of \citet[Table 2]{Rodrigues2017-zl} with our new DB-PB-LB lower bound algorithm in raising the trivial LB (\ref{eq:trivial_lower_bound}) of the problem. The results reported in Table \ref{tab:lower_bound_results} for the iterative LB algorithm were copied from \cite[Table 2]{Rodrigues2017-zl}, whilst our new DB-PB-LB algorithm was evaluated herein with a time limit of 600 seconds. The column named \emph{Optimality proof} in Table \ref{tab:lower_bound_results} indicates when our DB-PB-LB algorithm found the optimal solution.  Figure \ref{fig:rotation_example} compares the optimal solution obtained by our DB-PB algorithm for the same problem instance, either without considering rotations or considering two discrete rotations for two distinct pieces, to show the correctness of the new mathematical models and algorithms considering rotations.

Table \ref{tab:first_results} reports the Upper Bound (UB), Lower Bound (LB), and terminating gap in \%, as well as the preprocessing time, resolution time, and total running time in seconds obtained by the DB-CC \cite{Rodrigues2017-zl} model solved with Gurobi and our new exact DB-PB algorithm in the evaluation of all problem instances detailed in Table \ref{tab:features_dataset} with a Time Limit (TL) of 1 hour, whilst Table \ref{tab:second_results} reports the results comparing our new 0-1 DB-CC model solved with Google CPSat with our new exact DB-PB algorithm. The preprocessing times include the running times to compute the vertex and edge clique coverings, and the time to build the models in memory before solving them. The best total solution times or UB values are shown in bold red, whilst the ’X’ symbol in Tables \ref{tab:first_results} and \ref{tab:second_results} denotes that no feasible solution was found. The problem instances solved up to optimality within 1 hour for the first time in the literature are marked with an asterisk (*).

Figures \ref{fig:performance_profile_solution_time} and \ref{fig:performance_profile_total_time} show the performance profile curves comparing the performance ratio  (\ref{eq:performance_profile}) for resolution time and total time obtained by all methods in the evaluation of all problem instances. Table \ref{tab:binary_variables} reports the number of binary variables and constraints required by each mathematical programming model and our exact algorithm for solving all problem instances. Figure \ref{fig:solved_open_instances} shows the optimal solutions for seventeen open problem instances solved up to optimality within one hour for the first time in the literature. Figure \ref{fig:improved_UB_open_instances} shows the new best solutions obtained within one hour for ten open problem instances by reducing the current upper bounds reported in \cite{Rodrigues2017-zl}. Fifty-six problem instances are solved up to optimality within one hour. 

\begin{table}[h!]
\centering
\begin{minipage}[t]{0.3\textwidth}
\caption{Features of the problem instances.}
\label{tab:features_dataset}
\tiny
\begin{tabular}{lcccc}
& \multicolumn{2}{c}{\textbf{Pieces}} & \multicolumn{2}{c}{\textbf{Board}} \\
\textbf{Instance} & \textbf{Total} & \textbf{Types} & \textbf{Width} & \textbf{Length} \\ 
\hline
three &    3 &    3 &    7 &    7 \\ 
threep2 &    6 &    3 &    7 &   11 \\ 
threep2w9 &    6 &    3 &    9 &    9 \\ 
threep3 &    9 &    3 &    7 &   16 \\ 
threep3w9 &    9 &    3 &    9 &   13 \\ 
shapes4 &    4 &    4 &   13 &   24 \\ 
shapes8 &    8 &    4 &   20 &   28 \\ 
SHAPES-2 &    8 &    4 &   40 &   16 \\ 
SHAPES-4 &   16 &    4 &   40 &   28 \\ 
SHAPES-5 &   20 &    4 &   40 &   35 \\ 
SHAPES-7 &   28 &    4 &   40 &   48 \\ 
SHAPES-9 &   34 &    4 &   40 &   54 \\ 
SHAPES-15 &   43 &    4 &   40 &   67 \\ 
blasz2 &   20 &    4 &   15 &   27 \\ 
BLAZEWICZ1 &    7 &    7 &   15 &    8 \\ 
BLAZEWICZ2 &   14 &    7 &   15 &   16 \\ 
BLAZEWICZ3 &   21 &    7 &   15 &   22 \\ 
BLAZEWICZ4 &   28 &    7 &   15 &   29 \\ 
BLAZEWICZ5 &   35 &    7 &   15 &   36 \\ 
RCO1 &    7 &    7 &   15 &    8 \\ 
RCO2 &   14 &    7 &   15 &   17 \\ 
RCO3 &   21 &    7 &   15 &   25 \\ 
RCO4 &   28 &    7 &   15 &   29 \\ 
RCO5 &   35 &    7 &   15 &   41 \\ 
artif1\_2 &   13 &    8 &   27 &    8 \\ 
artif2\_4 &   26 &    8 &   27 &   14 \\ 
artif3\_6 &   39 &    8 &   27 &   20 \\ 
artif4\_8 &   52 &    8 &   27 &   28 \\ 
artif5\_10 &   65 &    8 &   27 &   34 \\ 
artif6\_12 &   78 &    8 &   27 &   41 \\ 
artif7\_14 &   91 &    8 &   27 &   48 \\ 
artif &   99 &    8 &   27 &   53 \\ 
shirts1\_2 &   13 &    8 &   40 &   13 \\ 
shirts2\_4 &   26 &    8 &   40 &   20 \\ 
shirts3\_6 &   39 &    8 &   40 &   26 \\ 
shirts4\_8 &   52 &    8 &   40 &   35 \\ 
shirts5\_10 &   65 &    8 &   40 &   42 \\ 
dagli1 &   10 &   10 &   60 &   25 \\ 
fu5 &    5 &    4 &   38 &   18 \\ 
fu6 &    6 &    5 &   38 &   24 \\ 
fu7 &    7 &    6 &   38 &   24 \\ 
fu8 &    8 &    7 &   38 &   24 \\ 
fu9 &    9 &    8 &   38 &   29 \\ 
fu10 &   10 &    9 &   38 &   34 \\ 
fu &   12 &   11 &   38 &   38 \\ 
J1-10-10-0 &   10 &    9 &   10 &   20 \\ 
J1-10-10-1 &   10 &   10 &   10 &   19 \\ 
J1-10-10-2 &   10 &    9 &   10 &   21 \\ 
J1-10-10-3 &   10 &   10 &   10 &   23 \\ 
J1-10-10-4 &   10 &   10 &   10 &   15 \\ 
J1-12-20-0 &   12 &   11 &   20 &   12 \\ 
J1-12-20-1 &   12 &   11 &   20 &   11 \\ 
J1-12-20-2 &   12 &   12 &   20 &   14 \\ 
J1-12-20-3 &   12 &   12 &   20 &   10 \\ 
J1-12-20-4 &   12 &   12 &   20 &   16 \\ 
J1-14-20-0 &   14 &   13 &   20 &   14 \\ 
J1-14-20-1 &   14 &   13 &   20 &   14 \\ 
J1-14-20-2 &   14 &   14 &   20 &   16 \\ 
J1-14-20-3 &   14 &   13 &   20 &   12 \\ 
J1-14-20-4 &   14 &   13 &   20 &   16 \\ 
J2-10-35-0 &   10 &    9 &   35 &   28 \\ 
J2-10-35-1 &   10 &    9 &   35 &   28 \\ 
J2-10-35-2 &   10 &   10 &   35 &   27 \\ 
J2-10-35-3 &   10 &    9 &   35 &   25 \\ 
J2-10-35-4 &   10 &   10 &   35 &   22 \\ 
J2-12-35-0 &   12 &   11 &   35 &   31 \\ 
J2-12-35-1 &   12 &   11 &   35 &   29 \\ 
J2-12-35-2 &   12 &   10 &   35 &   30 \\ 
J2-12-35-3 &   12 &   11 &   35 &   25 \\ 
J2-12-35-4 &   12 &   11 &   35 &   29 \\ 
J2-14-35-0 &   14 &   13 &   35 &   34 \\ 
J2-14-35-1 &   14 &   13 &   35 &   33 \\ 
J2-14-35-2 &   14 &   11 &   35 &   33 \\ 
J2-14-35-3 &   14 &   13 &   35 &   29 \\ 
J2-14-35-4 &   14 &   13 &   35 &   31 \\ 
Poly1a &   15 &   15 &   40 &   18 \\ 
Poly1b &   15 &   15 &   40 &   21 \\ 
Poly1c &   15 &   15 &   40 &   14 \\ 
Poly1d &   15 &   15 &   40 &   14 \\ 
Poly1e &   15 &   15 &   40 &   13 \\ 
Jakobs1 &   25 &   22 &   40 &   13 
\end{tabular}
\end{minipage}
\hfill
\begin{minipage}[t]{0.63\textwidth}
\caption{Comparison of the iterative and DB-PB-LB lower bound algorithms.\\Tightest or optimal-proven$^*$ LB values and solutions are shown in \textcolor{red}{red}.}
\label{tab:lower_bound_results}
\tiny
\begin{tabular}{lc|ccc|ccc}
& \textbf{Trivial} & \multicolumn{3}{c}{\textbf{Iterative lower bound} \cite{Rodrigues2017-zl}} & \multicolumn{3}{|c}{\textbf{DB-PB-LB (this work)}} \\
\textbf{Instance} & \textbf{LB} & \textbf{LB} & \textbf{Impr.(\%)} & \textbf{Time (secs)} & \textbf{LB} & \textbf{Impr. (\%)} & \textbf{Time (secs)} \\ 
\hline
three &   4 &   6 & 50.00 & 0.00 &   \textcolor{red}{\textbf{$6^*$}} & 50.00 & 0.01 \\ 
threep2 &   7 &  10 & 42.86 & 0.00 &  \textcolor{red}{\textbf{$10^*$}} & 42.86 & 0.00 \\ 
threep2w9 &   6 &   8 & 33.33 & 0.01 &   \textcolor{red}{\textbf{$8^*$}} & 33.33 & 0.00 \\ 
threep3 &  10 &  14 & 40.00 & 0.01 &  \textcolor{red}{\textbf{$14^*$}} & 40.00 & 0.05 \\ 
threep3w9 &   8 &  11 & 37.50 & 0.01 &  \textcolor{red}{\textbf{$12^*$}} & 50.00 & 0.11 \\ 
Shapes4 &  14 &  23 & 64.29 & 0.16 &  \textcolor{red}{\textbf{$24^*$}} & 71.43 & 0.00 \\ 
Shapes8 &  16 &  23 & 43.75 & 7.57 &  \textcolor{red}{\textbf{$26^*$}} & 62.50 & 2.24 \\ 
SHAPES-2 &  14 &  14 & 0.00 & 0.08 &  \textcolor{red}{\textbf{$14^*$}} & 0.00 & 0.00 \\ 
SHAPES-4 &  16 &  \textcolor{red}{21} & 31.25 & 127.73 &  19 & 18.75 & 600.00 \\ 
SHAPES-5 &  20 &  \textcolor{red}{25} & 25.00 & 28.49 &  22 & 10.00 & 600.00 \\ 
SHAPES-7 &  28 &  \textcolor{red}{34} & 21.43 & 108.49 &  28 & 0.00 & 600.00 \\ 
SHAPES-9 &  33 &  \textcolor{red}{38} & 15.15 & 133.25 &  33 & 0.00 & 600.00 \\ 
SHAPES-15 &  40 &  \textcolor{red}{47} & 17.50 & 281.59 &  40 & 0.00 & 600.00 \\ 
Blasz2 &  19 & \textcolor{red}{24} & 26.32 & 1.00 &  19 & 0.00 & 600.00 \\ 
BLAZEWICZ1 &   6 &   8 & 33.33 & 0.07 &   \textcolor{red}{\textbf{$8^*$}} & 33.33 & 0.01 \\ 
BLAZEWICZ2 &  11 &  13 & 18.18 & 1.46 &  \textcolor{red}{\textbf{$14^*$}} & 27.27 & 15.65 \\ 
BLAZEWICZ3 &  17 &  \textcolor{red}{19} & 11.76 & 5.69 &  17 & 0.00 & 600.00 \\ 
BLAZEWICZ4 &  22 &  \textcolor{red}{25} & 13.64 & 18.41 &  22 & 0.00 & 600.00 \\ 
BLAZEWICZ5 &  27 &  \textcolor{red}{31} & 14.81 & 48.46 &  27 & 0.00 & 600.00 \\ 
RCO1 &   7 &   8 & 14.29 & 0.02 &   \textcolor{red}{\textbf{$8^*$}} & 14.29 & 0.00 \\ 
RCO2 &  13 &  14 & 7.69 & 0.33 &  \textcolor{red}{\textbf{$15^*$}} & 15.38 & 35.35 \\ 
RCO3 &  19 &  \textcolor{red}{21} & 10.53 & 1.72 &  19 & 0.00 & 600.00 \\ 
RCO4 &  26 &  \textcolor{red}{27} & 3.85 & 2.26 &  26 & 0.00 & 600.00 \\ 
RC05 &  32 &  \textcolor{red}{34} & 6.25 & 0.97 &  32 & 0.00 & 600.00 \\ 
artif1\_2 &   7 &   7 & 0.00 & 0.02 &   \textcolor{red}{\textbf{$7^*$}} & 0.00 & 0.00 \\ 
artif2\_4 &  10 &  \textcolor{red}{12} & 20.00 & 9.57 &  \textcolor{red}{12} & 20.00 & 600.00 \\ 
artif3\_6 &  14 &  \textcolor{red}{17} & 21.43 & 6.35 &  14 & 0.00 & 600.00 \\ 
artif4\_8 &  19 &  \textcolor{red}{22} & 15.79 & 14.74 &  19 & 0.00 & 600.00 \\ 
artif5\_10 &  23 &  \textcolor{red}{27} & 17.39 & 25.49 &  23 & 0.00 & 600.00 \\ 
artif6\_12 &  28 &  \textcolor{red}{33} & 17.86 & 41.78 &  28 & 0.00 & 600.00 \\ 
artif7\_14 &  33 &  \textcolor{red}{38} & 15.15 & 65.05 &  33 & 0.00 & 600.00 \\ 
artif &  36 &  \textcolor{red}{41} & 13.89 & 96.36 &  36 & 0.00 & 600.00 \\ 
shirts1\_2 &  13 &  13 & 0.00 & 0.01 &  \textcolor{red}{\textbf{$13^*$}} & 0.00 & 0.00 \\ 
shirts2\_4 &  14 &  17 & 21.43 & 3.50 &  \textcolor{red}{\textbf{$17^*$}} & 21.43 & 5.01 \\ 
shirts3\_6 &  21 &  \textcolor{red}{24} & 14.29 & 14.47 &  \textcolor{red}{24} & 14.29 & 600.00 \\ 
shirts4\_8 &  28 &  \textcolor{red}{31} & 10.71 & 58.80 &  28 & 0.00 & 600.00 \\ 
shirts5\_10 &  35 &  \textcolor{red}{37} & 5.71 & 80.54 &  35 & 0.00 & 600.00 \\ 
Dagli1 &  23 &  23 & 0.00 & 32.07 &  \textcolor{red}{\textbf{$23^*$}} & 0.00 & 0.00 \\ 
fu5 &  14 &  18 & 28.57 & 0.04 &  \textcolor{red}{\textbf{$18^*$}} & 28.57 & 0.03 \\ 
fu6 &  17 &  23 & 35.29 & 3.80 &  \textcolor{red}{\textbf{$23^*$}} & 35.29 & 0.11 \\ 
fu7 &  19 &  24 & 26.32 & 4.90 &  \textcolor{red}{\textbf{$24^*$}} & 26.32 & 0.22 \\ 
fu8 &  20 &  24 & 20.00 & 15.90 &  \textcolor{red}{\textbf{$24^*$}} & 20.00 & 0.19 \\ 
fu9 &  23 &  24 & 4.35 & 10.16 &  \textcolor{red}{\textbf{$25^*$}} & 8.70 & 2.54 \\ 
fu10 &  26 &  28 & 7.69 & 21.58 &  \textcolor{red}{\textbf{$29^*$}} & 11.54 & 357.37 \\ 
fu &  29 &  \textcolor{red}{\textbf{31}} & 6.90 & 182.71 &  30 & 3.45 & 600.00 \\ 
J1-10-10-0 &  16 &  18 & 12.50 & 0.27 &  \textcolor{red}{\textbf{$18^*$}} & 12.50 & 0.00 \\ 
J1-10-10-1 &  16 &  16 & 0.00 & 0.14 &  \textcolor{red}{\textbf{$17^*$}} & 6.25 & 0.00 \\ 
J1-10-10-2 &  18 &  19 & 5.56 & 0.16 &  \textcolor{red}{\textbf{$20^*$}} & 11.11 & 0.01 \\ 
J1-10-10-3 &  18 &  20 & 11.11 & 0.82 &  \textcolor{red}{\textbf{$21^*$}} & 16.67 & 0.07 \\ 
J1-10-10-4 &  11 &  12 & 9.09 & 0.08 &  \textcolor{red}{\textbf{$13^*$}} & 18.18 & 0.01 \\ 
J1-12-20-0 &  10 &  11 & 10.00 & 0.53 &  \textcolor{red}{\textbf{$12^*$}} & 20.00 & 0.00 \\ 
J1-12-20-1 &   9 &  10 & 11.11 & 0.53 &  \textcolor{red}{\textbf{$10^*$}} & 11.11 & 0.01 \\ 
J1-12-20-2 &  11 &  12 & 9.09 & 0.72 &  \textcolor{red}{\textbf{$12^*$}} & 9.09 & 0.01 \\ 
J1-12-20-3 &   7 &   8 & 14.29 & 0.34 &   \textcolor{red}{\textbf{$8^*$}} & 14.29 & 0.01 \\ 
J1-12-20-4 &  12 &  12 & 0.00 & 0.47 &  \textcolor{red}{\textbf{$13^*$}} & 8.33 & 0.12 \\ 
J1-14-20-0 &  11 &  12 & 9.09 & 1.00 &  \textcolor{red}{\textbf{$12^*$}} & 9.09 & 0.03 \\ 
J1-14-20-1 &  11 &  11 & 0.00 & 0.21 &  \textcolor{red}{\textbf{$12^*$}} & 9.09 & 0.79 \\ 
J1-14-20-2 &  13 &  13 & 0.00 & 0.65 &  \textcolor{red}{\textbf{$14^*$}} & 7.69 & 2.17 \\ 
J1-14-20-3 &   9 &  10 & 11.11 & 0.48 &  \textcolor{red}{\textbf{$10^*$}} & 11.11 & 0.02 \\ 
J1-14-20-4 &  12 &  13 & 8.33 & 5.66 &  \textcolor{red}{\textbf{$14^*$}} & 16.67 & 0.26 \\ 
J2-10-35-0 &  19 &  22 & 15.79 & 102.08 &  \textcolor{red}{\textbf{$24^*$}} & 26.32 & 10.28 \\ 
J2-10-35-1 &  16 &  22 & 37.50 & 107.56 &  \textcolor{red}{\textbf{$24^*$}} & 50.00 & 19.90 \\ 
J2-10-35-2 &  17 &  21 & 23.53 & 76.72 &  \textcolor{red}{\textbf{$23^*$}} & 35.29 & 3.77 \\ 
J2-10-35-3 &  16 &  18 & 12.50 & 29.21 &  \textcolor{red}{\textbf{$20^*$}} & 25.00 & 2.49 \\ 
J2-10-35-4 &  14 &  17 & 21.43 & 34.77 &  \textcolor{red}{\textbf{$19^*$}} & 35.71 & 1.15 \\ 
J2-12-35-0 &  22 &  25 & 13.64 & 249.24 &  \textcolor{red}{\textbf{$29^*$}} & 31.82 & 542.15 \\ 
J2-12-35-1 &  18 &  23 & 27.78 & 179.60 &  \textcolor{red}{\textbf{$25^*$}} & 38.89 & 243.83 \\ 
J2-12-35-2 &  18 &  21 & 16.67 & 148.62 &  \textcolor{red}{\textbf{$24^*$}} & 33.33 & 245.41 \\ 
J2-12-35-3 &  16 &  19 & 18.75 & 75.31 &  \textcolor{red}{\textbf{$21^*$}} & 31.25 & 191.36 \\ 
J2-12-35-4 &  20 &  24 & 20.00 & 114.18 &  \textcolor{red}{\textbf{$26^*$}} & 30.00 & 351.03 \\ 
J2-14-35-0 &  24 &  26 & 8.33 & 254.99 &  \textcolor{red}{28} & 16.67 & 600.00 \\ 
J2-14-35-1 &  22 &  25 & 13.64 & 238.19 &  \textcolor{red}{26} & 18.18 & 600.00 \\ 
J2-14-35-2 &  21 &  23 & 9.52 & 145.02 &  \textcolor{red}{24} & 14.29 & 600.00 \\ 
J2-14-35-3 &  21 &  22 & 4.76 & 112.94 &  \textcolor{red}{24} & 14.29 & 600.00 \\ 
J2-14-35-4 &  21 &  24 & 14.29 & 191.63 &  \textcolor{red}{25} & 19.05 & 600.00 \\ 
Poly1a &  13 &  14 & 7.69 & 24.36 &  \textcolor{red}{15} & 15.38 & 600.00 \\ 
Poly1b &  13 &  16 & 23.08 & 71.23 &  \textcolor{red}{\textbf{$18^*$}} & 38.46 & 125.82 \\ 
Poly1c &  13 &  13 & 0.00 & 2.72 &  \textcolor{red}{\textbf{$13^*$}} & 0.00 & 0.00 \\ 
Poly1d &  11 &  11 & 0.00 & 3.93 &  \textcolor{red}{\textbf{$12^*$}} & 9.09 & 29.28 \\ 
Poly1e &  10 &  11 & 10.00 & 5.09 &  \textcolor{red}{\textbf{$12^*$}} & 20.00 & 452.47 \\ 
Jakobs1 &  10 &  \textcolor{red}{11} & 10.00 & 2.36 &  10 & 0.00 & 600.00 \\ 
\hline
\# Improvements &  &  & \textcolor{red}{\textbf{71}} &  &  & 57 &   \\ 
\multicolumn{2}{l|}{\# Optimal} & 0  &  &  & \textcolor{red}{\textbf{51}} \\
\multicolumn{2}{l|}{\# Best LB/ties}   & 24 &  &  &  \textcolor{red}{\textbf{59}}
\end{tabular}
\end{minipage}
\end{table}
\clearpage

\section{Discussion}
\label{sec_discussion}

\subsection{Comparison of lower bound algorithms}

\emph{Our new DB-PB-LB lower algorithm finds the optimal solutions in 51 of 81 problem instances, and obtains the best LB value for 59 of them}, as shown in Table \ref{tab:lower_bound_results}, whilst the iterative LB algorithm of \citet{Rodrigues2017-dy} only obtains the best LB values 24 times. However, the iterative LB algorithm improves the trivial LB (\ref{eq:trivial_lower_bound}) value for 71 instances, whilst DB-PB-LB can only do it for 57 problem instances, respectively.

\emph{Our new DB-PB-LB algorithm solves 51 problem instances up to optimality within 10 minutes, whilst the DB-CC model \cite{Rodrigues2017-zl} only solves 47 ones within 1 hour}, as drawn from the results shown in Tables \ref{tab:lower_bound_results} and \ref{tab:first_results}. Thus, our new lower bound algorithm outperforms both the iterative LB algorithm and the DB-CC model of \citet{Rodrigues2017-zl} in obtaining tighter LB values and solving problem instances.

\subsection{Comparison of the exact methods}

\emph{Our new exact DP-PB algorithm significantly outperforms the DB-CC and 0-1 DB-CC models in solving faster 53 and 55 problem instances, respectively}, as shown in Tables \ref{tab:first_results} and \ref{tab:second_results}. Among the instances solved up to optimality, the DB-CC \cite{Rodrigues2017-zl} model and our 0-1 DB-CC reformulation only outperform DB-PB in 6 and 1 problem instances, respectively. In addition, the performance profiles shown in Figures \ref{fig:performance_profile_solution_time} and \ref{fig:performance_profile_total_time} show that our DB-PB algorithm outperforms the two former ILP models by two orders of magnitude in most of the problem instances, even if the preprocessing is not considered. 

\emph{Our new exact DP-PB algorithm solves nine problem instances up to optimality more than the DB-CC model and fourteen instances more than our 0-1 DB-CC reformulation}, as shown in Tables \ref{tab:first_results} and \ref{tab:second_results}, respectively. However, our DB-PB algorithm cannot solve up to optimality four problem instances solved by the DB-CC \cite{Rodrigues2017-zl} model, among which we have \emph{Blasz2}, \emph{BLAZEWICZ3}, \emph{RCO3}, and \emph{shirt3\_6}. Likewise, the DB-CC \cite{Rodrigues2017-zl} model obtains tighter upper bounds than our DB-PB algorithm in most of the unsolved problem instances and two more feasible solutions, as shown in Table \ref{tab:first_results}. On the other hand, our 0-1 DB-CC reformulation can find feasible solutions for all problem instances, as shown in Table \ref{tab:second_results}. We attribute the outperformance of the DB-CC and 0-1 DB-CC models in obtaining tighter upper bounds than DB-PB among the unsolved problem instances to the set of feasibility pump and primal heuristics implemented by Gurobi and Google CPSat, as well as the lack of these type of primal heuristics in our new DB-PB algorithm, which we plan to tackle as forthcoming work.

\emph{Our new exact DB-PB algorithm finds the optimal solution for seventeen open instances, whilst DB-CC only solves four of them, as shown in Figure} \ref{fig:solved_open_instances}. However, the new open instances solved by the DB-CC model in comparison with the results reported in its introductory paper \cite[Table 4] {Rodrigues2017-zl} can only be attributed to the technological advances in Gurobi, as mentioned in Section \ref{sec:experimental_setup}.

\emph{The DB-CC model solves five problem instances more than our new 0-1 DB-CC reformulation, and DB-CC is slightly faster than the last one}, as shown in Tables \ref{tab:first_results} and \ref{tab:second_results}, and the performance profile shown in Figure \ref{fig:performance_profile_total_time}. However, we cannot discard the potential performance gain of Gurobi regarding Google CPSat in solving the two former models, respectively. Anyway, the main aim of evaluating our new 0-1 DB-CC reformulation was not to introduce a tighter formulation of the DB-CC model \cite{Rodrigues2017-zl} else to test our main hypothesis concerning the option of efficiently solving the DB model \cite{Toledo2013-oi} with an ad-hoc exact algorithm instead of using traditional Branch-and-Cut and Constraint Programming frameworks.

\emph{Our experimental results positively confirm our main hypothesis claiming that the preprocessing time required by the DB-CC is unnecessary and counterproductive}, and that \emph{ an ad-hoc exact algorithm might be a better option} to solve the discrete nesting problem in comparison with costly general-purpose Branch-and-Cut and CP-SAT methods.

%\subsection{Comparing the dimensions of the models}
\emph{The DB-CC model and our DB-PB algorithm share exactly the same number of binary of variables in all problem instances, but both significantly differ in the number of constraints, as shown in Table} \ref{tab:binary_variables}. As expected, the clique covering formulation of \citet{Rodrigues2017-zl} significantly reduces the number of constraints required to solve the original DB model \cite{Toledo2013-oi}. However, the results for our exact DB-PB algorithm shown in Tables \ref{tab:first_results} and \ref{tab:second_results}, and Figures \ref{fig:performance_profile_total_time} and \ref{fig:performance_profile_total_time}, positively confirm our main hypothesis claiming that the preprocessing time required by the DB-CC model \cite{Rodrigues2017-zl} is unnecessary and counterproductive.

\emph{Our 0-1 DB-CC model requires a lower number of constraints than DB-CC in all problem instances, excepting shirts1\_2, fu5, and RCO1, as detailed in Table \ref{tab:binary_variables}}. We attribute this reduction in the number of constraints to the new Constraints (\ref{ineq:BDB-CC_4}) that substitute the  vertex clique Constraints (\ref{ineq:DB-CC_innerfit}) of the DB-CC model. On the other hand, Table \ref{tab:binary_variables} shows that our new 0-1 DB-CC requires a few more binary variables than DB-CC to define the variables $z_m$ used to represent the set of feasible length upper bounds $\mathcal{L}^+$. 

%\subsection{Impact of technological advances and confirmation of previous results}
%\label{sec:technology}

Our replication of the DB-CC mode \cite{Rodrigues2017-zl}l allows us to confirm independently the results reported by \citet{Rodrigues2017-zl} and to draw some conclusions on the impact of the technological advances of MIP solvers. The DB-CC model implemented with Gurobi 11.0.3 herein solves two instances more within one hour than the same model in its introductory paper \cite{Rodrigues2017-zl} and also finds two feasible solutions more, as shown in Table \ref{tab:first_results}. Likewise, it also improves eight upper bounds among the unsolved problem instances, as shown in Figure \ref{fig:improved_UB_open_instances}. However, these improvements can only be attributed to the technological advances of Gurobi 11.0.3 compared to CPLEX 12.6 \cite{Rodrigues2017-zl}, as mentioned in Section \ref{sec:experimental_setup}. On the other hand, our DB--CC replication could not find a feasible solution for SHAPES-9 problem instance as done in \cite[Table 4]{Rodrigues2017-zl}, 

\subsection{The new state of the art}

\begin{figure}[t]
\centering
\begin{subfigure}[t]{0.48\textwidth}
\includegraphics[scale=0.45]{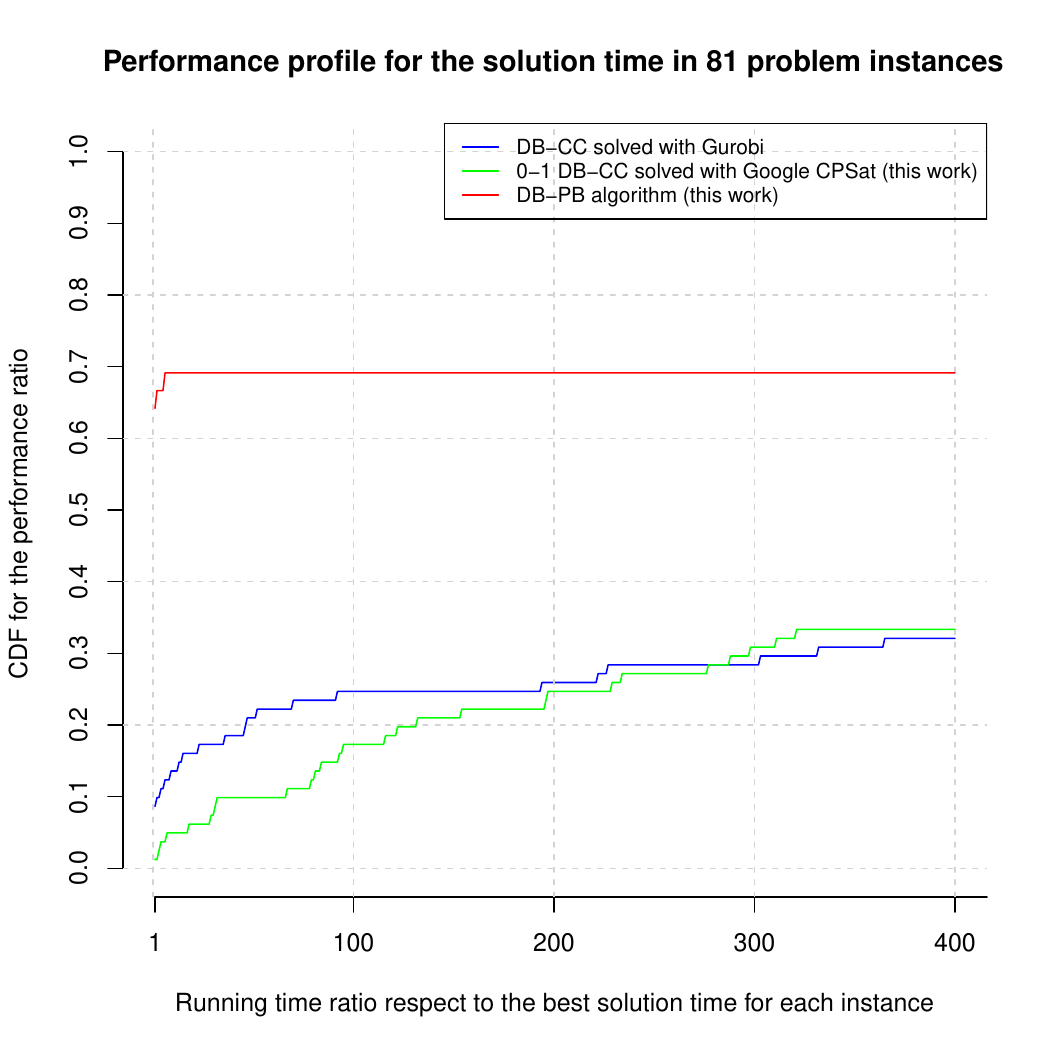}
\caption{Time for solving all problem instances.}
\label{fig:performance_profile_solution_time}
\end{subfigure}
\hfill
\begin{subfigure}[t]{0.48\textwidth}
\includegraphics[scale=0.45]{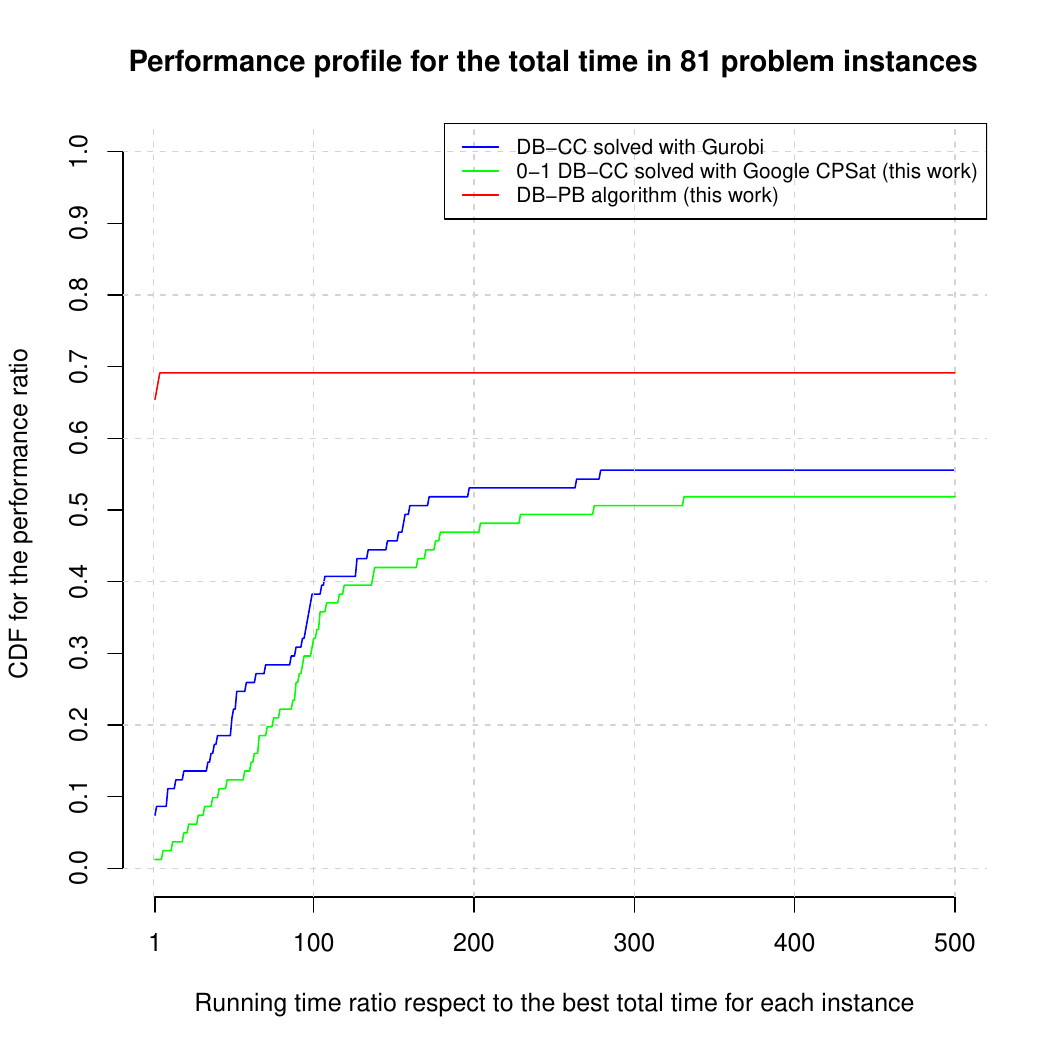}
\caption{Total time for building and solving all problem instances.}
\label{fig:performance_profile_total_time}
\end{subfigure}
\caption{Performance profiles \cite{Dolan2002-wd} showing the Cumulative Distribution Function (CDF) for the performance ratio $r_{\phi,m}$ comparing the performance of all exact DB methods in the evaluation of all problem instances.}
\end{figure}

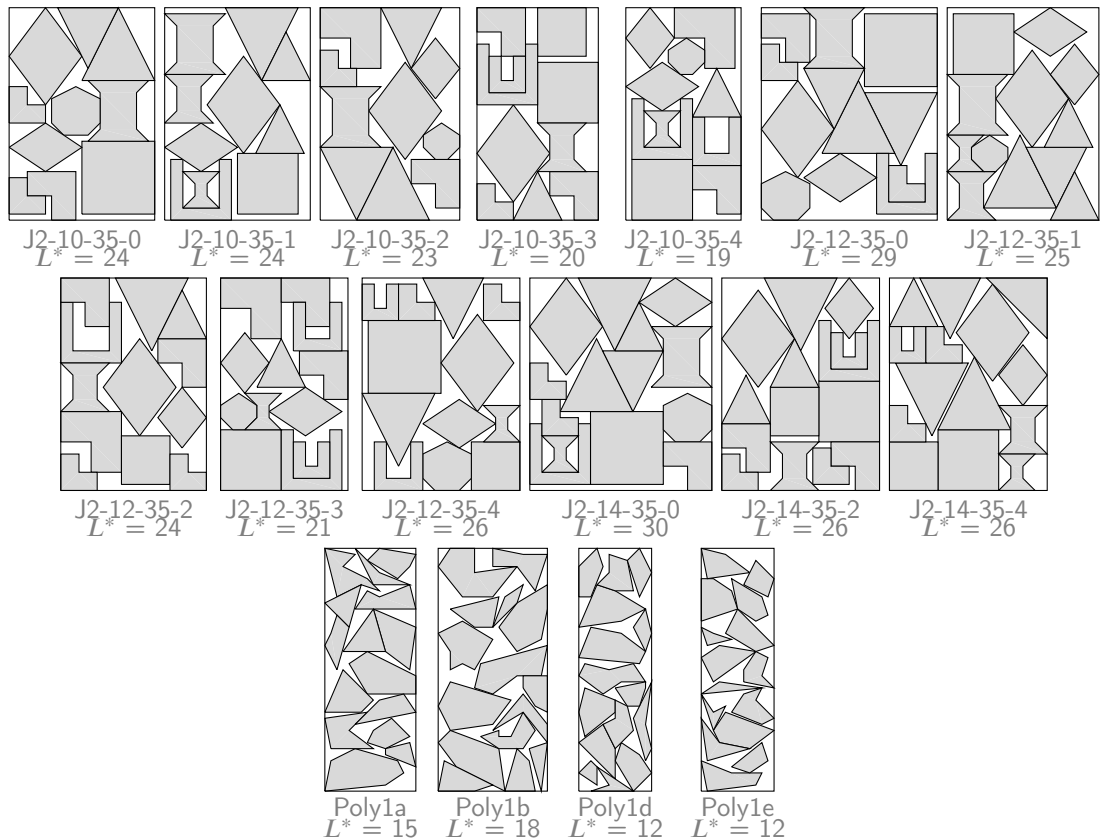
\begin{figure}[t]
\centering
% SOLUTION for problem instance = J2-10-35-0
% Scaling factor for coordinates = 1
\begin{tikzpicture}[scale = 0.08]
% Draw the convex parts of piece9
\fill[black!15](12,1) -- (24,1) -- (24,13) -- (12,13) -- (12,1);
% Draw the boundary of piece9
\draw[black](12,1) -- (24,1) -- (24,13) -- (12,13) -- (12,1);
% Draw the convex parts of piece4
\fill[black!15](6,19) -- (12,27) -- (6,35) -- (0,27) -- (6,19);
% Draw the boundary of piece4
\draw[black](6,19) -- (12,27) -- (6,35) -- (0,27) -- (6,19);
% Draw the convex parts of piece3
\fill[black!15](12,23) -- (24,23) -- (18,35) -- (12,23);
% Draw the boundary of piece3
\draw[black](12,23) -- (24,23) -- (18,35) -- (12,23);
% Draw the convex parts of piece2
\fill[black!15](15,15) -- (13,13) -- (23,13) -- (15,15);
\fill[black!15](15,15) -- (23,13) -- (21,15) -- (15,15);
\fill[black!15](15,21) -- (15,15) -- (21,15) -- (15,21);
\fill[black!15](15,21) -- (21,15) -- (21,21) -- (15,21);
\fill[black!15](21,21) -- (23,23) -- (13,23) -- (15,21) -- (21,21);
% Draw the boundary of piece2
\draw[black](13,13) -- (23,13) -- (21,15) -- (21,21) -- (23,23) -- (13,23) -- (15,21) -- (15,15) -- (13,13);
% Draw the convex parts of piece6
\fill[black!15](9,14) -- (13,14) -- (15,16) -- (15,20) -- (11,22) -- (7,20) -- (7,16) -- (9,14);
% Draw the boundary of piece6
\draw[black](9,14) -- (13,14) -- (15,16) -- (15,20) -- (11,22) -- (7,20) -- (7,16) -- (9,14);
% Draw the convex parts of piece5
\fill[black!15](13,25) -- (18,35) -- (8,35) -- (13,25);
% Draw the boundary of piece5
\draw[black](13,25) -- (18,35) -- (8,35) -- (13,25);
% Draw the convex parts of piece8
\fill[black!15](3,8) -- (3,4) -- (7,4) -- (3,8);
\fill[black!15](7,4) -- (7,0) -- (11,0) -- (7,4);
\fill[black!15](11,0) -- (11,8) -- (3,8) -- (7,4) -- (11,0);
% Draw the boundary of piece8
\draw[black](7,0) -- (11,0) -- (11,8) -- (3,8) -- (3,4) -- (7,4) -- (7,0);
% Draw the convex parts of piece1
\fill[black!15](0,12) -- (6,8) -- (12,12) -- (6,16) -- (0,12);
% Draw the boundary of piece1
\draw[black](0,12) -- (6,8) -- (12,12) -- (6,16) -- (0,12);
% Draw the convex parts of piece7
\fill[black!15](6,1) -- (6,4) -- (3,4) -- (6,1);
\fill[black!15](0,1) -- (6,1) -- (3,4) -- (0,1);
\fill[black!15](0,1) -- (3,4) -- (3,7) -- (0,7) -- (0,1);
% Draw the boundary of piece7
\draw[black](0,1) -- (6,1) -- (6,4) -- (3,4) -- (3,7) -- (0,7) -- (0,1);
% Draw the convex parts of piece7
\fill[black!15](6,16) -- (6,19) -- (3,19) -- (6,16);
\fill[black!15](0,16) -- (6,16) -- (3,19) -- (0,16);
\fill[black!15](0,16) -- (3,19) -- (3,22) -- (0,22) -- (0,16);
% Draw the boundary of piece7
\draw[black](0,16) -- (6,16) -- (6,19) -- (3,19) -- (3,22) -- (0,22) -- (0,16);
% Draw the bounding box of the solution
\draw[black](0,0) -- (24,0) -- (24,35)-- (0,35) -- (0,0);
\draw (12,-3) node {J2-10-35-0};
\draw (12,-6) node {$L^*$ = 24};
\end{tikzpicture}
% SOLUTION for problem instance = J2-10-35-1
% Scaling factor for coordinates = 1
\begin{tikzpicture}[scale = 0.08]
% Draw the convex parts of piece9
\fill[black!15](12,1) -- (22,1) -- (22,11) -- (12,11) -- (12,1);
% Draw the boundary of piece9
\draw[black](12,1) -- (22,1) -- (22,11) -- (12,11) -- (12,1);
% Draw the convex parts of piece7
\fill[black!15](13,11) -- (19,19) -- (13,27) -- (7,19) -- (13,11);
% Draw the boundary of piece7
\draw[black](13,11) -- (19,19) -- (13,27) -- (7,19) -- (13,11);
% Draw the convex parts of piece8
\fill[black!15](16,23) -- (22,35) -- (10,35) -- (16,23);
% Draw the boundary of piece8
\draw[black](16,23) -- (22,35) -- (10,35) -- (16,23);
% Draw the convex parts of piece4
\fill[black!15](2,26) -- (0,24) -- (10,24) -- (2,26);
\fill[black!15](2,26) -- (10,24) -- (8,26) -- (2,26);
\fill[black!15](2,32) -- (2,26) -- (8,26) -- (2,32);
\fill[black!15](2,32) -- (8,26) -- (8,32) -- (2,32);
\fill[black!15](8,32) -- (10,34) -- (0,34) -- (2,32) -- (8,32);
% Draw the boundary of piece4
\draw[black](0,24) -- (10,24) -- (8,26) -- (8,32) -- (10,34) -- (0,34) -- (2,32) -- (2,26) -- (0,24);
% Draw the convex parts of piece5
\fill[black!15](9,2) -- (3,2) -- (11,0) -- (9,2);
\fill[black!15](11,0) -- (11,10) -- (9,10) -- (9,2) -- (11,0);
\fill[black!15](1,0) -- (11,0) -- (3,2) -- (1,0);
\fill[black!15](1,0) -- (3,2) -- (3,10) -- (1,10) -- (1,0);
% Draw the boundary of piece5
\draw[black](1,0) -- (11,0) -- (11,10) -- (9,10) -- (9,2) -- (3,2) -- (3,10) -- (1,10) -- (1,0);
% Draw the convex parts of piece1
\fill[black!15](0,12) -- (6,8) -- (12,12) -- (6,16) -- (0,12);
% Draw the boundary of piece1
\draw[black](0,12) -- (6,8) -- (12,12) -- (6,16) -- (0,12);
% Draw the convex parts of piece3
\fill[black!15](2,18) -- (0,16) -- (8,16) -- (2,18);
\fill[black!15](2,18) -- (8,16) -- (6,18) -- (2,18);
\fill[black!15](2,22) -- (2,18) -- (6,18) -- (2,22);
\fill[black!15](2,22) -- (6,18) -- (6,22) -- (2,22);
\fill[black!15](6,22) -- (8,24) -- (0,24) -- (2,22) -- (6,22);
% Draw the boundary of piece3
\draw[black](0,16) -- (8,16) -- (6,18) -- (6,22) -- (8,24) -- (0,24) -- (2,22) -- (2,18) -- (0,16);
% Draw the convex parts of piece6
\fill[black!15](16,23) -- (24,23) -- (20,31) -- (16,23);
% Draw the boundary of piece6
\draw[black](16,23) -- (24,23) -- (20,31) -- (16,23);
% Draw the convex parts of piece6
\fill[black!15](15,11) -- (23,11) -- (19,19) -- (15,11);
% Draw the boundary of piece6
\draw[black](15,11) -- (23,11) -- (19,19) -- (15,11);
% Draw the convex parts of piece2
\fill[black!15](5,4) -- (3,2) -- (9,2) -- (5,4);
\fill[black!15](5,4) -- (9,2) -- (7,4) -- (5,4);
\fill[black!15](5,6) -- (5,4) -- (7,4) -- (5,6);
\fill[black!15](5,6) -- (7,4) -- (7,6) -- (5,6);
\fill[black!15](7,6) -- (9,8) -- (3,8) -- (5,6) -- (7,6);
% Draw the boundary of piece2
\draw[black](3,2) -- (9,2) -- (7,4) -- (7,6) -- (9,8) -- (3,8) -- (5,6) -- (5,4) -- (3,2);
% Draw the bounding box of the solution
\draw[black](0,0) -- (24,0) -- (24,35)-- (0,35) -- (0,0);
\draw (12,-3) node {J2-10-35-1};
\draw (12,-6) node {$L^*$ = 24};
\end{tikzpicture}
% SOLUTION for problem instance = J2-10-35-2
% Scaling factor for coordinates = 1
\begin{tikzpicture}[scale = 0.08]
% Draw the convex parts of piece4
\fill[black!15](14,10) -- (20,18) -- (14,26) -- (8,18) -- (14,10);
% Draw the boundary of piece4
\draw[black](14,10) -- (20,18) -- (14,26) -- (8,18) -- (14,10);
% Draw the convex parts of piece1
\fill[black!15](0,35) -- (0,30) -- (5,30) -- (0,35);
\fill[black!15](5,30) -- (5,25) -- (10,25) -- (5,30);
\fill[black!15](10,25) -- (10,35) -- (0,35) -- (5,30) -- (10,25);
% Draw the boundary of piece1
\draw[black](5,25) -- (10,25) -- (10,35) -- (0,35) -- (0,30) -- (5,30) -- (5,25);
% Draw the convex parts of piece5
\fill[black!15](6,0) -- (12,12) -- (0,12) -- (6,0);
% Draw the boundary of piece5
\draw[black](6,0) -- (12,12) -- (0,12) -- (6,0);
% Draw the convex parts of piece3
\fill[black!15](6,0) -- (18,0) -- (12,12) -- (6,0);
% Draw the boundary of piece3
\draw[black](6,0) -- (18,0) -- (12,12) -- (6,0);
% Draw the convex parts of piece2
\fill[black!15](2,14) -- (0,12) -- (10,12) -- (2,14);
\fill[black!15](2,14) -- (10,12) -- (8,14) -- (2,14);
\fill[black!15](2,20) -- (2,14) -- (8,14) -- (2,20);
\fill[black!15](2,20) -- (8,14) -- (8,20) -- (2,20);
\fill[black!15](8,20) -- (10,22) -- (0,22) -- (2,20) -- (8,20);
% Draw the boundary of piece2
\draw[black](0,12) -- (10,12) -- (8,14) -- (8,20) -- (10,22) -- (0,22) -- (2,20) -- (2,14) -- (0,12);
% Draw the convex parts of piece6
\fill[black!15](15,25) -- (20,35) -- (10,35) -- (15,25);
% Draw the boundary of piece6
\draw[black](15,25) -- (20,35) -- (10,35) -- (15,25);
% Draw the convex parts of piece9
\fill[black!15](15,10) -- (15,6) -- (19,6) -- (15,10);
\fill[black!15](19,6) -- (19,2) -- (23,2) -- (19,6);
\fill[black!15](23,2) -- (23,10) -- (15,10) -- (19,6) -- (23,2);
% Draw the boundary of piece9
\draw[black](19,2) -- (23,2) -- (23,10) -- (15,10) -- (15,6) -- (19,6) -- (19,2);
% Draw the convex parts of piece7
\fill[black!15](19,20) -- (23,25) -- (19,30) -- (15,25) -- (19,20);
% Draw the boundary of piece7
\draw[black](19,20) -- (23,25) -- (19,30) -- (15,25) -- (19,20);
% Draw the convex parts of piece10
\fill[black!15](6,22) -- (6,25) -- (3,25) -- (6,22);
\fill[black!15](0,22) -- (6,22) -- (3,25) -- (0,22);
\fill[black!15](0,22) -- (3,25) -- (3,28) -- (0,28) -- (0,22);
% Draw the boundary of piece10
\draw[black](0,22) -- (6,22) -- (6,25) -- (3,25) -- (3,28) -- (0,28) -- (0,22);
% Draw the convex parts of piece8
\fill[black!15](19,10) -- (21,10) -- (23,12) -- (23,14) -- (20,16) -- (17,14) -- (17,12) -- (19,10);
% Draw the boundary of piece8
\draw[black](19,10) -- (21,10) -- (23,12) -- (23,14) -- (20,16) -- (17,14) -- (17,12) -- (19,10);
% Draw the bounding box of the solution
\draw[black](0,0) -- (23,0) -- (23,35)-- (0,35) -- (0,0);
\draw (11.5,-3) node {J2-10-35-2};
\draw (11.5,-6) node {$L^*$ = 23};
\end{tikzpicture}
% SOLUTION for problem instance = J2-10-35-3
% Scaling factor for coordinates = 1
\begin{tikzpicture}[scale = 0.08]
% Draw the convex parts of piece9
\fill[black!15](10,16) -- (20,16) -- (20,26) -- (10,26) -- (10,16);
% Draw the boundary of piece9
\draw[black](10,16) -- (20,16) -- (20,26) -- (10,26) -- (10,16);
% Draw the convex parts of piece4
\fill[black!15](6,3) -- (12,11) -- (6,19) -- (0,11) -- (6,3);
% Draw the boundary of piece4
\draw[black](6,3) -- (12,11) -- (6,19) -- (0,11) -- (6,3);
% Draw the convex parts of piece8
\fill[black!15](10,27) -- (18,27) -- (18,35) -- (10,35) -- (10,27);
% Draw the boundary of piece8
\draw[black](10,27) -- (18,27) -- (18,35) -- (10,35) -- (10,27);
% Draw the convex parts of piece3
\fill[black!15](8,21) -- (2,21) -- (10,19) -- (8,21);
\fill[black!15](10,19) -- (10,29) -- (8,29) -- (8,21) -- (10,19);
\fill[black!15](0,19) -- (10,19) -- (2,21) -- (0,19);
\fill[black!15](0,19) -- (2,21) -- (2,29) -- (0,29) -- (0,19);
% Draw the boundary of piece3
\draw[black](0,19) -- (10,19) -- (10,29) -- (8,29) -- (8,21) -- (2,21) -- (2,29) -- (0,29) -- (0,19);
% Draw the convex parts of piece7
\fill[black!15](12,8) -- (12,4) -- (16,4) -- (12,8);
\fill[black!15](16,4) -- (16,0) -- (20,0) -- (16,4);
\fill[black!15](20,0) -- (20,8) -- (12,8) -- (16,4) -- (20,0);
% Draw the boundary of piece7
\draw[black](16,0) -- (20,0) -- (20,8) -- (12,8) -- (12,4) -- (16,4) -- (16,0);
% Draw the convex parts of piece7
\fill[black!15](0,35) -- (0,31) -- (4,31) -- (0,35);
\fill[black!15](4,31) -- (4,27) -- (8,27) -- (4,31);
\fill[black!15](8,27) -- (8,35) -- (0,35) -- (4,31) -- (8,27);
% Draw the boundary of piece7
\draw[black](4,27) -- (8,27) -- (8,35) -- (0,35) -- (0,31) -- (4,31) -- (4,27);
% Draw the convex parts of piece1
\fill[black!15](12,10) -- (10,8) -- (18,8) -- (12,10);
\fill[black!15](12,10) -- (18,8) -- (16,10) -- (12,10);
\fill[black!15](12,14) -- (12,10) -- (16,10) -- (12,14);
\fill[black!15](12,14) -- (16,10) -- (16,14) -- (12,14);
\fill[black!15](16,14) -- (18,16) -- (10,16) -- (12,14) -- (16,14);
% Draw the boundary of piece1
\draw[black](10,8) -- (18,8) -- (16,10) -- (16,14) -- (18,16) -- (10,16) -- (12,14) -- (12,10) -- (10,8);
% Draw the convex parts of piece5
\fill[black!15](6,0) -- (14,0) -- (10,8) -- (6,0);
% Draw the boundary of piece5
\draw[black](6,0) -- (14,0) -- (10,8) -- (6,0);
% Draw the convex parts of piece2
\fill[black!15](6,23) -- (4,23) -- (8,21) -- (6,23);
\fill[black!15](8,21) -- (8,27) -- (6,27) -- (6,23) -- (8,21);
\fill[black!15](2,21) -- (8,21) -- (4,23) -- (2,21);
\fill[black!15](2,21) -- (4,23) -- (4,27) -- (2,27) -- (2,21);
% Draw the boundary of piece2
\draw[black](2,21) -- (8,21) -- (8,27) -- (6,27) -- (6,23) -- (4,23) -- (4,27) -- (2,27) -- (2,21);
% Draw the convex parts of piece6
\fill[black!15](6,0) -- (6,3) -- (3,3) -- (6,0);
\fill[black!15](0,0) -- (6,0) -- (3,3) -- (0,0);
\fill[black!15](0,0) -- (3,3) -- (3,6) -- (0,6) -- (0,0);
% Draw the boundary of piece6
\draw[black](0,0) -- (6,0) -- (6,3) -- (3,3) -- (3,6) -- (0,6) -- (0,0);
% Draw the bounding box of the solution
\draw[black](0,0) -- (20,0) -- (20,35)-- (0,35) -- (0,0);
\draw (10,-3) node {J2-10-35-3};
\draw (10,-6) node {$L^*$ = 20};
\end{tikzpicture}
% SOLUTION for problem instance = J2-10-35-4
% Scaling factor for coordinates = 1
\begin{tikzpicture}[scale = 0.08]
% Draw the convex parts of piece10
\fill[black!15](1,0) -- (11,0) -- (11,10) -- (1,10) -- (1,0);
% Draw the boundary of piece10
\draw[black](1,0) -- (11,0) -- (11,10) -- (1,10) -- (1,0);
% Draw the convex parts of piece2
\fill[black!15](8,35) -- (8,30) -- (13,30) -- (8,35);
\fill[black!15](13,30) -- (13,25) -- (18,25) -- (13,30);
\fill[black!15](18,25) -- (18,35) -- (8,35) -- (13,30) -- (18,25);
% Draw the boundary of piece2
\draw[black](13,25) -- (18,25) -- (18,35) -- (8,35) -- (8,30) -- (13,30) -- (13,25);
% Draw the convex parts of piece5
\fill[black!15](9,12) -- (3,12) -- (11,10) -- (9,12);
\fill[black!15](11,10) -- (11,20) -- (9,20) -- (9,12) -- (11,10);
\fill[black!15](1,10) -- (11,10) -- (3,12) -- (1,10);
\fill[black!15](1,10) -- (3,12) -- (3,20) -- (1,20) -- (1,10);
% Draw the boundary of piece5
\draw[black](1,10) -- (11,10) -- (11,20) -- (9,20) -- (9,12) -- (3,12) -- (3,20) -- (1,20) -- (1,10);
% Draw the convex parts of piece9
\fill[black!15](11,9) -- (11,5) -- (15,5) -- (11,9);
\fill[black!15](15,5) -- (15,1) -- (19,1) -- (15,5);
\fill[black!15](19,1) -- (19,9) -- (11,9) -- (15,5) -- (19,1);
% Draw the boundary of piece9
\draw[black](15,1) -- (19,1) -- (19,9) -- (11,9) -- (11,5) -- (15,5) -- (15,1);
% Draw the convex parts of piece1
\fill[black!15](0,22) -- (6,18) -- (12,22) -- (6,26) -- (0,22);
% Draw the boundary of piece1
\draw[black](0,22) -- (6,18) -- (12,22) -- (6,26) -- (0,22);
% Draw the convex parts of piece7
\fill[black!15](4,25) -- (8,30) -- (4,35) -- (0,30) -- (4,25);
% Draw the boundary of piece7
\draw[black](4,25) -- (8,30) -- (4,35) -- (0,30) -- (4,25);
% Draw the convex parts of piece4
\fill[black!15](17,11) -- (13,11) -- (19,9) -- (17,11);
\fill[black!15](19,9) -- (19,17) -- (17,17) -- (17,11) -- (19,9);
\fill[black!15](11,9) -- (19,9) -- (13,11) -- (11,9);
\fill[black!15](11,9) -- (13,11) -- (13,17) -- (11,17) -- (11,9);
% Draw the boundary of piece4
\draw[black](11,9) -- (19,9) -- (19,17) -- (17,17) -- (17,11) -- (13,11) -- (13,17) -- (11,17) -- (11,9);
% Draw the convex parts of piece6
\fill[black!15](11,17) -- (19,17) -- (15,25) -- (11,17);
% Draw the boundary of piece6
\draw[black](11,17) -- (19,17) -- (15,25) -- (11,17);
% Draw the convex parts of piece8
\fill[black!15](9,24) -- (11,24) -- (13,26) -- (13,28) -- (10,30) -- (7,28) -- (7,26) -- (9,24);
% Draw the boundary of piece8
\draw[black](9,24) -- (11,24) -- (13,26) -- (13,28) -- (10,30) -- (7,28) -- (7,26) -- (9,24);
% Draw the convex parts of piece3
\fill[black!15](5,14) -- (3,12) -- (9,12) -- (5,14);
\fill[black!15](5,14) -- (9,12) -- (7,14) -- (5,14);
\fill[black!15](5,16) -- (5,14) -- (7,14) -- (5,16);
\fill[black!15](5,16) -- (7,14) -- (7,16) -- (5,16);
\fill[black!15](7,16) -- (9,18) -- (3,18) -- (5,16) -- (7,16);
% Draw the boundary of piece3
\draw[black](3,12) -- (9,12) -- (7,14) -- (7,16) -- (9,18) -- (3,18) -- (5,16) -- (5,14) -- (3,12);
% Draw the bounding box of the solution
\draw[black](0,0) -- (19,0) -- (19,35)-- (0,35) -- (0,0);
\draw (9.5,-3) node {J2-10-35-4};
\draw (9.5,-6) node {$L^*$ = 19};
\end{tikzpicture}
% SOLUTION for problem instance = J2-12-35-0
% Scaling factor for coordinates = 1
\begin{tikzpicture}[scale = 0.08]
% Draw the convex parts of piece9
\fill[black!15](17,22) -- (29,22) -- (29,34) -- (17,34) -- (17,22);
% Draw the boundary of piece9
\draw[black](17,22) -- (29,22) -- (29,34) -- (17,34) -- (17,22);
% Draw the convex parts of piece4
\fill[black!15](6,7) -- (12,15) -- (6,23) -- (0,15) -- (6,7);
% Draw the boundary of piece4
\draw[black](6,7) -- (12,15) -- (6,23) -- (0,15) -- (6,7);
% Draw the convex parts of piece11
\fill[black!15](23,9) -- (29,21) -- (17,21) -- (23,9);
% Draw the boundary of piece11
\draw[black](23,9) -- (29,21) -- (17,21) -- (23,9);
% Draw the convex parts of piece3
\fill[black!15](10,11) -- (22,11) -- (16,23) -- (10,11);
% Draw the boundary of piece3
\draw[black](10,11) -- (22,11) -- (16,23) -- (10,11);
% Draw the convex parts of piece2
\fill[black!15](9,27) -- (7,25) -- (17,25) -- (9,27);
\fill[black!15](9,27) -- (17,25) -- (15,27) -- (9,27);
\fill[black!15](9,33) -- (9,27) -- (15,27) -- (9,33);
\fill[black!15](9,33) -- (15,27) -- (15,33) -- (9,33);
\fill[black!15](15,33) -- (17,35) -- (7,35) -- (9,33) -- (15,33);
% Draw the boundary of piece2
\draw[black](7,25) -- (17,25) -- (15,27) -- (15,33) -- (17,35) -- (7,35) -- (9,33) -- (9,27) -- (7,25);
% Draw the convex parts of piece10
\fill[black!15](27,3) -- (21,3) -- (29,1) -- (27,3);
\fill[black!15](29,1) -- (29,11) -- (27,11) -- (27,3) -- (29,1);
\fill[black!15](19,1) -- (29,1) -- (21,3) -- (19,1);
\fill[black!15](19,1) -- (21,3) -- (21,11) -- (19,11) -- (19,1);
% Draw the boundary of piece10
\draw[black](19,1) -- (29,1) -- (29,11) -- (27,11) -- (27,3) -- (21,3) -- (21,11) -- (19,11) -- (19,1);
% Draw the convex parts of piece6
\fill[black!15](2,0) -- (6,0) -- (8,2) -- (8,6) -- (4,8) -- (0,6) -- (0,2) -- (2,0);
% Draw the boundary of piece6
\draw[black](2,0) -- (6,0) -- (8,2) -- (8,6) -- (4,8) -- (0,6) -- (0,2) -- (2,0);
% Draw the convex parts of piece5
\fill[black!15](12,15) -- (17,25) -- (7,25) -- (12,15);
% Draw the boundary of piece5
\draw[black](12,15) -- (17,25) -- (7,25) -- (12,15);
% Draw the convex parts of piece8
\fill[black!15](0,34) -- (0,30) -- (4,30) -- (0,34);
\fill[black!15](4,30) -- (4,26) -- (8,26) -- (4,30);
\fill[black!15](8,26) -- (8,34) -- (0,34) -- (4,30) -- (8,26);
% Draw the boundary of piece8
\draw[black](4,26) -- (8,26) -- (8,34) -- (0,34) -- (0,30) -- (4,30) -- (4,26);
% Draw the convex parts of piece1
\fill[black!15](7,7) -- (13,3) -- (19,7) -- (13,11) -- (7,7);
% Draw the boundary of piece1
\draw[black](7,7) -- (13,3) -- (19,7) -- (13,11) -- (7,7);
% Draw the convex parts of piece7
\fill[black!15](27,3) -- (27,6) -- (24,6) -- (27,3);
\fill[black!15](21,3) -- (27,3) -- (24,6) -- (21,3);
\fill[black!15](21,3) -- (24,6) -- (24,9) -- (21,9) -- (21,3);
% Draw the boundary of piece7
\draw[black](21,3) -- (27,3) -- (27,6) -- (24,6) -- (24,9) -- (21,9) -- (21,3);
% Draw the convex parts of piece7
\fill[black!15](6,23) -- (6,26) -- (3,26) -- (6,23);
\fill[black!15](0,23) -- (6,23) -- (3,26) -- (0,23);
\fill[black!15](0,23) -- (3,26) -- (3,29) -- (0,29) -- (0,23);
% Draw the boundary of piece7
\draw[black](0,23) -- (6,23) -- (6,26) -- (3,26) -- (3,29) -- (0,29) -- (0,23);
% Draw the bounding box of the solution
\draw[black](0,0) -- (29,0) -- (29,35)-- (0,35) -- (0,0);
\draw (14.5,-3) node {J2-12-35-0};
\draw (14.5,-6) node {$L^*$ = 29};
\end{tikzpicture}
% SOLUTION for problem instance = J2-12-35-1
% Scaling factor for coordinates = 1
\begin{tikzpicture}[scale = 0.08]
% Draw the convex parts of piece7
\fill[black!15](1,24) -- (11,24) -- (11,34) -- (1,34) -- (1,24);
% Draw the boundary of piece7
\draw[black](1,24) -- (11,24) -- (11,34) -- (1,34) -- (1,24);
% Draw the convex parts of piece6
\fill[black!15](14,12) -- (20,20) -- (14,28) -- (8,20) -- (14,12);
% Draw the boundary of piece6
\draw[black](14,12) -- (20,20) -- (14,28) -- (8,20) -- (14,12);
% Draw the convex parts of piece8
\fill[black!15](6,2) -- (18,2) -- (12,14) -- (6,2);
% Draw the boundary of piece8
\draw[black](6,2) -- (18,2) -- (12,14) -- (6,2);
% Draw the convex parts of piece4
\fill[black!15](2,16) -- (0,14) -- (10,14) -- (2,16);
\fill[black!15](2,16) -- (10,14) -- (8,16) -- (2,16);
\fill[black!15](2,22) -- (2,16) -- (8,16) -- (2,22);
\fill[black!15](2,22) -- (8,16) -- (8,22) -- (2,22);
\fill[black!15](8,22) -- (10,24) -- (0,24) -- (2,22) -- (8,22);
% Draw the boundary of piece4
\draw[black](0,14) -- (10,14) -- (8,16) -- (8,22) -- (10,24) -- (0,24) -- (2,22) -- (2,16) -- (0,14);
% Draw the convex parts of piece11
\fill[black!15](18,2) -- (23,12) -- (13,12) -- (18,2);
% Draw the boundary of piece11
\draw[black](18,2) -- (23,12) -- (13,12) -- (18,2);
% Draw the convex parts of piece1
\fill[black!15](11,31) -- (17,27) -- (23,31) -- (17,35) -- (11,31);
% Draw the boundary of piece1
\draw[black](11,31) -- (17,27) -- (23,31) -- (17,35) -- (11,31);
% Draw the convex parts of piece9
\fill[black!15](21,19) -- (25,24) -- (21,29) -- (17,24) -- (21,19);
% Draw the boundary of piece9
\draw[black](21,19) -- (25,24) -- (21,29) -- (17,24) -- (21,19);
% Draw the convex parts of piece3
\fill[black!15](2,2) -- (0,0) -- (8,0) -- (2,2);
\fill[black!15](2,2) -- (8,0) -- (6,2) -- (2,2);
\fill[black!15](2,6) -- (2,2) -- (6,2) -- (2,6);
\fill[black!15](2,6) -- (6,2) -- (6,6) -- (2,6);
\fill[black!15](6,6) -- (8,8) -- (0,8) -- (2,6) -- (6,6);
% Draw the boundary of piece3
\draw[black](0,0) -- (8,0) -- (6,2) -- (6,6) -- (8,8) -- (0,8) -- (2,6) -- (2,2) -- (0,0);
% Draw the convex parts of piece5
\fill[black!15](17,0) -- (25,0) -- (21,8) -- (17,0);
% Draw the boundary of piece5
\draw[black](17,0) -- (25,0) -- (21,8) -- (17,0);
% Draw the convex parts of piece5
\fill[black!15](16,12) -- (24,12) -- (20,20) -- (16,12);
% Draw the boundary of piece5
\draw[black](16,12) -- (24,12) -- (20,20) -- (16,12);
% Draw the convex parts of piece10
\fill[black!15](6,8) -- (8,8) -- (10,10) -- (10,12) -- (7,14) -- (4,12) -- (4,10) -- (6,8);
% Draw the boundary of piece10
\draw[black](6,8) -- (8,8) -- (10,10) -- (10,12) -- (7,14) -- (4,12) -- (4,10) -- (6,8);
% Draw the convex parts of piece2
\fill[black!15](2,10) -- (0,8) -- (6,8) -- (2,10);
\fill[black!15](2,10) -- (6,8) -- (4,10) -- (2,10);
\fill[black!15](2,12) -- (2,10) -- (4,10) -- (2,12);
\fill[black!15](2,12) -- (4,10) -- (4,12) -- (2,12);
\fill[black!15](4,12) -- (6,14) -- (0,14) -- (2,12) -- (4,12);
% Draw the boundary of piece2
\draw[black](0,8) -- (6,8) -- (4,10) -- (4,12) -- (6,14) -- (0,14) -- (2,12) -- (2,10) -- (0,8);
% Draw the bounding box of the solution
\draw[black](0,0) -- (25,0) -- (25,35)-- (0,35) -- (0,0);
\draw (12.5,-3) node {J2-12-35-1};
\draw (12.5,-6) node {$L^*$ = 25};
\end{tikzpicture}
\\
% SOLUTION for problem instance = J2-12-35-2
% Scaling factor for coordinates = 1
\begin{tikzpicture}[scale = 0.08]
% Draw the convex parts of piece9
\fill[black!15](13,9) -- (19,17) -- (13,25) -- (7,17) -- (13,9);
% Draw the boundary of piece9
\draw[black](13,9) -- (19,17) -- (13,25) -- (7,17) -- (13,9);
% Draw the convex parts of piece5
\fill[black!15](0,13) -- (0,8) -- (5,8) -- (0,13);
\fill[black!15](5,8) -- (5,3) -- (10,3) -- (5,8);
\fill[black!15](10,3) -- (10,13) -- (0,13) -- (5,8) -- (10,3);
% Draw the boundary of piece5
\draw[black](5,3) -- (10,3) -- (10,13) -- (0,13) -- (0,8) -- (5,8) -- (5,3);
% Draw the convex parts of piece3
\fill[black!15](15,23) -- (21,35) -- (9,35) -- (15,23);
% Draw the boundary of piece3
\draw[black](15,23) -- (21,35) -- (9,35) -- (15,23);
% Draw the convex parts of piece4
\fill[black!15](10,1) -- (18,1) -- (18,9) -- (10,9) -- (10,1);
% Draw the boundary of piece4
\draw[black](10,1) -- (18,1) -- (18,9) -- (10,9) -- (10,1);
% Draw the convex parts of piece10
\fill[black!15](8,23) -- (2,23) -- (10,21) -- (8,23);
\fill[black!15](10,21) -- (10,31) -- (8,31) -- (8,23) -- (10,21);
\fill[black!15](0,21) -- (10,21) -- (2,23) -- (0,21);
\fill[black!15](0,21) -- (2,23) -- (2,31) -- (0,31) -- (0,21);
% Draw the boundary of piece10
\draw[black](0,21) -- (10,21) -- (10,31) -- (8,31) -- (8,23) -- (2,23) -- (2,31) -- (0,31) -- (0,21);
% Draw the convex parts of piece2
\fill[black!15](16,25) -- (16,21) -- (20,21) -- (16,25);
\fill[black!15](20,21) -- (20,17) -- (24,17) -- (20,21);
\fill[black!15](24,17) -- (24,25) -- (16,25) -- (20,21) -- (24,17);
% Draw the boundary of piece2
\draw[black](20,17) -- (24,17) -- (24,25) -- (16,25) -- (16,21) -- (20,21) -- (20,17);
% Draw the convex parts of piece2
\fill[black!15](0,35) -- (0,31) -- (4,31) -- (0,35);
\fill[black!15](4,31) -- (4,27) -- (8,27) -- (4,31);
\fill[black!15](8,27) -- (8,35) -- (0,35) -- (4,31) -- (8,27);
% Draw the boundary of piece2
\draw[black](4,27) -- (8,27) -- (8,35) -- (0,35) -- (0,31) -- (4,31) -- (4,27);
% Draw the convex parts of piece8
\fill[black!15](20,7) -- (24,12) -- (20,17) -- (16,12) -- (20,7);
% Draw the boundary of piece8
\draw[black](20,7) -- (24,12) -- (20,17) -- (16,12) -- (20,7);
% Draw the convex parts of piece6
\fill[black!15](2,15) -- (0,13) -- (8,13) -- (2,15);
\fill[black!15](2,15) -- (8,13) -- (6,15) -- (2,15);
\fill[black!15](2,19) -- (2,15) -- (6,15) -- (2,19);
\fill[black!15](2,19) -- (6,15) -- (6,19) -- (2,19);
\fill[black!15](6,19) -- (8,21) -- (0,21) -- (2,19) -- (6,19);
% Draw the boundary of piece6
\draw[black](0,13) -- (8,13) -- (6,15) -- (6,19) -- (8,21) -- (0,21) -- (2,19) -- (2,15) -- (0,13);
% Draw the convex parts of piece7
\fill[black!15](16,25) -- (24,25) -- (20,33) -- (16,25);
% Draw the boundary of piece7
\draw[black](16,25) -- (24,25) -- (20,33) -- (16,25);
% Draw the convex parts of piece1
\fill[black!15](24,0) -- (24,3) -- (21,3) -- (24,0);
\fill[black!15](18,0) -- (24,0) -- (21,3) -- (18,0);
\fill[black!15](18,0) -- (21,3) -- (21,6) -- (18,6) -- (18,0);
% Draw the boundary of piece1
\draw[black](18,0) -- (24,0) -- (24,3) -- (21,3) -- (21,6) -- (18,6) -- (18,0);
% Draw the convex parts of piece1
\fill[black!15](6,0) -- (6,3) -- (3,3) -- (6,0);
\fill[black!15](0,0) -- (6,0) -- (3,3) -- (0,0);
\fill[black!15](0,0) -- (3,3) -- (3,6) -- (0,6) -- (0,0);
% Draw the boundary of piece1
\draw[black](0,0) -- (6,0) -- (6,3) -- (3,3) -- (3,6) -- (0,6) -- (0,0);
% Draw the bounding box of the solution
\draw[black](0,0) -- (24,0) -- (24,35)-- (0,35) -- (0,0);
\draw (12,-3) node {J2-12-35-2};
\draw (12,-6) node {$L^*$ = 24};
\end{tikzpicture}
% SOLUTION for problem instance = J2-12-35-3
% Scaling factor for coordinates = 1
\begin{tikzpicture}[scale = 0.08]
% Draw the convex parts of piece9
\fill[black!15](0,0) -- (10,0) -- (10,10) -- (0,10) -- (0,0);
% Draw the boundary of piece9
\draw[black](0,0) -- (10,0) -- (10,10) -- (0,10) -- (0,0);
% Draw the convex parts of piece7
\fill[black!15](0,35) -- (0,30) -- (5,30) -- (0,35);
\fill[black!15](5,30) -- (5,25) -- (10,25) -- (5,30);
\fill[black!15](10,25) -- (10,35) -- (0,35) -- (5,30) -- (10,25);
% Draw the boundary of piece7
\draw[black](5,25) -- (10,25) -- (10,35) -- (0,35) -- (0,30) -- (5,30) -- (5,25);
% Draw the convex parts of piece8
\fill[black!15](18,2) -- (12,2) -- (20,0) -- (18,2);
\fill[black!15](20,0) -- (20,10) -- (18,10) -- (18,2) -- (20,0);
\fill[black!15](10,0) -- (20,0) -- (12,2) -- (10,0);
\fill[black!15](10,0) -- (12,2) -- (12,10) -- (10,10) -- (10,0);
% Draw the boundary of piece8
\draw[black](10,0) -- (20,0) -- (20,10) -- (18,10) -- (18,2) -- (12,2) -- (12,10) -- (10,10) -- (10,0);
% Draw the convex parts of piece6
\fill[black!15](13,23) -- (13,19) -- (17,19) -- (13,23);
\fill[black!15](17,19) -- (17,15) -- (21,15) -- (17,19);
\fill[black!15](21,15) -- (21,23) -- (13,23) -- (17,19) -- (21,15);
% Draw the boundary of piece6
\draw[black](17,15) -- (21,15) -- (21,23) -- (13,23) -- (13,19) -- (17,19) -- (17,15);
% Draw the convex parts of piece6
\fill[black!15](10,35) -- (10,31) -- (14,31) -- (10,35);
\fill[black!15](14,31) -- (14,27) -- (18,27) -- (14,31);
\fill[black!15](18,27) -- (18,35) -- (10,35) -- (14,31) -- (18,27);
% Draw the boundary of piece6
\draw[black](14,27) -- (18,27) -- (18,35) -- (10,35) -- (10,31) -- (14,31) -- (14,27);
% Draw the convex parts of piece2
\fill[black!15](8,13) -- (14,9) -- (20,13) -- (14,17) -- (8,13);
% Draw the boundary of piece2
\draw[black](8,13) -- (14,9) -- (20,13) -- (14,17) -- (8,13);
% Draw the convex parts of piece11
\fill[black!15](18,25) -- (14,25) -- (20,23) -- (18,25);
\fill[black!15](20,23) -- (20,31) -- (18,31) -- (18,25) -- (20,23);
\fill[black!15](12,23) -- (20,23) -- (14,25) -- (12,23);
\fill[black!15](12,23) -- (14,25) -- (14,31) -- (12,31) -- (12,23);
% Draw the boundary of piece11
\draw[black](12,23) -- (20,23) -- (20,31) -- (18,31) -- (18,25) -- (14,25) -- (14,31) -- (12,31) -- (12,23);
% Draw the convex parts of piece5
\fill[black!15](4,16) -- (8,21) -- (4,26) -- (0,21) -- (4,16);
% Draw the boundary of piece5
\draw[black](4,16) -- (8,21) -- (4,26) -- (0,21) -- (4,16);
% Draw the convex parts of piece3
\fill[black!15](6,17) -- (14,17) -- (10,25) -- (6,17);
% Draw the boundary of piece3
\draw[black](6,17) -- (14,17) -- (10,25) -- (6,17);
% Draw the convex parts of piece10
\fill[black!15](16,4) -- (14,4) -- (18,2) -- (16,4);
\fill[black!15](18,2) -- (18,8) -- (16,8) -- (16,4) -- (18,2);
\fill[black!15](12,2) -- (18,2) -- (14,4) -- (12,2);
\fill[black!15](12,2) -- (14,4) -- (14,8) -- (12,8) -- (12,2);
% Draw the boundary of piece10
\draw[black](12,2) -- (18,2) -- (18,8) -- (16,8) -- (16,4) -- (14,4) -- (14,8) -- (12,8) -- (12,2);
% Draw the convex parts of piece1
\fill[black!15](2,10) -- (4,10) -- (6,12) -- (6,14) -- (3,16) -- (0,14) -- (0,12) -- (2,10);
% Draw the boundary of piece1
\draw[black](2,10) -- (4,10) -- (6,12) -- (6,14) -- (3,16) -- (0,14) -- (0,12) -- (2,10);
% Draw the convex parts of piece4
\fill[black!15](6,12) -- (4,10) -- (10,10) -- (6,12);
\fill[black!15](6,12) -- (10,10) -- (8,12) -- (6,12);
\fill[black!15](6,14) -- (6,12) -- (8,12) -- (6,14);
\fill[black!15](6,14) -- (8,12) -- (8,14) -- (6,14);
\fill[black!15](8,14) -- (10,16) -- (4,16) -- (6,14) -- (8,14);
% Draw the boundary of piece4
\draw[black](4,10) -- (10,10) -- (8,12) -- (8,14) -- (10,16) -- (4,16) -- (6,14) -- (6,12) -- (4,10);
% Draw the bounding box of the solution
\draw[black](0,0) -- (21,0) -- (21,35)-- (0,35) -- (0,0);
\draw (10.5,-3) node {J2-12-35-3};
\draw (10.5,-6) node {$L^*$ = 21};
\end{tikzpicture}
% SOLUTION for problem instance = J2-12-35-4
% Scaling factor for coordinates = 1
\begin{tikzpicture}[scale = 0.08]
% Draw the convex parts of piece7
\fill[black!15](1,16) -- (13,16) -- (13,28) -- (1,28) -- (1,16);
% Draw the boundary of piece7
\draw[black](1,16) -- (13,16) -- (13,28) -- (1,28) -- (1,16);
% Draw the convex parts of piece11
\fill[black!15](19,13) -- (25,21) -- (19,29) -- (13,21) -- (19,13);
% Draw the boundary of piece11
\draw[black](19,13) -- (25,21) -- (19,29) -- (13,21) -- (19,13);
% Draw the convex parts of piece4
\fill[black!15](6,4) -- (12,16) -- (0,16) -- (6,4);
% Draw the boundary of piece4
\draw[black](6,4) -- (12,16) -- (0,16) -- (6,4);
% Draw the convex parts of piece10
\fill[black!15](18,0) -- (26,0) -- (26,8) -- (18,8) -- (18,0);
% Draw the boundary of piece10
\draw[black](18,0) -- (26,0) -- (26,8) -- (18,8) -- (18,0);
% Draw the convex parts of piece6
\fill[black!15](12,0) -- (16,0) -- (18,2) -- (18,6) -- (14,8) -- (10,6) -- (10,2) -- (12,0);
% Draw the boundary of piece6
\draw[black](12,0) -- (16,0) -- (18,2) -- (18,6) -- (14,8) -- (10,6) -- (10,2) -- (12,0);
% Draw the convex parts of piece2
\fill[black!15](15,25) -- (20,35) -- (10,35) -- (15,25);
% Draw the boundary of piece2
\draw[black](15,25) -- (20,35) -- (10,35) -- (15,25);
% Draw the convex parts of piece1
\fill[black!15](10,11) -- (16,7) -- (22,11) -- (16,15) -- (10,11);
% Draw the boundary of piece1
\draw[black](10,11) -- (16,7) -- (22,11) -- (16,15) -- (10,11);
% Draw the convex parts of piece5
\fill[black!15](8,2) -- (4,2) -- (10,0) -- (8,2);
\fill[black!15](10,0) -- (10,8) -- (8,8) -- (8,2) -- (10,0);
\fill[black!15](2,0) -- (10,0) -- (4,2) -- (2,0);
\fill[black!15](2,0) -- (4,2) -- (4,8) -- (2,8) -- (2,0);
% Draw the boundary of piece5
\draw[black](2,0) -- (10,0) -- (10,8) -- (8,8) -- (8,2) -- (4,2) -- (4,8) -- (2,8) -- (2,0);
% Draw the convex parts of piece9
\fill[black!15](4,30) -- (2,30) -- (6,28) -- (4,30);
\fill[black!15](6,28) -- (6,34) -- (4,34) -- (4,30) -- (6,28);
\fill[black!15](0,28) -- (6,28) -- (2,30) -- (0,28);
\fill[black!15](0,28) -- (2,30) -- (2,34) -- (0,34) -- (0,28);
% Draw the boundary of piece9
\draw[black](0,28) -- (6,28) -- (6,34) -- (4,34) -- (4,30) -- (2,30) -- (2,34) -- (0,34) -- (0,28);
% Draw the convex parts of piece8
\fill[black!15](26,28) -- (26,31) -- (23,31) -- (26,28);
\fill[black!15](20,28) -- (26,28) -- (23,31) -- (20,28);
\fill[black!15](20,28) -- (23,31) -- (23,34) -- (20,34) -- (20,28);
% Draw the boundary of piece8
\draw[black](20,28) -- (26,28) -- (26,31) -- (23,31) -- (23,34) -- (20,34) -- (20,28);
% Draw the convex parts of piece8
\fill[black!15](12,28) -- (12,31) -- (9,31) -- (12,28);
\fill[black!15](6,28) -- (12,28) -- (9,31) -- (6,28);
\fill[black!15](6,28) -- (9,31) -- (9,34) -- (6,34) -- (6,28);
% Draw the boundary of piece8
\draw[black](6,28) -- (12,28) -- (12,31) -- (9,31) -- (9,34) -- (6,34) -- (6,28);
% Draw the convex parts of piece3
\fill[black!15](22,10) -- (20,8) -- (26,8) -- (22,10);
\fill[black!15](22,10) -- (26,8) -- (24,10) -- (22,10);
\fill[black!15](22,12) -- (22,10) -- (24,10) -- (22,12);
\fill[black!15](22,12) -- (24,10) -- (24,12) -- (22,12);
\fill[black!15](24,12) -- (26,14) -- (20,14) -- (22,12) -- (24,12);
% Draw the boundary of piece3
\draw[black](20,8) -- (26,8) -- (24,10) -- (24,12) -- (26,14) -- (20,14) -- (22,12) -- (22,10) -- (20,8);
% Draw the bounding box of the solution
\draw[black](0,0) -- (26,0) -- (26,35)-- (0,35) -- (0,0);
\draw (13,-3) node {J2-12-35-4};
\draw (13,-6) node {$L^*$ = 26};
\end{tikzpicture}
% SOLUTION for problem instance = J2-14-35-0
% Scaling factor for coordinates = 1
\begin{tikzpicture}[scale = 0.08]
% Draw the convex parts of piece7
\fill[black!15](10,1) -- (22,1) -- (22,13) -- (10,13) -- (10,1);
% Draw the boundary of piece7
\draw[black](10,1) -- (22,1) -- (22,13) -- (10,13) -- (10,1);
% Draw the convex parts of piece4
\fill[black!15](6,19) -- (12,27) -- (6,35) -- (0,27) -- (6,19);
% Draw the boundary of piece4
\draw[black](6,19) -- (12,27) -- (6,35) -- (0,27) -- (6,19);
% Draw the convex parts of piece9
\fill[black!15](14,23) -- (20,35) -- (8,35) -- (14,23);
% Draw the boundary of piece9
\draw[black](14,23) -- (20,35) -- (8,35) -- (14,23);
% Draw the convex parts of piece1
\fill[black!15](5,13) -- (17,13) -- (11,25) -- (5,13);
% Draw the boundary of piece1
\draw[black](5,13) -- (17,13) -- (11,25) -- (5,13);
% Draw the convex parts of piece2
\fill[black!15](22,19) -- (20,17) -- (30,17) -- (22,19);
\fill[black!15](22,19) -- (30,17) -- (28,19) -- (22,19);
\fill[black!15](22,25) -- (22,19) -- (28,19) -- (22,25);
\fill[black!15](22,25) -- (28,19) -- (28,25) -- (22,25);
\fill[black!15](28,25) -- (30,27) -- (20,27) -- (22,25) -- (28,25);
% Draw the boundary of piece2
\draw[black](20,17) -- (30,17) -- (28,19) -- (28,25) -- (30,27) -- (20,27) -- (22,25) -- (22,19) -- (20,17);
% Draw the convex parts of piece8
\fill[black!15](8,3) -- (2,3) -- (10,1) -- (8,3);
\fill[black!15](10,1) -- (10,11) -- (8,11) -- (8,3) -- (10,1);
\fill[black!15](0,1) -- (10,1) -- (2,3) -- (0,1);
\fill[black!15](0,1) -- (2,3) -- (2,11) -- (0,11) -- (0,1);
% Draw the boundary of piece8
\draw[black](0,1) -- (10,1) -- (10,11) -- (8,11) -- (8,3) -- (2,3) -- (2,11) -- (0,11) -- (0,1);
% Draw the convex parts of piece5
\fill[black!15](24,8) -- (28,8) -- (30,10) -- (30,14) -- (26,16) -- (22,14) -- (22,10) -- (24,8);
% Draw the boundary of piece5
\draw[black](24,8) -- (28,8) -- (30,10) -- (30,14) -- (26,16) -- (22,14) -- (22,10) -- (24,8);
% Draw the convex parts of piece13
\fill[black!15](17,13) -- (22,23) -- (12,23) -- (17,13);
% Draw the boundary of piece13
\draw[black](17,13) -- (22,23) -- (12,23) -- (17,13);
% Draw the convex parts of piece10
\fill[black!15](22,8) -- (22,4) -- (26,4) -- (22,8);
\fill[black!15](26,4) -- (26,0) -- (30,0) -- (26,4);
\fill[black!15](30,0) -- (30,8) -- (22,8) -- (26,4) -- (30,0);
% Draw the boundary of piece10
\draw[black](26,0) -- (30,0) -- (30,8) -- (22,8) -- (22,4) -- (26,4) -- (26,0);
% Draw the convex parts of piece3
\fill[black!15](18,31) -- (24,27) -- (30,31) -- (24,35) -- (18,31);
% Draw the boundary of piece3
\draw[black](18,31) -- (24,27) -- (30,31) -- (24,35) -- (18,31);
% Draw the convex parts of piece12
\fill[black!15](14,23) -- (22,23) -- (18,31) -- (14,23);
% Draw the boundary of piece12
\draw[black](14,23) -- (22,23) -- (18,31) -- (14,23);
% Draw the convex parts of piece6
\fill[black!15](8,9) -- (8,12) -- (5,12) -- (8,9);
\fill[black!15](2,9) -- (8,9) -- (5,12) -- (2,9);
\fill[black!15](2,9) -- (5,12) -- (5,15) -- (2,15) -- (2,9);
% Draw the boundary of piece6
\draw[black](2,9) -- (8,9) -- (8,12) -- (5,12) -- (5,15) -- (2,15) -- (2,9);
% Draw the convex parts of piece6
\fill[black!15](6,15) -- (6,18) -- (3,18) -- (6,15);
\fill[black!15](0,15) -- (6,15) -- (3,18) -- (0,15);
\fill[black!15](0,15) -- (3,18) -- (3,21) -- (0,21) -- (0,15);
% Draw the boundary of piece6
\draw[black](0,15) -- (6,15) -- (6,18) -- (3,18) -- (3,21) -- (0,21) -- (0,15);
% Draw the convex parts of piece11
\fill[black!15](4,5) -- (2,3) -- (8,3) -- (4,5);
\fill[black!15](4,5) -- (8,3) -- (6,5) -- (4,5);
\fill[black!15](4,7) -- (4,5) -- (6,5) -- (4,7);
\fill[black!15](4,7) -- (6,5) -- (6,7) -- (4,7);
\fill[black!15](6,7) -- (8,9) -- (2,9) -- (4,7) -- (6,7);
% Draw the boundary of piece11
\draw[black](2,3) -- (8,3) -- (6,5) -- (6,7) -- (8,9) -- (2,9) -- (4,7) -- (4,5) -- (2,3);
% Draw the bounding box of the solution
\draw[black](0,0) -- (30,0) -- (30,35)-- (0,35) -- (0,0);
\draw (15,-3) node {J2-14-35-0};
\draw (15,-6) node {$L^*$ = 30};
\end{tikzpicture}
% SOLUTION for problem instance = J2-14-35-2
% Scaling factor for coordinates = 1
\begin{tikzpicture}[scale = 0.08]
% Draw the convex parts of piece10
\fill[black!15](16,8) -- (26,8) -- (26,18) -- (16,18) -- (16,8);
% Draw the boundary of piece10
\draw[black](16,8) -- (26,8) -- (26,18) -- (16,18) -- (16,8);
% Draw the convex parts of piece4
\fill[black!15](6,17) -- (12,25) -- (6,33) -- (0,25) -- (6,17);
% Draw the boundary of piece4
\draw[black](6,17) -- (12,25) -- (6,33) -- (0,25) -- (6,17);
% Draw the convex parts of piece8
\fill[black!15](13,23) -- (19,35) -- (7,35) -- (13,23);
% Draw the boundary of piece8
\draw[black](13,23) -- (19,35) -- (7,35) -- (13,23);
% Draw the convex parts of piece9
\fill[black!15](8,9) -- (16,9) -- (16,17) -- (8,17) -- (8,9);
% Draw the boundary of piece9
\draw[black](8,9) -- (16,9) -- (16,17) -- (8,17) -- (8,9);
% Draw the convex parts of piece2
\fill[black!15](24,20) -- (18,20) -- (26,18) -- (24,20);
\fill[black!15](26,18) -- (26,28) -- (24,28) -- (24,20) -- (26,18);
\fill[black!15](16,18) -- (26,18) -- (18,20) -- (16,18);
\fill[black!15](16,18) -- (18,20) -- (18,28) -- (16,28) -- (16,18);
% Draw the boundary of piece2
\draw[black](16,18) -- (26,18) -- (26,28) -- (24,28) -- (24,20) -- (18,20) -- (18,28) -- (16,28) -- (16,18);
% Draw the convex parts of piece7
\fill[black!15](18,8) -- (18,4) -- (22,4) -- (18,8);
\fill[black!15](22,4) -- (22,0) -- (26,0) -- (22,4);
\fill[black!15](26,0) -- (26,8) -- (18,8) -- (22,4) -- (26,0);
% Draw the boundary of piece7
\draw[black](22,0) -- (26,0) -- (26,8) -- (18,8) -- (18,4) -- (22,4) -- (22,0);
% Draw the convex parts of piece7
\fill[black!15](0,11) -- (0,7) -- (4,7) -- (0,11);
\fill[black!15](4,7) -- (4,3) -- (8,3) -- (4,7);
\fill[black!15](8,3) -- (8,11) -- (0,11) -- (4,7) -- (8,3);
% Draw the boundary of piece7
\draw[black](4,3) -- (8,3) -- (8,11) -- (0,11) -- (0,7) -- (4,7) -- (4,3);
% Draw the convex parts of piece6
\fill[black!15](10,2) -- (8,0) -- (16,0) -- (10,2);
\fill[black!15](10,2) -- (16,0) -- (14,2) -- (10,2);
\fill[black!15](10,6) -- (10,2) -- (14,2) -- (10,6);
\fill[black!15](10,6) -- (14,2) -- (14,6) -- (10,6);
\fill[black!15](14,6) -- (16,8) -- (8,8) -- (10,6) -- (14,6);
% Draw the boundary of piece6
\draw[black](8,0) -- (16,0) -- (14,2) -- (14,6) -- (16,8) -- (8,8) -- (10,6) -- (10,2) -- (8,0);
% Draw the convex parts of piece5
\fill[black!15](21,25) -- (25,30) -- (21,35) -- (17,30) -- (21,25);
% Draw the boundary of piece5
\draw[black](21,25) -- (25,30) -- (21,35) -- (17,30) -- (21,25);
% Draw the convex parts of piece3
\fill[black!15](8,17) -- (16,17) -- (12,25) -- (8,17);
% Draw the boundary of piece3
\draw[black](8,17) -- (16,17) -- (12,25) -- (8,17);
% Draw the convex parts of piece3
\fill[black!15](0,11) -- (8,11) -- (4,19) -- (0,11);
% Draw the boundary of piece3
\draw[black](0,11) -- (8,11) -- (4,19) -- (0,11);
% Draw the convex parts of piece11
\fill[black!15](22,22) -- (20,22) -- (24,20) -- (22,22);
\fill[black!15](24,20) -- (24,26) -- (22,26) -- (22,22) -- (24,20);
\fill[black!15](18,20) -- (24,20) -- (20,22) -- (18,20);
\fill[black!15](18,20) -- (20,22) -- (20,26) -- (18,26) -- (18,20);
% Draw the boundary of piece11
\draw[black](18,20) -- (24,20) -- (24,26) -- (22,26) -- (22,22) -- (20,22) -- (20,26) -- (18,26) -- (18,20);
% Draw the convex parts of piece1
\fill[black!15](21,1) -- (21,4) -- (18,4) -- (21,1);
\fill[black!15](15,1) -- (21,1) -- (18,4) -- (15,1);
\fill[black!15](15,1) -- (18,4) -- (18,7) -- (15,7) -- (15,1);
% Draw the boundary of piece1
\draw[black](15,1) -- (21,1) -- (21,4) -- (18,4) -- (18,7) -- (15,7) -- (15,1);
% Draw the convex parts of piece1
\fill[black!15](6,0) -- (6,3) -- (3,3) -- (6,0);
\fill[black!15](0,0) -- (6,0) -- (3,3) -- (0,0);
\fill[black!15](0,0) -- (3,3) -- (3,6) -- (0,6) -- (0,0);
% Draw the boundary of piece1
\draw[black](0,0) -- (6,0) -- (6,3) -- (3,3) -- (3,6) -- (0,6) -- (0,0);
% Draw the bounding box of the solution
\draw[black](0,0) -- (26,0) -- (26,35)-- (0,35) -- (0,0);
\draw (13,-3) node {J2-14-35-2};
\draw (13,-6) node {$L^*$ = 26};
\end{tikzpicture}
% SOLUTION for problem instance = J2-14-35-4
% Scaling factor for coordinates = 1
\begin{tikzpicture}[scale = 0.08]
% Draw the convex parts of piece5
\fill[black!15](8,0) -- (18,0) -- (18,10) -- (8,10) -- (8,0);
% Draw the boundary of piece5
\draw[black](8,0) -- (18,0) -- (18,10) -- (8,10) -- (8,0);
% Draw the convex parts of piece10
\fill[black!15](17,18) -- (23,26) -- (17,34) -- (11,26) -- (17,18);
% Draw the boundary of piece10
\draw[black](17,18) -- (23,26) -- (17,34) -- (11,26) -- (17,18);
% Draw the convex parts of piece4
\fill[black!15](8,10) -- (20,10) -- (14,22) -- (8,10);
% Draw the boundary of piece4
\draw[black](8,10) -- (20,10) -- (14,22) -- (8,10);
% Draw the convex parts of piece3
\fill[black!15](7,9) -- (13,21) -- (1,21) -- (7,9);
% Draw the boundary of piece3
\draw[black](7,9) -- (13,21) -- (1,21) -- (7,9);
% Draw the convex parts of piece2
\fill[black!15](0,0) -- (8,0) -- (8,8) -- (0,8) -- (0,0);
% Draw the boundary of piece2
\draw[black](0,0) -- (8,0) -- (8,8) -- (0,8) -- (0,0);
% Draw the convex parts of piece13
\fill[black!15](26,25) -- (26,35) -- (16,35) -- (26,25);
% Draw the boundary of piece13
\draw[black](26,25) -- (26,35) -- (16,35) -- (26,25);
% Draw the convex parts of piece12
\fill[black!15](10,25) -- (15,35) -- (5,35) -- (10,25);
% Draw the boundary of piece12
\draw[black](10,25) -- (15,35) -- (5,35) -- (10,25);
% Draw the convex parts of piece9
\fill[black!15](20,8) -- (18,6) -- (26,6) -- (20,8);
\fill[black!15](20,8) -- (26,6) -- (24,8) -- (20,8);
\fill[black!15](20,12) -- (20,8) -- (24,8) -- (20,12);
\fill[black!15](20,12) -- (24,8) -- (24,12) -- (20,12);
\fill[black!15](24,12) -- (26,14) -- (18,14) -- (20,12) -- (24,12);
% Draw the boundary of piece9
\draw[black](18,6) -- (26,6) -- (24,8) -- (24,12) -- (26,14) -- (18,14) -- (20,12) -- (20,8) -- (18,6);
% Draw the convex parts of piece7
\fill[black!15](22,14) -- (26,19) -- (22,24) -- (18,19) -- (22,14);
% Draw the boundary of piece7
\draw[black](22,14) -- (26,19) -- (22,24) -- (18,19) -- (22,14);
% Draw the convex parts of piece6
\fill[black!15](0,27) -- (8,27) -- (4,35) -- (0,27);
% Draw the boundary of piece6
\draw[black](0,27) -- (8,27) -- (4,35) -- (0,27);
% Draw the convex parts of piece11
\fill[black!15](4,23) -- (2,23) -- (6,21) -- (4,23);
\fill[black!15](6,21) -- (6,27) -- (4,27) -- (4,23) -- (6,21);
\fill[black!15](0,21) -- (6,21) -- (2,23) -- (0,21);
\fill[black!15](0,21) -- (2,23) -- (2,27) -- (0,27) -- (0,21);
% Draw the boundary of piece11
\draw[black](0,21) -- (6,21) -- (6,27) -- (4,27) -- (4,23) -- (2,23) -- (2,27) -- (0,27) -- (0,21);
% Draw the convex parts of piece1
\fill[black!15](12,21) -- (12,24) -- (9,24) -- (12,21);
\fill[black!15](6,21) -- (12,21) -- (9,24) -- (6,21);
\fill[black!15](6,21) -- (9,24) -- (9,27) -- (6,27) -- (6,21);
% Draw the boundary of piece1
\draw[black](6,21) -- (12,21) -- (12,24) -- (9,24) -- (9,27) -- (6,27) -- (6,21);
% Draw the convex parts of piece1
\fill[black!15](6,8) -- (6,11) -- (3,11) -- (6,8);
\fill[black!15](0,8) -- (6,8) -- (3,11) -- (0,8);
\fill[black!15](0,8) -- (3,11) -- (3,14) -- (0,14) -- (0,8);
% Draw the boundary of piece1
\draw[black](0,8) -- (6,8) -- (6,11) -- (3,11) -- (3,14) -- (0,14) -- (0,8);
% Draw the convex parts of piece8
\fill[black!15](20,2) -- (18,0) -- (24,0) -- (20,2);
\fill[black!15](20,2) -- (24,0) -- (22,2) -- (20,2);
\fill[black!15](20,4) -- (20,2) -- (22,2) -- (20,4);
\fill[black!15](20,4) -- (22,2) -- (22,4) -- (20,4);
\fill[black!15](22,4) -- (24,6) -- (18,6) -- (20,4) -- (22,4);
% Draw the boundary of piece8
\draw[black](18,0) -- (24,0) -- (22,2) -- (22,4) -- (24,6) -- (18,6) -- (20,4) -- (20,2) -- (18,0);
% Draw the bounding box of the solution
\draw[black](0,0) -- (26,0) -- (26,35)-- (0,35) -- (0,0);
\draw (13,-3) node {J2-14-35-4};
\draw (13,-6) node {$L^*$ = 26};
\end{tikzpicture}
\\
% SOLUTION for problem instance = Poly1a
% Scaling factor for coordinates = 1
\begin{tikzpicture}[scale = 0.08]
% Draw the convex parts of piece4
\fill[black!15](0,0) -- (10,1) -- (13,3) -- (12,5) -- (5,7) -- (1,5) -- (0,0);
% Draw the boundary of piece4
\draw[black](0,0) -- (10,1) -- (13,3) -- (12,5) -- (5,7) -- (1,5) -- (0,0);
% Draw the convex parts of piece9
\fill[black!15](7,16) -- (15,16) -- (15,19) -- (7,23) -- (4,21) -- (7,16);
% Draw the boundary of piece9
\draw[black](7,16) -- (15,16) -- (15,19) -- (7,23) -- (4,21) -- (7,16);
% Draw the convex parts of piece3
\fill[black!15](10,23) -- (14,20) -- (15,27) -- (8,30) -- (10,23);
% Draw the boundary of piece3
\draw[black](10,23) -- (14,20) -- (15,27) -- (8,30) -- (10,23);
% Draw the convex parts of piece8
\fill[black!15](0,13) -- (6,13) -- (8,15) -- (3,20) -- (0,13);
% Draw the boundary of piece8
\draw[black](0,13) -- (6,13) -- (8,15) -- (3,20) -- (0,13);
% Draw the convex parts of piece11
\fill[black!15](6,34) -- (9,33) -- (7,35) -- (6,34);
\fill[black!15](2,32) -- (6,34) -- (7,35) -- (3,37) -- (2,32);
\fill[black!15](0,40) -- (2,32) -- (3,37) -- (0,40);
\fill[black!15](3,37) -- (5,39) -- (0,40) -- (3,37);
% Draw the boundary of piece11
\draw[black](2,32) -- (6,34) -- (9,33) -- (7,35) -- (3,37) -- (5,39) -- (0,40) -- (2,32);
% Draw the convex parts of piece10
\fill[black!15](1,7) -- (8,7) -- (4,10) -- (1,7);
\fill[black!15](1,7) -- (4,10) -- (7,13) -- (0,12) -- (1,7);
% Draw the boundary of piece10
\draw[black](1,7) -- (8,7) -- (4,10) -- (7,13) -- (0,12) -- (1,7);
% Draw the convex parts of piece13
\fill[black!15](15,30) -- (14,34) -- (9,32) -- (15,30);
\fill[black!15](5,30) -- (15,30) -- (9,32) -- (5,30);
\fill[black!15](5,30) -- (9,32) -- (7,33) -- (5,30);
% Draw the boundary of piece13
\draw[black](5,30) -- (15,30) -- (14,34) -- (9,32) -- (7,33) -- (5,30);
% Draw the convex parts of piece1
\fill[black!15](3,23) -- (10,23) -- (8,30) -- (3,23);
% Draw the boundary of piece1
\draw[black](3,23) -- (10,23) -- (8,30) -- (3,23);
% Draw the convex parts of piece0
\fill[black!15](5,11) -- (15,13) -- (13,16) -- (9,15) -- (5,11);
% Draw the boundary of piece0
\draw[black](5,11) -- (15,13) -- (13,16) -- (9,15) -- (5,11);
% Draw the convex parts of piece14
\fill[black!15](7,38) -- (12,35) -- (14,37) -- (15,39) -- (9,40) -- (7,38);
% Draw the boundary of piece14
\draw[black](7,38) -- (12,35) -- (14,37) -- (15,39) -- (9,40) -- (7,38);
% Draw the convex parts of piece2
\fill[black!15](7,8) -- (10,7) -- (14,9) -- (11,12) -- (8,10) -- (7,8);
% Draw the boundary of piece2
\draw[black](7,8) -- (10,7) -- (14,9) -- (11,12) -- (8,10) -- (7,8);
% Draw the convex parts of piece6
\fill[black!15](4,27) -- (6,34) -- (0,31) -- (4,27);
% Draw the boundary of piece6
\draw[black](4,27) -- (6,34) -- (0,31) -- (4,27);
% Draw the convex parts of piece7
\fill[black!15](1,26) -- (0,24) -- (2,20) -- (1,26);
\fill[black!15](2,20) -- (4,27) -- (2,29) -- (1,26) -- (2,20);
% Draw the boundary of piece7
\draw[black](2,20) -- (4,27) -- (2,29) -- (1,26) -- (0,24) -- (2,20);
% Draw the convex parts of piece12
\fill[black!15](8,34) -- (14,34) -- (9,36) -- (8,34);
\fill[black!15](5,40) -- (8,34) -- (9,36) -- (5,40);
% Draw the boundary of piece12
\draw[black](5,40) -- (8,34) -- (14,34) -- (9,36) -- (5,40);
% Draw the convex parts of piece5
\fill[black!15](10,6) -- (15,4) -- (14,8) -- (10,7) -- (10,6);
% Draw the boundary of piece5
\draw[black](10,6) -- (15,4) -- (14,8) -- (10,7) -- (10,6);
% Draw the bounding box of the solution
\draw[black](0,0) -- (15,0) -- (15,40)-- (0,40) -- (0,0);
\draw (7.5,-3) node {Poly1a};
\draw (7.5,-6) node {$L^*$ = 15};
\end{tikzpicture}
% SOLUTION for problem instance = Poly1b
% Scaling factor for coordinates = 1
\begin{tikzpicture}[scale = 0.08]
% Draw the convex parts of piece7
\fill[black!15](0,14) -- (2,10) -- (9,11) -- (12,15) -- (2,18) -- (0,14);
% Draw the boundary of piece7
\draw[black](0,14) -- (2,10) -- (9,11) -- (12,15) -- (2,18) -- (0,14);
% Draw the convex parts of piece5
\fill[black!15](9,32) -- (13,37) -- (6,36) -- (9,32);
\fill[black!15](4,32) -- (9,32) -- (6,36) -- (4,32);
\fill[black!15](0,37) -- (4,32) -- (6,36) -- (6,40) -- (2,40) -- (0,37);
% Draw the boundary of piece5
\draw[black](0,37) -- (4,32) -- (9,32) -- (13,37) -- (6,36) -- (6,40) -- (2,40) -- (0,37);
% Draw the convex parts of piece8
\fill[black!15](0,8) -- (3,4) -- (11,0) -- (10,7) -- (2,10) -- (0,8);
% Draw the boundary of piece8
\draw[black](0,8) -- (3,4) -- (11,0) -- (10,7) -- (2,10) -- (0,8);
% Draw the convex parts of piece1
\fill[black!15](10,27) -- (13,23) -- (17,25) -- (18,30) -- (18,34) -- (13,31) -- (10,27);
% Draw the boundary of piece1
\draw[black](10,27) -- (13,23) -- (17,25) -- (18,30) -- (18,34) -- (13,31) -- (10,27);
% Draw the convex parts of piece4
\fill[black!15](6,17) -- (9,16) -- (13,19) -- (6,17);
\fill[black!15](18,24) -- (7,21) -- (6,17) -- (13,19) -- (18,24);
\fill[black!15](13,19) -- (17,19) -- (18,24) -- (13,19);
% Draw the boundary of piece4
\draw[black](6,17) -- (9,16) -- (13,19) -- (17,19) -- (18,24) -- (7,21) -- (6,17);
% Draw the convex parts of piece15
\fill[black!15](2,23) -- (2,20) -- (4,21) -- (2,23);
\fill[black!15](0,26) -- (2,23) -- (4,21) -- (6,20) -- (9,25) -- (4,28) -- (0,26);
% Draw the boundary of piece15
\draw[black](0,26) -- (2,23) -- (2,20) -- (4,21) -- (6,20) -- (9,25) -- (4,28) -- (0,26);
% Draw the convex parts of piece2
\fill[black!15](13,37) -- (16,33) -- (18,39) -- (14,39) -- (13,37);
\fill[black!15](9,37) -- (13,37) -- (14,39) -- (9,40) -- (9,37);
% Draw the boundary of piece2
\draw[black](9,37) -- (13,37) -- (16,33) -- (18,39) -- (14,39) -- (9,40) -- (9,37);
% Draw the convex parts of piece3
\fill[black!15](0,0) -- (6,0) -- (9,1) -- (1,5) -- (0,0);
% Draw the boundary of piece3
\draw[black](0,0) -- (6,0) -- (9,1) -- (1,5) -- (0,0);
% Draw the convex parts of piece6
\fill[black!15](2,30) -- (6,27) -- (9,27) -- (10,32) -- (2,30);
% Draw the boundary of piece6
\draw[black](2,30) -- (6,27) -- (9,27) -- (10,32) -- (2,30);
% Draw the convex parts of piece14
\fill[black!15](12,9) -- (10,9) -- (7,9) -- (10,7) -- (14,7) -- (12,9);
\fill[black!15](14,7) -- (16,10) -- (14,13) -- (12,9) -- (14,7);
\fill[black!15](7,9) -- (10,9) -- (10,10) -- (7,9);
% Draw the boundary of piece14
\draw[black](7,9) -- (10,7) -- (14,7) -- (16,10) -- (14,13) -- (12,9) -- (10,9) -- (10,10) -- (7,9);
% Draw the convex parts of piece12
\fill[black!15](11,0) -- (16,0) -- (15,5) -- (13,7) -- (11,0);
% Draw the boundary of piece12
\draw[black](11,0) -- (16,0) -- (15,5) -- (13,7) -- (11,0);
% Draw the convex parts of piece13
\fill[black!15](14,13) -- (18,8) -- (18,11) -- (14,16) -- (14,13);
\fill[black!15](9,10) -- (14,13) -- (14,16) -- (9,10);
% Draw the boundary of piece13
\draw[black](9,10) -- (14,13) -- (18,8) -- (18,11) -- (14,16) -- (9,10);
% Draw the convex parts of piece9
\fill[black!15](14,17) -- (16,14) -- (18,13) -- (18,19) -- (14,19) -- (14,17);
% Draw the boundary of piece9
\draw[black](14,17) -- (16,14) -- (18,13) -- (18,19) -- (14,19) -- (14,17);
% Draw the convex parts of piece10
\fill[black!15](15,7) -- (17,0) -- (18,8) -- (15,7);
% Draw the boundary of piece10
\draw[black](15,7) -- (17,0) -- (18,8) -- (15,7);
% Draw the convex parts of piece11
\fill[black!15](9,27) -- (13,31) -- (14,34) -- (11,32) -- (9,27);
% Draw the boundary of piece11
\draw[black](9,27) -- (13,31) -- (14,34) -- (11,32) -- (9,27);
% Draw the bounding box of the solution
\draw[black](0,0) -- (18,0) -- (18,40)-- (0,40) -- (0,0);
\draw (9,-3) node {Poly1b};
\draw (9,-6) node {$L^*$ = 18};
\end{tikzpicture}
% SOLUTION for problem instance = Poly1d
% Scaling factor for coordinates = 1
\begin{tikzpicture}[scale = 0.08]
% Draw the convex parts of piece2
\fill[black!15](0,22) -- (10,21) -- (11,23) -- (7,26) -- (2,27) -- (0,22);
% Draw the boundary of piece2
\draw[black](0,22) -- (10,21) -- (11,23) -- (7,26) -- (2,27) -- (0,22);
% Draw the convex parts of piece11
\fill[black!15](0,27) -- (7,28) -- (11,30) -- (6,33) -- (3,34) -- (0,27);
% Draw the boundary of piece11
\draw[black](0,27) -- (7,28) -- (11,30) -- (6,33) -- (3,34) -- (0,27);
% Draw the convex parts of piece13
\fill[black!15](5,20) -- (0,19) -- (4,16) -- (7,18) -- (5,20);
\fill[black!15](7,18) -- (11,18) -- (9,21) -- (5,20) -- (7,18);
\fill[black!15](0,19) -- (5,20) -- (1,21) -- (0,19);
% Draw the boundary of piece13
\draw[black](0,19) -- (4,16) -- (7,18) -- (11,18) -- (9,21) -- (5,20) -- (1,21) -- (0,19);
% Draw the convex parts of piece7
\fill[black!15](0,13) -- (1,9) -- (6,13) -- (6,15) -- (1,17) -- (0,13);
% Draw the boundary of piece7
\draw[black](0,13) -- (1,9) -- (6,13) -- (6,15) -- (1,17) -- (0,13);
% Draw the convex parts of piece1
\fill[black!15](0,8) -- (6,2) -- (6,7) -- (4,11) -- (0,8);
% Draw the boundary of piece1
\draw[black](0,8) -- (6,2) -- (6,7) -- (4,11) -- (0,8);
% Draw the convex parts of piece14
\fill[black!15](6,5) -- (9,0) -- (12,3) -- (10,6) -- (7,8) -- (6,5);
% Draw the boundary of piece14
\draw[black](6,5) -- (9,0) -- (12,3) -- (10,6) -- (7,8) -- (6,5);
% Draw the convex parts of piece10
\fill[black!15](5,36) -- (3,34) -- (6,33) -- (9,33) -- (5,36);
\fill[black!15](9,33) -- (8,38) -- (5,39) -- (5,36) -- (9,33);
% Draw the boundary of piece10
\draw[black](3,34) -- (6,33) -- (9,33) -- (8,38) -- (5,39) -- (5,36) -- (3,34);
% Draw the convex parts of piece6
\fill[black!15](0,32) -- (2,32) -- (4,37) -- (3,40) -- (0,37) -- (0,32);
% Draw the boundary of piece6
\draw[black](0,32) -- (2,32) -- (4,37) -- (3,40) -- (0,37) -- (0,32);
% Draw the convex parts of piece8
\fill[black!15](6,12) -- (4,11) -- (6,7) -- (9,10) -- (6,12);
\fill[black!15](9,10) -- (8,14) -- (6,15) -- (6,12) -- (9,10);
% Draw the boundary of piece8
\draw[black](4,11) -- (6,7) -- (9,10) -- (8,14) -- (6,15) -- (6,12) -- (4,11);
% Draw the convex parts of piece5
\fill[black!15](11,5) -- (12,7) -- (11,11) -- (11,5);
\fill[black!15](9,12) -- (11,5) -- (11,11) -- (9,12);
\fill[black!15](11,11) -- (12,18) -- (9,12) -- (11,11);
% Draw the boundary of piece5
\draw[black](11,5) -- (12,7) -- (11,11) -- (12,18) -- (9,12) -- (11,5);
% Draw the convex parts of piece9
\fill[black!15](8,38) -- (10,32) -- (12,36) -- (10,40) -- (8,38);
% Draw the boundary of piece9
\draw[black](8,38) -- (10,32) -- (12,36) -- (10,40) -- (8,38);
% Draw the convex parts of piece3
\fill[black!15](6,17) -- (8,14) -- (10,14) -- (11,17) -- (11,18) -- (6,17);
% Draw the boundary of piece3
\draw[black](6,17) -- (8,14) -- (10,14) -- (11,17) -- (11,18) -- (6,17);
% Draw the convex parts of piece4
\fill[black!15](11,23) -- (12,26) -- (11,30) -- (9,27) -- (11,23);
% Draw the boundary of piece4
\draw[black](11,23) -- (12,26) -- (11,30) -- (9,27) -- (11,23);
% Draw the convex parts of piece12
\fill[black!15](5,1) -- (2,0) -- (9,0) -- (5,1);
\fill[black!15](9,0) -- (3,4) -- (5,1) -- (9,0);
% Draw the boundary of piece12
\draw[black](2,0) -- (9,0) -- (3,4) -- (5,1) -- (2,0);
% Draw the convex parts of piece15
\fill[black!15](0,3) -- (3,1) -- (3,4) -- (2,5) -- (0,3);
% Draw the boundary of piece15
\draw[black](0,3) -- (3,1) -- (3,4) -- (2,5) -- (0,3);
% Draw the bounding box of the solution
\draw[black](0,0) -- (12,0) -- (12,40)-- (0,40) -- (0,0);
\draw (6,-3) node {Poly1d};
\draw (6,-6) node {$L^*$ = 12};
\end{tikzpicture}
% SOLUTION for problem instance = Poly1e
% Scaling factor for coordinates = 1
\begin{tikzpicture}[scale = 0.08]
% Draw the convex parts of piece14
\fill[black!15](2,10) -- (3,8) -- (8,7) -- (12,10) -- (4,13) -- (2,10);
% Draw the boundary of piece14
\draw[black](2,10) -- (3,8) -- (8,7) -- (12,10) -- (4,13) -- (2,10);
% Draw the convex parts of piece11
\fill[black!15](0,28) -- (6,30) -- (5,33) -- (2,35) -- (0,35) -- (0,28);
% Draw the boundary of piece11
\draw[black](0,28) -- (6,30) -- (5,33) -- (2,35) -- (0,35) -- (0,28);
% Draw the convex parts of piece8
\fill[black!15](9,1) -- (10,3) -- (5,3) -- (9,1);
\fill[black!15](0,5) -- (1,0) -- (9,1) -- (5,3) -- (0,5);
% Draw the boundary of piece8
\draw[black](0,5) -- (1,0) -- (9,1) -- (10,3) -- (5,3) -- (0,5);
% Draw the convex parts of piece7
\fill[black!15](0,22) -- (2,18) -- (9,23) -- (3,24) -- (0,22);
% Draw the boundary of piece7
\draw[black](0,22) -- (2,18) -- (9,23) -- (3,24) -- (0,22);
% Draw the convex parts of piece3
\fill[black!15](4,26) -- (9,23) -- (12,24) -- (6,30) -- (4,26);
% Draw the boundary of piece3
\draw[black](4,26) -- (9,23) -- (12,24) -- (6,30) -- (4,26);
% Draw the convex parts of piece9
\fill[black!15](0,40) -- (1,35) -- (8,36) -- (4,39) -- (0,40);
% Draw the boundary of piece9
\draw[black](0,40) -- (1,35) -- (8,36) -- (4,39) -- (0,40);
% Draw the convex parts of piece5
\fill[black!15](6,30) -- (9,28) -- (11,29) -- (10,32) -- (5,35) -- (6,30);
% Draw the boundary of piece5
\draw[black](6,30) -- (9,28) -- (11,29) -- (10,32) -- (5,35) -- (6,30);
% Draw the convex parts of piece13
\fill[black!15](10,11) -- (12,12) -- (9,15) -- (10,11);
\fill[black!15](7,16) -- (5,13) -- (10,11) -- (9,15) -- (7,16);
\fill[black!15](9,15) -- (10,17) -- (7,16) -- (9,15);
% Draw the boundary of piece13
\draw[black](5,13) -- (10,11) -- (12,12) -- (9,15) -- (10,17) -- (7,16) -- (5,13);
% Draw the convex parts of piece1
\fill[black!15](1,5) -- (8,7) -- (3,8) -- (1,5);
\fill[black!15](1,5) -- (3,8) -- (0,12) -- (1,5);
% Draw the boundary of piece1
\draw[black](1,5) -- (8,7) -- (3,8) -- (0,12) -- (1,5);
% Draw the convex parts of piece6
\fill[black!15](7,34) -- (11,32) -- (12,35) -- (9,38) -- (7,34);
% Draw the boundary of piece6
\draw[black](7,34) -- (11,32) -- (12,35) -- (9,38) -- (7,34);
% Draw the convex parts of piece2
\fill[black!15](7,18) -- (12,17) -- (10,19) -- (7,18);
\fill[black!15](7,18) -- (10,19) -- (11,22) -- (9,23) -- (7,20) -- (7,18);
% Draw the boundary of piece2
\draw[black](7,18) -- (12,17) -- (10,19) -- (11,22) -- (9,23) -- (7,20) -- (7,18);
% Draw the convex parts of piece15
\fill[black!15](5,4) -- (8,3) -- (11,5) -- (10,8) -- (5,4);
% Draw the boundary of piece15
\draw[black](5,4) -- (8,3) -- (11,5) -- (10,8) -- (5,4);
% Draw the convex parts of piece4
\fill[black!15](0,16) -- (10,17) -- (3,18) -- (0,16);
% Draw the boundary of piece4
\draw[black](0,16) -- (10,17) -- (3,18) -- (0,16);
% Draw the convex parts of piece12
\fill[black!15](2,10) -- (4,14) -- (2,13) -- (2,10);
\fill[black!15](0,16) -- (2,10) -- (2,13) -- (0,16);
% Draw the boundary of piece12
\draw[black](0,16) -- (2,10) -- (4,14) -- (2,13) -- (0,16);
% Draw the convex parts of piece10
\fill[black!15](1,24) -- (5,25) -- (2,26) -- (1,24);
\fill[black!15](0,27) -- (1,24) -- (2,26) -- (0,27);
% Draw the boundary of piece10
\draw[black](0,27) -- (1,24) -- (5,25) -- (2,26) -- (0,27);
% Draw the bounding box of the solution
\draw[black](0,0) -- (12,0) -- (12,40)-- (0,40) -- (0,0);
\draw (6,-3) node {Poly1e};
\draw (6,-6) node {$L^*$ = 12};
\end{tikzpicture}
\caption{Seventeen open problem instances solved up to optimality within one hour for the first time in the literature. \citet{Rodrigues2017-zl} reported a wrong optimal solution for J2-10-35-3 with $L=21$ instead of the optimal solution with $L^*=20$ shown above, which we attribute to a minor tolerance error in their experiments. All the new optimal solutions shown above are found by our DB-PB algorithm, although the DB-CC model \cite{Rodrigues_de_Souza_Queiroz2020-tp} also finds four of them. However, this latter improvement can be attributed only to the advances in Gurobi 11.0.3 compared to CPLEX 12.6 \cite{Rodrigues2017-zl}.}
\label{fig:solved_open_instances}
\end{figure}

\begin{figure}
\centering
% SOLUTION for problem instance = SHAPES-5
% Scaling factor for coordinates = 1
\begin{tikzpicture}[scale = 0.07]
% Draw the convex parts of piece2
\fill[black!15](0,20) -- (6,14) -- (12,20) -- (6,26) -- (0,20);
% Draw the boundary of piece2
\draw[black](0,20) -- (6,14) -- (12,20) -- (6,26) -- (0,20);
% Draw the convex parts of piece2
\fill[black!15](6,14) -- (12,8) -- (18,14) -- (12,20) -- (6,14);
% Draw the boundary of piece2
\draw[black](6,14) -- (12,8) -- (18,14) -- (12,20) -- (6,14);
% Draw the convex parts of piece2
\fill[black!15](13,19) -- (19,13) -- (25,19) -- (19,25) -- (13,19);
% Draw the boundary of piece2
\draw[black](13,19) -- (19,13) -- (25,19) -- (19,25) -- (13,19);
% Draw the convex parts of piece2
\fill[black!15](16,10) -- (22,4) -- (28,10) -- (22,16) -- (16,10);
% Draw the boundary of piece2
\draw[black](16,10) -- (22,4) -- (28,10) -- (22,16) -- (16,10);
% Draw the convex parts of piece2
\fill[black!15](18,26) -- (24,20) -- (30,26) -- (24,32) -- (18,26);
% Draw the boundary of piece2
\draw[black](18,26) -- (24,20) -- (30,26) -- (24,32) -- (18,26);
% Draw the convex parts of piece4
\fill[black!15](2,12) -- (2,10) -- (4,10) -- (4,12) -- (2,12);
\fill[black!15](4,14) -- (4,16) -- (2,16) -- (2,14) -- (4,14);
\fill[black!15](0,12) -- (2,12) -- (4,12) -- (6,12) -- (6,14) -- (4,14) -- (2,14) -- (0,14) -- (0,12);
% Draw the boundary of piece4
\draw[black](0,12) -- (2,12) -- (2,10) -- (4,10) -- (4,12) -- (6,12) -- (6,14) -- (4,14) -- (4,16) -- (2,16) -- (2,14) -- (0,14) -- (0,12);
% Draw the convex parts of piece4
\fill[black!15](12,22) -- (12,20) -- (14,20) -- (14,22) -- (12,22);
\fill[black!15](14,24) -- (14,26) -- (12,26) -- (12,24) -- (14,24);
\fill[black!15](10,22) -- (12,22) -- (14,22) -- (16,22) -- (16,24) -- (14,24) -- (12,24) -- (10,24) -- (10,22);
% Draw the boundary of piece4
\draw[black](10,22) -- (12,22) -- (12,20) -- (14,20) -- (14,22) -- (16,22) -- (16,24) -- (14,24) -- (14,26) -- (12,26) -- (12,24) -- (10,24) -- (10,22);
% Draw the convex parts of piece4
\fill[black!15](16,28) -- (16,26) -- (18,26) -- (18,28) -- (16,28);
\fill[black!15](18,30) -- (18,32) -- (16,32) -- (16,30) -- (18,30);
\fill[black!15](14,28) -- (16,28) -- (18,28) -- (20,28) -- (20,30) -- (18,30) -- (16,30) -- (14,30) -- (14,28);
% Draw the boundary of piece4
\draw[black](14,28) -- (16,28) -- (16,26) -- (18,26) -- (18,28) -- (20,28) -- (20,30) -- (18,30) -- (18,32) -- (16,32) -- (16,30) -- (14,30) -- (14,28);
% Draw the convex parts of piece4
\fill[black!15](24,36) -- (24,34) -- (26,34) -- (26,36) -- (24,36);
\fill[black!15](26,38) -- (26,40) -- (24,40) -- (24,38) -- (26,38);
\fill[black!15](22,36) -- (24,36) -- (26,36) -- (28,36) -- (28,38) -- (26,38) -- (24,38) -- (22,38) -- (22,36);
% Draw the boundary of piece4
\draw[black](22,36) -- (24,36) -- (24,34) -- (26,34) -- (26,36) -- (28,36) -- (28,38) -- (26,38) -- (26,40) -- (24,40) -- (24,38) -- (22,38) -- (22,36);
% Draw the convex parts of piece4
\fill[black!15](26,14) -- (26,12) -- (28,12) -- (28,14) -- (26,14);
\fill[black!15](28,16) -- (28,18) -- (26,18) -- (26,16) -- (28,16);
\fill[black!15](24,14) -- (26,14) -- (28,14) -- (30,14) -- (30,16) -- (28,16) -- (26,16) -- (24,16) -- (24,14);
% Draw the boundary of piece4
\draw[black](24,14) -- (26,14) -- (26,12) -- (28,12) -- (28,14) -- (30,14) -- (30,16) -- (28,16) -- (28,18) -- (26,18) -- (26,16) -- (24,16) -- (24,14);
% Draw the convex parts of piece3
\fill[black!15](8,5) -- (8,2) -- (11,2) -- (8,5);
\fill[black!15](6,6) -- (0,6) -- (8,5) -- (6,6);
\fill[black!15](8,5) -- (11,2) -- (11,4) -- (7,8) -- (8,5);
\fill[black!15](8,5) -- (7,8) -- (6,8) -- (6,6) -- (8,5);
\fill[black!15](0,2) -- (2,2) -- (2,5) -- (0,2);
\fill[black!15](0,6) -- (0,2) -- (2,5) -- (0,6);
\fill[black!15](2,5) -- (8,5) -- (0,6) -- (2,5);
% Draw the boundary of piece3
\draw[black](0,2) -- (2,2) -- (2,5) -- (8,5) -- (8,2) -- (11,2) -- (11,4) -- (7,8) -- (6,8) -- (6,6) -- (0,6) -- (0,2);
% Draw the convex parts of piece3
\fill[black!15](8,37) -- (8,34) -- (11,34) -- (8,37);
\fill[black!15](6,38) -- (0,38) -- (8,37) -- (6,38);
\fill[black!15](8,37) -- (11,34) -- (11,36) -- (7,40) -- (8,37);
\fill[black!15](8,37) -- (7,40) -- (6,40) -- (6,38) -- (8,37);
\fill[black!15](0,34) -- (2,34) -- (2,37) -- (0,34);
\fill[black!15](0,38) -- (0,34) -- (2,37) -- (0,38);
\fill[black!15](2,37) -- (8,37) -- (0,38) -- (2,37);
% Draw the boundary of piece3
\draw[black](0,34) -- (2,34) -- (2,37) -- (8,37) -- (8,34) -- (11,34) -- (11,36) -- (7,40) -- (6,40) -- (6,38) -- (0,38) -- (0,34);
% Draw the convex parts of piece3
\fill[black!15](9,9) -- (9,6) -- (12,6) -- (9,9);
\fill[black!15](7,10) -- (1,10) -- (9,9) -- (7,10);
\fill[black!15](9,9) -- (12,6) -- (12,8) -- (8,12) -- (9,9);
\fill[black!15](9,9) -- (8,12) -- (7,12) -- (7,10) -- (9,9);
\fill[black!15](1,6) -- (3,6) -- (3,9) -- (1,6);
\fill[black!15](1,10) -- (1,6) -- (3,9) -- (1,10);
\fill[black!15](3,9) -- (9,9) -- (1,10) -- (3,9);
% Draw the boundary of piece3
\draw[black](1,6) -- (3,6) -- (3,9) -- (9,9) -- (9,6) -- (12,6) -- (12,8) -- (8,12) -- (7,12) -- (7,10) -- (1,10) -- (1,6);
% Draw the convex parts of piece3
\fill[black!15](19,5) -- (19,2) -- (22,2) -- (19,5);
\fill[black!15](17,6) -- (11,6) -- (19,5) -- (17,6);
\fill[black!15](19,5) -- (22,2) -- (22,4) -- (18,8) -- (19,5);
\fill[black!15](19,5) -- (18,8) -- (17,8) -- (17,6) -- (19,5);
\fill[black!15](11,2) -- (13,2) -- (13,5) -- (11,2);
\fill[black!15](11,6) -- (11,2) -- (13,5) -- (11,6);
\fill[black!15](13,5) -- (19,5) -- (11,6) -- (13,5);
% Draw the boundary of piece3
\draw[black](11,2) -- (13,2) -- (13,5) -- (19,5) -- (19,2) -- (22,2) -- (22,4) -- (18,8) -- (17,8) -- (17,6) -- (11,6) -- (11,2);
% Draw the convex parts of piece3
\fill[black!15](19,37) -- (19,34) -- (22,34) -- (19,37);
\fill[black!15](17,38) -- (11,38) -- (19,37) -- (17,38);
\fill[black!15](19,37) -- (22,34) -- (22,36) -- (18,40) -- (19,37);
\fill[black!15](19,37) -- (18,40) -- (17,40) -- (17,38) -- (19,37);
\fill[black!15](11,34) -- (13,34) -- (13,37) -- (11,34);
\fill[black!15](11,38) -- (11,34) -- (13,37) -- (11,38);
\fill[black!15](13,37) -- (19,37) -- (11,38) -- (13,37);
% Draw the boundary of piece3
\draw[black](11,34) -- (13,34) -- (13,37) -- (19,37) -- (19,34) -- (22,34) -- (22,36) -- (18,40) -- (17,40) -- (17,38) -- (11,38) -- (11,34);
% Draw the convex parts of piece1
\fill[black!15](12,28) -- (2,28) -- (14,26) -- (12,28);
\fill[black!15](14,26) -- (14,31) -- (12,31) -- (12,28) -- (14,26);
\fill[black!15](0,26) -- (14,26) -- (2,28) -- (0,26);
\fill[black!15](0,26) -- (2,28) -- (2,31) -- (0,31) -- (0,26);
% Draw the boundary of piece1
\draw[black](0,26) -- (14,26) -- (14,31) -- (12,31) -- (12,28) -- (2,28) -- (2,31) -- (0,31) -- (0,26);
% Draw the convex parts of piece1
\fill[black!15](14,2) -- (4,2) -- (16,0) -- (14,2);
\fill[black!15](16,0) -- (16,5) -- (14,5) -- (14,2) -- (16,0);
\fill[black!15](2,0) -- (16,0) -- (4,2) -- (2,0);
\fill[black!15](2,0) -- (4,2) -- (4,5) -- (2,5) -- (2,0);
% Draw the boundary of piece1
\draw[black](2,0) -- (16,0) -- (16,5) -- (14,5) -- (14,2) -- (4,2) -- (4,5) -- (2,5) -- (2,0);
% Draw the convex parts of piece1
\fill[black!15](14,34) -- (4,34) -- (16,32) -- (14,34);
\fill[black!15](16,32) -- (16,37) -- (14,37) -- (14,34) -- (16,32);
\fill[black!15](2,32) -- (16,32) -- (4,34) -- (2,32);
\fill[black!15](2,32) -- (4,34) -- (4,37) -- (2,37) -- (2,32);
% Draw the boundary of piece1
\draw[black](2,32) -- (16,32) -- (16,37) -- (14,37) -- (14,34) -- (4,34) -- (4,37) -- (2,37) -- (2,32);
% Draw the convex parts of piece1
\fill[black!15](28,2) -- (18,2) -- (30,0) -- (28,2);
\fill[black!15](30,0) -- (30,5) -- (28,5) -- (28,2) -- (30,0);
\fill[black!15](16,0) -- (30,0) -- (18,2) -- (16,0);
\fill[black!15](16,0) -- (18,2) -- (18,5) -- (16,5) -- (16,0);
% Draw the boundary of piece1
\draw[black](16,0) -- (30,0) -- (30,5) -- (28,5) -- (28,2) -- (18,2) -- (18,5) -- (16,5) -- (16,0);
% Draw the convex parts of piece1
\fill[black!15](28,34) -- (18,34) -- (30,32) -- (28,34);
\fill[black!15](30,32) -- (30,37) -- (28,37) -- (28,34) -- (30,32);
\fill[black!15](16,32) -- (30,32) -- (18,34) -- (16,32);
\fill[black!15](16,32) -- (18,34) -- (18,37) -- (16,37) -- (16,32);
% Draw the boundary of piece1
\draw[black](16,32) -- (30,32) -- (30,37) -- (28,37) -- (28,34) -- (18,34) -- (18,37) -- (16,37) -- (16,32);
% Draw the bounding box of the solution
\draw[black](0,0) -- (30,0) -- (30,40)-- (0,40) -- (0,0);
\draw (15,-2) node {SHAPES-5};
\draw (15,-5) node {$L$ = 30};
% Draw the bounding box of the solution
\end{tikzpicture}
% SOLUTION for problem instance = SHAPES-7
% Scaling factor for coordinates = 1
\begin{tikzpicture}[scale = 0.07]
% Draw the convex parts of piece2
\fill[black!15](0,30) -- (6,24) -- (12,30) -- (6,36) -- (0,30);
% Draw the boundary of piece2
\draw[black](0,30) -- (6,24) -- (12,30) -- (6,36) -- (0,30);
% Draw the convex parts of piece2
\fill[black!15](6,24) -- (12,18) -- (18,24) -- (12,30) -- (6,24);
% Draw the boundary of piece2
\draw[black](6,24) -- (12,18) -- (18,24) -- (12,30) -- (6,24);
% Draw the convex parts of piece2
\fill[black!15](8,34) -- (14,28) -- (20,34) -- (14,40) -- (8,34);
% Draw the boundary of piece2
\draw[black](8,34) -- (14,28) -- (20,34) -- (14,40) -- (8,34);
% Draw the convex parts of piece2
\fill[black!15](13,8) -- (19,2) -- (25,8) -- (19,14) -- (13,8);
% Draw the boundary of piece2
\draw[black](13,8) -- (19,2) -- (25,8) -- (19,14) -- (13,8);
% Draw the convex parts of piece2
\fill[black!15](17,29) -- (23,23) -- (29,29) -- (23,35) -- (17,29);
% Draw the boundary of piece2
\draw[black](17,29) -- (23,23) -- (29,29) -- (23,35) -- (17,29);
% Draw the convex parts of piece2
\fill[black!15](22,22) -- (28,16) -- (34,22) -- (28,28) -- (22,22);
% Draw the boundary of piece2
\draw[black](22,22) -- (28,16) -- (34,22) -- (28,28) -- (22,22);
% Draw the convex parts of piece2
\fill[black!15](30,26) -- (36,20) -- (42,26) -- (36,32) -- (30,26);
% Draw the boundary of piece2
\draw[black](30,26) -- (36,20) -- (42,26) -- (36,32) -- (30,26);
% Draw the convex parts of piece4
\fill[black!15](2,14) -- (2,12) -- (4,12) -- (4,14) -- (2,14);
\fill[black!15](4,16) -- (4,18) -- (2,18) -- (2,16) -- (4,16);
\fill[black!15](0,14) -- (2,14) -- (4,14) -- (6,14) -- (6,16) -- (4,16) -- (2,16) -- (0,16) -- (0,14);
% Draw the boundary of piece4
\draw[black](0,14) -- (2,14) -- (2,12) -- (4,12) -- (4,14) -- (6,14) -- (6,16) -- (4,16) -- (4,18) -- (2,18) -- (2,16) -- (0,16) -- (0,14);
% Draw the convex parts of piece4
\fill[black!15](2,22) -- (2,20) -- (4,20) -- (4,22) -- (2,22);
\fill[black!15](4,24) -- (4,26) -- (2,26) -- (2,24) -- (4,24);
\fill[black!15](0,22) -- (2,22) -- (4,22) -- (6,22) -- (6,24) -- (4,24) -- (2,24) -- (0,24) -- (0,22);
% Draw the boundary of piece4
\draw[black](0,22) -- (2,22) -- (2,20) -- (4,20) -- (4,22) -- (6,22) -- (6,24) -- (4,24) -- (4,26) -- (2,26) -- (2,24) -- (0,24) -- (0,22);
% Draw the convex parts of piece4
\fill[black!15](2,36) -- (2,34) -- (4,34) -- (4,36) -- (2,36);
\fill[black!15](4,38) -- (4,40) -- (2,40) -- (2,38) -- (4,38);
\fill[black!15](0,36) -- (2,36) -- (4,36) -- (6,36) -- (6,38) -- (4,38) -- (2,38) -- (0,38) -- (0,36);
% Draw the boundary of piece4
\draw[black](0,36) -- (2,36) -- (2,34) -- (4,34) -- (4,36) -- (6,36) -- (6,38) -- (4,38) -- (4,40) -- (2,40) -- (2,38) -- (0,38) -- (0,36);
% Draw the convex parts of piece4
\fill[black!15](3,5) -- (3,3) -- (5,3) -- (5,5) -- (3,5);
\fill[black!15](5,7) -- (5,9) -- (3,9) -- (3,7) -- (5,7);
\fill[black!15](1,5) -- (3,5) -- (5,5) -- (7,5) -- (7,7) -- (5,7) -- (3,7) -- (1,7) -- (1,5);
% Draw the boundary of piece4
\draw[black](1,5) -- (3,5) -- (3,3) -- (5,3) -- (5,5) -- (7,5) -- (7,7) -- (5,7) -- (5,9) -- (3,9) -- (3,7) -- (1,7) -- (1,5);
% Draw the convex parts of piece4
\fill[black!15](9,5) -- (9,3) -- (11,3) -- (11,5) -- (9,5);
\fill[black!15](11,7) -- (11,9) -- (9,9) -- (9,7) -- (11,7);
\fill[black!15](7,5) -- (9,5) -- (11,5) -- (13,5) -- (13,7) -- (11,7) -- (9,7) -- (7,7) -- (7,5);
% Draw the boundary of piece4
\draw[black](7,5) -- (9,5) -- (9,3) -- (11,3) -- (11,5) -- (13,5) -- (13,7) -- (11,7) -- (11,9) -- (9,9) -- (9,7) -- (7,7) -- (7,5);
% Draw the convex parts of piece4
\fill[black!15](36,4) -- (36,2) -- (38,2) -- (38,4) -- (36,4);
\fill[black!15](38,6) -- (38,8) -- (36,8) -- (36,6) -- (38,6);
\fill[black!15](34,4) -- (36,4) -- (38,4) -- (40,4) -- (40,6) -- (38,6) -- (36,6) -- (34,6) -- (34,4);
% Draw the boundary of piece4
\draw[black](34,4) -- (36,4) -- (36,2) -- (38,2) -- (38,4) -- (40,4) -- (40,6) -- (38,6) -- (38,8) -- (36,8) -- (36,6) -- (34,6) -- (34,4);
% Draw the convex parts of piece4
\fill[black!15](38,18) -- (38,16) -- (40,16) -- (40,18) -- (38,18);
\fill[black!15](40,20) -- (40,22) -- (38,22) -- (38,20) -- (40,20);
\fill[black!15](36,18) -- (38,18) -- (40,18) -- (42,18) -- (42,20) -- (40,20) -- (38,20) -- (36,20) -- (36,18);
% Draw the boundary of piece4
\draw[black](36,18) -- (38,18) -- (38,16) -- (40,16) -- (40,18) -- (42,18) -- (42,20) -- (40,20) -- (40,22) -- (38,22) -- (38,20) -- (36,20) -- (36,18);
% Draw the convex parts of piece3
\fill[black!15](8,19) -- (8,16) -- (11,16) -- (8,19);
\fill[black!15](6,20) -- (0,20) -- (8,19) -- (6,20);
\fill[black!15](8,19) -- (11,16) -- (11,18) -- (7,22) -- (8,19);
\fill[black!15](8,19) -- (7,22) -- (6,22) -- (6,20) -- (8,19);
\fill[black!15](0,16) -- (2,16) -- (2,19) -- (0,16);
\fill[black!15](0,20) -- (0,16) -- (2,19) -- (0,20);
\fill[black!15](2,19) -- (8,19) -- (0,20) -- (2,19);
% Draw the boundary of piece3
\draw[black](0,16) -- (2,16) -- (2,19) -- (8,19) -- (8,16) -- (11,16) -- (11,18) -- (7,22) -- (6,22) -- (6,20) -- (0,20) -- (0,16);
% Draw the convex parts of piece3
\fill[black!15](23,19) -- (23,16) -- (26,16) -- (23,19);
\fill[black!15](21,20) -- (15,20) -- (23,19) -- (21,20);
\fill[black!15](23,19) -- (26,16) -- (26,18) -- (22,22) -- (23,19);
\fill[black!15](23,19) -- (22,22) -- (21,22) -- (21,20) -- (23,19);
\fill[black!15](15,16) -- (17,16) -- (17,19) -- (15,16);
\fill[black!15](15,20) -- (15,16) -- (17,19) -- (15,20);
\fill[black!15](17,19) -- (23,19) -- (15,20) -- (17,19);
% Draw the boundary of piece3
\draw[black](15,16) -- (17,16) -- (17,19) -- (23,19) -- (23,16) -- (26,16) -- (26,18) -- (22,22) -- (21,22) -- (21,20) -- (15,20) -- (15,16);
% Draw the convex parts of piece3
\fill[black!15](28,37) -- (28,34) -- (31,34) -- (28,37);
\fill[black!15](26,38) -- (20,38) -- (28,37) -- (26,38);
\fill[black!15](28,37) -- (31,34) -- (31,36) -- (27,40) -- (28,37);
\fill[black!15](28,37) -- (27,40) -- (26,40) -- (26,38) -- (28,37);
\fill[black!15](20,34) -- (22,34) -- (22,37) -- (20,34);
\fill[black!15](20,38) -- (20,34) -- (22,37) -- (20,38);
\fill[black!15](22,37) -- (28,37) -- (20,38) -- (22,37);
% Draw the boundary of piece3
\draw[black](20,34) -- (22,34) -- (22,37) -- (28,37) -- (28,34) -- (31,34) -- (31,36) -- (27,40) -- (26,40) -- (26,38) -- (20,38) -- (20,34);
% Draw the convex parts of piece3
\fill[black!15](31,5) -- (31,2) -- (34,2) -- (31,5);
\fill[black!15](29,6) -- (23,6) -- (31,5) -- (29,6);
\fill[black!15](31,5) -- (34,2) -- (34,4) -- (30,8) -- (31,5);
\fill[black!15](31,5) -- (30,8) -- (29,8) -- (29,6) -- (31,5);
\fill[black!15](23,2) -- (25,2) -- (25,5) -- (23,2);
\fill[black!15](23,6) -- (23,2) -- (25,5) -- (23,6);
\fill[black!15](25,5) -- (31,5) -- (23,6) -- (25,5);
% Draw the boundary of piece3
\draw[black](23,2) -- (25,2) -- (25,5) -- (31,5) -- (31,2) -- (34,2) -- (34,4) -- (30,8) -- (29,8) -- (29,6) -- (23,6) -- (23,2);
% Draw the convex parts of piece3
\fill[black!15](39,13) -- (39,10) -- (42,10) -- (39,13);
\fill[black!15](37,14) -- (31,14) -- (39,13) -- (37,14);
\fill[black!15](39,13) -- (42,10) -- (42,12) -- (38,16) -- (39,13);
\fill[black!15](39,13) -- (38,16) -- (37,16) -- (37,14) -- (39,13);
\fill[black!15](31,10) -- (33,10) -- (33,13) -- (31,10);
\fill[black!15](31,14) -- (31,10) -- (33,13) -- (31,14);
\fill[black!15](33,13) -- (39,13) -- (31,14) -- (33,13);
% Draw the boundary of piece3
\draw[black](31,10) -- (33,10) -- (33,13) -- (39,13) -- (39,10) -- (42,10) -- (42,12) -- (38,16) -- (37,16) -- (37,14) -- (31,14) -- (31,10);
% Draw the convex parts of piece3
\fill[black!15](39,32) -- (39,29) -- (42,29) -- (39,32);
\fill[black!15](37,33) -- (31,33) -- (39,32) -- (37,33);
\fill[black!15](39,32) -- (42,29) -- (42,31) -- (38,35) -- (39,32);
\fill[black!15](39,32) -- (38,35) -- (37,35) -- (37,33) -- (39,32);
\fill[black!15](31,29) -- (33,29) -- (33,32) -- (31,29);
\fill[black!15](31,33) -- (31,29) -- (33,32) -- (31,33);
\fill[black!15](33,32) -- (39,32) -- (31,33) -- (33,32);
% Draw the boundary of piece3
\draw[black](31,29) -- (33,29) -- (33,32) -- (39,32) -- (39,29) -- (42,29) -- (42,31) -- (38,35) -- (37,35) -- (37,33) -- (31,33) -- (31,29);
% Draw the convex parts of piece3
\fill[black!15](39,37) -- (39,34) -- (42,34) -- (39,37);
\fill[black!15](37,38) -- (31,38) -- (39,37) -- (37,38);
\fill[black!15](39,37) -- (42,34) -- (42,36) -- (38,40) -- (39,37);
\fill[black!15](39,37) -- (38,40) -- (37,40) -- (37,38) -- (39,37);
\fill[black!15](31,34) -- (33,34) -- (33,37) -- (31,34);
\fill[black!15](31,38) -- (31,34) -- (33,37) -- (31,38);
\fill[black!15](33,37) -- (39,37) -- (31,38) -- (33,37);
% Draw the boundary of piece3
\draw[black](31,34) -- (33,34) -- (33,37) -- (39,37) -- (39,34) -- (42,34) -- (42,36) -- (38,40) -- (37,40) -- (37,38) -- (31,38) -- (31,34);
% Draw the convex parts of piece1
\fill[black!15](12,2) -- (2,2) -- (14,0) -- (12,2);
\fill[black!15](14,0) -- (14,5) -- (12,5) -- (12,2) -- (14,0);
\fill[black!15](0,0) -- (14,0) -- (2,2) -- (0,0);
\fill[black!15](0,0) -- (2,2) -- (2,5) -- (0,5) -- (0,0);
% Draw the boundary of piece1
\draw[black](0,0) -- (14,0) -- (14,5) -- (12,5) -- (12,2) -- (2,2) -- (2,5) -- (0,5) -- (0,0);
% Draw the convex parts of piece1
\fill[black!15](12,11) -- (2,11) -- (14,9) -- (12,11);
\fill[black!15](14,9) -- (14,14) -- (12,14) -- (12,11) -- (14,9);
\fill[black!15](0,9) -- (14,9) -- (2,11) -- (0,9);
\fill[black!15](0,9) -- (2,11) -- (2,14) -- (0,14) -- (0,9);
% Draw the boundary of piece1
\draw[black](0,9) -- (14,9) -- (14,14) -- (12,14) -- (12,11) -- (2,11) -- (2,14) -- (0,14) -- (0,9);
% Draw the convex parts of piece1
\fill[black!15](18,16) -- (8,16) -- (20,14) -- (18,16);
\fill[black!15](20,14) -- (20,19) -- (18,19) -- (18,16) -- (20,14);
\fill[black!15](6,14) -- (20,14) -- (8,16) -- (6,14);
\fill[black!15](6,14) -- (8,16) -- (8,19) -- (6,19) -- (6,14);
% Draw the boundary of piece1
\draw[black](6,14) -- (20,14) -- (20,19) -- (18,19) -- (18,16) -- (8,16) -- (8,19) -- (6,19) -- (6,14);
% Draw the convex parts of piece1
\fill[black!15](26,2) -- (16,2) -- (28,0) -- (26,2);
\fill[black!15](28,0) -- (28,5) -- (26,5) -- (26,2) -- (28,0);
\fill[black!15](14,0) -- (28,0) -- (16,2) -- (14,0);
\fill[black!15](14,0) -- (16,2) -- (16,5) -- (14,5) -- (14,0);
% Draw the boundary of piece1
\draw[black](14,0) -- (28,0) -- (28,5) -- (26,5) -- (26,2) -- (16,2) -- (16,5) -- (14,5) -- (14,0);
% Draw the convex parts of piece1
\fill[black!15](32,16) -- (22,16) -- (34,14) -- (32,16);
\fill[black!15](34,14) -- (34,19) -- (32,19) -- (32,16) -- (34,14);
\fill[black!15](20,14) -- (34,14) -- (22,16) -- (20,14);
\fill[black!15](20,14) -- (22,16) -- (22,19) -- (20,19) -- (20,14);
% Draw the boundary of piece1
\draw[black](20,14) -- (34,14) -- (34,19) -- (32,19) -- (32,16) -- (22,16) -- (22,19) -- (20,19) -- (20,14);
% Draw the convex parts of piece1
\fill[black!15](37,10) -- (27,10) -- (39,8) -- (37,10);
\fill[black!15](39,8) -- (39,13) -- (37,13) -- (37,10) -- (39,8);
\fill[black!15](25,8) -- (39,8) -- (27,10) -- (25,8);
\fill[black!15](25,8) -- (27,10) -- (27,13) -- (25,13) -- (25,8);
% Draw the boundary of piece1
\draw[black](25,8) -- (39,8) -- (39,13) -- (37,13) -- (37,10) -- (27,10) -- (27,13) -- (25,13) -- (25,8);
% Draw the convex parts of piece1
\fill[black!15](40,2) -- (30,2) -- (42,0) -- (40,2);
\fill[black!15](42,0) -- (42,5) -- (40,5) -- (40,2) -- (42,0);
\fill[black!15](28,0) -- (42,0) -- (30,2) -- (28,0);
\fill[black!15](28,0) -- (30,2) -- (30,5) -- (28,5) -- (28,0);
% Draw the boundary of piece1
\draw[black](28,0) -- (42,0) -- (42,5) -- (40,5) -- (40,2) -- (30,2) -- (30,5) -- (28,5) -- (28,0);
% Draw the bounding box of the solution
\draw[black](0,0) -- (42,0) -- (42,40)-- (0,40) -- (0,0);
\draw (21,-2) node {SHAPES-7};
\draw (21,-5) node {$L$ = 42};
% Draw the bounding box of the solution
\end{tikzpicture}
% SOLUTION for problem instance = SHAPES-9
% Scaling factor for coordinates = 1
\begin{tikzpicture}[scale = 0.07]
% Draw the convex parts of piece2
\fill[black!15](1,11) -- (7,5) -- (13,11) -- (7,17) -- (1,11);
% Draw the boundary of piece2
\draw[black](1,11) -- (7,5) -- (13,11) -- (7,17) -- (1,11);
% Draw the convex parts of piece2
\fill[black!15](19,21) -- (25,15) -- (31,21) -- (25,27) -- (19,21);
% Draw the boundary of piece2
\draw[black](19,21) -- (25,15) -- (31,21) -- (25,27) -- (19,21);
% Draw the convex parts of piece2
\fill[black!15](20,34) -- (26,28) -- (32,34) -- (26,40) -- (20,34);
% Draw the boundary of piece2
\draw[black](20,34) -- (26,28) -- (32,34) -- (26,40) -- (20,34);
% Draw the convex parts of piece2
\fill[black!15](25,27) -- (31,21) -- (37,27) -- (31,33) -- (25,27);
% Draw the boundary of piece2
\draw[black](25,27) -- (31,21) -- (37,27) -- (31,33) -- (25,27);
% Draw the convex parts of piece2
\fill[black!15](26,16) -- (32,10) -- (38,16) -- (32,22) -- (26,16);
% Draw the boundary of piece2
\draw[black](26,16) -- (32,10) -- (38,16) -- (32,22) -- (26,16);
% Draw the convex parts of piece2
\fill[black!15](32,22) -- (38,16) -- (44,22) -- (38,28) -- (32,22);
% Draw the boundary of piece2
\draw[black](32,22) -- (38,16) -- (44,22) -- (38,28) -- (32,22);
% Draw the convex parts of piece2
\fill[black!15](40,18) -- (46,12) -- (52,18) -- (46,24) -- (40,18);
% Draw the boundary of piece2
\draw[black](40,18) -- (46,12) -- (52,18) -- (46,24) -- (40,18);
% Draw the convex parts of piece4
\fill[black!15](2,36) -- (2,34) -- (4,34) -- (4,36) -- (2,36);
\fill[black!15](4,38) -- (4,40) -- (2,40) -- (2,38) -- (4,38);
\fill[black!15](0,36) -- (2,36) -- (4,36) -- (6,36) -- (6,38) -- (4,38) -- (2,38) -- (0,38) -- (0,36);
% Draw the boundary of piece4
\draw[black](0,36) -- (2,36) -- (2,34) -- (4,34) -- (4,36) -- (6,36) -- (6,38) -- (4,38) -- (4,40) -- (2,40) -- (2,38) -- (0,38) -- (0,36);
% Draw the convex parts of piece4
\fill[black!15](13,23) -- (13,21) -- (15,21) -- (15,23) -- (13,23);
\fill[black!15](15,25) -- (15,27) -- (13,27) -- (13,25) -- (15,25);
\fill[black!15](11,23) -- (13,23) -- (15,23) -- (17,23) -- (17,25) -- (15,25) -- (13,25) -- (11,25) -- (11,23);
% Draw the boundary of piece4
\draw[black](11,23) -- (13,23) -- (13,21) -- (15,21) -- (15,23) -- (17,23) -- (17,25) -- (15,25) -- (15,27) -- (13,27) -- (13,25) -- (11,25) -- (11,23);
% Draw the convex parts of piece4
\fill[black!15](15,12) -- (15,10) -- (17,10) -- (17,12) -- (15,12);
\fill[black!15](17,14) -- (17,16) -- (15,16) -- (15,14) -- (17,14);
\fill[black!15](13,12) -- (15,12) -- (17,12) -- (19,12) -- (19,14) -- (17,14) -- (15,14) -- (13,14) -- (13,12);
% Draw the boundary of piece4
\draw[black](13,12) -- (15,12) -- (15,10) -- (17,10) -- (17,12) -- (19,12) -- (19,14) -- (17,14) -- (17,16) -- (15,16) -- (15,14) -- (13,14) -- (13,12);
% Draw the convex parts of piece4
\fill[black!15](17,2) -- (17,0) -- (19,0) -- (19,2) -- (17,2);
\fill[black!15](19,4) -- (19,6) -- (17,6) -- (17,4) -- (19,4);
\fill[black!15](15,2) -- (17,2) -- (19,2) -- (21,2) -- (21,4) -- (19,4) -- (17,4) -- (15,4) -- (15,2);
% Draw the boundary of piece4
\draw[black](15,2) -- (17,2) -- (17,0) -- (19,0) -- (19,2) -- (21,2) -- (21,4) -- (19,4) -- (19,6) -- (17,6) -- (17,4) -- (15,4) -- (15,2);
% Draw the convex parts of piece4
\fill[black!15](17,25) -- (17,23) -- (19,23) -- (19,25) -- (17,25);
\fill[black!15](19,27) -- (19,29) -- (17,29) -- (17,27) -- (19,27);
\fill[black!15](15,25) -- (17,25) -- (19,25) -- (21,25) -- (21,27) -- (19,27) -- (17,27) -- (15,27) -- (15,25);
% Draw the boundary of piece4
\draw[black](15,25) -- (17,25) -- (17,23) -- (19,23) -- (19,25) -- (21,25) -- (21,27) -- (19,27) -- (19,29) -- (17,29) -- (17,27) -- (15,27) -- (15,25);
% Draw the convex parts of piece4
\fill[black!15](32,6) -- (32,4) -- (34,4) -- (34,6) -- (32,6);
\fill[black!15](34,8) -- (34,10) -- (32,10) -- (32,8) -- (34,8);
\fill[black!15](30,6) -- (32,6) -- (34,6) -- (36,6) -- (36,8) -- (34,8) -- (32,8) -- (30,8) -- (30,6);
% Draw the boundary of piece4
\draw[black](30,6) -- (32,6) -- (32,4) -- (34,4) -- (34,6) -- (36,6) -- (36,8) -- (34,8) -- (34,10) -- (32,10) -- (32,8) -- (30,8) -- (30,6);
% Draw the convex parts of piece4
\fill[black!15](37,36) -- (37,34) -- (39,34) -- (39,36) -- (37,36);
\fill[black!15](39,38) -- (39,40) -- (37,40) -- (37,38) -- (39,38);
\fill[black!15](35,36) -- (37,36) -- (39,36) -- (41,36) -- (41,38) -- (39,38) -- (37,38) -- (35,38) -- (35,36);
% Draw the boundary of piece4
\draw[black](35,36) -- (37,36) -- (37,34) -- (39,34) -- (39,36) -- (41,36) -- (41,38) -- (39,38) -- (39,40) -- (37,40) -- (37,38) -- (35,38) -- (35,36);
% Draw the convex parts of piece4
\fill[black!15](43,26) -- (43,24) -- (45,24) -- (45,26) -- (43,26);
\fill[black!15](45,28) -- (45,30) -- (43,30) -- (43,28) -- (45,28);
\fill[black!15](41,26) -- (43,26) -- (45,26) -- (47,26) -- (47,28) -- (45,28) -- (43,28) -- (41,28) -- (41,26);
% Draw the boundary of piece4
\draw[black](41,26) -- (43,26) -- (43,24) -- (45,24) -- (45,26) -- (47,26) -- (47,28) -- (45,28) -- (45,30) -- (43,30) -- (43,28) -- (41,28) -- (41,26);
% Draw the convex parts of piece4
\fill[black!15](48,24) -- (48,22) -- (50,22) -- (50,24) -- (48,24);
\fill[black!15](50,26) -- (50,28) -- (48,28) -- (48,26) -- (50,26);
\fill[black!15](46,24) -- (48,24) -- (50,24) -- (52,24) -- (52,26) -- (50,26) -- (48,26) -- (46,26) -- (46,24);
% Draw the boundary of piece4
\draw[black](46,24) -- (48,24) -- (48,22) -- (50,22) -- (50,24) -- (52,24) -- (52,26) -- (50,26) -- (50,28) -- (48,28) -- (48,26) -- (46,26) -- (46,24);
% Draw the convex parts of piece3
\fill[black!15](8,22) -- (8,19) -- (11,19) -- (8,22);
\fill[black!15](6,23) -- (0,23) -- (8,22) -- (6,23);
\fill[black!15](8,22) -- (11,19) -- (11,21) -- (7,25) -- (8,22);
\fill[black!15](8,22) -- (7,25) -- (6,25) -- (6,23) -- (8,22);
\fill[black!15](0,19) -- (2,19) -- (2,22) -- (0,19);
\fill[black!15](0,23) -- (0,19) -- (2,22) -- (0,23);
\fill[black!15](2,22) -- (8,22) -- (0,23) -- (2,22);
% Draw the boundary of piece3
\draw[black](0,19) -- (2,19) -- (2,22) -- (8,22) -- (8,19) -- (11,19) -- (11,21) -- (7,25) -- (6,25) -- (6,23) -- (0,23) -- (0,19);
% Draw the convex parts of piece3
\fill[black!15](8,32) -- (8,29) -- (11,29) -- (8,32);
\fill[black!15](6,33) -- (0,33) -- (8,32) -- (6,33);
\fill[black!15](8,32) -- (11,29) -- (11,31) -- (7,35) -- (8,32);
\fill[black!15](8,32) -- (7,35) -- (6,35) -- (6,33) -- (8,32);
\fill[black!15](0,29) -- (2,29) -- (2,32) -- (0,29);
\fill[black!15](0,33) -- (0,29) -- (2,32) -- (0,33);
\fill[black!15](2,32) -- (8,32) -- (0,33) -- (2,32);
% Draw the boundary of piece3
\draw[black](0,29) -- (2,29) -- (2,32) -- (8,32) -- (8,29) -- (11,29) -- (11,31) -- (7,35) -- (6,35) -- (6,33) -- (0,33) -- (0,29);
% Draw the convex parts of piece3
\fill[black!15](10,3) -- (10,0) -- (13,0) -- (10,3);
\fill[black!15](8,4) -- (2,4) -- (10,3) -- (8,4);
\fill[black!15](10,3) -- (13,0) -- (13,2) -- (9,6) -- (10,3);
\fill[black!15](10,3) -- (9,6) -- (8,6) -- (8,4) -- (10,3);
\fill[black!15](2,0) -- (4,0) -- (4,3) -- (2,0);
\fill[black!15](2,4) -- (2,0) -- (4,3) -- (2,4);
\fill[black!15](4,3) -- (10,3) -- (2,4) -- (4,3);
% Draw the boundary of piece3
\draw[black](2,0) -- (4,0) -- (4,3) -- (10,3) -- (10,0) -- (13,0) -- (13,2) -- (9,6) -- (8,6) -- (8,4) -- (2,4) -- (2,0);
% Draw the convex parts of piece3
\fill[black!15](19,16) -- (19,13) -- (22,13) -- (19,16);
\fill[black!15](17,17) -- (11,17) -- (19,16) -- (17,17);
\fill[black!15](19,16) -- (22,13) -- (22,15) -- (18,19) -- (19,16);
\fill[black!15](19,16) -- (18,19) -- (17,19) -- (17,17) -- (19,16);
\fill[black!15](11,13) -- (13,13) -- (13,16) -- (11,13);
\fill[black!15](11,17) -- (11,13) -- (13,16) -- (11,17);
\fill[black!15](13,16) -- (19,16) -- (11,17) -- (13,16);
% Draw the boundary of piece3
\draw[black](11,13) -- (13,13) -- (13,16) -- (19,16) -- (19,13) -- (22,13) -- (22,15) -- (18,19) -- (17,19) -- (17,17) -- (11,17) -- (11,13);
% Draw the convex parts of piece3
\fill[black!15](20,32) -- (20,29) -- (23,29) -- (20,32);
\fill[black!15](18,33) -- (12,33) -- (20,32) -- (18,33);
\fill[black!15](20,32) -- (23,29) -- (23,31) -- (19,35) -- (20,32);
\fill[black!15](20,32) -- (19,35) -- (18,35) -- (18,33) -- (20,32);
\fill[black!15](12,29) -- (14,29) -- (14,32) -- (12,29);
\fill[black!15](12,33) -- (12,29) -- (14,32) -- (12,33);
\fill[black!15](14,32) -- (20,32) -- (12,33) -- (14,32);
% Draw the boundary of piece3
\draw[black](12,29) -- (14,29) -- (14,32) -- (20,32) -- (20,29) -- (23,29) -- (23,31) -- (19,35) -- (18,35) -- (18,33) -- (12,33) -- (12,29);
% Draw the convex parts of piece3
\fill[black!15](27,12) -- (27,9) -- (30,9) -- (27,12);
\fill[black!15](25,13) -- (19,13) -- (27,12) -- (25,13);
\fill[black!15](27,12) -- (30,9) -- (30,11) -- (26,15) -- (27,12);
\fill[black!15](27,12) -- (26,15) -- (25,15) -- (25,13) -- (27,12);
\fill[black!15](19,9) -- (21,9) -- (21,12) -- (19,9);
\fill[black!15](19,13) -- (19,9) -- (21,12) -- (19,13);
\fill[black!15](21,12) -- (27,12) -- (19,13) -- (21,12);
% Draw the boundary of piece3
\draw[black](19,9) -- (21,9) -- (21,12) -- (27,12) -- (27,9) -- (30,9) -- (30,11) -- (26,15) -- (25,15) -- (25,13) -- (19,13) -- (19,9);
% Draw the convex parts of piece3
\fill[black!15](29,5) -- (29,2) -- (32,2) -- (29,5);
\fill[black!15](27,6) -- (21,6) -- (29,5) -- (27,6);
\fill[black!15](29,5) -- (32,2) -- (32,4) -- (28,8) -- (29,5);
\fill[black!15](29,5) -- (28,8) -- (27,8) -- (27,6) -- (29,5);
\fill[black!15](21,2) -- (23,2) -- (23,5) -- (21,2);
\fill[black!15](21,6) -- (21,2) -- (23,5) -- (21,6);
\fill[black!15](23,5) -- (29,5) -- (21,6) -- (23,5);
% Draw the boundary of piece3
\draw[black](21,2) -- (23,2) -- (23,5) -- (29,5) -- (29,2) -- (32,2) -- (32,4) -- (28,8) -- (27,8) -- (27,6) -- (21,6) -- (21,2);
% Draw the convex parts of piece3
\fill[black!15](49,31) -- (49,28) -- (52,28) -- (49,31);
\fill[black!15](47,32) -- (41,32) -- (49,31) -- (47,32);
\fill[black!15](49,31) -- (52,28) -- (52,30) -- (48,34) -- (49,31);
\fill[black!15](49,31) -- (48,34) -- (47,34) -- (47,32) -- (49,31);
\fill[black!15](41,28) -- (43,28) -- (43,31) -- (41,28);
\fill[black!15](41,32) -- (41,28) -- (43,31) -- (41,32);
\fill[black!15](43,31) -- (49,31) -- (41,32) -- (43,31);
% Draw the boundary of piece3
\draw[black](41,28) -- (43,28) -- (43,31) -- (49,31) -- (49,28) -- (52,28) -- (52,30) -- (48,34) -- (47,34) -- (47,32) -- (41,32) -- (41,28);
% Draw the convex parts of piece3
\fill[black!15](49,37) -- (49,34) -- (52,34) -- (49,37);
\fill[black!15](47,38) -- (41,38) -- (49,37) -- (47,38);
\fill[black!15](49,37) -- (52,34) -- (52,36) -- (48,40) -- (49,37);
\fill[black!15](49,37) -- (48,40) -- (47,40) -- (47,38) -- (49,37);
\fill[black!15](41,34) -- (43,34) -- (43,37) -- (41,34);
\fill[black!15](41,38) -- (41,34) -- (43,37) -- (41,38);
\fill[black!15](43,37) -- (49,37) -- (41,38) -- (43,37);
% Draw the boundary of piece3
\draw[black](41,34) -- (43,34) -- (43,37) -- (49,37) -- (49,34) -- (52,34) -- (52,36) -- (48,40) -- (47,40) -- (47,38) -- (41,38) -- (41,34);
% Draw the convex parts of piece1
\fill[black!15](15,19) -- (5,19) -- (17,17) -- (15,19);
\fill[black!15](17,17) -- (17,22) -- (15,22) -- (15,19) -- (17,17);
\fill[black!15](3,17) -- (17,17) -- (5,19) -- (3,17);
\fill[black!15](3,17) -- (5,19) -- (5,22) -- (3,22) -- (3,17);
% Draw the boundary of piece1
\draw[black](3,17) -- (17,17) -- (17,22) -- (15,22) -- (15,19) -- (5,19) -- (5,22) -- (3,22) -- (3,17);
% Draw the convex parts of piece1
\fill[black!15](15,29) -- (5,29) -- (17,27) -- (15,29);
\fill[black!15](17,27) -- (17,32) -- (15,32) -- (15,29) -- (17,27);
\fill[black!15](3,27) -- (17,27) -- (5,29) -- (3,27);
\fill[black!15](3,27) -- (5,29) -- (5,32) -- (3,32) -- (3,27);
% Draw the boundary of piece1
\draw[black](3,27) -- (17,27) -- (17,32) -- (15,32) -- (15,29) -- (5,29) -- (5,32) -- (3,32) -- (3,27);
% Draw the convex parts of piece1
\fill[black!15](19,37) -- (9,37) -- (21,35) -- (19,37);
\fill[black!15](21,35) -- (21,40) -- (19,40) -- (19,37) -- (21,35);
\fill[black!15](7,35) -- (21,35) -- (9,37) -- (7,35);
\fill[black!15](7,35) -- (9,37) -- (9,40) -- (7,40) -- (7,35);
% Draw the boundary of piece1
\draw[black](7,35) -- (21,35) -- (21,40) -- (19,40) -- (19,37) -- (9,37) -- (9,40) -- (7,40) -- (7,35);
% Draw the convex parts of piece1
\fill[black!15](25,8) -- (15,8) -- (27,6) -- (25,8);
\fill[black!15](27,6) -- (27,11) -- (25,11) -- (25,8) -- (27,6);
\fill[black!15](13,6) -- (27,6) -- (15,8) -- (13,6);
\fill[black!15](13,6) -- (15,8) -- (15,11) -- (13,11) -- (13,6);
% Draw the boundary of piece1
\draw[black](13,6) -- (27,6) -- (27,11) -- (25,11) -- (25,8) -- (15,8) -- (15,11) -- (13,11) -- (13,6);
% Draw the convex parts of piece1
\fill[black!15](36,2) -- (26,2) -- (38,0) -- (36,2);
\fill[black!15](38,0) -- (38,5) -- (36,5) -- (36,2) -- (38,0);
\fill[black!15](24,0) -- (38,0) -- (26,2) -- (24,0);
\fill[black!15](24,0) -- (26,2) -- (26,5) -- (24,5) -- (24,0);
% Draw the boundary of piece1
\draw[black](24,0) -- (38,0) -- (38,5) -- (36,5) -- (36,2) -- (26,2) -- (26,5) -- (24,5) -- (24,0);
% Draw the convex parts of piece1
\fill[black!15](45,34) -- (35,34) -- (47,32) -- (45,34);
\fill[black!15](47,32) -- (47,37) -- (45,37) -- (45,34) -- (47,32);
\fill[black!15](33,32) -- (47,32) -- (35,34) -- (33,32);
\fill[black!15](33,32) -- (35,34) -- (35,37) -- (33,37) -- (33,32);
% Draw the boundary of piece1
\draw[black](33,32) -- (47,32) -- (47,37) -- (45,37) -- (45,34) -- (35,34) -- (35,37) -- (33,37) -- (33,32);
% Draw the convex parts of piece1
\fill[black!15](49,7) -- (39,7) -- (51,5) -- (49,7);
\fill[black!15](51,5) -- (51,10) -- (49,10) -- (49,7) -- (51,5);
\fill[black!15](37,5) -- (51,5) -- (39,7) -- (37,5);
\fill[black!15](37,5) -- (39,7) -- (39,10) -- (37,10) -- (37,5);
% Draw the boundary of piece1
\draw[black](37,5) -- (51,5) -- (51,10) -- (49,10) -- (49,7) -- (39,7) -- (39,10) -- (37,10) -- (37,5);
% Draw the convex parts of piece1
\fill[black!15](50,2) -- (40,2) -- (52,0) -- (50,2);
\fill[black!15](52,0) -- (52,5) -- (50,5) -- (50,2) -- (52,0);
\fill[black!15](38,0) -- (52,0) -- (40,2) -- (38,0);
\fill[black!15](38,0) -- (40,2) -- (40,5) -- (38,5) -- (38,0);
% Draw the boundary of piece1
\draw[black](38,0) -- (52,0) -- (52,5) -- (50,5) -- (50,2) -- (40,2) -- (40,5) -- (38,5) -- (38,0);
% Draw the convex parts of piece1
\fill[black!15](50,12) -- (40,12) -- (52,10) -- (50,12);
\fill[black!15](52,10) -- (52,15) -- (50,15) -- (50,12) -- (52,10);
\fill[black!15](38,10) -- (52,10) -- (40,12) -- (38,10);
\fill[black!15](38,10) -- (40,12) -- (40,15) -- (38,15) -- (38,10);
% Draw the boundary of piece1
\draw[black](38,10) -- (52,10) -- (52,15) -- (50,15) -- (50,12) -- (40,12) -- (40,15) -- (38,15) -- (38,10);
% Draw the bounding box of the solution
\draw[black](0,0) -- (52,0) -- (52,40)-- (0,40) -- (0,0);
\draw (26,-2) node {SHAPES-9};
\draw (26,-5) node {L = 52};
\end{tikzpicture}   
% SOLUTION for problem instance = SHAPES-15
% Scaling factor for coordinates = 1
\begin{tikzpicture}[scale = 0.07]
% Draw the convex parts of piece2
\fill[black!15](0,25) -- (6,19) -- (12,25) -- (6,31) -- (0,25);
% Draw the boundary of piece2
\draw[black](0,25) -- (6,19) -- (12,25) -- (6,31) -- (0,25);
% Draw the convex parts of piece2
\fill[black!15](4,33) -- (10,27) -- (16,33) -- (10,39) -- (4,33);
% Draw the boundary of piece2
\draw[black](4,33) -- (10,27) -- (16,33) -- (10,39) -- (4,33);
% Draw the convex parts of piece2
\fill[black!15](11,24) -- (17,18) -- (23,24) -- (17,30) -- (11,24);
% Draw the boundary of piece2
\draw[black](11,24) -- (17,18) -- (23,24) -- (17,30) -- (11,24);
% Draw the convex parts of piece2
\fill[black!15](18,29) -- (24,23) -- (30,29) -- (24,35) -- (18,29);
% Draw the boundary of piece2
\draw[black](18,29) -- (24,23) -- (30,29) -- (24,35) -- (18,29);
% Draw the convex parts of piece2
\fill[black!15](26,33) -- (32,27) -- (38,33) -- (32,39) -- (26,33);
% Draw the boundary of piece2
\draw[black](26,33) -- (32,27) -- (38,33) -- (32,39) -- (26,33);
% Draw the convex parts of piece2
\fill[black!15](35,30) -- (41,24) -- (47,30) -- (41,36) -- (35,30);
% Draw the boundary of piece2
\draw[black](35,30) -- (41,24) -- (47,30) -- (41,36) -- (35,30);
% Draw the convex parts of piece2
\fill[black!15](52,9) -- (58,3) -- (64,9) -- (58,15) -- (52,9);
% Draw the boundary of piece2
\draw[black](52,9) -- (58,3) -- (64,9) -- (58,15) -- (52,9);
% Draw the convex parts of piece4
\fill[black!15](2,17) -- (2,15) -- (4,15) -- (4,17) -- (2,17);
\fill[black!15](4,19) -- (4,21) -- (2,21) -- (2,19) -- (4,19);
\fill[black!15](0,17) -- (2,17) -- (4,17) -- (6,17) -- (6,19) -- (4,19) -- (2,19) -- (0,19) -- (0,17);
% Draw the boundary of piece4
\draw[black](0,17) -- (2,17) -- (2,15) -- (4,15) -- (4,17) -- (6,17) -- (6,19) -- (4,19) -- (4,21) -- (2,21) -- (2,19) -- (0,19) -- (0,17);
% Draw the convex parts of piece4
\fill[black!15](2,35) -- (2,33) -- (4,33) -- (4,35) -- (2,35);
\fill[black!15](4,37) -- (4,39) -- (2,39) -- (2,37) -- (4,37);
\fill[black!15](0,35) -- (2,35) -- (4,35) -- (6,35) -- (6,37) -- (4,37) -- (2,37) -- (0,37) -- (0,35);
% Draw the boundary of piece4
\draw[black](0,35) -- (2,35) -- (2,33) -- (4,33) -- (4,35) -- (6,35) -- (6,37) -- (4,37) -- (4,39) -- (2,39) -- (2,37) -- (0,37) -- (0,35);
% Draw the convex parts of piece4
\fill[black!15](15,36) -- (15,34) -- (17,34) -- (17,36) -- (15,36);
\fill[black!15](17,38) -- (17,40) -- (15,40) -- (15,38) -- (17,38);
\fill[black!15](13,36) -- (15,36) -- (17,36) -- (19,36) -- (19,38) -- (17,38) -- (15,38) -- (13,38) -- (13,36);
% Draw the boundary of piece4
\draw[black](13,36) -- (15,36) -- (15,34) -- (17,34) -- (17,36) -- (19,36) -- (19,38) -- (17,38) -- (17,40) -- (15,40) -- (15,38) -- (13,38) -- (13,36);
% Draw the convex parts of piece4
\fill[black!15](19,34) -- (19,32) -- (21,32) -- (21,34) -- (19,34);
\fill[black!15](21,36) -- (21,38) -- (19,38) -- (19,36) -- (21,36);
\fill[black!15](17,34) -- (19,34) -- (21,34) -- (23,34) -- (23,36) -- (21,36) -- (19,36) -- (17,36) -- (17,34);
% Draw the boundary of piece4
\draw[black](17,34) -- (19,34) -- (19,32) -- (21,32) -- (21,34) -- (23,34) -- (23,36) -- (21,36) -- (21,38) -- (19,38) -- (19,36) -- (17,36) -- (17,34);
% Draw the convex parts of piece4
\fill[black!15](25,36) -- (25,34) -- (27,34) -- (27,36) -- (25,36);
\fill[black!15](27,38) -- (27,40) -- (25,40) -- (25,38) -- (27,38);
\fill[black!15](23,36) -- (25,36) -- (27,36) -- (29,36) -- (29,38) -- (27,38) -- (25,38) -- (23,38) -- (23,36);
% Draw the boundary of piece4
\draw[black](23,36) -- (25,36) -- (25,34) -- (27,34) -- (27,36) -- (29,36) -- (29,38) -- (27,38) -- (27,40) -- (25,40) -- (25,38) -- (23,38) -- (23,36);
% Draw the convex parts of piece4
\fill[black!15](27,21) -- (27,19) -- (29,19) -- (29,21) -- (27,21);
\fill[black!15](29,23) -- (29,25) -- (27,25) -- (27,23) -- (29,23);
\fill[black!15](25,21) -- (27,21) -- (29,21) -- (31,21) -- (31,23) -- (29,23) -- (27,23) -- (25,23) -- (25,21);
% Draw the boundary of piece4
\draw[black](25,21) -- (27,21) -- (27,19) -- (29,19) -- (29,21) -- (31,21) -- (31,23) -- (29,23) -- (29,25) -- (27,25) -- (27,23) -- (25,23) -- (25,21);
% Draw the convex parts of piece4
\fill[black!15](32,2) -- (32,0) -- (34,0) -- (34,2) -- (32,2);
\fill[black!15](34,4) -- (34,6) -- (32,6) -- (32,4) -- (34,4);
\fill[black!15](30,2) -- (32,2) -- (34,2) -- (36,2) -- (36,4) -- (34,4) -- (32,4) -- (30,4) -- (30,2);
% Draw the boundary of piece4
\draw[black](30,2) -- (32,2) -- (32,0) -- (34,0) -- (34,2) -- (36,2) -- (36,4) -- (34,4) -- (34,6) -- (32,6) -- (32,4) -- (30,4) -- (30,2);
% Draw the convex parts of piece4
\fill[black!15](37,36) -- (37,34) -- (39,34) -- (39,36) -- (37,36);
\fill[black!15](39,38) -- (39,40) -- (37,40) -- (37,38) -- (39,38);
\fill[black!15](35,36) -- (37,36) -- (39,36) -- (41,36) -- (41,38) -- (39,38) -- (37,38) -- (35,38) -- (35,36);
% Draw the boundary of piece4
\draw[black](35,36) -- (37,36) -- (37,34) -- (39,34) -- (39,36) -- (41,36) -- (41,38) -- (39,38) -- (39,40) -- (37,40) -- (37,38) -- (35,38) -- (35,36);
% Draw the convex parts of piece4
\fill[black!15](41,15) -- (41,13) -- (43,13) -- (43,15) -- (41,15);
\fill[black!15](43,17) -- (43,19) -- (41,19) -- (41,17) -- (43,17);
\fill[black!15](39,15) -- (41,15) -- (43,15) -- (45,15) -- (45,17) -- (43,17) -- (41,17) -- (39,17) -- (39,15);
% Draw the boundary of piece4
\draw[black](39,15) -- (41,15) -- (41,13) -- (43,13) -- (43,15) -- (45,15) -- (45,17) -- (43,17) -- (43,19) -- (41,19) -- (41,17) -- (39,17) -- (39,15);
% Draw the convex parts of piece4
\fill[black!15](44,36) -- (44,34) -- (46,34) -- (46,36) -- (44,36);
\fill[black!15](46,38) -- (46,40) -- (44,40) -- (44,38) -- (46,38);
\fill[black!15](42,36) -- (44,36) -- (46,36) -- (48,36) -- (48,38) -- (46,38) -- (44,38) -- (42,38) -- (42,36);
% Draw the boundary of piece4
\draw[black](42,36) -- (44,36) -- (44,34) -- (46,34) -- (46,36) -- (48,36) -- (48,38) -- (46,38) -- (46,40) -- (44,40) -- (44,38) -- (42,38) -- (42,36);
% Draw the convex parts of piece4
\fill[black!15](58,26) -- (58,24) -- (60,24) -- (60,26) -- (58,26);
\fill[black!15](60,28) -- (60,30) -- (58,30) -- (58,28) -- (60,28);
\fill[black!15](56,26) -- (58,26) -- (60,26) -- (62,26) -- (62,28) -- (60,28) -- (58,28) -- (56,28) -- (56,26);
% Draw the boundary of piece4
\draw[black](56,26) -- (58,26) -- (58,24) -- (60,24) -- (60,26) -- (62,26) -- (62,28) -- (60,28) -- (60,30) -- (58,30) -- (58,28) -- (56,28) -- (56,26);
% Draw the convex parts of piece4
\fill[black!15](60,30) -- (60,28) -- (62,28) -- (62,30) -- (60,30);
\fill[black!15](62,32) -- (62,34) -- (60,34) -- (60,32) -- (62,32);
\fill[black!15](58,30) -- (60,30) -- (62,30) -- (64,30) -- (64,32) -- (62,32) -- (60,32) -- (58,32) -- (58,30);
% Draw the boundary of piece4
\draw[black](58,30) -- (60,30) -- (60,28) -- (62,28) -- (62,30) -- (64,30) -- (64,32) -- (62,32) -- (62,34) -- (60,34) -- (60,32) -- (58,32) -- (58,30);
% Draw the convex parts of piece3
\fill[black!15](8,5) -- (8,2) -- (11,2) -- (8,5);
\fill[black!15](6,6) -- (0,6) -- (8,5) -- (6,6);
\fill[black!15](8,5) -- (11,2) -- (11,4) -- (7,8) -- (8,5);
\fill[black!15](8,5) -- (7,8) -- (6,8) -- (6,6) -- (8,5);
\fill[black!15](0,2) -- (2,2) -- (2,5) -- (0,2);
\fill[black!15](0,6) -- (0,2) -- (2,5) -- (0,6);
\fill[black!15](2,5) -- (8,5) -- (0,6) -- (2,5);
% Draw the boundary of piece3
\draw[black](0,2) -- (2,2) -- (2,5) -- (8,5) -- (8,2) -- (11,2) -- (11,4) -- (7,8) -- (6,8) -- (6,6) -- (0,6) -- (0,2);
% Draw the convex parts of piece3
\fill[black!15](8,14) -- (8,11) -- (11,11) -- (8,14);
\fill[black!15](6,15) -- (0,15) -- (8,14) -- (6,15);
\fill[black!15](8,14) -- (11,11) -- (11,13) -- (7,17) -- (8,14);
\fill[black!15](8,14) -- (7,17) -- (6,17) -- (6,15) -- (8,14);
\fill[black!15](0,11) -- (2,11) -- (2,14) -- (0,11);
\fill[black!15](0,15) -- (0,11) -- (2,14) -- (0,15);
\fill[black!15](2,14) -- (8,14) -- (0,15) -- (2,14);
% Draw the boundary of piece3
\draw[black](0,11) -- (2,11) -- (2,14) -- (8,14) -- (8,11) -- (11,11) -- (11,13) -- (7,17) -- (6,17) -- (6,15) -- (0,15) -- (0,11);
% Draw the convex parts of piece3
\fill[black!15](19,5) -- (19,2) -- (22,2) -- (19,5);
\fill[black!15](17,6) -- (11,6) -- (19,5) -- (17,6);
\fill[black!15](19,5) -- (22,2) -- (22,4) -- (18,8) -- (19,5);
\fill[black!15](19,5) -- (18,8) -- (17,8) -- (17,6) -- (19,5);
\fill[black!15](11,2) -- (13,2) -- (13,5) -- (11,2);
\fill[black!15](11,6) -- (11,2) -- (13,5) -- (11,6);
\fill[black!15](13,5) -- (19,5) -- (11,6) -- (13,5);
% Draw the boundary of piece3
\draw[black](11,2) -- (13,2) -- (13,5) -- (19,5) -- (19,2) -- (22,2) -- (22,4) -- (18,8) -- (17,8) -- (17,6) -- (11,6) -- (11,2);
% Draw the convex parts of piece3
\fill[black!15](19,13) -- (19,10) -- (22,10) -- (19,13);
\fill[black!15](17,14) -- (11,14) -- (19,13) -- (17,14);
\fill[black!15](19,13) -- (22,10) -- (22,12) -- (18,16) -- (19,13);
\fill[black!15](19,13) -- (18,16) -- (17,16) -- (17,14) -- (19,13);
\fill[black!15](11,10) -- (13,10) -- (13,13) -- (11,10);
\fill[black!15](11,14) -- (11,10) -- (13,13) -- (11,14);
\fill[black!15](13,13) -- (19,13) -- (11,14) -- (13,13);
% Draw the boundary of piece3
\draw[black](11,10) -- (13,10) -- (13,13) -- (19,13) -- (19,10) -- (22,10) -- (22,12) -- (18,16) -- (17,16) -- (17,14) -- (11,14) -- (11,10);
% Draw the convex parts of piece3
\fill[black!15](37,26) -- (37,23) -- (40,23) -- (37,26);
\fill[black!15](35,27) -- (29,27) -- (37,26) -- (35,27);
\fill[black!15](37,26) -- (40,23) -- (40,25) -- (36,29) -- (37,26);
\fill[black!15](37,26) -- (36,29) -- (35,29) -- (35,27) -- (37,26);
\fill[black!15](29,23) -- (31,23) -- (31,26) -- (29,23);
\fill[black!15](29,27) -- (29,23) -- (31,26) -- (29,27);
\fill[black!15](31,26) -- (37,26) -- (29,27) -- (31,26);
% Draw the boundary of piece3
\draw[black](29,23) -- (31,23) -- (31,26) -- (37,26) -- (37,23) -- (40,23) -- (40,25) -- (36,29) -- (35,29) -- (35,27) -- (29,27) -- (29,23);
% Draw the convex parts of piece3
\fill[black!15](54,5) -- (54,2) -- (57,2) -- (54,5);
\fill[black!15](52,6) -- (46,6) -- (54,5) -- (52,6);
\fill[black!15](54,5) -- (57,2) -- (57,4) -- (53,8) -- (54,5);
\fill[black!15](54,5) -- (53,8) -- (52,8) -- (52,6) -- (54,5);
\fill[black!15](46,2) -- (48,2) -- (48,5) -- (46,2);
\fill[black!15](46,6) -- (46,2) -- (48,5) -- (46,6);
\fill[black!15](48,5) -- (54,5) -- (46,6) -- (48,5);
% Draw the boundary of piece3
\draw[black](46,2) -- (48,2) -- (48,5) -- (54,5) -- (54,2) -- (57,2) -- (57,4) -- (53,8) -- (52,8) -- (52,6) -- (46,6) -- (46,2);
% Draw the convex parts of piece3
\fill[black!15](54,26) -- (54,23) -- (57,23) -- (54,26);
\fill[black!15](52,27) -- (46,27) -- (54,26) -- (52,27);
\fill[black!15](54,26) -- (57,23) -- (57,25) -- (53,29) -- (54,26);
\fill[black!15](54,26) -- (53,29) -- (52,29) -- (52,27) -- (54,26);
\fill[black!15](46,23) -- (48,23) -- (48,26) -- (46,23);
\fill[black!15](46,27) -- (46,23) -- (48,26) -- (46,27);
\fill[black!15](48,26) -- (54,26) -- (46,27) -- (48,26);
% Draw the boundary of piece3
\draw[black](46,23) -- (48,23) -- (48,26) -- (54,26) -- (54,23) -- (57,23) -- (57,25) -- (53,29) -- (52,29) -- (52,27) -- (46,27) -- (46,23);
% Draw the convex parts of piece3
\fill[black!15](55,31) -- (55,28) -- (58,28) -- (55,31);
\fill[black!15](53,32) -- (47,32) -- (55,31) -- (53,32);
\fill[black!15](55,31) -- (58,28) -- (58,30) -- (54,34) -- (55,31);
\fill[black!15](55,31) -- (54,34) -- (53,34) -- (53,32) -- (55,31);
\fill[black!15](47,28) -- (49,28) -- (49,31) -- (47,28);
\fill[black!15](47,32) -- (47,28) -- (49,31) -- (47,32);
\fill[black!15](49,31) -- (55,31) -- (47,32) -- (49,31);
% Draw the boundary of piece3
\draw[black](47,28) -- (49,28) -- (49,31) -- (55,31) -- (55,28) -- (58,28) -- (58,30) -- (54,34) -- (53,34) -- (53,32) -- (47,32) -- (47,28);
% Draw the convex parts of piece3
\fill[black!15](61,15) -- (61,12) -- (64,12) -- (61,15);
\fill[black!15](59,16) -- (53,16) -- (61,15) -- (59,16);
\fill[black!15](61,15) -- (64,12) -- (64,14) -- (60,18) -- (61,15);
\fill[black!15](61,15) -- (60,18) -- (59,18) -- (59,16) -- (61,15);
\fill[black!15](53,12) -- (55,12) -- (55,15) -- (53,12);
\fill[black!15](53,16) -- (53,12) -- (55,15) -- (53,16);
\fill[black!15](55,15) -- (61,15) -- (53,16) -- (55,15);
% Draw the boundary of piece3
\draw[black](53,12) -- (55,12) -- (55,15) -- (61,15) -- (61,12) -- (64,12) -- (64,14) -- (60,18) -- (59,18) -- (59,16) -- (53,16) -- (53,12);
% Draw the convex parts of piece1
\fill[black!15](14,2) -- (4,2) -- (16,0) -- (14,2);
\fill[black!15](16,0) -- (16,5) -- (14,5) -- (14,2) -- (16,0);
\fill[black!15](2,0) -- (16,0) -- (4,2) -- (2,0);
\fill[black!15](2,0) -- (4,2) -- (4,5) -- (2,5) -- (2,0);
% Draw the boundary of piece1
\draw[black](2,0) -- (16,0) -- (16,5) -- (14,5) -- (14,2) -- (4,2) -- (4,5) -- (2,5) -- (2,0);
% Draw the convex parts of piece1
\fill[black!15](17,10) -- (7,10) -- (19,8) -- (17,10);
\fill[black!15](19,8) -- (19,13) -- (17,13) -- (17,10) -- (19,8);
\fill[black!15](5,8) -- (19,8) -- (7,10) -- (5,8);
\fill[black!15](5,8) -- (7,10) -- (7,13) -- (5,13) -- (5,8);
% Draw the boundary of piece1
\draw[black](5,8) -- (19,8) -- (19,13) -- (17,13) -- (17,10) -- (7,10) -- (7,13) -- (5,13) -- (5,8);
% Draw the convex parts of piece1
\fill[black!15](21,18) -- (11,18) -- (23,16) -- (21,18);
\fill[black!15](23,16) -- (23,21) -- (21,21) -- (21,18) -- (23,16);
\fill[black!15](9,16) -- (23,16) -- (11,18) -- (9,16);
\fill[black!15](9,16) -- (11,18) -- (11,21) -- (9,21) -- (9,16);
% Draw the boundary of piece1
\draw[black](9,16) -- (23,16) -- (23,21) -- (21,21) -- (21,18) -- (11,18) -- (11,21) -- (9,21) -- (9,16);
% Draw the convex parts of piece1
\fill[black!15](28,2) -- (18,2) -- (30,0) -- (28,2);
\fill[black!15](30,0) -- (30,5) -- (28,5) -- (28,2) -- (30,0);
\fill[black!15](16,0) -- (30,0) -- (18,2) -- (16,0);
\fill[black!15](16,0) -- (18,2) -- (18,5) -- (16,5) -- (16,0);
% Draw the boundary of piece1
\draw[black](16,0) -- (30,0) -- (30,5) -- (28,5) -- (28,2) -- (18,2) -- (18,5) -- (16,5) -- (16,0);
% Draw the convex parts of piece1
\fill[black!15](35,8) -- (25,8) -- (37,6) -- (35,8);
\fill[black!15](37,6) -- (37,11) -- (35,11) -- (35,8) -- (37,6);
\fill[black!15](23,6) -- (37,6) -- (25,8) -- (23,6);
\fill[black!15](23,6) -- (25,8) -- (25,11) -- (23,11) -- (23,6);
% Draw the boundary of piece1
\draw[black](23,6) -- (37,6) -- (37,11) -- (35,11) -- (35,8) -- (25,8) -- (25,11) -- (23,11) -- (23,6);
% Draw the convex parts of piece1
\fill[black!15](35,13) -- (25,13) -- (37,11) -- (35,13);
\fill[black!15](37,11) -- (37,16) -- (35,16) -- (35,13) -- (37,11);
\fill[black!15](23,11) -- (37,11) -- (25,13) -- (23,11);
\fill[black!15](23,11) -- (25,13) -- (25,16) -- (23,16) -- (23,11);
% Draw the boundary of piece1
\draw[black](23,11) -- (37,11) -- (37,16) -- (35,16) -- (35,13) -- (25,13) -- (25,16) -- (23,16) -- (23,11);
% Draw the convex parts of piece1
\fill[black!15](35,18) -- (25,18) -- (37,16) -- (35,18);
\fill[black!15](37,16) -- (37,21) -- (35,21) -- (35,18) -- (37,16);
\fill[black!15](23,16) -- (37,16) -- (25,18) -- (23,16);
\fill[black!15](23,16) -- (25,18) -- (25,21) -- (23,21) -- (23,16);
% Draw the boundary of piece1
\draw[black](23,16) -- (37,16) -- (37,21) -- (35,21) -- (35,18) -- (25,18) -- (25,21) -- (23,21) -- (23,16);
% Draw the convex parts of piece1
\fill[black!15](43,23) -- (33,23) -- (45,21) -- (43,23);
\fill[black!15](45,21) -- (45,26) -- (43,26) -- (43,23) -- (45,21);
\fill[black!15](31,21) -- (45,21) -- (33,23) -- (31,21);
\fill[black!15](31,21) -- (33,23) -- (33,26) -- (31,26) -- (31,21);
% Draw the boundary of piece1
\draw[black](31,21) -- (45,21) -- (45,26) -- (43,26) -- (43,23) -- (33,23) -- (33,26) -- (31,26) -- (31,21);
% Draw the convex parts of piece1
\fill[black!15](48,2) -- (38,2) -- (50,0) -- (48,2);
\fill[black!15](50,0) -- (50,5) -- (48,5) -- (48,2) -- (50,0);
\fill[black!15](36,0) -- (50,0) -- (38,2) -- (36,0);
\fill[black!15](36,0) -- (38,2) -- (38,5) -- (36,5) -- (36,0);
% Draw the boundary of piece1
\draw[black](36,0) -- (50,0) -- (50,5) -- (48,5) -- (48,2) -- (38,2) -- (38,5) -- (36,5) -- (36,0);
% Draw the convex parts of piece1
\fill[black!15](49,8) -- (39,8) -- (51,6) -- (49,8);
\fill[black!15](51,6) -- (51,11) -- (49,11) -- (49,8) -- (51,6);
\fill[black!15](37,6) -- (51,6) -- (39,8) -- (37,6);
\fill[black!15](37,6) -- (39,8) -- (39,11) -- (37,11) -- (37,6);
% Draw the boundary of piece1
\draw[black](37,6) -- (51,6) -- (51,11) -- (49,11) -- (49,8) -- (39,8) -- (39,11) -- (37,11) -- (37,6);
% Draw the convex parts of piece1
\fill[black!15](49,13) -- (39,13) -- (51,11) -- (49,13);
\fill[black!15](51,11) -- (51,16) -- (49,16) -- (49,13) -- (51,11);
\fill[black!15](37,11) -- (51,11) -- (39,13) -- (37,11);
\fill[black!15](37,11) -- (39,13) -- (39,16) -- (37,16) -- (37,11);
% Draw the boundary of piece1
\draw[black](37,11) -- (51,11) -- (51,16) -- (49,16) -- (49,13) -- (39,13) -- (39,16) -- (37,16) -- (37,11);
% Draw the convex parts of piece1
\fill[black!15](57,18) -- (47,18) -- (59,16) -- (57,18);
\fill[black!15](59,16) -- (59,21) -- (57,21) -- (57,18) -- (59,16);
\fill[black!15](45,16) -- (59,16) -- (47,18) -- (45,16);
\fill[black!15](45,16) -- (47,18) -- (47,21) -- (45,21) -- (45,16);
% Draw the boundary of piece1
\draw[black](45,16) -- (59,16) -- (59,21) -- (57,21) -- (57,18) -- (47,18) -- (47,21) -- (45,21) -- (45,16);
% Draw the convex parts of piece1
\fill[black!15](60,37) -- (50,37) -- (62,35) -- (60,37);
\fill[black!15](62,35) -- (62,40) -- (60,40) -- (60,37) -- (62,35);
\fill[black!15](48,35) -- (62,35) -- (50,37) -- (48,35);
\fill[black!15](48,35) -- (50,37) -- (50,40) -- (48,40) -- (48,35);
% Draw the boundary of piece1
\draw[black](48,35) -- (62,35) -- (62,40) -- (60,40) -- (60,37) -- (50,37) -- (50,40) -- (48,40) -- (48,35);
% Draw the convex parts of piece1
\fill[black!15](62,2) -- (52,2) -- (64,0) -- (62,2);
\fill[black!15](64,0) -- (64,5) -- (62,5) -- (62,2) -- (64,0);
\fill[black!15](50,0) -- (64,0) -- (52,2) -- (50,0);
\fill[black!15](50,0) -- (52,2) -- (52,5) -- (50,5) -- (50,0);
% Draw the boundary of piece1
\draw[black](50,0) -- (64,0) -- (64,5) -- (62,5) -- (62,2) -- (52,2) -- (52,5) -- (50,5) -- (50,0);
% Draw the convex parts of piece1
\fill[black!15](62,23) -- (52,23) -- (64,21) -- (62,23);
\fill[black!15](64,21) -- (64,26) -- (62,26) -- (62,23) -- (64,21);
\fill[black!15](50,21) -- (64,21) -- (52,23) -- (50,21);
\fill[black!15](50,21) -- (52,23) -- (52,26) -- (50,26) -- (50,21);
% Draw the boundary of piece1
\draw[black](50,21) -- (64,21) -- (64,26) -- (62,26) -- (62,23) -- (52,23) -- (52,26) -- (50,26) -- (50,21);
% Draw the bounding box of the solution
\draw[black](0,0) -- (64,0) -- (64,40)-- (0,40) -- (0,0);
\draw (32,-2) node {SHAPES-15};
\draw (32,-5) node {L = 64};
\end{tikzpicture}
\\
% SOLUTION for problem instance = artif5_10
% Scaling factor for coordinates = 1
\begin{tikzpicture}[scale = 0.10]
% Draw the convex parts of piece3
\fill[black!15](4,0) -- (7,2) -- (6,4) -- (3,5) -- (0,3) -- (4,0);
% Draw the boundary of piece3
\draw[black](4,0) -- (7,2) -- (6,4) -- (3,5) -- (0,3) -- (4,0);
% Draw the convex parts of piece3
\fill[black!15](4,15) -- (7,17) -- (6,19) -- (3,20) -- (0,18) -- (4,15);
% Draw the boundary of piece3
\draw[black](4,15) -- (7,17) -- (6,19) -- (3,20) -- (0,18) -- (4,15);
% Draw the convex parts of piece3
\fill[black!15](4,20) -- (7,22) -- (6,24) -- (3,25) -- (0,23) -- (4,20);
% Draw the boundary of piece3
\draw[black](4,20) -- (7,22) -- (6,24) -- (3,25) -- (0,23) -- (4,20);
% Draw the convex parts of piece3
\fill[black!15](16,22) -- (19,24) -- (18,26) -- (15,27) -- (12,25) -- (16,22);
% Draw the boundary of piece3
\draw[black](16,22) -- (19,24) -- (18,26) -- (15,27) -- (12,25) -- (16,22);
% Draw the convex parts of piece3
\fill[black!15](25,14) -- (28,16) -- (27,18) -- (24,19) -- (21,17) -- (25,14);
% Draw the boundary of piece3
\draw[black](25,14) -- (28,16) -- (27,18) -- (24,19) -- (21,17) -- (25,14);
% Draw the convex parts of piece2
\fill[black!15](10,3) -- (8,3) -- (11,0) -- (10,3);
\fill[black!15](11,0) -- (13,3) -- (12,5) -- (10,3) -- (11,0);
\fill[black!15](7,0) -- (11,0) -- (8,3) -- (7,4) -- (7,0);
% Draw the boundary of piece2
\draw[black](7,0) -- (11,0) -- (13,3) -- (12,5) -- (10,3) -- (8,3) -- (7,4) -- (7,0);
% Draw the convex parts of piece2
\fill[black!15](10,7) -- (8,7) -- (11,4) -- (10,7);
\fill[black!15](11,4) -- (13,7) -- (12,9) -- (10,7) -- (11,4);
\fill[black!15](7,4) -- (11,4) -- (8,7) -- (7,8) -- (7,4);
% Draw the boundary of piece2
\draw[black](7,4) -- (11,4) -- (13,7) -- (12,9) -- (10,7) -- (8,7) -- (7,8) -- (7,4);
% Draw the convex parts of piece2
\fill[black!15](10,11) -- (8,11) -- (11,8) -- (10,11);
\fill[black!15](11,8) -- (13,11) -- (12,13) -- (10,11) -- (11,8);
\fill[black!15](7,8) -- (11,8) -- (8,11) -- (7,12) -- (7,8);
% Draw the boundary of piece2
\draw[black](7,8) -- (11,8) -- (13,11) -- (12,13) -- (10,11) -- (8,11) -- (7,12) -- (7,8);
% Draw the convex parts of piece2
\fill[black!15](10,15) -- (8,15) -- (11,12) -- (10,15);
\fill[black!15](11,12) -- (13,15) -- (12,17) -- (10,15) -- (11,12);
\fill[black!15](7,12) -- (11,12) -- (8,15) -- (7,16) -- (7,12);
% Draw the boundary of piece2
\draw[black](7,12) -- (11,12) -- (13,15) -- (12,17) -- (10,15) -- (8,15) -- (7,16) -- (7,12);
% Draw the convex parts of piece2
\fill[black!15](10,19) -- (8,19) -- (11,16) -- (10,19);
\fill[black!15](11,16) -- (13,19) -- (12,21) -- (10,19) -- (11,16);
\fill[black!15](7,16) -- (11,16) -- (8,19) -- (7,20) -- (7,16);
% Draw the boundary of piece2
\draw[black](7,16) -- (11,16) -- (13,19) -- (12,21) -- (10,19) -- (8,19) -- (7,20) -- (7,16);
% Draw the convex parts of piece2
\fill[black!15](10,23) -- (8,23) -- (11,20) -- (10,23);
\fill[black!15](11,20) -- (13,23) -- (12,25) -- (10,23) -- (11,20);
\fill[black!15](7,20) -- (11,20) -- (8,23) -- (7,24) -- (7,20);
% Draw the boundary of piece2
\draw[black](7,20) -- (11,20) -- (13,23) -- (12,25) -- (10,23) -- (8,23) -- (7,24) -- (7,20);
% Draw the convex parts of piece2
\fill[black!15](16,3) -- (14,3) -- (17,0) -- (16,3);
\fill[black!15](17,0) -- (19,3) -- (18,5) -- (16,3) -- (17,0);
\fill[black!15](13,0) -- (17,0) -- (14,3) -- (13,4) -- (13,0);
% Draw the boundary of piece2
\draw[black](13,0) -- (17,0) -- (19,3) -- (18,5) -- (16,3) -- (14,3) -- (13,4) -- (13,0);
% Draw the convex parts of piece2
\fill[black!15](16,7) -- (14,7) -- (17,4) -- (16,7);
\fill[black!15](17,4) -- (19,7) -- (18,9) -- (16,7) -- (17,4);
\fill[black!15](13,4) -- (17,4) -- (14,7) -- (13,8) -- (13,4);
% Draw the boundary of piece2
\draw[black](13,4) -- (17,4) -- (19,7) -- (18,9) -- (16,7) -- (14,7) -- (13,8) -- (13,4);
% Draw the convex parts of piece2
\fill[black!15](16,11) -- (14,11) -- (17,8) -- (16,11);
\fill[black!15](17,8) -- (19,11) -- (18,13) -- (16,11) -- (17,8);
\fill[black!15](13,8) -- (17,8) -- (14,11) -- (13,12) -- (13,8);
% Draw the boundary of piece2
\draw[black](13,8) -- (17,8) -- (19,11) -- (18,13) -- (16,11) -- (14,11) -- (13,12) -- (13,8);
% Draw the convex parts of piece2
\fill[black!15](16,15) -- (14,15) -- (17,12) -- (16,15);
\fill[black!15](17,12) -- (19,15) -- (18,17) -- (16,15) -- (17,12);
\fill[black!15](13,12) -- (17,12) -- (14,15) -- (13,16) -- (13,12);
% Draw the boundary of piece2
\draw[black](13,12) -- (17,12) -- (19,15) -- (18,17) -- (16,15) -- (14,15) -- (13,16) -- (13,12);
% Draw the convex parts of piece7
\fill[black!15](17,19) -- (15,18) -- (16,17) -- (17,19);
\fill[black!15](17,19) -- (16,17) -- (18,18) -- (18,19) -- (17,19);
\fill[black!15](18,19) -- (19,20) -- (17,20) -- (17,19) -- (18,19);
\fill[black!15](14,17) -- (13,16) -- (16,16) -- (14,17);
\fill[black!15](14,17) -- (16,16) -- (16,17) -- (15,18) -- (14,17);
\fill[black!15](15,18) -- (15,20) -- (13,20) -- (14,17) -- (15,18);
% Draw the boundary of piece7
\draw[black](13,16) -- (16,16) -- (16,17) -- (18,18) -- (18,19) -- (19,20) -- (17,20) -- (17,19) -- (15,18) -- (15,20) -- (13,20) -- (14,17) -- (13,16);
% Draw the convex parts of piece7
\fill[black!15](22,22) -- (20,21) -- (21,20) -- (22,22);
\fill[black!15](22,22) -- (21,20) -- (23,21) -- (23,22) -- (22,22);
\fill[black!15](23,22) -- (24,23) -- (22,23) -- (22,22) -- (23,22);
\fill[black!15](19,20) -- (18,19) -- (21,19) -- (19,20);
\fill[black!15](19,20) -- (21,19) -- (21,20) -- (20,21) -- (19,20);
\fill[black!15](20,21) -- (20,23) -- (18,23) -- (19,20) -- (20,21);
% Draw the boundary of piece7
\draw[black](18,19) -- (21,19) -- (21,20) -- (23,21) -- (23,22) -- (24,23) -- (22,23) -- (22,22) -- (20,21) -- (20,23) -- (18,23) -- (19,20) -- (18,19);
% Draw the convex parts of piece7
\fill[black!15](23,3) -- (21,2) -- (22,1) -- (23,3);
\fill[black!15](23,3) -- (22,1) -- (24,2) -- (24,3) -- (23,3);
\fill[black!15](24,3) -- (25,4) -- (23,4) -- (23,3) -- (24,3);
\fill[black!15](20,1) -- (19,0) -- (22,0) -- (20,1);
\fill[black!15](20,1) -- (22,0) -- (22,1) -- (21,2) -- (20,1);
\fill[black!15](21,2) -- (21,4) -- (19,4) -- (20,1) -- (21,2);
% Draw the boundary of piece7
\draw[black](19,0) -- (22,0) -- (22,1) -- (24,2) -- (24,3) -- (25,4) -- (23,4) -- (23,3) -- (21,2) -- (21,4) -- (19,4) -- (20,1) -- (19,0);
% Draw the convex parts of piece7
\fill[black!15](25,26) -- (23,25) -- (24,24) -- (25,26);
\fill[black!15](25,26) -- (24,24) -- (26,25) -- (26,26) -- (25,26);
\fill[black!15](26,26) -- (27,27) -- (25,27) -- (25,26) -- (26,26);
\fill[black!15](22,24) -- (21,23) -- (24,23) -- (22,24);
\fill[black!15](22,24) -- (24,23) -- (24,24) -- (23,25) -- (22,24);
\fill[black!15](23,25) -- (23,27) -- (21,27) -- (22,24) -- (23,25);
% Draw the boundary of piece7
\draw[black](21,23) -- (24,23) -- (24,24) -- (26,25) -- (26,26) -- (27,27) -- (25,27) -- (25,26) -- (23,25) -- (23,27) -- (21,27) -- (22,24) -- (21,23);
% Draw the convex parts of piece7
\fill[black!15](31,12) -- (29,11) -- (30,10) -- (31,12);
\fill[black!15](31,12) -- (30,10) -- (32,11) -- (32,12) -- (31,12);
\fill[black!15](32,12) -- (33,13) -- (31,13) -- (31,12) -- (32,12);
\fill[black!15](28,10) -- (27,9) -- (30,9) -- (28,10);
\fill[black!15](28,10) -- (30,9) -- (30,10) -- (29,11) -- (28,10);
\fill[black!15](29,11) -- (29,13) -- (27,13) -- (28,10) -- (29,11);
% Draw the boundary of piece7
\draw[black](27,9) -- (30,9) -- (30,10) -- (32,11) -- (32,12) -- (33,13) -- (31,13) -- (31,12) -- (29,11) -- (29,13) -- (27,13) -- (28,10) -- (27,9);
% Draw the convex parts of piece5
\fill[black!15](0,5) -- (3,5) -- (4,6) -- (1,7) -- (0,5);
\fill[black!15](0,5) -- (1,7) -- (1,8) -- (0,5);
\fill[black!15](0,10) -- (0,5) -- (1,8) -- (0,10);
\fill[black!15](2,10) -- (0,10) -- (1,8) -- (2,9) -- (2,10);
\fill[black!15](2,9) -- (3,9) -- (2,10) -- (2,9);
% Draw the boundary of piece5
\draw[black](0,5) -- (3,5) -- (4,6) -- (1,7) -- (1,8) -- (2,9) -- (3,9) -- (2,10) -- (0,10) -- (0,5);
% Draw the convex parts of piece5
\fill[black!15](0,10) -- (3,10) -- (4,11) -- (1,12) -- (0,10);
\fill[black!15](0,10) -- (1,12) -- (1,13) -- (0,10);
\fill[black!15](0,15) -- (0,10) -- (1,13) -- (0,15);
\fill[black!15](2,15) -- (0,15) -- (1,13) -- (2,14) -- (2,15);
\fill[black!15](2,14) -- (3,14) -- (2,15) -- (2,14);
% Draw the boundary of piece5
\draw[black](0,10) -- (3,10) -- (4,11) -- (1,12) -- (1,13) -- (2,14) -- (3,14) -- (2,15) -- (0,15) -- (0,10);
% Draw the convex parts of piece5
\fill[black!15](19,4) -- (22,4) -- (23,5) -- (20,6) -- (19,4);
\fill[black!15](19,4) -- (20,6) -- (20,7) -- (19,4);
\fill[black!15](19,9) -- (19,4) -- (20,7) -- (19,9);
\fill[black!15](21,9) -- (19,9) -- (20,7) -- (21,8) -- (21,9);
\fill[black!15](21,8) -- (22,8) -- (21,9) -- (21,8);
% Draw the boundary of piece5
\draw[black](19,4) -- (22,4) -- (23,5) -- (20,6) -- (20,7) -- (21,8) -- (22,8) -- (21,9) -- (19,9) -- (19,4);
% Draw the convex parts of piece5
\fill[black!15](19,9) -- (22,9) -- (23,10) -- (20,11) -- (19,9);
\fill[black!15](19,9) -- (20,11) -- (20,12) -- (19,9);
\fill[black!15](19,14) -- (19,9) -- (20,12) -- (19,14);
\fill[black!15](21,14) -- (19,14) -- (20,12) -- (21,13) -- (21,14);
\fill[black!15](21,13) -- (22,13) -- (21,14) -- (21,13);
% Draw the boundary of piece5
\draw[black](19,9) -- (22,9) -- (23,10) -- (20,11) -- (20,12) -- (21,13) -- (22,13) -- (21,14) -- (19,14) -- (19,9);
% Draw the convex parts of piece5
\fill[black!15](19,14) -- (22,14) -- (23,15) -- (20,16) -- (19,14);
\fill[black!15](19,14) -- (20,16) -- (20,17) -- (19,14);
\fill[black!15](19,19) -- (19,14) -- (20,17) -- (19,19);
\fill[black!15](21,19) -- (19,19) -- (20,17) -- (21,18) -- (21,19);
\fill[black!15](21,18) -- (22,18) -- (21,19) -- (21,18);
% Draw the boundary of piece5
\draw[black](19,14) -- (22,14) -- (23,15) -- (20,16) -- (20,17) -- (21,18) -- (22,18) -- (21,19) -- (19,19) -- (19,14);
% Draw the convex parts of piece5
\fill[black!15](23,4) -- (26,4) -- (27,5) -- (24,6) -- (23,4);
\fill[black!15](23,4) -- (24,6) -- (24,7) -- (23,4);
\fill[black!15](23,9) -- (23,4) -- (24,7) -- (23,9);
\fill[black!15](25,9) -- (23,9) -- (24,7) -- (25,8) -- (25,9);
\fill[black!15](25,8) -- (26,8) -- (25,9) -- (25,8);
% Draw the boundary of piece5
\draw[black](23,4) -- (26,4) -- (27,5) -- (24,6) -- (24,7) -- (25,8) -- (26,8) -- (25,9) -- (23,9) -- (23,4);
% Draw the convex parts of piece5
\fill[black!15](23,9) -- (26,9) -- (27,10) -- (24,11) -- (23,9);
\fill[black!15](23,9) -- (24,11) -- (24,12) -- (23,9);
\fill[black!15](23,14) -- (23,9) -- (24,12) -- (23,14);
\fill[black!15](25,14) -- (23,14) -- (24,12) -- (25,13) -- (25,14);
\fill[black!15](25,13) -- (26,13) -- (25,14) -- (25,13);
% Draw the boundary of piece5
\draw[black](23,9) -- (26,9) -- (27,10) -- (24,11) -- (24,12) -- (25,13) -- (26,13) -- (25,14) -- (23,14) -- (23,9);
% Draw the convex parts of piece5
\fill[black!15](24,19) -- (27,19) -- (28,20) -- (25,21) -- (24,19);
\fill[black!15](24,19) -- (25,21) -- (25,22) -- (24,19);
\fill[black!15](24,24) -- (24,19) -- (25,22) -- (24,24);
\fill[black!15](26,24) -- (24,24) -- (25,22) -- (26,23) -- (26,24);
\fill[black!15](26,23) -- (27,23) -- (26,24) -- (26,23);
% Draw the boundary of piece5
\draw[black](24,19) -- (27,19) -- (28,20) -- (25,21) -- (25,22) -- (26,23) -- (27,23) -- (26,24) -- (24,24) -- (24,19);
% Draw the convex parts of piece5
\fill[black!15](28,15) -- (31,15) -- (32,16) -- (29,17) -- (28,15);
\fill[black!15](28,15) -- (29,17) -- (29,18) -- (28,15);
\fill[black!15](28,20) -- (28,15) -- (29,18) -- (28,20);
\fill[black!15](30,20) -- (28,20) -- (29,18) -- (30,19) -- (30,20);
\fill[black!15](30,19) -- (31,19) -- (30,20) -- (30,19);
% Draw the boundary of piece5
\draw[black](28,15) -- (31,15) -- (32,16) -- (29,17) -- (29,18) -- (30,19) -- (31,19) -- (30,20) -- (28,20) -- (28,15);
% Draw the convex parts of piece5
\fill[black!15](29,0) -- (32,0) -- (33,1) -- (30,2) -- (29,0);
\fill[black!15](29,0) -- (30,2) -- (30,3) -- (29,0);
\fill[black!15](29,5) -- (29,0) -- (30,3) -- (29,5);
\fill[black!15](31,5) -- (29,5) -- (30,3) -- (31,4) -- (31,5);
\fill[black!15](31,4) -- (32,4) -- (31,5) -- (31,4);
% Draw the boundary of piece5
\draw[black](29,0) -- (32,0) -- (33,1) -- (30,2) -- (30,3) -- (31,4) -- (32,4) -- (31,5) -- (29,5) -- (29,0);
% Draw the convex parts of piece8
\fill[black!15](5,9) -- (4,9) -- (6,8) -- (5,9);
\fill[black!15](6,10) -- (5,9) -- (6,8) -- (6,10);
\fill[black!15](6,8) -- (7,11) -- (6,11) -- (6,10) -- (6,8);
\fill[black!15](4,6) -- (6,6) -- (5,7) -- (4,6);
\fill[black!15](4,6) -- (5,7) -- (5,8) -- (4,6);
\fill[black!15](5,8) -- (6,8) -- (4,9) -- (5,8);
\fill[black!15](4,6) -- (5,8) -- (4,9) -- (4,6);
\fill[black!15](4,6) -- (4,9) -- (4,10) -- (3,10) -- (4,6);
% Draw the boundary of piece8
\draw[black](4,6) -- (6,6) -- (5,7) -- (5,8) -- (6,8) -- (7,11) -- (6,11) -- (6,10) -- (5,9) -- (4,9) -- (4,10) -- (3,10) -- (4,6);
% Draw the convex parts of piece8
\fill[black!15](5,14) -- (4,14) -- (6,13) -- (5,14);
\fill[black!15](6,15) -- (5,14) -- (6,13) -- (6,15);
\fill[black!15](6,13) -- (7,16) -- (6,16) -- (6,15) -- (6,13);
\fill[black!15](4,11) -- (6,11) -- (5,12) -- (4,11);
\fill[black!15](4,11) -- (5,12) -- (5,13) -- (4,11);
\fill[black!15](5,13) -- (6,13) -- (4,14) -- (5,13);
\fill[black!15](4,11) -- (5,13) -- (4,14) -- (4,11);
\fill[black!15](4,11) -- (4,14) -- (4,15) -- (3,15) -- (4,11);
% Draw the boundary of piece8
\draw[black](4,11) -- (6,11) -- (5,12) -- (5,13) -- (6,13) -- (7,16) -- (6,16) -- (6,15) -- (5,14) -- (4,14) -- (4,15) -- (3,15) -- (4,11);
% Draw the convex parts of piece8
\fill[black!15](27,3) -- (26,3) -- (28,2) -- (27,3);
\fill[black!15](28,4) -- (27,3) -- (28,2) -- (28,4);
\fill[black!15](28,2) -- (29,5) -- (28,5) -- (28,4) -- (28,2);
\fill[black!15](26,0) -- (28,0) -- (27,1) -- (26,0);
\fill[black!15](26,0) -- (27,1) -- (27,2) -- (26,0);
\fill[black!15](27,2) -- (28,2) -- (26,3) -- (27,2);
\fill[black!15](26,0) -- (27,2) -- (26,3) -- (26,0);
\fill[black!15](26,0) -- (26,3) -- (26,4) -- (25,4) -- (26,0);
% Draw the boundary of piece8
\draw[black](26,0) -- (28,0) -- (27,1) -- (27,2) -- (28,2) -- (29,5) -- (28,5) -- (28,4) -- (27,3) -- (26,3) -- (26,4) -- (25,4) -- (26,0);
% Draw the convex parts of piece8
\fill[black!15](29,8) -- (28,8) -- (30,7) -- (29,8);
\fill[black!15](30,9) -- (29,8) -- (30,7) -- (30,9);
\fill[black!15](30,7) -- (31,10) -- (30,10) -- (30,9) -- (30,7);
\fill[black!15](28,5) -- (30,5) -- (29,6) -- (28,5);
\fill[black!15](28,5) -- (29,6) -- (29,7) -- (28,5);
\fill[black!15](29,7) -- (30,7) -- (28,8) -- (29,7);
\fill[black!15](28,5) -- (29,7) -- (28,8) -- (28,5);
\fill[black!15](28,5) -- (28,8) -- (28,9) -- (27,9) -- (28,5);
% Draw the boundary of piece8
\draw[black](28,5) -- (30,5) -- (29,6) -- (29,7) -- (30,7) -- (31,10) -- (30,10) -- (30,9) -- (29,8) -- (28,8) -- (28,9) -- (27,9) -- (28,5);
% Draw the convex parts of piece8
\fill[black!15](29,23) -- (28,23) -- (30,22) -- (29,23);
\fill[black!15](30,24) -- (29,23) -- (30,22) -- (30,24);
\fill[black!15](30,22) -- (31,25) -- (30,25) -- (30,24) -- (30,22);
\fill[black!15](28,20) -- (30,20) -- (29,21) -- (28,20);
\fill[black!15](28,20) -- (29,21) -- (29,22) -- (28,20);
\fill[black!15](29,22) -- (30,22) -- (28,23) -- (29,22);
\fill[black!15](28,20) -- (29,22) -- (28,23) -- (28,20);
\fill[black!15](28,20) -- (28,23) -- (28,24) -- (27,24) -- (28,20);
% Draw the boundary of piece8
\draw[black](28,20) -- (30,20) -- (29,21) -- (29,22) -- (30,22) -- (31,25) -- (30,25) -- (30,24) -- (29,23) -- (28,23) -- (28,24) -- (27,24) -- (28,20);
% Draw the convex parts of piece4
\fill[black!15](2,26) -- (1,26) -- (3,25) -- (2,26);
\fill[black!15](3,25) -- (3,27) -- (2,27) -- (2,26) -- (3,25);
\fill[black!15](0,25) -- (3,25) -- (1,26) -- (0,25);
\fill[black!15](0,25) -- (1,26) -- (1,27) -- (0,27) -- (0,25);
% Draw the boundary of piece4
\draw[black](0,25) -- (3,25) -- (3,27) -- (2,27) -- (2,26) -- (1,26) -- (1,27) -- (0,27) -- (0,25);
% Draw the convex parts of piece4
\fill[black!15](5,26) -- (4,26) -- (6,25) -- (5,26);
\fill[black!15](6,25) -- (6,27) -- (5,27) -- (5,26) -- (6,25);
\fill[black!15](3,25) -- (6,25) -- (4,26) -- (3,25);
\fill[black!15](3,25) -- (4,26) -- (4,27) -- (3,27) -- (3,25);
% Draw the boundary of piece4
\draw[black](3,25) -- (6,25) -- (6,27) -- (5,27) -- (5,26) -- (4,26) -- (4,27) -- (3,27) -- (3,25);
% Draw the convex parts of piece4
\fill[black!15](8,26) -- (7,26) -- (9,25) -- (8,26);
\fill[black!15](9,25) -- (9,27) -- (8,27) -- (8,26) -- (9,25);
\fill[black!15](6,25) -- (9,25) -- (7,26) -- (6,25);
\fill[black!15](6,25) -- (7,26) -- (7,27) -- (6,27) -- (6,25);
% Draw the boundary of piece4
\draw[black](6,25) -- (9,25) -- (9,27) -- (8,27) -- (8,26) -- (7,26) -- (7,27) -- (6,27) -- (6,25);
% Draw the convex parts of piece4
\fill[black!15](11,26) -- (10,26) -- (12,25) -- (11,26);
\fill[black!15](12,25) -- (12,27) -- (11,27) -- (11,26) -- (12,25);
\fill[black!15](9,25) -- (12,25) -- (10,26) -- (9,25);
\fill[black!15](9,25) -- (10,26) -- (10,27) -- (9,27) -- (9,25);
% Draw the boundary of piece4
\draw[black](9,25) -- (12,25) -- (12,27) -- (11,27) -- (11,26) -- (10,26) -- (10,27) -- (9,27) -- (9,25);
% Draw the convex parts of piece4
\fill[black!15](15,21) -- (14,21) -- (16,20) -- (15,21);
\fill[black!15](16,20) -- (16,22) -- (15,22) -- (15,21) -- (16,20);
\fill[black!15](13,20) -- (16,20) -- (14,21) -- (13,20);
\fill[black!15](13,20) -- (14,21) -- (14,22) -- (13,22) -- (13,20);
% Draw the boundary of piece4
\draw[black](13,20) -- (16,20) -- (16,22) -- (15,22) -- (15,21) -- (14,21) -- (14,22) -- (13,22) -- (13,20);
% Draw the convex parts of piece4
\fill[black!15](28,25) -- (27,25) -- (29,24) -- (28,25);
\fill[black!15](29,24) -- (29,26) -- (28,26) -- (28,25) -- (29,24);
\fill[black!15](26,24) -- (29,24) -- (27,25) -- (26,24);
\fill[black!15](26,24) -- (27,25) -- (27,26) -- (26,26) -- (26,24);
% Draw the boundary of piece4
\draw[black](26,24) -- (29,24) -- (29,26) -- (28,26) -- (28,25) -- (27,25) -- (27,26) -- (26,26) -- (26,24);
% Draw the convex parts of piece4
\fill[black!15](29,14) -- (28,14) -- (30,13) -- (29,14);
\fill[black!15](30,13) -- (30,15) -- (29,15) -- (29,14) -- (30,13);
\fill[black!15](27,13) -- (30,13) -- (28,14) -- (27,13);
\fill[black!15](27,13) -- (28,14) -- (28,15) -- (27,15) -- (27,13);
% Draw the boundary of piece4
\draw[black](27,13) -- (30,13) -- (30,15) -- (29,15) -- (29,14) -- (28,14) -- (28,15) -- (27,15) -- (27,13);
% Draw the convex parts of piece4
\fill[black!15](31,26) -- (30,26) -- (32,25) -- (31,26);
\fill[black!15](32,25) -- (32,27) -- (31,27) -- (31,26) -- (32,25);
\fill[black!15](29,25) -- (32,25) -- (30,26) -- (29,25);
\fill[black!15](29,25) -- (30,26) -- (30,27) -- (29,27) -- (29,25);
% Draw the boundary of piece4
\draw[black](29,25) -- (32,25) -- (32,27) -- (31,27) -- (31,26) -- (30,26) -- (30,27) -- (29,27) -- (29,25);
% Draw the convex parts of piece4
\fill[black!15](32,6) -- (31,6) -- (33,5) -- (32,6);
\fill[black!15](33,5) -- (33,7) -- (32,7) -- (32,6) -- (33,5);
\fill[black!15](30,5) -- (33,5) -- (31,6) -- (30,5);
\fill[black!15](30,5) -- (31,6) -- (31,7) -- (30,7) -- (30,5);
% Draw the boundary of piece4
\draw[black](30,5) -- (33,5) -- (33,7) -- (32,7) -- (32,6) -- (31,6) -- (31,7) -- (30,7) -- (30,5);
% Draw the convex parts of piece4
\fill[black!15](32,14) -- (31,14) -- (33,13) -- (32,14);
\fill[black!15](33,13) -- (33,15) -- (32,15) -- (32,14) -- (33,13);
\fill[black!15](30,13) -- (33,13) -- (31,14) -- (30,13);
\fill[black!15](30,13) -- (31,14) -- (31,15) -- (30,15) -- (30,13);
% Draw the boundary of piece4
\draw[black](30,13) -- (33,13) -- (33,15) -- (32,15) -- (32,14) -- (31,14) -- (31,15) -- (30,15) -- (30,13);
% Draw the convex parts of piece1
\fill[black!15](8,23) -- (10,23) -- (10,25) -- (8,25) -- (8,23);
% Draw the boundary of piece1
\draw[black](8,23) -- (10,23) -- (10,25) -- (8,25) -- (8,23);
% Draw the convex parts of piece1
\fill[black!15](16,20) -- (18,20) -- (18,22) -- (16,22) -- (16,20);
% Draw the boundary of piece1
\draw[black](16,20) -- (18,20) -- (18,22) -- (16,22) -- (16,20);
% Draw the convex parts of piece1
\fill[black!15](19,23) -- (21,23) -- (21,25) -- (19,25) -- (19,23);
% Draw the boundary of piece1
\draw[black](19,23) -- (21,23) -- (21,25) -- (19,25) -- (19,23);
% Draw the convex parts of piece1
\fill[black!15](19,25) -- (21,25) -- (21,27) -- (19,27) -- (19,25);
% Draw the boundary of piece1
\draw[black](19,25) -- (21,25) -- (21,27) -- (19,27) -- (19,25);
% Draw the convex parts of piece1
\fill[black!15](21,6) -- (23,6) -- (23,8) -- (21,8) -- (21,6);
% Draw the boundary of piece1
\draw[black](21,6) -- (23,6) -- (23,8) -- (21,8) -- (21,6);
% Draw the convex parts of piece1
\fill[black!15](21,11) -- (23,11) -- (23,13) -- (21,13) -- (21,11);
% Draw the boundary of piece1
\draw[black](21,11) -- (23,11) -- (23,13) -- (21,13) -- (21,11);
% Draw the convex parts of piece1
\fill[black!15](25,11) -- (27,11) -- (27,13) -- (25,13) -- (25,11);
% Draw the boundary of piece1
\draw[black](25,11) -- (27,11) -- (27,13) -- (25,13) -- (25,11);
% Draw the convex parts of piece1
\fill[black!15](31,2) -- (33,2) -- (33,4) -- (31,4) -- (31,2);
% Draw the boundary of piece1
\draw[black](31,2) -- (33,2) -- (33,4) -- (31,4) -- (31,2);
% Draw the convex parts of piece1
\fill[black!15](31,7) -- (33,7) -- (33,9) -- (31,9) -- (31,7);
% Draw the boundary of piece1
\draw[black](31,7) -- (33,7) -- (33,9) -- (31,9) -- (31,7);
% Draw the convex parts of piece1
\fill[black!15](31,21) -- (33,21) -- (33,23) -- (31,23) -- (31,21);
% Draw the boundary of piece1
\draw[black](31,21) -- (33,21) -- (33,23) -- (31,23) -- (31,21);
% Draw the convex parts of piece6
\fill[black!15](7,3) -- (7,7) -- (6,6) -- (6,4) -- (7,3);
% Draw the boundary of piece6
\draw[black](7,3) -- (7,7) -- (6,6) -- (6,4) -- (7,3);
% Draw the convex parts of piece6
\fill[black!15](7,18) -- (7,22) -- (6,21) -- (6,19) -- (7,18);
% Draw the boundary of piece6
\draw[black](7,18) -- (7,22) -- (6,21) -- (6,19) -- (7,18);
% Draw the convex parts of piece6
\fill[black!15](24,19) -- (24,23) -- (23,22) -- (23,20) -- (24,19);
% Draw the boundary of piece6
\draw[black](24,19) -- (24,23) -- (23,22) -- (23,20) -- (24,19);
% Draw the convex parts of piece6
\fill[black!15](25,0) -- (25,4) -- (24,3) -- (24,1) -- (25,0);
% Draw the boundary of piece6
\draw[black](25,0) -- (25,4) -- (24,3) -- (24,1) -- (25,0);
% Draw the convex parts of piece6
\fill[black!15](27,6) -- (27,10) -- (26,9) -- (26,7) -- (27,6);
% Draw the boundary of piece6
\draw[black](27,6) -- (27,10) -- (26,9) -- (26,7) -- (27,6);
% Draw the convex parts of piece6
\fill[black!15](31,19) -- (31,23) -- (30,22) -- (30,20) -- (31,19);
% Draw the boundary of piece6
\draw[black](31,19) -- (31,23) -- (30,22) -- (30,20) -- (31,19);
% Draw the convex parts of piece6
\fill[black!15](32,17) -- (32,21) -- (31,20) -- (31,18) -- (32,17);
% Draw the boundary of piece6
\draw[black](32,17) -- (32,21) -- (31,20) -- (31,18) -- (32,17);
% Draw the convex parts of piece6
\fill[black!15](33,9) -- (33,13) -- (32,12) -- (32,10) -- (33,9);
% Draw the boundary of piece6
\draw[black](33,9) -- (33,13) -- (32,12) -- (32,10) -- (33,9);
% Draw the convex parts of piece6
\fill[black!15](33,16) -- (33,20) -- (32,19) -- (32,17) -- (33,16);
% Draw the boundary of piece6
\draw[black](33,16) -- (33,20) -- (32,19) -- (32,17) -- (33,16);
% Draw the convex parts of piece6
\fill[black!15](33,23) -- (33,27) -- (32,26) -- (32,24) -- (33,23);
% Draw the boundary of piece6
\draw[black](33,23) -- (33,27) -- (32,26) -- (32,24) -- (33,23);
% Draw the bounding box of the solution
\draw[black](0,0) -- (33,0) -- (33,27)-- (0,27) -- (0,0);
\draw (16.5,-2) node {artif5\_10};
\draw (16.5,-4.5) node {L = 33};
% Draw the bounding box of the solution
\end{tikzpicture}
% SOLUTION for problem instance = artif6_12
% Scaling factor for coordinates = 1
\begin{tikzpicture}[scale = 0.10]
% Draw the convex parts of piece3
\fill[black!15](4,15) -- (7,17) -- (6,19) -- (3,20) -- (0,18) -- (4,15);
% Draw the boundary of piece3
\draw[black](4,15) -- (7,17) -- (6,19) -- (3,20) -- (0,18) -- (4,15);
% Draw the convex parts of piece3
\fill[black!15](9,22) -- (12,24) -- (11,26) -- (8,27) -- (5,25) -- (9,22);
% Draw the boundary of piece3
\draw[black](9,22) -- (12,24) -- (11,26) -- (8,27) -- (5,25) -- (9,22);
% Draw the convex parts of piece3
\fill[black!15](14,3) -- (17,5) -- (16,7) -- (13,8) -- (10,6) -- (14,3);
% Draw the boundary of piece3
\draw[black](14,3) -- (17,5) -- (16,7) -- (13,8) -- (10,6) -- (14,3);
% Draw the convex parts of piece3
\fill[black!15](19,22) -- (22,24) -- (21,26) -- (18,27) -- (15,25) -- (19,22);
% Draw the boundary of piece3
\draw[black](19,22) -- (22,24) -- (21,26) -- (18,27) -- (15,25) -- (19,22);
% Draw the convex parts of piece3
\fill[black!15](32,14) -- (35,16) -- (34,18) -- (31,19) -- (28,17) -- (32,14);
% Draw the boundary of piece3
\draw[black](32,14) -- (35,16) -- (34,18) -- (31,19) -- (28,17) -- (32,14);
% Draw the convex parts of piece3
\fill[black!15](34,4) -- (37,6) -- (36,8) -- (33,9) -- (30,7) -- (34,4);
% Draw the boundary of piece3
\draw[black](34,4) -- (37,6) -- (36,8) -- (33,9) -- (30,7) -- (34,4);
% Draw the convex parts of piece2
\fill[black!15](10,10) -- (8,10) -- (11,7) -- (10,10);
\fill[black!15](11,7) -- (13,10) -- (12,12) -- (10,10) -- (11,7);
\fill[black!15](7,7) -- (11,7) -- (8,10) -- (7,11) -- (7,7);
% Draw the boundary of piece2
\draw[black](7,7) -- (11,7) -- (13,10) -- (12,12) -- (10,10) -- (8,10) -- (7,11) -- (7,7);
% Draw the convex parts of piece2
\fill[black!15](10,14) -- (8,14) -- (11,11) -- (10,14);
\fill[black!15](11,11) -- (13,14) -- (12,16) -- (10,14) -- (11,11);
\fill[black!15](7,11) -- (11,11) -- (8,14) -- (7,15) -- (7,11);
% Draw the boundary of piece2
\draw[black](7,11) -- (11,11) -- (13,14) -- (12,16) -- (10,14) -- (8,14) -- (7,15) -- (7,11);
% Draw the convex parts of piece2
\fill[black!15](10,18) -- (8,18) -- (11,15) -- (10,18);
\fill[black!15](11,15) -- (13,18) -- (12,20) -- (10,18) -- (11,15);
\fill[black!15](7,15) -- (11,15) -- (8,18) -- (7,19) -- (7,15);
% Draw the boundary of piece2
\draw[black](7,15) -- (11,15) -- (13,18) -- (12,20) -- (10,18) -- (8,18) -- (7,19) -- (7,15);
% Draw the convex parts of piece2
\fill[black!15](10,22) -- (8,22) -- (11,19) -- (10,22);
\fill[black!15](11,19) -- (13,22) -- (12,24) -- (10,22) -- (11,19);
\fill[black!15](7,19) -- (11,19) -- (8,22) -- (7,23) -- (7,19);
% Draw the boundary of piece2
\draw[black](7,19) -- (11,19) -- (13,22) -- (12,24) -- (10,22) -- (8,22) -- (7,23) -- (7,19);
% Draw the convex parts of piece2
\fill[black!15](15,3) -- (13,3) -- (16,0) -- (15,3);
\fill[black!15](16,0) -- (18,3) -- (17,5) -- (15,3) -- (16,0);
\fill[black!15](12,0) -- (16,0) -- (13,3) -- (12,4) -- (12,0);
% Draw the boundary of piece2
\draw[black](12,0) -- (16,0) -- (18,3) -- (17,5) -- (15,3) -- (13,3) -- (12,4) -- (12,0);
% Draw the convex parts of piece2
\fill[black!15](16,16) -- (14,16) -- (17,13) -- (16,16);
\fill[black!15](17,13) -- (19,16) -- (18,18) -- (16,16) -- (17,13);
\fill[black!15](13,13) -- (17,13) -- (14,16) -- (13,17) -- (13,13);
% Draw the boundary of piece2
\draw[black](13,13) -- (17,13) -- (19,16) -- (18,18) -- (16,16) -- (14,16) -- (13,17) -- (13,13);
% Draw the convex parts of piece2
\fill[black!15](16,20) -- (14,20) -- (17,17) -- (16,20);
\fill[black!15](17,17) -- (19,20) -- (18,22) -- (16,20) -- (17,17);
\fill[black!15](13,17) -- (17,17) -- (14,20) -- (13,21) -- (13,17);
% Draw the boundary of piece2
\draw[black](13,17) -- (17,17) -- (19,20) -- (18,22) -- (16,20) -- (14,20) -- (13,21) -- (13,17);
% Draw the convex parts of piece2
\fill[black!15](21,3) -- (19,3) -- (22,0) -- (21,3);
\fill[black!15](22,0) -- (24,3) -- (23,5) -- (21,3) -- (22,0);
\fill[black!15](18,0) -- (22,0) -- (19,3) -- (18,4) -- (18,0);
% Draw the boundary of piece2
\draw[black](18,0) -- (22,0) -- (24,3) -- (23,5) -- (21,3) -- (19,3) -- (18,4) -- (18,0);
% Draw the convex parts of piece2
\fill[black!15](21,7) -- (19,7) -- (22,4) -- (21,7);
\fill[black!15](22,4) -- (24,7) -- (23,9) -- (21,7) -- (22,4);
\fill[black!15](18,4) -- (22,4) -- (19,7) -- (18,8) -- (18,4);
% Draw the boundary of piece2
\draw[black](18,4) -- (22,4) -- (24,7) -- (23,9) -- (21,7) -- (19,7) -- (18,8) -- (18,4);
% Draw the convex parts of piece2
\fill[black!15](27,14) -- (25,14) -- (28,11) -- (27,14);
\fill[black!15](28,11) -- (30,14) -- (29,16) -- (27,14) -- (28,11);
\fill[black!15](24,11) -- (28,11) -- (25,14) -- (24,15) -- (24,11);
% Draw the boundary of piece2
\draw[black](24,11) -- (28,11) -- (30,14) -- (29,16) -- (27,14) -- (25,14) -- (24,15) -- (24,11);
% Draw the convex parts of piece2
\fill[black!15](30,3) -- (28,3) -- (31,0) -- (30,3);
\fill[black!15](31,0) -- (33,3) -- (32,5) -- (30,3) -- (31,0);
\fill[black!15](27,0) -- (31,0) -- (28,3) -- (27,4) -- (27,0);
% Draw the boundary of piece2
\draw[black](27,0) -- (31,0) -- (33,3) -- (32,5) -- (30,3) -- (28,3) -- (27,4) -- (27,0);
% Draw the convex parts of piece2
\fill[black!15](33,14) -- (31,14) -- (34,11) -- (33,14);
\fill[black!15](34,11) -- (36,14) -- (35,16) -- (33,14) -- (34,11);
\fill[black!15](30,11) -- (34,11) -- (31,14) -- (30,15) -- (30,11);
% Draw the boundary of piece2
\draw[black](30,11) -- (34,11) -- (36,14) -- (35,16) -- (33,14) -- (31,14) -- (30,15) -- (30,11);
% Draw the convex parts of piece7
\fill[black!15](20,11) -- (18,10) -- (19,9) -- (20,11);
\fill[black!15](20,11) -- (19,9) -- (21,10) -- (21,11) -- (20,11);
\fill[black!15](21,11) -- (22,12) -- (20,12) -- (20,11) -- (21,11);
\fill[black!15](17,9) -- (16,8) -- (19,8) -- (17,9);
\fill[black!15](17,9) -- (19,8) -- (19,9) -- (18,10) -- (17,9);
\fill[black!15](18,10) -- (18,12) -- (16,12) -- (17,9) -- (18,10);
% Draw the boundary of piece7
\draw[black](16,8) -- (19,8) -- (19,9) -- (21,10) -- (21,11) -- (22,12) -- (20,12) -- (20,11) -- (18,10) -- (18,12) -- (16,12) -- (17,9) -- (16,8);
% Draw the convex parts of piece7
\fill[black!15](23,15) -- (21,14) -- (22,13) -- (23,15);
\fill[black!15](23,15) -- (22,13) -- (24,14) -- (24,15) -- (23,15);
\fill[black!15](24,15) -- (25,16) -- (23,16) -- (23,15) -- (24,15);
\fill[black!15](20,13) -- (19,12) -- (22,12) -- (20,13);
\fill[black!15](20,13) -- (22,12) -- (22,13) -- (21,14) -- (20,13);
\fill[black!15](21,14) -- (21,16) -- (19,16) -- (20,13) -- (21,14);
% Draw the boundary of piece7
\draw[black](19,12) -- (22,12) -- (22,13) -- (24,14) -- (24,15) -- (25,16) -- (23,16) -- (23,15) -- (21,14) -- (21,16) -- (19,16) -- (20,13) -- (19,12);
% Draw the convex parts of piece7
\fill[black!15](23,19) -- (21,18) -- (22,17) -- (23,19);
\fill[black!15](23,19) -- (22,17) -- (24,18) -- (24,19) -- (23,19);
\fill[black!15](24,19) -- (25,20) -- (23,20) -- (23,19) -- (24,19);
\fill[black!15](20,17) -- (19,16) -- (22,16) -- (20,17);
\fill[black!15](20,17) -- (22,16) -- (22,17) -- (21,18) -- (20,17);
\fill[black!15](21,18) -- (21,20) -- (19,20) -- (20,17) -- (21,18);
% Draw the boundary of piece7
\draw[black](19,16) -- (22,16) -- (22,17) -- (24,18) -- (24,19) -- (25,20) -- (23,20) -- (23,19) -- (21,18) -- (21,20) -- (19,20) -- (20,17) -- (19,16);
% Draw the convex parts of piece7
\fill[black!15](25,26) -- (23,25) -- (24,24) -- (25,26);
\fill[black!15](25,26) -- (24,24) -- (26,25) -- (26,26) -- (25,26);
\fill[black!15](26,26) -- (27,27) -- (25,27) -- (25,26) -- (26,26);
\fill[black!15](22,24) -- (21,23) -- (24,23) -- (22,24);
\fill[black!15](22,24) -- (24,23) -- (24,24) -- (23,25) -- (22,24);
\fill[black!15](23,25) -- (23,27) -- (21,27) -- (22,24) -- (23,25);
% Draw the boundary of piece7
\draw[black](21,23) -- (24,23) -- (24,24) -- (26,25) -- (26,26) -- (27,27) -- (25,27) -- (25,26) -- (23,25) -- (23,27) -- (21,27) -- (22,24) -- (21,23);
% Draw the convex parts of piece7
\fill[black!15](28,23) -- (26,22) -- (27,21) -- (28,23);
\fill[black!15](28,23) -- (27,21) -- (29,22) -- (29,23) -- (28,23);
\fill[black!15](29,23) -- (30,24) -- (28,24) -- (28,23) -- (29,23);
\fill[black!15](25,21) -- (24,20) -- (27,20) -- (25,21);
\fill[black!15](25,21) -- (27,20) -- (27,21) -- (26,22) -- (25,21);
\fill[black!15](26,22) -- (26,24) -- (24,24) -- (25,21) -- (26,22);
% Draw the boundary of piece7
\draw[black](24,20) -- (27,20) -- (27,21) -- (29,22) -- (29,23) -- (30,24) -- (28,24) -- (28,23) -- (26,22) -- (26,24) -- (24,24) -- (25,21) -- (24,20);
% Draw the convex parts of piece7
\fill[black!15](33,26) -- (31,25) -- (32,24) -- (33,26);
\fill[black!15](33,26) -- (32,24) -- (34,25) -- (34,26) -- (33,26);
\fill[black!15](34,26) -- (35,27) -- (33,27) -- (33,26) -- (34,26);
\fill[black!15](30,24) -- (29,23) -- (32,23) -- (30,24);
\fill[black!15](30,24) -- (32,23) -- (32,24) -- (31,25) -- (30,24);
\fill[black!15](31,25) -- (31,27) -- (29,27) -- (30,24) -- (31,25);
% Draw the boundary of piece7
\draw[black](29,23) -- (32,23) -- (32,24) -- (34,25) -- (34,26) -- (35,27) -- (33,27) -- (33,26) -- (31,25) -- (31,27) -- (29,27) -- (30,24) -- (29,23);
% Draw the convex parts of piece5
\fill[black!15](0,0) -- (3,0) -- (4,1) -- (1,2) -- (0,0);
\fill[black!15](0,0) -- (1,2) -- (1,3) -- (0,0);
\fill[black!15](0,5) -- (0,0) -- (1,3) -- (0,5);
\fill[black!15](2,5) -- (0,5) -- (1,3) -- (2,4) -- (2,5);
\fill[black!15](2,4) -- (3,4) -- (2,5) -- (2,4);
% Draw the boundary of piece5
\draw[black](0,0) -- (3,0) -- (4,1) -- (1,2) -- (1,3) -- (2,4) -- (3,4) -- (2,5) -- (0,5) -- (0,0);
% Draw the convex parts of piece5
\fill[black!15](0,5) -- (3,5) -- (4,6) -- (1,7) -- (0,5);
\fill[black!15](0,5) -- (1,7) -- (1,8) -- (0,5);
\fill[black!15](0,10) -- (0,5) -- (1,8) -- (0,10);
\fill[black!15](2,10) -- (0,10) -- (1,8) -- (2,9) -- (2,10);
\fill[black!15](2,9) -- (3,9) -- (2,10) -- (2,9);
% Draw the boundary of piece5
\draw[black](0,5) -- (3,5) -- (4,6) -- (1,7) -- (1,8) -- (2,9) -- (3,9) -- (2,10) -- (0,10) -- (0,5);
% Draw the convex parts of piece5
\fill[black!15](0,10) -- (3,10) -- (4,11) -- (1,12) -- (0,10);
\fill[black!15](0,10) -- (1,12) -- (1,13) -- (0,10);
\fill[black!15](0,15) -- (0,10) -- (1,13) -- (0,15);
\fill[black!15](2,15) -- (0,15) -- (1,13) -- (2,14) -- (2,15);
\fill[black!15](2,14) -- (3,14) -- (2,15) -- (2,14);
% Draw the boundary of piece5
\draw[black](0,10) -- (3,10) -- (4,11) -- (1,12) -- (1,13) -- (2,14) -- (3,14) -- (2,15) -- (0,15) -- (0,10);
% Draw the convex parts of piece5
\fill[black!15](0,20) -- (3,20) -- (4,21) -- (1,22) -- (0,20);
\fill[black!15](0,20) -- (1,22) -- (1,23) -- (0,20);
\fill[black!15](0,25) -- (0,20) -- (1,23) -- (0,25);
\fill[black!15](2,25) -- (0,25) -- (1,23) -- (2,24) -- (2,25);
\fill[black!15](2,24) -- (3,24) -- (2,25) -- (2,24);
% Draw the boundary of piece5
\draw[black](0,20) -- (3,20) -- (4,21) -- (1,22) -- (1,23) -- (2,24) -- (3,24) -- (2,25) -- (0,25) -- (0,20);
% Draw the convex parts of piece5
\fill[black!15](3,22) -- (6,22) -- (7,23) -- (4,24) -- (3,22);
\fill[black!15](3,22) -- (4,24) -- (4,25) -- (3,22);
\fill[black!15](3,27) -- (3,22) -- (4,25) -- (3,27);
\fill[black!15](5,27) -- (3,27) -- (4,25) -- (5,26) -- (5,27);
\fill[black!15](5,26) -- (6,26) -- (5,27) -- (5,26);
% Draw the boundary of piece5
\draw[black](3,22) -- (6,22) -- (7,23) -- (4,24) -- (4,25) -- (5,26) -- (6,26) -- (5,27) -- (3,27) -- (3,22);
% Draw the convex parts of piece5
\fill[black!15](4,0) -- (7,0) -- (8,1) -- (5,2) -- (4,0);
\fill[black!15](4,0) -- (5,2) -- (5,3) -- (4,0);
\fill[black!15](4,5) -- (4,0) -- (5,3) -- (4,5);
\fill[black!15](6,5) -- (4,5) -- (5,3) -- (6,4) -- (6,5);
\fill[black!15](6,4) -- (7,4) -- (6,5) -- (6,4);
% Draw the boundary of piece5
\draw[black](4,0) -- (7,0) -- (8,1) -- (5,2) -- (5,3) -- (6,4) -- (7,4) -- (6,5) -- (4,5) -- (4,0);
% Draw the convex parts of piece5
\fill[black!15](4,5) -- (7,5) -- (8,6) -- (5,7) -- (4,5);
\fill[black!15](4,5) -- (5,7) -- (5,8) -- (4,5);
\fill[black!15](4,10) -- (4,5) -- (5,8) -- (4,10);
\fill[black!15](6,10) -- (4,10) -- (5,8) -- (6,9) -- (6,10);
\fill[black!15](6,9) -- (7,9) -- (6,10) -- (6,9);
% Draw the boundary of piece5
\draw[black](4,5) -- (7,5) -- (8,6) -- (5,7) -- (5,8) -- (6,9) -- (7,9) -- (6,10) -- (4,10) -- (4,5);
% Draw the convex parts of piece5
\fill[black!15](8,0) -- (11,0) -- (12,1) -- (9,2) -- (8,0);
\fill[black!15](8,0) -- (9,2) -- (9,3) -- (8,0);
\fill[black!15](8,5) -- (8,0) -- (9,3) -- (8,5);
\fill[black!15](10,5) -- (8,5) -- (9,3) -- (10,4) -- (10,5);
\fill[black!15](10,4) -- (11,4) -- (10,5) -- (10,4);
% Draw the boundary of piece5
\draw[black](8,0) -- (11,0) -- (12,1) -- (9,2) -- (9,3) -- (10,4) -- (11,4) -- (10,5) -- (8,5) -- (8,0);
% Draw the convex parts of piece5
\fill[black!15](13,8) -- (16,8) -- (17,9) -- (14,10) -- (13,8);
\fill[black!15](13,8) -- (14,10) -- (14,11) -- (13,8);
\fill[black!15](13,13) -- (13,8) -- (14,11) -- (13,13);
\fill[black!15](15,13) -- (13,13) -- (14,11) -- (15,12) -- (15,13);
\fill[black!15](15,12) -- (16,12) -- (15,13) -- (15,12);
% Draw the boundary of piece5
\draw[black](13,8) -- (16,8) -- (17,9) -- (14,10) -- (14,11) -- (15,12) -- (16,12) -- (15,13) -- (13,13) -- (13,8);
% Draw the convex parts of piece5
\fill[black!15](13,22) -- (16,22) -- (17,23) -- (14,24) -- (13,22);
\fill[black!15](13,22) -- (14,24) -- (14,25) -- (13,22);
\fill[black!15](13,27) -- (13,22) -- (14,25) -- (13,27);
\fill[black!15](15,27) -- (13,27) -- (14,25) -- (15,26) -- (15,27);
\fill[black!15](15,26) -- (16,26) -- (15,27) -- (15,26);
% Draw the boundary of piece5
\draw[black](13,22) -- (16,22) -- (17,23) -- (14,24) -- (14,25) -- (15,26) -- (16,26) -- (15,27) -- (13,27) -- (13,22);
% Draw the convex parts of piece5
\fill[black!15](25,15) -- (28,15) -- (29,16) -- (26,17) -- (25,15);
\fill[black!15](25,15) -- (26,17) -- (26,18) -- (25,15);
\fill[black!15](25,20) -- (25,15) -- (26,18) -- (25,20);
\fill[black!15](27,20) -- (25,20) -- (26,18) -- (27,19) -- (27,20);
\fill[black!15](27,19) -- (28,19) -- (27,20) -- (27,19);
% Draw the boundary of piece5
\draw[black](25,15) -- (28,15) -- (29,16) -- (26,17) -- (26,18) -- (27,19) -- (28,19) -- (27,20) -- (25,20) -- (25,15);
% Draw the convex parts of piece5
\fill[black!15](28,4) -- (31,4) -- (32,5) -- (29,6) -- (28,4);
\fill[black!15](28,4) -- (29,6) -- (29,7) -- (28,4);
\fill[black!15](28,9) -- (28,4) -- (29,7) -- (28,9);
\fill[black!15](30,9) -- (28,9) -- (29,7) -- (30,8) -- (30,9);
\fill[black!15](30,8) -- (31,8) -- (30,9) -- (30,8);
% Draw the boundary of piece5
\draw[black](28,4) -- (31,4) -- (32,5) -- (29,6) -- (29,7) -- (30,8) -- (31,8) -- (30,9) -- (28,9) -- (28,4);
% Draw the convex parts of piece8
\fill[black!15](5,14) -- (4,14) -- (6,13) -- (5,14);
\fill[black!15](6,15) -- (5,14) -- (6,13) -- (6,15);
\fill[black!15](6,13) -- (7,16) -- (6,16) -- (6,15) -- (6,13);
\fill[black!15](4,11) -- (6,11) -- (5,12) -- (4,11);
\fill[black!15](4,11) -- (5,12) -- (5,13) -- (4,11);
\fill[black!15](5,13) -- (6,13) -- (4,14) -- (5,13);
\fill[black!15](4,11) -- (5,13) -- (4,14) -- (4,11);
\fill[black!15](4,11) -- (4,14) -- (4,15) -- (3,15) -- (4,11);
% Draw the boundary of piece8
\draw[black](4,11) -- (6,11) -- (5,12) -- (5,13) -- (6,13) -- (7,16) -- (6,16) -- (6,15) -- (5,14) -- (4,14) -- (4,15) -- (3,15) -- (4,11);
% Draw the convex parts of piece8
\fill[black!15](26,5) -- (25,5) -- (27,4) -- (26,5);
\fill[black!15](27,6) -- (26,5) -- (27,4) -- (27,6);
\fill[black!15](27,4) -- (28,7) -- (27,7) -- (27,6) -- (27,4);
\fill[black!15](25,2) -- (27,2) -- (26,3) -- (25,2);
\fill[black!15](25,2) -- (26,3) -- (26,4) -- (25,2);
\fill[black!15](26,4) -- (27,4) -- (25,5) -- (26,4);
\fill[black!15](25,2) -- (26,4) -- (25,5) -- (25,2);
\fill[black!15](25,2) -- (25,5) -- (25,6) -- (24,6) -- (25,2);
% Draw the boundary of piece8
\draw[black](25,2) -- (27,2) -- (26,3) -- (26,4) -- (27,4) -- (28,7) -- (27,7) -- (27,6) -- (26,5) -- (25,5) -- (25,6) -- (24,6) -- (25,2);
% Draw the convex parts of piece8
\fill[black!15](26,9) -- (25,9) -- (27,8) -- (26,9);
\fill[black!15](27,10) -- (26,9) -- (27,8) -- (27,10);
\fill[black!15](27,8) -- (28,11) -- (27,11) -- (27,10) -- (27,8);
\fill[black!15](25,6) -- (27,6) -- (26,7) -- (25,6);
\fill[black!15](25,6) -- (26,7) -- (26,8) -- (25,6);
\fill[black!15](26,8) -- (27,8) -- (25,9) -- (26,8);
\fill[black!15](25,6) -- (26,8) -- (25,9) -- (25,6);
\fill[black!15](25,6) -- (25,9) -- (25,10) -- (24,10) -- (25,6);
% Draw the boundary of piece8
\draw[black](25,6) -- (27,6) -- (26,7) -- (26,8) -- (27,8) -- (28,11) -- (27,11) -- (27,10) -- (26,9) -- (25,9) -- (25,10) -- (24,10) -- (25,6);
% Draw the convex parts of piece8
\fill[black!15](31,22) -- (30,22) -- (32,21) -- (31,22);
\fill[black!15](32,23) -- (31,22) -- (32,21) -- (32,23);
\fill[black!15](32,21) -- (33,24) -- (32,24) -- (32,23) -- (32,21);
\fill[black!15](30,19) -- (32,19) -- (31,20) -- (30,19);
\fill[black!15](30,19) -- (31,20) -- (31,21) -- (30,19);
\fill[black!15](31,21) -- (32,21) -- (30,22) -- (31,21);
\fill[black!15](30,19) -- (31,21) -- (30,22) -- (30,19);
\fill[black!15](30,19) -- (30,22) -- (30,23) -- (29,23) -- (30,19);
% Draw the boundary of piece8
\draw[black](30,19) -- (32,19) -- (31,20) -- (31,21) -- (32,21) -- (33,24) -- (32,24) -- (32,23) -- (31,22) -- (30,22) -- (30,23) -- (29,23) -- (30,19);
% Draw the convex parts of piece8
\fill[black!15](36,19) -- (35,19) -- (37,18) -- (36,19);
\fill[black!15](37,20) -- (36,19) -- (37,18) -- (37,20);
\fill[black!15](37,18) -- (38,21) -- (37,21) -- (37,20) -- (37,18);
\fill[black!15](35,16) -- (37,16) -- (36,17) -- (35,16);
\fill[black!15](35,16) -- (36,17) -- (36,18) -- (35,16);
\fill[black!15](36,18) -- (37,18) -- (35,19) -- (36,18);
\fill[black!15](35,16) -- (36,18) -- (35,19) -- (35,16);
\fill[black!15](35,16) -- (35,19) -- (35,20) -- (34,20) -- (35,16);
% Draw the boundary of piece8
\draw[black](35,16) -- (37,16) -- (36,17) -- (36,18) -- (37,18) -- (38,21) -- (37,21) -- (37,20) -- (36,19) -- (35,19) -- (35,20) -- (34,20) -- (35,16);
% Draw the convex parts of piece8
\fill[black!15](36,25) -- (35,25) -- (37,24) -- (36,25);
\fill[black!15](37,26) -- (36,25) -- (37,24) -- (37,26);
\fill[black!15](37,24) -- (38,27) -- (37,27) -- (37,26) -- (37,24);
\fill[black!15](35,22) -- (37,22) -- (36,23) -- (35,22);
\fill[black!15](35,22) -- (36,23) -- (36,24) -- (35,22);
\fill[black!15](36,24) -- (37,24) -- (35,25) -- (36,24);
\fill[black!15](35,22) -- (36,24) -- (35,25) -- (35,22);
\fill[black!15](35,22) -- (35,25) -- (35,26) -- (34,26) -- (35,22);
% Draw the boundary of piece8
\draw[black](35,22) -- (37,22) -- (36,23) -- (36,24) -- (37,24) -- (38,27) -- (37,27) -- (37,26) -- (36,25) -- (35,25) -- (35,26) -- (34,26) -- (35,22);
% Draw the convex parts of piece4
\fill[black!15](2,26) -- (1,26) -- (3,25) -- (2,26);
\fill[black!15](3,25) -- (3,27) -- (2,27) -- (2,26) -- (3,25);
\fill[black!15](0,25) -- (3,25) -- (1,26) -- (0,25);
\fill[black!15](0,25) -- (1,26) -- (1,27) -- (0,27) -- (0,25);
% Draw the boundary of piece4
\draw[black](0,25) -- (3,25) -- (3,27) -- (2,27) -- (2,26) -- (1,26) -- (1,27) -- (0,27) -- (0,25);
% Draw the convex parts of piece4
\fill[black!15](6,21) -- (5,21) -- (7,20) -- (6,21);
\fill[black!15](7,20) -- (7,22) -- (6,22) -- (6,21) -- (7,20);
\fill[black!15](4,20) -- (7,20) -- (5,21) -- (4,20);
\fill[black!15](4,20) -- (5,21) -- (5,22) -- (4,22) -- (4,20);
% Draw the boundary of piece4
\draw[black](4,20) -- (7,20) -- (7,22) -- (6,22) -- (6,21) -- (5,21) -- (5,22) -- (4,22) -- (4,20);
% Draw the convex parts of piece4
\fill[black!15](23,10) -- (22,10) -- (24,9) -- (23,10);
\fill[black!15](24,9) -- (24,11) -- (23,11) -- (23,10) -- (24,9);
\fill[black!15](21,9) -- (24,9) -- (22,10) -- (21,9);
\fill[black!15](21,9) -- (22,10) -- (22,11) -- (21,11) -- (21,9);
% Draw the boundary of piece4
\draw[black](21,9) -- (24,9) -- (24,11) -- (23,11) -- (23,10) -- (22,10) -- (22,11) -- (21,11) -- (21,9);
% Draw the convex parts of piece4
\fill[black!15](23,22) -- (22,22) -- (24,21) -- (23,22);
\fill[black!15](24,21) -- (24,23) -- (23,23) -- (23,22) -- (24,21);
\fill[black!15](21,21) -- (24,21) -- (22,22) -- (21,21);
\fill[black!15](21,21) -- (22,22) -- (22,23) -- (21,23) -- (21,21);
% Draw the boundary of piece4
\draw[black](21,21) -- (24,21) -- (24,23) -- (23,23) -- (23,22) -- (22,22) -- (22,23) -- (21,23) -- (21,21);
% Draw the convex parts of piece4
\fill[black!15](26,1) -- (25,1) -- (27,0) -- (26,1);
\fill[black!15](27,0) -- (27,2) -- (26,2) -- (26,1) -- (27,0);
\fill[black!15](24,0) -- (27,0) -- (25,1) -- (24,0);
\fill[black!15](24,0) -- (25,1) -- (25,2) -- (24,2) -- (24,0);
% Draw the boundary of piece4
\draw[black](24,0) -- (27,0) -- (27,2) -- (26,2) -- (26,1) -- (25,1) -- (25,2) -- (24,2) -- (24,0);
% Draw the convex parts of piece4
\fill[black!15](28,25) -- (27,25) -- (29,24) -- (28,25);
\fill[black!15](29,24) -- (29,26) -- (28,26) -- (28,25) -- (29,24);
\fill[black!15](26,24) -- (29,24) -- (27,25) -- (26,24);
\fill[black!15](26,24) -- (27,25) -- (27,26) -- (26,26) -- (26,24);
% Draw the boundary of piece4
\draw[black](26,24) -- (29,24) -- (29,26) -- (28,26) -- (28,25) -- (27,25) -- (27,26) -- (26,26) -- (26,24);
% Draw the convex parts of piece4
\fill[black!15](30,10) -- (29,10) -- (31,9) -- (30,10);
\fill[black!15](31,9) -- (31,11) -- (30,11) -- (30,10) -- (31,9);
\fill[black!15](28,9) -- (31,9) -- (29,10) -- (28,9);
\fill[black!15](28,9) -- (29,10) -- (29,11) -- (28,11) -- (28,9);
% Draw the boundary of piece4
\draw[black](28,9) -- (31,9) -- (31,11) -- (30,11) -- (30,10) -- (29,10) -- (29,11) -- (28,11) -- (28,9);
% Draw the convex parts of piece4
\fill[black!15](34,10) -- (33,10) -- (35,9) -- (34,10);
\fill[black!15](35,9) -- (35,11) -- (34,11) -- (34,10) -- (35,9);
\fill[black!15](32,9) -- (35,9) -- (33,10) -- (32,9);
\fill[black!15](32,9) -- (33,10) -- (33,11) -- (32,11) -- (32,9);
% Draw the boundary of piece4
\draw[black](32,9) -- (35,9) -- (35,11) -- (34,11) -- (34,10) -- (33,10) -- (33,11) -- (32,11) -- (32,9);
% Draw the convex parts of piece4
\fill[black!15](35,1) -- (34,1) -- (36,0) -- (35,1);
\fill[black!15](36,0) -- (36,2) -- (35,2) -- (35,1) -- (36,0);
\fill[black!15](33,0) -- (36,0) -- (34,1) -- (33,0);
\fill[black!15](33,0) -- (34,1) -- (34,2) -- (33,2) -- (33,0);
% Draw the boundary of piece4
\draw[black](33,0) -- (36,0) -- (36,2) -- (35,2) -- (35,1) -- (34,1) -- (34,2) -- (33,2) -- (33,0);
% Draw the convex parts of piece4
\fill[black!15](36,21) -- (35,21) -- (37,20) -- (36,21);
\fill[black!15](37,20) -- (37,22) -- (36,22) -- (36,21) -- (37,20);
\fill[black!15](34,20) -- (37,20) -- (35,21) -- (34,20);
\fill[black!15](34,20) -- (35,21) -- (35,22) -- (34,22) -- (34,20);
% Draw the boundary of piece4
\draw[black](34,20) -- (37,20) -- (37,22) -- (36,22) -- (36,21) -- (35,21) -- (35,22) -- (34,22) -- (34,20);
% Draw the convex parts of piece4
\fill[black!15](37,3) -- (36,3) -- (38,2) -- (37,3);
\fill[black!15](38,2) -- (38,4) -- (37,4) -- (37,3) -- (38,2);
\fill[black!15](35,2) -- (38,2) -- (36,3) -- (35,2);
\fill[black!15](35,2) -- (36,3) -- (36,4) -- (35,4) -- (35,2);
% Draw the boundary of piece4
\draw[black](35,2) -- (38,2) -- (38,4) -- (37,4) -- (37,3) -- (36,3) -- (36,4) -- (35,4) -- (35,2);
% Draw the convex parts of piece4
\fill[black!15](37,11) -- (36,11) -- (38,10) -- (37,11);
\fill[black!15](38,10) -- (38,12) -- (37,12) -- (37,11) -- (38,10);
\fill[black!15](35,10) -- (38,10) -- (36,11) -- (35,10);
\fill[black!15](35,10) -- (36,11) -- (36,12) -- (35,12) -- (35,10);
% Draw the boundary of piece4
\draw[black](35,10) -- (38,10) -- (38,12) -- (37,12) -- (37,11) -- (36,11) -- (36,12) -- (35,12) -- (35,10);
% Draw the convex parts of piece1
\fill[black!15](2,7) -- (4,7) -- (4,9) -- (2,9) -- (2,7);
% Draw the boundary of piece1
\draw[black](2,7) -- (4,7) -- (4,9) -- (2,9) -- (2,7);
% Draw the convex parts of piece1
\fill[black!15](6,2) -- (8,2) -- (8,4) -- (6,4) -- (6,2);
% Draw the boundary of piece1
\draw[black](6,2) -- (8,2) -- (8,4) -- (6,4) -- (6,2);
% Draw the convex parts of piece1
\fill[black!15](8,5) -- (10,5) -- (10,7) -- (8,7) -- (8,5);
% Draw the boundary of piece1
\draw[black](8,5) -- (10,5) -- (10,7) -- (8,7) -- (8,5);
% Draw the convex parts of piece1
\fill[black!15](10,2) -- (12,2) -- (12,4) -- (10,4) -- (10,2);
% Draw the boundary of piece1
\draw[black](10,2) -- (12,2) -- (12,4) -- (10,4) -- (10,2);
% Draw the convex parts of piece1
\fill[black!15](14,20) -- (16,20) -- (16,22) -- (14,22) -- (14,20);
% Draw the boundary of piece1
\draw[black](14,20) -- (16,20) -- (16,22) -- (14,22) -- (14,20);
% Draw the convex parts of piece1
\fill[black!15](19,7) -- (21,7) -- (21,9) -- (19,9) -- (19,7);
% Draw the boundary of piece1
\draw[black](19,7) -- (21,7) -- (21,9) -- (19,9) -- (19,7);
% Draw the convex parts of piece1
\fill[black!15](19,20) -- (21,20) -- (21,22) -- (19,22) -- (19,20);
% Draw the boundary of piece1
\draw[black](19,20) -- (21,20) -- (21,22) -- (19,22) -- (19,20);
% Draw the convex parts of piece1
\fill[black!15](21,19) -- (23,19) -- (23,21) -- (21,21) -- (21,19);
% Draw the boundary of piece1
\draw[black](21,19) -- (23,19) -- (23,21) -- (21,21) -- (21,19);
% Draw the convex parts of piece1
\fill[black!15](22,11) -- (24,11) -- (24,13) -- (22,13) -- (22,11);
% Draw the boundary of piece1
\draw[black](22,11) -- (24,11) -- (24,13) -- (22,13) -- (22,11);
% Draw the convex parts of piece1
\fill[black!15](32,19) -- (34,19) -- (34,21) -- (32,21) -- (32,19);
% Draw the boundary of piece1
\draw[black](32,19) -- (34,19) -- (34,21) -- (32,21) -- (32,19);
% Draw the convex parts of piece1
\fill[black!15](33,2) -- (35,2) -- (35,4) -- (33,4) -- (33,2);
% Draw the boundary of piece1
\draw[black](33,2) -- (35,2) -- (35,4) -- (33,4) -- (33,2);
% Draw the convex parts of piece1
\fill[black!15](36,0) -- (38,0) -- (38,2) -- (36,2) -- (36,0);
% Draw the boundary of piece1
\draw[black](36,0) -- (38,0) -- (38,2) -- (36,2) -- (36,0);
% Draw the convex parts of piece6
\fill[black!15](4,2) -- (4,6) -- (3,5) -- (3,3) -- (4,2);
% Draw the boundary of piece6
\draw[black](4,2) -- (4,6) -- (3,5) -- (3,3) -- (4,2);
% Draw the convex parts of piece6
\fill[black!15](7,10) -- (7,14) -- (6,13) -- (6,11) -- (7,10);
% Draw the boundary of piece6
\draw[black](7,10) -- (7,14) -- (6,13) -- (6,11) -- (7,10);
% Draw the convex parts of piece6
\fill[black!15](13,23) -- (13,27) -- (12,26) -- (12,24) -- (13,23);
% Draw the boundary of piece6
\draw[black](13,23) -- (13,27) -- (12,26) -- (12,24) -- (13,23);
% Draw the convex parts of piece6
\fill[black!15](18,4) -- (18,8) -- (17,7) -- (17,5) -- (18,4);
% Draw the boundary of piece6
\draw[black](18,4) -- (18,8) -- (17,7) -- (17,5) -- (18,4);
% Draw the convex parts of piece6
\fill[black!15](19,11) -- (19,15) -- (18,14) -- (18,12) -- (19,11);
% Draw the boundary of piece6
\draw[black](19,11) -- (19,15) -- (18,14) -- (18,12) -- (19,11);
% Draw the convex parts of piece6
\fill[black!15](25,16) -- (25,20) -- (24,19) -- (24,17) -- (25,16);
% Draw the boundary of piece6
\draw[black](25,16) -- (25,20) -- (24,19) -- (24,17) -- (25,16);
% Draw the convex parts of piece6
\fill[black!15](29,18) -- (29,22) -- (28,21) -- (28,19) -- (29,18);
% Draw the boundary of piece6
\draw[black](29,18) -- (29,22) -- (28,21) -- (28,19) -- (29,18);
% Draw the convex parts of piece6
\fill[black!15](34,21) -- (34,25) -- (33,24) -- (33,22) -- (34,21);
% Draw the boundary of piece6
\draw[black](34,21) -- (34,25) -- (33,24) -- (33,22) -- (34,21);
% Draw the convex parts of piece6
\fill[black!15](37,12) -- (37,16) -- (36,15) -- (36,13) -- (37,12);
% Draw the boundary of piece6
\draw[black](37,12) -- (37,16) -- (36,15) -- (36,13) -- (37,12);
% Draw the convex parts of piece6
\fill[black!15](38,6) -- (38,10) -- (37,9) -- (37,7) -- (38,6);
% Draw the boundary of piece6
\draw[black](38,6) -- (38,10) -- (37,9) -- (37,7) -- (38,6);
% Draw the convex parts of piece6
\fill[black!15](38,12) -- (38,16) -- (37,15) -- (37,13) -- (38,12);
% Draw the boundary of piece6
\draw[black](38,12) -- (38,16) -- (37,15) -- (37,13) -- (38,12);
% Draw the convex parts of piece6
\fill[black!15](38,21) -- (38,25) -- (37,24) -- (37,22) -- (38,21);
% Draw the boundary of piece6
\draw[black](38,21) -- (38,25) -- (37,24) -- (37,22) -- (38,21);
% Draw the bounding box of the solution
\draw[black](0,0) -- (38,0) -- (38,27)-- (0,27) -- (0,0);
\draw (19,-2) node {artif6\_12};
\draw (19,-4.5) node {L = 38};
\end{tikzpicture}
% SOLUTION for problem instance = artif7_14
% Scaling factor for coordinates = 1
\begin{tikzpicture}[scale = 0.10]
% Draw the convex parts of piece3
\fill[black!15](4,5) -- (7,7) -- (6,9) -- (3,10) -- (0,8) -- (4,5);
% Draw the boundary of piece3
\draw[black](4,5) -- (7,7) -- (6,9) -- (3,10) -- (0,8) -- (4,5);
% Draw the convex parts of piece3
\fill[black!15](4,15) -- (7,17) -- (6,19) -- (3,20) -- (0,18) -- (4,15);
% Draw the boundary of piece3
\draw[black](4,15) -- (7,17) -- (6,19) -- (3,20) -- (0,18) -- (4,15);
% Draw the convex parts of piece3
\fill[black!15](9,2) -- (12,4) -- (11,6) -- (8,7) -- (5,5) -- (9,2);
% Draw the boundary of piece3
\draw[black](9,2) -- (12,4) -- (11,6) -- (8,7) -- (5,5) -- (9,2);
% Draw the convex parts of piece3
\fill[black!15](9,7) -- (12,9) -- (11,11) -- (8,12) -- (5,10) -- (9,7);
% Draw the boundary of piece3
\draw[black](9,7) -- (12,9) -- (11,11) -- (8,12) -- (5,10) -- (9,7);
% Draw the convex parts of piece3
\fill[black!15](37,15) -- (40,17) -- (39,19) -- (36,20) -- (33,18) -- (37,15);
% Draw the boundary of piece3
\draw[black](37,15) -- (40,17) -- (39,19) -- (36,20) -- (33,18) -- (37,15);
% Draw the convex parts of piece3
\fill[black!15](39,0) -- (42,2) -- (41,4) -- (38,5) -- (35,3) -- (39,0);
% Draw the boundary of piece3
\draw[black](39,0) -- (42,2) -- (41,4) -- (38,5) -- (35,3) -- (39,0);
% Draw the convex parts of piece3
\fill[black!15](42,4) -- (45,6) -- (44,8) -- (41,9) -- (38,7) -- (42,4);
% Draw the boundary of piece3
\draw[black](42,4) -- (45,6) -- (44,8) -- (41,9) -- (38,7) -- (42,4);
% Draw the convex parts of piece2
\fill[black!15](10,15) -- (8,15) -- (11,12) -- (10,15);
\fill[black!15](11,12) -- (13,15) -- (12,17) -- (10,15) -- (11,12);
\fill[black!15](7,12) -- (11,12) -- (8,15) -- (7,16) -- (7,12);
% Draw the boundary of piece2
\draw[black](7,12) -- (11,12) -- (13,15) -- (12,17) -- (10,15) -- (8,15) -- (7,16) -- (7,12);
% Draw the convex parts of piece2
\fill[black!15](10,19) -- (8,19) -- (11,16) -- (10,19);
\fill[black!15](11,16) -- (13,19) -- (12,21) -- (10,19) -- (11,16);
\fill[black!15](7,16) -- (11,16) -- (8,19) -- (7,20) -- (7,16);
% Draw the boundary of piece2
\draw[black](7,16) -- (11,16) -- (13,19) -- (12,21) -- (10,19) -- (8,19) -- (7,20) -- (7,16);
% Draw the convex parts of piece2
\fill[black!15](10,23) -- (8,23) -- (11,20) -- (10,23);
\fill[black!15](11,20) -- (13,23) -- (12,25) -- (10,23) -- (11,20);
\fill[black!15](7,20) -- (11,20) -- (8,23) -- (7,24) -- (7,20);
% Draw the boundary of piece2
\draw[black](7,20) -- (11,20) -- (13,23) -- (12,25) -- (10,23) -- (8,23) -- (7,24) -- (7,20);
% Draw the convex parts of piece2
\fill[black!15](15,3) -- (13,3) -- (16,0) -- (15,3);
\fill[black!15](16,0) -- (18,3) -- (17,5) -- (15,3) -- (16,0);
\fill[black!15](12,0) -- (16,0) -- (13,3) -- (12,4) -- (12,0);
% Draw the boundary of piece2
\draw[black](12,0) -- (16,0) -- (18,3) -- (17,5) -- (15,3) -- (13,3) -- (12,4) -- (12,0);
% Draw the convex parts of piece2
\fill[black!15](15,7) -- (13,7) -- (16,4) -- (15,7);
\fill[black!15](16,4) -- (18,7) -- (17,9) -- (15,7) -- (16,4);
\fill[black!15](12,4) -- (16,4) -- (13,7) -- (12,8) -- (12,4);
% Draw the boundary of piece2
\draw[black](12,4) -- (16,4) -- (18,7) -- (17,9) -- (15,7) -- (13,7) -- (12,8) -- (12,4);
% Draw the convex parts of piece2
\fill[black!15](15,11) -- (13,11) -- (16,8) -- (15,11);
\fill[black!15](16,8) -- (18,11) -- (17,13) -- (15,11) -- (16,8);
\fill[black!15](12,8) -- (16,8) -- (13,11) -- (12,12) -- (12,8);
% Draw the boundary of piece2
\draw[black](12,8) -- (16,8) -- (18,11) -- (17,13) -- (15,11) -- (13,11) -- (12,12) -- (12,8);
% Draw the convex parts of piece2
\fill[black!15](16,16) -- (14,16) -- (17,13) -- (16,16);
\fill[black!15](17,13) -- (19,16) -- (18,18) -- (16,16) -- (17,13);
\fill[black!15](13,13) -- (17,13) -- (14,16) -- (13,17) -- (13,13);
% Draw the boundary of piece2
\draw[black](13,13) -- (17,13) -- (19,16) -- (18,18) -- (16,16) -- (14,16) -- (13,17) -- (13,13);
% Draw the convex parts of piece2
\fill[black!15](16,20) -- (14,20) -- (17,17) -- (16,20);
\fill[black!15](17,17) -- (19,20) -- (18,22) -- (16,20) -- (17,17);
\fill[black!15](13,17) -- (17,17) -- (14,20) -- (13,21) -- (13,17);
% Draw the boundary of piece2
\draw[black](13,17) -- (17,17) -- (19,20) -- (18,22) -- (16,20) -- (14,20) -- (13,21) -- (13,17);
% Draw the convex parts of piece2
\fill[black!15](16,24) -- (14,24) -- (17,21) -- (16,24);
\fill[black!15](17,21) -- (19,24) -- (18,26) -- (16,24) -- (17,21);
\fill[black!15](13,21) -- (17,21) -- (14,24) -- (13,25) -- (13,21);
% Draw the boundary of piece2
\draw[black](13,21) -- (17,21) -- (19,24) -- (18,26) -- (16,24) -- (14,24) -- (13,25) -- (13,21);
% Draw the convex parts of piece2
\fill[black!15](21,3) -- (19,3) -- (22,0) -- (21,3);
\fill[black!15](22,0) -- (24,3) -- (23,5) -- (21,3) -- (22,0);
\fill[black!15](18,0) -- (22,0) -- (19,3) -- (18,4) -- (18,0);
% Draw the boundary of piece2
\draw[black](18,0) -- (22,0) -- (24,3) -- (23,5) -- (21,3) -- (19,3) -- (18,4) -- (18,0);
% Draw the convex parts of piece2
\fill[black!15](21,7) -- (19,7) -- (22,4) -- (21,7);
\fill[black!15](22,4) -- (24,7) -- (23,9) -- (21,7) -- (22,4);
\fill[black!15](18,4) -- (22,4) -- (19,7) -- (18,8) -- (18,4);
% Draw the boundary of piece2
\draw[black](18,4) -- (22,4) -- (24,7) -- (23,9) -- (21,7) -- (19,7) -- (18,8) -- (18,4);
% Draw the convex parts of piece2
\fill[black!15](21,11) -- (19,11) -- (22,8) -- (21,11);
\fill[black!15](22,8) -- (24,11) -- (23,13) -- (21,11) -- (22,8);
\fill[black!15](18,8) -- (22,8) -- (19,11) -- (18,12) -- (18,8);
% Draw the boundary of piece2
\draw[black](18,8) -- (22,8) -- (24,11) -- (23,13) -- (21,11) -- (19,11) -- (18,12) -- (18,8);
% Draw the convex parts of piece2
\fill[black!15](22,16) -- (20,16) -- (23,13) -- (22,16);
\fill[black!15](23,13) -- (25,16) -- (24,18) -- (22,16) -- (23,13);
\fill[black!15](19,13) -- (23,13) -- (20,16) -- (19,17) -- (19,13);
% Draw the boundary of piece2
\draw[black](19,13) -- (23,13) -- (25,16) -- (24,18) -- (22,16) -- (20,16) -- (19,17) -- (19,13);
% Draw the convex parts of piece2
\fill[black!15](22,20) -- (20,20) -- (23,17) -- (22,20);
\fill[black!15](23,17) -- (25,20) -- (24,22) -- (22,20) -- (23,17);
\fill[black!15](19,17) -- (23,17) -- (20,20) -- (19,21) -- (19,17);
% Draw the boundary of piece2
\draw[black](19,17) -- (23,17) -- (25,20) -- (24,22) -- (22,20) -- (20,20) -- (19,21) -- (19,17);
% Draw the convex parts of piece7
\fill[black!15](26,26) -- (24,25) -- (25,24) -- (26,26);
\fill[black!15](26,26) -- (25,24) -- (27,25) -- (27,26) -- (26,26);
\fill[black!15](27,26) -- (28,27) -- (26,27) -- (26,26) -- (27,26);
\fill[black!15](23,24) -- (22,23) -- (25,23) -- (23,24);
\fill[black!15](23,24) -- (25,23) -- (25,24) -- (24,25) -- (23,24);
\fill[black!15](24,25) -- (24,27) -- (22,27) -- (23,24) -- (24,25);
% Draw the boundary of piece7
\draw[black](22,23) -- (25,23) -- (25,24) -- (27,25) -- (27,26) -- (28,27) -- (26,27) -- (26,26) -- (24,25) -- (24,27) -- (22,27) -- (23,24) -- (22,23);
% Draw the convex parts of piece7
\fill[black!15](28,7) -- (26,6) -- (27,5) -- (28,7);
\fill[black!15](28,7) -- (27,5) -- (29,6) -- (29,7) -- (28,7);
\fill[black!15](29,7) -- (30,8) -- (28,8) -- (28,7) -- (29,7);
\fill[black!15](25,5) -- (24,4) -- (27,4) -- (25,5);
\fill[black!15](25,5) -- (27,4) -- (27,5) -- (26,6) -- (25,5);
\fill[black!15](26,6) -- (26,8) -- (24,8) -- (25,5) -- (26,6);
% Draw the boundary of piece7
\draw[black](24,4) -- (27,4) -- (27,5) -- (29,6) -- (29,7) -- (30,8) -- (28,8) -- (28,7) -- (26,6) -- (26,8) -- (24,8) -- (25,5) -- (24,4);
% Draw the convex parts of piece7
\fill[black!15](28,11) -- (26,10) -- (27,9) -- (28,11);
\fill[black!15](28,11) -- (27,9) -- (29,10) -- (29,11) -- (28,11);
\fill[black!15](29,11) -- (30,12) -- (28,12) -- (28,11) -- (29,11);
\fill[black!15](25,9) -- (24,8) -- (27,8) -- (25,9);
\fill[black!15](25,9) -- (27,8) -- (27,9) -- (26,10) -- (25,9);
\fill[black!15](26,10) -- (26,12) -- (24,12) -- (25,9) -- (26,10);
% Draw the boundary of piece7
\draw[black](24,8) -- (27,8) -- (27,9) -- (29,10) -- (29,11) -- (30,12) -- (28,12) -- (28,11) -- (26,10) -- (26,12) -- (24,12) -- (25,9) -- (24,8);
% Draw the convex parts of piece7
\fill[black!15](29,15) -- (27,14) -- (28,13) -- (29,15);
\fill[black!15](29,15) -- (28,13) -- (30,14) -- (30,15) -- (29,15);
\fill[black!15](30,15) -- (31,16) -- (29,16) -- (29,15) -- (30,15);
\fill[black!15](26,13) -- (25,12) -- (28,12) -- (26,13);
\fill[black!15](26,13) -- (28,12) -- (28,13) -- (27,14) -- (26,13);
\fill[black!15](27,14) -- (27,16) -- (25,16) -- (26,13) -- (27,14);
% Draw the boundary of piece7
\draw[black](25,12) -- (28,12) -- (28,13) -- (30,14) -- (30,15) -- (31,16) -- (29,16) -- (29,15) -- (27,14) -- (27,16) -- (25,16) -- (26,13) -- (25,12);
% Draw the convex parts of piece7
\fill[black!15](29,23) -- (27,22) -- (28,21) -- (29,23);
\fill[black!15](29,23) -- (28,21) -- (30,22) -- (30,23) -- (29,23);
\fill[black!15](30,23) -- (31,24) -- (29,24) -- (29,23) -- (30,23);
\fill[black!15](26,21) -- (25,20) -- (28,20) -- (26,21);
\fill[black!15](26,21) -- (28,20) -- (28,21) -- (27,22) -- (26,21);
\fill[black!15](27,22) -- (27,24) -- (25,24) -- (26,21) -- (27,22);
% Draw the boundary of piece7
\draw[black](25,20) -- (28,20) -- (28,21) -- (30,22) -- (30,23) -- (31,24) -- (29,24) -- (29,23) -- (27,22) -- (27,24) -- (25,24) -- (26,21) -- (25,20);
% Draw the convex parts of piece7
\fill[black!15](42,26) -- (40,25) -- (41,24) -- (42,26);
\fill[black!15](42,26) -- (41,24) -- (43,25) -- (43,26) -- (42,26);
\fill[black!15](43,26) -- (44,27) -- (42,27) -- (42,26) -- (43,26);
\fill[black!15](39,24) -- (38,23) -- (41,23) -- (39,24);
\fill[black!15](39,24) -- (41,23) -- (41,24) -- (40,25) -- (39,24);
\fill[black!15](40,25) -- (40,27) -- (38,27) -- (39,24) -- (40,25);
% Draw the boundary of piece7
\draw[black](38,23) -- (41,23) -- (41,24) -- (43,25) -- (43,26) -- (44,27) -- (42,27) -- (42,26) -- (40,25) -- (40,27) -- (38,27) -- (39,24) -- (38,23);
% Draw the convex parts of piece7
\fill[black!15](44,17) -- (42,16) -- (43,15) -- (44,17);
\fill[black!15](44,17) -- (43,15) -- (45,16) -- (45,17) -- (44,17);
\fill[black!15](45,17) -- (46,18) -- (44,18) -- (44,17) -- (45,17);
\fill[black!15](41,15) -- (40,14) -- (43,14) -- (41,15);
\fill[black!15](41,15) -- (43,14) -- (43,15) -- (42,16) -- (41,15);
\fill[black!15](42,16) -- (42,18) -- (40,18) -- (41,15) -- (42,16);
% Draw the boundary of piece7
\draw[black](40,14) -- (43,14) -- (43,15) -- (45,16) -- (45,17) -- (46,18) -- (44,18) -- (44,17) -- (42,16) -- (42,18) -- (40,18) -- (41,15) -- (40,14);
% Draw the convex parts of piece5
\fill[black!15](0,10) -- (3,10) -- (4,11) -- (1,12) -- (0,10);
\fill[black!15](0,10) -- (1,12) -- (1,13) -- (0,10);
\fill[black!15](0,15) -- (0,10) -- (1,13) -- (0,15);
\fill[black!15](2,15) -- (0,15) -- (1,13) -- (2,14) -- (2,15);
\fill[black!15](2,14) -- (3,14) -- (2,15) -- (2,14);
% Draw the boundary of piece5
\draw[black](0,10) -- (3,10) -- (4,11) -- (1,12) -- (1,13) -- (2,14) -- (3,14) -- (2,15) -- (0,15) -- (0,10);
% Draw the convex parts of piece5
\fill[black!15](0,20) -- (3,20) -- (4,21) -- (1,22) -- (0,20);
\fill[black!15](0,20) -- (1,22) -- (1,23) -- (0,20);
\fill[black!15](0,25) -- (0,20) -- (1,23) -- (0,25);
\fill[black!15](2,25) -- (0,25) -- (1,23) -- (2,24) -- (2,25);
\fill[black!15](2,24) -- (3,24) -- (2,25) -- (2,24);
% Draw the boundary of piece5
\draw[black](0,20) -- (3,20) -- (4,21) -- (1,22) -- (1,23) -- (2,24) -- (3,24) -- (2,25) -- (0,25) -- (0,20);
% Draw the convex parts of piece5
\fill[black!15](3,0) -- (6,0) -- (7,1) -- (4,2) -- (3,0);
\fill[black!15](3,0) -- (4,2) -- (4,3) -- (3,0);
\fill[black!15](3,5) -- (3,0) -- (4,3) -- (3,5);
\fill[black!15](5,5) -- (3,5) -- (4,3) -- (5,4) -- (5,5);
\fill[black!15](5,4) -- (6,4) -- (5,5) -- (5,4);
% Draw the boundary of piece5
\draw[black](3,0) -- (6,0) -- (7,1) -- (4,2) -- (4,3) -- (5,4) -- (6,4) -- (5,5) -- (3,5) -- (3,0);
% Draw the convex parts of piece5
\fill[black!15](3,22) -- (6,22) -- (7,23) -- (4,24) -- (3,22);
\fill[black!15](3,22) -- (4,24) -- (4,25) -- (3,22);
\fill[black!15](3,27) -- (3,22) -- (4,25) -- (3,27);
\fill[black!15](5,27) -- (3,27) -- (4,25) -- (5,26) -- (5,27);
\fill[black!15](5,26) -- (6,26) -- (5,27) -- (5,26);
% Draw the boundary of piece5
\draw[black](3,22) -- (6,22) -- (7,23) -- (4,24) -- (4,25) -- (5,26) -- (6,26) -- (5,27) -- (3,27) -- (3,22);
% Draw the convex parts of piece5
\fill[black!15](19,21) -- (22,21) -- (23,22) -- (20,23) -- (19,21);
\fill[black!15](19,21) -- (20,23) -- (20,24) -- (19,21);
\fill[black!15](19,26) -- (19,21) -- (20,24) -- (19,26);
\fill[black!15](21,26) -- (19,26) -- (20,24) -- (21,25) -- (21,26);
\fill[black!15](21,25) -- (22,25) -- (21,26) -- (21,25);
% Draw the boundary of piece5
\draw[black](19,21) -- (22,21) -- (23,22) -- (20,23) -- (20,24) -- (21,25) -- (22,25) -- (21,26) -- (19,26) -- (19,21);
% Draw the convex parts of piece5
\fill[black!15](30,5) -- (33,5) -- (34,6) -- (31,7) -- (30,5);
\fill[black!15](30,5) -- (31,7) -- (31,8) -- (30,5);
\fill[black!15](30,10) -- (30,5) -- (31,8) -- (30,10);
\fill[black!15](32,10) -- (30,10) -- (31,8) -- (32,9) -- (32,10);
\fill[black!15](32,9) -- (33,9) -- (32,10) -- (32,9);
% Draw the boundary of piece5
\draw[black](30,5) -- (33,5) -- (34,6) -- (31,7) -- (31,8) -- (32,9) -- (33,9) -- (32,10) -- (30,10) -- (30,5);
% Draw the convex parts of piece5
\fill[black!15](30,10) -- (33,10) -- (34,11) -- (31,12) -- (30,10);
\fill[black!15](30,10) -- (31,12) -- (31,13) -- (30,10);
\fill[black!15](30,15) -- (30,10) -- (31,13) -- (30,15);
\fill[black!15](32,15) -- (30,15) -- (31,13) -- (32,14) -- (32,15);
\fill[black!15](32,14) -- (33,14) -- (32,15) -- (32,14);
% Draw the boundary of piece5
\draw[black](30,10) -- (33,10) -- (34,11) -- (31,12) -- (31,13) -- (32,14) -- (33,14) -- (32,15) -- (30,15) -- (30,10);
% Draw the convex parts of piece5
\fill[black!15](31,15) -- (34,15) -- (35,16) -- (32,17) -- (31,15);
\fill[black!15](31,15) -- (32,17) -- (32,18) -- (31,15);
\fill[black!15](31,20) -- (31,15) -- (32,18) -- (31,20);
\fill[black!15](33,20) -- (31,20) -- (32,18) -- (33,19) -- (33,20);
\fill[black!15](33,19) -- (34,19) -- (33,20) -- (33,19);
% Draw the boundary of piece5
\draw[black](31,15) -- (34,15) -- (35,16) -- (32,17) -- (32,18) -- (33,19) -- (34,19) -- (33,20) -- (31,20) -- (31,15);
% Draw the convex parts of piece5
\fill[black!15](31,20) -- (34,20) -- (35,21) -- (32,22) -- (31,20);
\fill[black!15](31,20) -- (32,22) -- (32,23) -- (31,20);
\fill[black!15](31,25) -- (31,20) -- (32,23) -- (31,25);
\fill[black!15](33,25) -- (31,25) -- (32,23) -- (33,24) -- (33,25);
\fill[black!15](33,24) -- (34,24) -- (33,25) -- (33,24);
% Draw the boundary of piece5
\draw[black](31,20) -- (34,20) -- (35,21) -- (32,22) -- (32,23) -- (33,24) -- (34,24) -- (33,25) -- (31,25) -- (31,20);
% Draw the convex parts of piece5
\fill[black!15](33,0) -- (36,0) -- (37,1) -- (34,2) -- (33,0);
\fill[black!15](33,0) -- (34,2) -- (34,3) -- (33,0);
\fill[black!15](33,5) -- (33,0) -- (34,3) -- (33,5);
\fill[black!15](35,5) -- (33,5) -- (34,3) -- (35,4) -- (35,5);
\fill[black!15](35,4) -- (36,4) -- (35,5) -- (35,4);
% Draw the boundary of piece5
\draw[black](33,0) -- (36,0) -- (37,1) -- (34,2) -- (34,3) -- (35,4) -- (36,4) -- (35,5) -- (33,5) -- (33,0);
% Draw the convex parts of piece5
\fill[black!15](34,5) -- (37,5) -- (38,6) -- (35,7) -- (34,5);
\fill[black!15](34,5) -- (35,7) -- (35,8) -- (34,5);
\fill[black!15](34,10) -- (34,5) -- (35,8) -- (34,10);
\fill[black!15](36,10) -- (34,10) -- (35,8) -- (36,9) -- (36,10);
\fill[black!15](36,9) -- (37,9) -- (36,10) -- (36,9);
% Draw the boundary of piece5
\draw[black](34,5) -- (37,5) -- (38,6) -- (35,7) -- (35,8) -- (36,9) -- (37,9) -- (36,10) -- (34,10) -- (34,5);
% Draw the convex parts of piece5
\fill[black!15](34,10) -- (37,10) -- (38,11) -- (35,12) -- (34,10);
\fill[black!15](34,10) -- (35,12) -- (35,13) -- (34,10);
\fill[black!15](34,15) -- (34,10) -- (35,13) -- (34,15);
\fill[black!15](36,15) -- (34,15) -- (35,13) -- (36,14) -- (36,15);
\fill[black!15](36,14) -- (37,14) -- (36,15) -- (36,14);
% Draw the boundary of piece5
\draw[black](34,10) -- (37,10) -- (38,11) -- (35,12) -- (35,13) -- (36,14) -- (37,14) -- (36,15) -- (34,15) -- (34,10);
% Draw the convex parts of piece5
\fill[black!15](34,22) -- (37,22) -- (38,23) -- (35,24) -- (34,22);
\fill[black!15](34,22) -- (35,24) -- (35,25) -- (34,22);
\fill[black!15](34,27) -- (34,22) -- (35,25) -- (34,27);
\fill[black!15](36,27) -- (34,27) -- (35,25) -- (36,26) -- (36,27);
\fill[black!15](36,26) -- (37,26) -- (36,27) -- (36,26);
% Draw the boundary of piece5
\draw[black](34,22) -- (37,22) -- (38,23) -- (35,24) -- (35,25) -- (36,26) -- (37,26) -- (36,27) -- (34,27) -- (34,22);
% Draw the convex parts of piece5
\fill[black!15](38,9) -- (41,9) -- (42,10) -- (39,11) -- (38,9);
\fill[black!15](38,9) -- (39,11) -- (39,12) -- (38,9);
\fill[black!15](38,14) -- (38,9) -- (39,12) -- (38,14);
\fill[black!15](40,14) -- (38,14) -- (39,12) -- (40,13) -- (40,14);
\fill[black!15](40,13) -- (41,13) -- (40,14) -- (40,13);
% Draw the boundary of piece5
\draw[black](38,9) -- (41,9) -- (42,10) -- (39,11) -- (39,12) -- (40,13) -- (41,13) -- (40,14) -- (38,14) -- (38,9);
% Draw the convex parts of piece8
\fill[black!15](5,14) -- (4,14) -- (6,13) -- (5,14);
\fill[black!15](6,15) -- (5,14) -- (6,13) -- (6,15);
\fill[black!15](6,13) -- (7,16) -- (6,16) -- (6,15) -- (6,13);
\fill[black!15](4,11) -- (6,11) -- (5,12) -- (4,11);
\fill[black!15](4,11) -- (5,12) -- (5,13) -- (4,11);
\fill[black!15](5,13) -- (6,13) -- (4,14) -- (5,13);
\fill[black!15](4,11) -- (5,13) -- (4,14) -- (4,11);
\fill[black!15](4,11) -- (4,14) -- (4,15) -- (3,15) -- (4,11);
% Draw the boundary of piece8
\draw[black](4,11) -- (6,11) -- (5,12) -- (5,13) -- (6,13) -- (7,16) -- (6,16) -- (6,15) -- (5,14) -- (4,14) -- (4,15) -- (3,15) -- (4,11);
% Draw the convex parts of piece8
\fill[black!15](26,3) -- (25,3) -- (27,2) -- (26,3);
\fill[black!15](27,4) -- (26,3) -- (27,2) -- (27,4);
\fill[black!15](27,2) -- (28,5) -- (27,5) -- (27,4) -- (27,2);
\fill[black!15](25,0) -- (27,0) -- (26,1) -- (25,0);
\fill[black!15](25,0) -- (26,1) -- (26,2) -- (25,0);
\fill[black!15](26,2) -- (27,2) -- (25,3) -- (26,2);
\fill[black!15](25,0) -- (26,2) -- (25,3) -- (25,0);
\fill[black!15](25,0) -- (25,3) -- (25,4) -- (24,4) -- (25,0);
% Draw the boundary of piece8
\draw[black](25,0) -- (27,0) -- (26,1) -- (26,2) -- (27,2) -- (28,5) -- (27,5) -- (27,4) -- (26,3) -- (25,3) -- (25,4) -- (24,4) -- (25,0);
% Draw the convex parts of piece8
\fill[black!15](27,19) -- (26,19) -- (28,18) -- (27,19);
\fill[black!15](28,20) -- (27,19) -- (28,18) -- (28,20);
\fill[black!15](28,18) -- (29,21) -- (28,21) -- (28,20) -- (28,18);
\fill[black!15](26,16) -- (28,16) -- (27,17) -- (26,16);
\fill[black!15](26,16) -- (27,17) -- (27,18) -- (26,16);
\fill[black!15](27,18) -- (28,18) -- (26,19) -- (27,18);
\fill[black!15](26,16) -- (27,18) -- (26,19) -- (26,16);
\fill[black!15](26,16) -- (26,19) -- (26,20) -- (25,20) -- (26,16);
% Draw the boundary of piece8
\draw[black](26,16) -- (28,16) -- (27,17) -- (27,18) -- (28,18) -- (29,21) -- (28,21) -- (28,20) -- (27,19) -- (26,19) -- (26,20) -- (25,20) -- (26,16);
% Draw the convex parts of piece8
\fill[black!15](30,3) -- (29,3) -- (31,2) -- (30,3);
\fill[black!15](31,4) -- (30,3) -- (31,2) -- (31,4);
\fill[black!15](31,2) -- (32,5) -- (31,5) -- (31,4) -- (31,2);
\fill[black!15](29,0) -- (31,0) -- (30,1) -- (29,0);
\fill[black!15](29,0) -- (30,1) -- (30,2) -- (29,0);
\fill[black!15](30,2) -- (31,2) -- (29,3) -- (30,2);
\fill[black!15](29,0) -- (30,2) -- (29,3) -- (29,0);
\fill[black!15](29,0) -- (29,3) -- (29,4) -- (28,4) -- (29,0);
% Draw the boundary of piece8
\draw[black](29,0) -- (31,0) -- (30,1) -- (30,2) -- (31,2) -- (32,5) -- (31,5) -- (31,4) -- (30,3) -- (29,3) -- (29,4) -- (28,4) -- (29,0);
% Draw the convex parts of piece8
\fill[black!15](41,22) -- (40,22) -- (42,21) -- (41,22);
\fill[black!15](42,23) -- (41,22) -- (42,21) -- (42,23);
\fill[black!15](42,21) -- (43,24) -- (42,24) -- (42,23) -- (42,21);
\fill[black!15](40,19) -- (42,19) -- (41,20) -- (40,19);
\fill[black!15](40,19) -- (41,20) -- (41,21) -- (40,19);
\fill[black!15](41,21) -- (42,21) -- (40,22) -- (41,21);
\fill[black!15](40,19) -- (41,21) -- (40,22) -- (40,19);
\fill[black!15](40,19) -- (40,22) -- (40,23) -- (39,23) -- (40,19);
% Draw the boundary of piece8
\draw[black](40,19) -- (42,19) -- (41,20) -- (41,21) -- (42,21) -- (43,24) -- (42,24) -- (42,23) -- (41,22) -- (40,22) -- (40,23) -- (39,23) -- (40,19);
% Draw the convex parts of piece8
\fill[black!15](44,3) -- (43,3) -- (45,2) -- (44,3);
\fill[black!15](45,4) -- (44,3) -- (45,2) -- (45,4);
\fill[black!15](45,2) -- (46,5) -- (45,5) -- (45,4) -- (45,2);
\fill[black!15](43,0) -- (45,0) -- (44,1) -- (43,0);
\fill[black!15](43,0) -- (44,1) -- (44,2) -- (43,0);
\fill[black!15](44,2) -- (45,2) -- (43,3) -- (44,2);
\fill[black!15](43,0) -- (44,2) -- (43,3) -- (43,0);
\fill[black!15](43,0) -- (43,3) -- (43,4) -- (42,4) -- (43,0);
% Draw the boundary of piece8
\draw[black](43,0) -- (45,0) -- (44,1) -- (44,2) -- (45,2) -- (46,5) -- (45,5) -- (45,4) -- (44,3) -- (43,3) -- (43,4) -- (42,4) -- (43,0);
% Draw the convex parts of piece8
\fill[black!15](44,13) -- (43,13) -- (45,12) -- (44,13);
\fill[black!15](45,14) -- (44,13) -- (45,12) -- (45,14);
\fill[black!15](45,12) -- (46,15) -- (45,15) -- (45,14) -- (45,12);
\fill[black!15](43,10) -- (45,10) -- (44,11) -- (43,10);
\fill[black!15](43,10) -- (44,11) -- (44,12) -- (43,10);
\fill[black!15](44,12) -- (45,12) -- (43,13) -- (44,12);
\fill[black!15](43,10) -- (44,12) -- (43,13) -- (43,10);
\fill[black!15](43,10) -- (43,13) -- (43,14) -- (42,14) -- (43,10);
% Draw the boundary of piece8
\draw[black](43,10) -- (45,10) -- (44,11) -- (44,12) -- (45,12) -- (46,15) -- (45,15) -- (45,14) -- (44,13) -- (43,13) -- (43,14) -- (42,14) -- (43,10);
% Draw the convex parts of piece4
\fill[black!15](2,1) -- (1,1) -- (3,0) -- (2,1);
\fill[black!15](3,0) -- (3,2) -- (2,2) -- (2,1) -- (3,0);
\fill[black!15](0,0) -- (3,0) -- (1,1) -- (0,0);
\fill[black!15](0,0) -- (1,1) -- (1,2) -- (0,2) -- (0,0);
% Draw the boundary of piece4
\draw[black](0,0) -- (3,0) -- (3,2) -- (2,2) -- (2,1) -- (1,1) -- (1,2) -- (0,2) -- (0,0);
% Draw the convex parts of piece4
\fill[black!15](2,3) -- (1,3) -- (3,2) -- (2,3);
\fill[black!15](3,2) -- (3,4) -- (2,4) -- (2,3) -- (3,2);
\fill[black!15](0,2) -- (3,2) -- (1,3) -- (0,2);
\fill[black!15](0,2) -- (1,3) -- (1,4) -- (0,4) -- (0,2);
% Draw the boundary of piece4
\draw[black](0,2) -- (3,2) -- (3,4) -- (2,4) -- (2,3) -- (1,3) -- (1,4) -- (0,4) -- (0,2);
% Draw the convex parts of piece4
\fill[black!15](2,26) -- (1,26) -- (3,25) -- (2,26);
\fill[black!15](3,25) -- (3,27) -- (2,27) -- (2,26) -- (3,25);
\fill[black!15](0,25) -- (3,25) -- (1,26) -- (0,25);
\fill[black!15](0,25) -- (1,26) -- (1,27) -- (0,27) -- (0,25);
% Draw the boundary of piece4
\draw[black](0,25) -- (3,25) -- (3,27) -- (2,27) -- (2,26) -- (1,26) -- (1,27) -- (0,27) -- (0,25);
% Draw the convex parts of piece4
\fill[black!15](8,26) -- (7,26) -- (9,25) -- (8,26);
\fill[black!15](9,25) -- (9,27) -- (8,27) -- (8,26) -- (9,25);
\fill[black!15](6,25) -- (9,25) -- (7,26) -- (6,25);
\fill[black!15](6,25) -- (7,26) -- (7,27) -- (6,27) -- (6,25);
% Draw the boundary of piece4
\draw[black](6,25) -- (9,25) -- (9,27) -- (8,27) -- (8,26) -- (7,26) -- (7,27) -- (6,27) -- (6,25);
% Draw the convex parts of piece4
\fill[black!15](9,1) -- (8,1) -- (10,0) -- (9,1);
\fill[black!15](10,0) -- (10,2) -- (9,2) -- (9,1) -- (10,0);
\fill[black!15](7,0) -- (10,0) -- (8,1) -- (7,0);
\fill[black!15](7,0) -- (8,1) -- (8,2) -- (7,2) -- (7,0);
% Draw the boundary of piece4
\draw[black](7,0) -- (10,0) -- (10,2) -- (9,2) -- (9,1) -- (8,1) -- (8,2) -- (7,2) -- (7,0);
% Draw the convex parts of piece4
\fill[black!15](11,26) -- (10,26) -- (12,25) -- (11,26);
\fill[black!15](12,25) -- (12,27) -- (11,27) -- (11,26) -- (12,25);
\fill[black!15](9,25) -- (12,25) -- (10,26) -- (9,25);
\fill[black!15](9,25) -- (10,26) -- (10,27) -- (9,27) -- (9,25);
% Draw the boundary of piece4
\draw[black](9,25) -- (12,25) -- (12,27) -- (11,27) -- (11,26) -- (10,26) -- (10,27) -- (9,27) -- (9,25);
% Draw the convex parts of piece4
\fill[black!15](16,26) -- (15,26) -- (17,25) -- (16,26);
\fill[black!15](17,25) -- (17,27) -- (16,27) -- (16,26) -- (17,25);
\fill[black!15](14,25) -- (17,25) -- (15,26) -- (14,25);
\fill[black!15](14,25) -- (15,26) -- (15,27) -- (14,27) -- (14,25);
% Draw the boundary of piece4
\draw[black](14,25) -- (17,25) -- (17,27) -- (16,27) -- (16,26) -- (15,26) -- (15,27) -- (14,27) -- (14,25);
% Draw the convex parts of piece4
\fill[black!15](30,26) -- (29,26) -- (31,25) -- (30,26);
\fill[black!15](31,25) -- (31,27) -- (30,27) -- (30,26) -- (31,25);
\fill[black!15](28,25) -- (31,25) -- (29,26) -- (28,25);
\fill[black!15](28,25) -- (29,26) -- (29,27) -- (28,27) -- (28,25);
% Draw the boundary of piece4
\draw[black](28,25) -- (31,25) -- (31,27) -- (30,27) -- (30,26) -- (29,26) -- (29,27) -- (28,27) -- (28,25);
% Draw the convex parts of piece4
\fill[black!15](33,26) -- (32,26) -- (34,25) -- (33,26);
\fill[black!15](34,25) -- (34,27) -- (33,27) -- (33,26) -- (34,25);
\fill[black!15](31,25) -- (34,25) -- (32,26) -- (31,25);
\fill[black!15](31,25) -- (32,26) -- (32,27) -- (31,27) -- (31,25);
% Draw the boundary of piece4
\draw[black](31,25) -- (34,25) -- (34,27) -- (33,27) -- (33,26) -- (32,26) -- (32,27) -- (31,27) -- (31,25);
% Draw the convex parts of piece4
\fill[black!15](37,21) -- (36,21) -- (38,20) -- (37,21);
\fill[black!15](38,20) -- (38,22) -- (37,22) -- (37,21) -- (38,20);
\fill[black!15](35,20) -- (38,20) -- (36,21) -- (35,20);
\fill[black!15](35,20) -- (36,21) -- (36,22) -- (35,22) -- (35,20);
% Draw the boundary of piece4
\draw[black](35,20) -- (38,20) -- (38,22) -- (37,22) -- (37,21) -- (36,21) -- (36,22) -- (35,22) -- (35,20);
% Draw the convex parts of piece4
\fill[black!15](45,19) -- (44,19) -- (46,18) -- (45,19);
\fill[black!15](46,18) -- (46,20) -- (45,20) -- (45,19) -- (46,18);
\fill[black!15](43,18) -- (46,18) -- (44,19) -- (43,18);
\fill[black!15](43,18) -- (44,19) -- (44,20) -- (43,20) -- (43,18);
% Draw the boundary of piece4
\draw[black](43,18) -- (46,18) -- (46,20) -- (45,20) -- (45,19) -- (44,19) -- (44,20) -- (43,20) -- (43,18);
% Draw the convex parts of piece4
\fill[black!15](45,21) -- (44,21) -- (46,20) -- (45,21);
\fill[black!15](46,20) -- (46,22) -- (45,22) -- (45,21) -- (46,20);
\fill[black!15](43,20) -- (46,20) -- (44,21) -- (43,20);
\fill[black!15](43,20) -- (44,21) -- (44,22) -- (43,22) -- (43,20);
% Draw the boundary of piece4
\draw[black](43,20) -- (46,20) -- (46,22) -- (45,22) -- (45,21) -- (44,21) -- (44,22) -- (43,22) -- (43,20);
% Draw the convex parts of piece4
\fill[black!15](45,23) -- (44,23) -- (46,22) -- (45,23);
\fill[black!15](46,22) -- (46,24) -- (45,24) -- (45,23) -- (46,22);
\fill[black!15](43,22) -- (46,22) -- (44,23) -- (43,22);
\fill[black!15](43,22) -- (44,23) -- (44,24) -- (43,24) -- (43,22);
% Draw the boundary of piece4
\draw[black](43,22) -- (46,22) -- (46,24) -- (45,24) -- (45,23) -- (44,23) -- (44,24) -- (43,24) -- (43,22);
% Draw the convex parts of piece4
\fill[black!15](45,25) -- (44,25) -- (46,24) -- (45,25);
\fill[black!15](46,24) -- (46,26) -- (45,26) -- (45,25) -- (46,24);
\fill[black!15](43,24) -- (46,24) -- (44,25) -- (43,24);
\fill[black!15](43,24) -- (44,25) -- (44,26) -- (43,26) -- (43,24);
% Draw the boundary of piece4
\draw[black](43,24) -- (46,24) -- (46,26) -- (45,26) -- (45,25) -- (44,25) -- (44,26) -- (43,26) -- (43,24);
% Draw the convex parts of piece1
\fill[black!15](0,4) -- (2,4) -- (2,6) -- (0,6) -- (0,4);
% Draw the boundary of piece1
\draw[black](0,4) -- (2,4) -- (2,6) -- (0,6) -- (0,4);
% Draw the convex parts of piece1
\fill[black!15](4,20) -- (6,20) -- (6,22) -- (4,22) -- (4,20);
% Draw the boundary of piece1
\draw[black](4,20) -- (6,20) -- (6,22) -- (4,22) -- (4,20);
% Draw the convex parts of piece1
\fill[black!15](8,23) -- (10,23) -- (10,25) -- (8,25) -- (8,23);
% Draw the boundary of piece1
\draw[black](8,23) -- (10,23) -- (10,25) -- (8,25) -- (8,23);
% Draw the convex parts of piece1
\fill[black!15](10,0) -- (12,0) -- (12,2) -- (10,2) -- (10,0);
% Draw the boundary of piece1
\draw[black](10,0) -- (12,0) -- (12,2) -- (10,2) -- (10,0);
% Draw the convex parts of piece1
\fill[black!15](12,25) -- (14,25) -- (14,27) -- (12,27) -- (12,25);
% Draw the boundary of piece1
\draw[black](12,25) -- (14,25) -- (14,27) -- (12,27) -- (12,25);
% Draw the convex parts of piece1
\fill[black!15](13,11) -- (15,11) -- (15,13) -- (13,13) -- (13,11);
% Draw the boundary of piece1
\draw[black](13,11) -- (15,11) -- (15,13) -- (13,13) -- (13,11);
% Draw the convex parts of piece1
\fill[black!15](19,11) -- (21,11) -- (21,13) -- (19,13) -- (19,11);
% Draw the boundary of piece1
\draw[black](19,11) -- (21,11) -- (21,13) -- (19,13) -- (19,11);
% Draw the convex parts of piece1
\fill[black!15](27,23) -- (29,23) -- (29,25) -- (27,25) -- (27,23);
% Draw the boundary of piece1
\draw[black](27,23) -- (29,23) -- (29,25) -- (27,25) -- (27,23);
% Draw the convex parts of piece1
\fill[black!15](31,0) -- (33,0) -- (33,2) -- (31,2) -- (31,0);
% Draw the boundary of piece1
\draw[black](31,0) -- (33,0) -- (33,2) -- (31,2) -- (31,0);
% Draw the convex parts of piece1
\fill[black!15](32,7) -- (34,7) -- (34,9) -- (32,9) -- (32,7);
% Draw the boundary of piece1
\draw[black](32,7) -- (34,7) -- (34,9) -- (32,9) -- (32,7);
% Draw the convex parts of piece1
\fill[black!15](32,12) -- (34,12) -- (34,14) -- (32,14) -- (32,12);
% Draw the boundary of piece1
\draw[black](32,12) -- (34,12) -- (34,14) -- (32,14) -- (32,12);
% Draw the convex parts of piece1
\fill[black!15](36,7) -- (38,7) -- (38,9) -- (36,9) -- (36,7);
% Draw the boundary of piece1
\draw[black](36,7) -- (38,7) -- (38,9) -- (36,9) -- (36,7);
% Draw the convex parts of piece1
\fill[black!15](36,12) -- (38,12) -- (38,14) -- (36,14) -- (36,12);
% Draw the boundary of piece1
\draw[black](36,12) -- (38,12) -- (38,14) -- (36,14) -- (36,12);
% Draw the convex parts of piece1
\fill[black!15](40,11) -- (42,11) -- (42,13) -- (40,13) -- (40,11);
% Draw the boundary of piece1
\draw[black](40,11) -- (42,11) -- (42,13) -- (40,13) -- (40,11);
% Draw the convex parts of piece6
\fill[black!15](7,18) -- (7,22) -- (6,21) -- (6,19) -- (7,18);
% Draw the boundary of piece6
\draw[black](7,18) -- (7,22) -- (6,21) -- (6,19) -- (7,18);
% Draw the convex parts of piece6
\fill[black!15](12,5) -- (12,9) -- (11,8) -- (11,6) -- (12,5);
% Draw the boundary of piece6
\draw[black](12,5) -- (12,9) -- (11,8) -- (11,6) -- (12,5);
% Draw the convex parts of piece6
\fill[black!15](19,11) -- (19,15) -- (18,14) -- (18,12) -- (19,11);
% Draw the boundary of piece6
\draw[black](19,11) -- (19,15) -- (18,14) -- (18,12) -- (19,11);
% Draw the convex parts of piece6
\fill[black!15](29,15) -- (29,19) -- (28,18) -- (28,16) -- (29,15);
% Draw the boundary of piece6
\draw[black](29,15) -- (29,19) -- (28,18) -- (28,16) -- (29,15);
% Draw the convex parts of piece6
\fill[black!15](30,4) -- (30,8) -- (29,7) -- (29,5) -- (30,4);
% Draw the boundary of piece6
\draw[black](30,4) -- (30,8) -- (29,7) -- (29,5) -- (30,4);
% Draw the convex parts of piece6
\fill[black!15](30,8) -- (30,12) -- (29,11) -- (29,9) -- (30,8);
% Draw the boundary of piece6
\draw[black](30,8) -- (30,12) -- (29,11) -- (29,9) -- (30,8);
% Draw the convex parts of piece6
\fill[black!15](30,18) -- (30,22) -- (29,21) -- (29,19) -- (30,18);
% Draw the boundary of piece6
\draw[black](30,18) -- (30,22) -- (29,21) -- (29,19) -- (30,18);
% Draw the convex parts of piece6
\fill[black!15](31,16) -- (31,20) -- (30,19) -- (30,17) -- (31,16);
% Draw the boundary of piece6
\draw[black](31,16) -- (31,20) -- (30,19) -- (30,17) -- (31,16);
% Draw the convex parts of piece6
\fill[black!15](31,20) -- (31,24) -- (30,23) -- (30,21) -- (31,20);
% Draw the boundary of piece6
\draw[black](31,20) -- (31,24) -- (30,23) -- (30,21) -- (31,20);
% Draw the convex parts of piece6
\fill[black!15](38,23) -- (38,27) -- (37,26) -- (37,24) -- (38,23);
% Draw the boundary of piece6
\draw[black](38,23) -- (38,27) -- (37,26) -- (37,24) -- (38,23);
% Draw the convex parts of piece6
\fill[black!15](39,19) -- (39,23) -- (38,22) -- (38,20) -- (39,19);
% Draw the boundary of piece6
\draw[black](39,19) -- (39,23) -- (38,22) -- (38,20) -- (39,19);
% Draw the convex parts of piece6
\fill[black!15](43,18) -- (43,22) -- (42,21) -- (42,19) -- (43,18);
% Draw the boundary of piece6
\draw[black](43,18) -- (43,22) -- (42,21) -- (42,19) -- (43,18);
% Draw the convex parts of piece6
\fill[black!15](46,5) -- (46,9) -- (45,8) -- (45,6) -- (46,5);
% Draw the boundary of piece6
\draw[black](46,5) -- (46,9) -- (45,8) -- (45,6) -- (46,5);
% Draw the convex parts of piece6
\fill[black!15](46,9) -- (46,13) -- (45,12) -- (45,10) -- (46,9);
% Draw the boundary of piece6
\draw[black](46,9) -- (46,13) -- (45,12) -- (45,10) -- (46,9);
% Draw the bounding box of the solution
\draw[black](0,0) -- (46,0) -- (46,27)-- (0,27) -- (0,0);
\draw (23,-2) node {artif7\_14};
\draw (23,-4.5) node {L = 46};
% Draw the bounding box of the solution
\end{tikzpicture}
% SOLUTION for problem instance = artif
% Scaling factor for coordinates = 1
\begin{tikzpicture}[scale = 0.10]
% Draw the convex parts of piece3
\fill[black!15](9,2) -- (12,4) -- (11,6) -- (8,7) -- (5,5) -- (9,2);
% Draw the boundary of piece3
\draw[black](9,2) -- (12,4) -- (11,6) -- (8,7) -- (5,5) -- (9,2);
% Draw the convex parts of piece3
\fill[black!15](9,7) -- (12,9) -- (11,11) -- (8,12) -- (5,10) -- (9,7);
% Draw the boundary of piece3
\draw[black](9,7) -- (12,9) -- (11,11) -- (8,12) -- (5,10) -- (9,7);
% Draw the convex parts of piece3
\fill[black!15](10,17) -- (13,19) -- (12,21) -- (9,22) -- (6,20) -- (10,17);
% Draw the boundary of piece3
\draw[black](10,17) -- (13,19) -- (12,21) -- (9,22) -- (6,20) -- (10,17);
% Draw the convex parts of piece3
\fill[black!15](31,22) -- (34,24) -- (33,26) -- (30,27) -- (27,25) -- (31,22);
% Draw the boundary of piece3
\draw[black](31,22) -- (34,24) -- (33,26) -- (30,27) -- (27,25) -- (31,22);
% Draw the convex parts of piece3
\fill[black!15](33,11) -- (36,13) -- (35,15) -- (32,16) -- (29,14) -- (33,11);
% Draw the boundary of piece3
\draw[black](33,11) -- (36,13) -- (35,15) -- (32,16) -- (29,14) -- (33,11);
% Draw the convex parts of piece3
\fill[black!15](41,22) -- (44,24) -- (43,26) -- (40,27) -- (37,25) -- (41,22);
% Draw the boundary of piece3
\draw[black](41,22) -- (44,24) -- (43,26) -- (40,27) -- (37,25) -- (41,22);
% Draw the convex parts of piece3
\fill[black!15](42,9) -- (45,11) -- (44,13) -- (41,14) -- (38,12) -- (42,9);
% Draw the boundary of piece3
\draw[black](42,9) -- (45,11) -- (44,13) -- (41,14) -- (38,12) -- (42,9);
% Draw the convex parts of piece3
\fill[black!15](46,0) -- (49,2) -- (48,4) -- (45,5) -- (42,3) -- (46,0);
% Draw the boundary of piece3
\draw[black](46,0) -- (49,2) -- (48,4) -- (45,5) -- (42,3) -- (46,0);
% Draw the convex parts of piece2
\fill[black!15](10,25) -- (8,25) -- (11,22) -- (10,25);
\fill[black!15](11,22) -- (13,25) -- (12,27) -- (10,25) -- (11,22);
\fill[black!15](7,22) -- (11,22) -- (8,25) -- (7,26) -- (7,22);
% Draw the boundary of piece2
\draw[black](7,22) -- (11,22) -- (13,25) -- (12,27) -- (10,25) -- (8,25) -- (7,26) -- (7,22);
% Draw the convex parts of piece2
\fill[black!15](15,3) -- (13,3) -- (16,0) -- (15,3);
\fill[black!15](16,0) -- (18,3) -- (17,5) -- (15,3) -- (16,0);
\fill[black!15](12,0) -- (16,0) -- (13,3) -- (12,4) -- (12,0);
% Draw the boundary of piece2
\draw[black](12,0) -- (16,0) -- (18,3) -- (17,5) -- (15,3) -- (13,3) -- (12,4) -- (12,0);
% Draw the convex parts of piece2
\fill[black!15](15,7) -- (13,7) -- (16,4) -- (15,7);
\fill[black!15](16,4) -- (18,7) -- (17,9) -- (15,7) -- (16,4);
\fill[black!15](12,4) -- (16,4) -- (13,7) -- (12,8) -- (12,4);
% Draw the boundary of piece2
\draw[black](12,4) -- (16,4) -- (18,7) -- (17,9) -- (15,7) -- (13,7) -- (12,8) -- (12,4);
% Draw the convex parts of piece2
\fill[black!15](15,11) -- (13,11) -- (16,8) -- (15,11);
\fill[black!15](16,8) -- (18,11) -- (17,13) -- (15,11) -- (16,8);
\fill[black!15](12,8) -- (16,8) -- (13,11) -- (12,12) -- (12,8);
% Draw the boundary of piece2
\draw[black](12,8) -- (16,8) -- (18,11) -- (17,13) -- (15,11) -- (13,11) -- (12,12) -- (12,8);
% Draw the convex parts of piece2
\fill[black!15](15,15) -- (13,15) -- (16,12) -- (15,15);
\fill[black!15](16,12) -- (18,15) -- (17,17) -- (15,15) -- (16,12);
\fill[black!15](12,12) -- (16,12) -- (13,15) -- (12,16) -- (12,12);
% Draw the boundary of piece2
\draw[black](12,12) -- (16,12) -- (18,15) -- (17,17) -- (15,15) -- (13,15) -- (12,16) -- (12,12);
% Draw the convex parts of piece2
\fill[black!15](16,20) -- (14,20) -- (17,17) -- (16,20);
\fill[black!15](17,17) -- (19,20) -- (18,22) -- (16,20) -- (17,17);
\fill[black!15](13,17) -- (17,17) -- (14,20) -- (13,21) -- (13,17);
% Draw the boundary of piece2
\draw[black](13,17) -- (17,17) -- (19,20) -- (18,22) -- (16,20) -- (14,20) -- (13,21) -- (13,17);
% Draw the convex parts of piece2
\fill[black!15](16,24) -- (14,24) -- (17,21) -- (16,24);
\fill[black!15](17,21) -- (19,24) -- (18,26) -- (16,24) -- (17,21);
\fill[black!15](13,21) -- (17,21) -- (14,24) -- (13,25) -- (13,21);
% Draw the boundary of piece2
\draw[black](13,21) -- (17,21) -- (19,24) -- (18,26) -- (16,24) -- (14,24) -- (13,25) -- (13,21);
% Draw the convex parts of piece2
\fill[black!15](21,3) -- (19,3) -- (22,0) -- (21,3);
\fill[black!15](22,0) -- (24,3) -- (23,5) -- (21,3) -- (22,0);
\fill[black!15](18,0) -- (22,0) -- (19,3) -- (18,4) -- (18,0);
% Draw the boundary of piece2
\draw[black](18,0) -- (22,0) -- (24,3) -- (23,5) -- (21,3) -- (19,3) -- (18,4) -- (18,0);
% Draw the convex parts of piece2
\fill[black!15](21,7) -- (19,7) -- (22,4) -- (21,7);
\fill[black!15](22,4) -- (24,7) -- (23,9) -- (21,7) -- (22,4);
\fill[black!15](18,4) -- (22,4) -- (19,7) -- (18,8) -- (18,4);
% Draw the boundary of piece2
\draw[black](18,4) -- (22,4) -- (24,7) -- (23,9) -- (21,7) -- (19,7) -- (18,8) -- (18,4);
% Draw the convex parts of piece2
\fill[black!15](21,11) -- (19,11) -- (22,8) -- (21,11);
\fill[black!15](22,8) -- (24,11) -- (23,13) -- (21,11) -- (22,8);
\fill[black!15](18,8) -- (22,8) -- (19,11) -- (18,12) -- (18,8);
% Draw the boundary of piece2
\draw[black](18,8) -- (22,8) -- (24,11) -- (23,13) -- (21,11) -- (19,11) -- (18,12) -- (18,8);
% Draw the convex parts of piece2
\fill[black!15](21,15) -- (19,15) -- (22,12) -- (21,15);
\fill[black!15](22,12) -- (24,15) -- (23,17) -- (21,15) -- (22,12);
\fill[black!15](18,12) -- (22,12) -- (19,15) -- (18,16) -- (18,12);
% Draw the boundary of piece2
\draw[black](18,12) -- (22,12) -- (24,15) -- (23,17) -- (21,15) -- (19,15) -- (18,16) -- (18,12);
% Draw the convex parts of piece2
\fill[black!15](27,3) -- (25,3) -- (28,0) -- (27,3);
\fill[black!15](28,0) -- (30,3) -- (29,5) -- (27,3) -- (28,0);
\fill[black!15](24,0) -- (28,0) -- (25,3) -- (24,4) -- (24,0);
% Draw the boundary of piece2
\draw[black](24,0) -- (28,0) -- (30,3) -- (29,5) -- (27,3) -- (25,3) -- (24,4) -- (24,0);
% Draw the convex parts of piece2
\fill[black!15](28,19) -- (26,19) -- (29,16) -- (28,19);
\fill[black!15](29,16) -- (31,19) -- (30,21) -- (28,19) -- (29,16);
\fill[black!15](25,16) -- (29,16) -- (26,19) -- (25,20) -- (25,16);
% Draw the boundary of piece2
\draw[black](25,16) -- (29,16) -- (31,19) -- (30,21) -- (28,19) -- (26,19) -- (25,20) -- (25,16);
% Draw the convex parts of piece2
\fill[black!15](38,8) -- (36,8) -- (39,5) -- (38,8);
\fill[black!15](39,5) -- (41,8) -- (40,10) -- (38,8) -- (39,5);
\fill[black!15](35,5) -- (39,5) -- (36,8) -- (35,9) -- (35,5);
% Draw the boundary of piece2
\draw[black](35,5) -- (39,5) -- (41,8) -- (40,10) -- (38,8) -- (36,8) -- (35,9) -- (35,5);
% Draw the convex parts of piece2
\fill[black!15](42,21) -- (40,21) -- (43,18) -- (42,21);
\fill[black!15](43,18) -- (45,21) -- (44,23) -- (42,21) -- (43,18);
\fill[black!15](39,18) -- (43,18) -- (40,21) -- (39,22) -- (39,18);
% Draw the boundary of piece2
\draw[black](39,18) -- (43,18) -- (45,21) -- (44,23) -- (42,21) -- (40,21) -- (39,22) -- (39,18);
% Draw the convex parts of piece7
\fill[black!15](7,16) -- (5,15) -- (6,14) -- (7,16);
\fill[black!15](7,16) -- (6,14) -- (8,15) -- (8,16) -- (7,16);
\fill[black!15](8,16) -- (9,17) -- (7,17) -- (7,16) -- (8,16);
\fill[black!15](4,14) -- (3,13) -- (6,13) -- (4,14);
\fill[black!15](4,14) -- (6,13) -- (6,14) -- (5,15) -- (4,14);
\fill[black!15](5,15) -- (5,17) -- (3,17) -- (4,14) -- (5,15);
% Draw the boundary of piece7
\draw[black](3,13) -- (6,13) -- (6,14) -- (8,15) -- (8,16) -- (9,17) -- (7,17) -- (7,16) -- (5,15) -- (5,17) -- (3,17) -- (4,14) -- (3,13);
% Draw the convex parts of piece7
\fill[black!15](23,19) -- (21,18) -- (22,17) -- (23,19);
\fill[black!15](23,19) -- (22,17) -- (24,18) -- (24,19) -- (23,19);
\fill[black!15](24,19) -- (25,20) -- (23,20) -- (23,19) -- (24,19);
\fill[black!15](20,17) -- (19,16) -- (22,16) -- (20,17);
\fill[black!15](20,17) -- (22,16) -- (22,17) -- (21,18) -- (20,17);
\fill[black!15](21,18) -- (21,20) -- (19,20) -- (20,17) -- (21,18);
% Draw the boundary of piece7
\draw[black](19,16) -- (22,16) -- (22,17) -- (24,18) -- (24,19) -- (25,20) -- (23,20) -- (23,19) -- (21,18) -- (21,20) -- (19,20) -- (20,17) -- (19,16);
% Draw the convex parts of piece7
\fill[black!15](26,26) -- (24,25) -- (25,24) -- (26,26);
\fill[black!15](26,26) -- (25,24) -- (27,25) -- (27,26) -- (26,26);
\fill[black!15](27,26) -- (28,27) -- (26,27) -- (26,26) -- (27,26);
\fill[black!15](23,24) -- (22,23) -- (25,23) -- (23,24);
\fill[black!15](23,24) -- (25,23) -- (25,24) -- (24,25) -- (23,24);
\fill[black!15](24,25) -- (24,27) -- (22,27) -- (23,24) -- (24,25);
% Draw the boundary of piece7
\draw[black](22,23) -- (25,23) -- (25,24) -- (27,25) -- (27,26) -- (28,27) -- (26,27) -- (26,26) -- (24,25) -- (24,27) -- (22,27) -- (23,24) -- (22,23);
% Draw the convex parts of piece7
\fill[black!15](28,7) -- (26,6) -- (27,5) -- (28,7);
\fill[black!15](28,7) -- (27,5) -- (29,6) -- (29,7) -- (28,7);
\fill[black!15](29,7) -- (30,8) -- (28,8) -- (28,7) -- (29,7);
\fill[black!15](25,5) -- (24,4) -- (27,4) -- (25,5);
\fill[black!15](25,5) -- (27,4) -- (27,5) -- (26,6) -- (25,5);
\fill[black!15](26,6) -- (26,8) -- (24,8) -- (25,5) -- (26,6);
% Draw the boundary of piece7
\draw[black](24,4) -- (27,4) -- (27,5) -- (29,6) -- (29,7) -- (30,8) -- (28,8) -- (28,7) -- (26,6) -- (26,8) -- (24,8) -- (25,5) -- (24,4);
% Draw the convex parts of piece7
\fill[black!15](28,15) -- (26,14) -- (27,13) -- (28,15);
\fill[black!15](28,15) -- (27,13) -- (29,14) -- (29,15) -- (28,15);
\fill[black!15](29,15) -- (30,16) -- (28,16) -- (28,15) -- (29,15);
\fill[black!15](25,13) -- (24,12) -- (27,12) -- (25,13);
\fill[black!15](25,13) -- (27,12) -- (27,13) -- (26,14) -- (25,13);
\fill[black!15](26,14) -- (26,16) -- (24,16) -- (25,13) -- (26,14);
% Draw the boundary of piece7
\draw[black](24,12) -- (27,12) -- (27,13) -- (29,14) -- (29,15) -- (30,16) -- (28,16) -- (28,15) -- (26,14) -- (26,16) -- (24,16) -- (25,13) -- (24,12);
% Draw the convex parts of piece7
\fill[black!15](33,10) -- (31,9) -- (32,8) -- (33,10);
\fill[black!15](33,10) -- (32,8) -- (34,9) -- (34,10) -- (33,10);
\fill[black!15](34,10) -- (35,11) -- (33,11) -- (33,10) -- (34,10);
\fill[black!15](30,8) -- (29,7) -- (32,7) -- (30,8);
\fill[black!15](30,8) -- (32,7) -- (32,8) -- (31,9) -- (30,8);
\fill[black!15](31,9) -- (31,11) -- (29,11) -- (30,8) -- (31,9);
% Draw the boundary of piece7
\draw[black](29,7) -- (32,7) -- (32,8) -- (34,9) -- (34,10) -- (35,11) -- (33,11) -- (33,10) -- (31,9) -- (31,11) -- (29,11) -- (30,8) -- (29,7);
% Draw the convex parts of piece7
\fill[black!15](46,17) -- (44,16) -- (45,15) -- (46,17);
\fill[black!15](46,17) -- (45,15) -- (47,16) -- (47,17) -- (46,17);
\fill[black!15](47,17) -- (48,18) -- (46,18) -- (46,17) -- (47,17);
\fill[black!15](43,15) -- (42,14) -- (45,14) -- (43,15);
\fill[black!15](43,15) -- (45,14) -- (45,15) -- (44,16) -- (43,15);
\fill[black!15](44,16) -- (44,18) -- (42,18) -- (43,15) -- (44,16);
% Draw the boundary of piece7
\draw[black](42,14) -- (45,14) -- (45,15) -- (47,16) -- (47,17) -- (48,18) -- (46,18) -- (46,17) -- (44,16) -- (44,18) -- (42,18) -- (43,15) -- (42,14);
% Draw the convex parts of piece7
\fill[black!15](47,26) -- (45,25) -- (46,24) -- (47,26);
\fill[black!15](47,26) -- (46,24) -- (48,25) -- (48,26) -- (47,26);
\fill[black!15](48,26) -- (49,27) -- (47,27) -- (47,26) -- (48,26);
\fill[black!15](44,24) -- (43,23) -- (46,23) -- (44,24);
\fill[black!15](44,24) -- (46,23) -- (46,24) -- (45,25) -- (44,24);
\fill[black!15](45,25) -- (45,27) -- (43,27) -- (44,24) -- (45,25);
% Draw the boundary of piece7
\draw[black](43,23) -- (46,23) -- (46,24) -- (48,25) -- (48,26) -- (49,27) -- (47,27) -- (47,26) -- (45,25) -- (45,27) -- (43,27) -- (44,24) -- (43,23);
% Draw the convex parts of piece5
\fill[black!15](0,8) -- (3,8) -- (4,9) -- (1,10) -- (0,8);
\fill[black!15](0,8) -- (1,10) -- (1,11) -- (0,8);
\fill[black!15](0,13) -- (0,8) -- (1,11) -- (0,13);
\fill[black!15](2,13) -- (0,13) -- (1,11) -- (2,12) -- (2,13);
\fill[black!15](2,12) -- (3,12) -- (2,13) -- (2,12);
% Draw the boundary of piece5
\draw[black](0,8) -- (3,8) -- (4,9) -- (1,10) -- (1,11) -- (2,12) -- (3,12) -- (2,13) -- (0,13) -- (0,8);
% Draw the convex parts of piece5
\fill[black!15](0,13) -- (3,13) -- (4,14) -- (1,15) -- (0,13);
\fill[black!15](0,13) -- (1,15) -- (1,16) -- (0,13);
\fill[black!15](0,18) -- (0,13) -- (1,16) -- (0,18);
\fill[black!15](2,18) -- (0,18) -- (1,16) -- (2,17) -- (2,18);
\fill[black!15](2,17) -- (3,17) -- (2,18) -- (2,17);
% Draw the boundary of piece5
\draw[black](0,13) -- (3,13) -- (4,14) -- (1,15) -- (1,16) -- (2,17) -- (3,17) -- (2,18) -- (0,18) -- (0,13);
% Draw the convex parts of piece5
\fill[black!15](0,18) -- (3,18) -- (4,19) -- (1,20) -- (0,18);
\fill[black!15](0,18) -- (1,20) -- (1,21) -- (0,18);
\fill[black!15](0,23) -- (0,18) -- (1,21) -- (0,23);
\fill[black!15](2,23) -- (0,23) -- (1,21) -- (2,22) -- (2,23);
\fill[black!15](2,22) -- (3,22) -- (2,23) -- (2,22);
% Draw the boundary of piece5
\draw[black](0,18) -- (3,18) -- (4,19) -- (1,20) -- (1,21) -- (2,22) -- (3,22) -- (2,23) -- (0,23) -- (0,18);
% Draw the convex parts of piece5
\fill[black!15](3,2) -- (6,2) -- (7,3) -- (4,4) -- (3,2);
\fill[black!15](3,2) -- (4,4) -- (4,5) -- (3,2);
\fill[black!15](3,7) -- (3,2) -- (4,5) -- (3,7);
\fill[black!15](5,7) -- (3,7) -- (4,5) -- (5,6) -- (5,7);
\fill[black!15](5,6) -- (6,6) -- (5,7) -- (5,6);
% Draw the boundary of piece5
\draw[black](3,2) -- (6,2) -- (7,3) -- (4,4) -- (4,5) -- (5,6) -- (6,6) -- (5,7) -- (3,7) -- (3,2);
% Draw the convex parts of piece5
\fill[black!15](3,22) -- (6,22) -- (7,23) -- (4,24) -- (3,22);
\fill[black!15](3,22) -- (4,24) -- (4,25) -- (3,22);
\fill[black!15](3,27) -- (3,22) -- (4,25) -- (3,27);
\fill[black!15](5,27) -- (3,27) -- (4,25) -- (5,26) -- (5,27);
\fill[black!15](5,26) -- (6,26) -- (5,27) -- (5,26);
% Draw the boundary of piece5
\draw[black](3,22) -- (6,22) -- (7,23) -- (4,24) -- (4,25) -- (5,26) -- (6,26) -- (5,27) -- (3,27) -- (3,22);
% Draw the convex parts of piece5
\fill[black!15](4,17) -- (7,17) -- (8,18) -- (5,19) -- (4,17);
\fill[black!15](4,17) -- (5,19) -- (5,20) -- (4,17);
\fill[black!15](4,22) -- (4,17) -- (5,20) -- (4,22);
\fill[black!15](6,22) -- (4,22) -- (5,20) -- (6,21) -- (6,22);
\fill[black!15](6,21) -- (7,21) -- (6,22) -- (6,21);
% Draw the boundary of piece5
\draw[black](4,17) -- (7,17) -- (8,18) -- (5,19) -- (5,20) -- (6,21) -- (7,21) -- (6,22) -- (4,22) -- (4,17);
% Draw the convex parts of piece5
\fill[black!15](19,20) -- (22,20) -- (23,21) -- (20,22) -- (19,20);
\fill[black!15](19,20) -- (20,22) -- (20,23) -- (19,20);
\fill[black!15](19,25) -- (19,20) -- (20,23) -- (19,25);
\fill[black!15](21,25) -- (19,25) -- (20,23) -- (21,24) -- (21,25);
\fill[black!15](21,24) -- (22,24) -- (21,25) -- (21,24);
% Draw the boundary of piece5
\draw[black](19,20) -- (22,20) -- (23,21) -- (20,22) -- (20,23) -- (21,24) -- (22,24) -- (21,25) -- (19,25) -- (19,20);
% Draw the convex parts of piece5
\fill[black!15](32,0) -- (35,0) -- (36,1) -- (33,2) -- (32,0);
\fill[black!15](32,0) -- (33,2) -- (33,3) -- (32,0);
\fill[black!15](32,5) -- (32,0) -- (33,3) -- (32,5);
\fill[black!15](34,5) -- (32,5) -- (33,3) -- (34,4) -- (34,5);
\fill[black!15](34,4) -- (35,4) -- (34,5) -- (34,4);
% Draw the boundary of piece5
\draw[black](32,0) -- (35,0) -- (36,1) -- (33,2) -- (33,3) -- (34,4) -- (35,4) -- (34,5) -- (32,5) -- (32,0);
% Draw the convex parts of piece5
\fill[black!15](34,22) -- (37,22) -- (38,23) -- (35,24) -- (34,22);
\fill[black!15](34,22) -- (35,24) -- (35,25) -- (34,22);
\fill[black!15](34,27) -- (34,22) -- (35,25) -- (34,27);
\fill[black!15](36,27) -- (34,27) -- (35,25) -- (36,26) -- (36,27);
\fill[black!15](36,26) -- (37,26) -- (36,27) -- (36,26);
% Draw the boundary of piece5
\draw[black](34,22) -- (37,22) -- (38,23) -- (35,24) -- (35,25) -- (36,26) -- (37,26) -- (36,27) -- (34,27) -- (34,22);
% Draw the convex parts of piece5
\fill[black!15](36,0) -- (39,0) -- (40,1) -- (37,2) -- (36,0);
\fill[black!15](36,0) -- (37,2) -- (37,3) -- (36,0);
\fill[black!15](36,5) -- (36,0) -- (37,3) -- (36,5);
\fill[black!15](38,5) -- (36,5) -- (37,3) -- (38,4) -- (38,5);
\fill[black!15](38,4) -- (39,4) -- (38,5) -- (38,4);
% Draw the boundary of piece5
\draw[black](36,0) -- (39,0) -- (40,1) -- (37,2) -- (37,3) -- (38,4) -- (39,4) -- (38,5) -- (36,5) -- (36,0);
% Draw the convex parts of piece5
\fill[black!15](36,9) -- (39,9) -- (40,10) -- (37,11) -- (36,9);
\fill[black!15](36,9) -- (37,11) -- (37,12) -- (36,9);
\fill[black!15](36,14) -- (36,9) -- (37,12) -- (36,14);
\fill[black!15](38,14) -- (36,14) -- (37,12) -- (38,13) -- (38,14);
\fill[black!15](38,13) -- (39,13) -- (38,14) -- (38,13);
% Draw the boundary of piece5
\draw[black](36,9) -- (39,9) -- (40,10) -- (37,11) -- (37,12) -- (38,13) -- (39,13) -- (38,14) -- (36,14) -- (36,9);
% Draw the convex parts of piece5
\fill[black!15](40,0) -- (43,0) -- (44,1) -- (41,2) -- (40,0);
\fill[black!15](40,0) -- (41,2) -- (41,3) -- (40,0);
\fill[black!15](40,5) -- (40,0) -- (41,3) -- (40,5);
\fill[black!15](42,5) -- (40,5) -- (41,3) -- (42,4) -- (42,5);
\fill[black!15](42,4) -- (43,4) -- (42,5) -- (42,4);
% Draw the boundary of piece5
\draw[black](40,0) -- (43,0) -- (44,1) -- (41,2) -- (41,3) -- (42,4) -- (43,4) -- (42,5) -- (40,5) -- (40,0);
% Draw the convex parts of piece5
\fill[black!15](45,5) -- (48,5) -- (49,6) -- (46,7) -- (45,5);
\fill[black!15](45,5) -- (46,7) -- (46,8) -- (45,5);
\fill[black!15](45,10) -- (45,5) -- (46,8) -- (45,10);
\fill[black!15](47,10) -- (45,10) -- (46,8) -- (47,9) -- (47,10);
\fill[black!15](47,9) -- (48,9) -- (47,10) -- (47,9);
% Draw the boundary of piece5
\draw[black](45,5) -- (48,5) -- (49,6) -- (46,7) -- (46,8) -- (47,9) -- (48,9) -- (47,10) -- (45,10) -- (45,5);
% Draw the convex parts of piece5
\fill[black!15](45,10) -- (48,10) -- (49,11) -- (46,12) -- (45,10);
\fill[black!15](45,10) -- (46,12) -- (46,13) -- (45,10);
\fill[black!15](45,15) -- (45,10) -- (46,13) -- (45,15);
\fill[black!15](47,15) -- (45,15) -- (46,13) -- (47,14) -- (47,15);
\fill[black!15](47,14) -- (48,14) -- (47,15) -- (47,14);
% Draw the boundary of piece5
\draw[black](45,10) -- (48,10) -- (49,11) -- (46,12) -- (46,13) -- (47,14) -- (48,14) -- (47,15) -- (45,15) -- (45,10);
% Draw the convex parts of piece5
\fill[black!15](45,18) -- (48,18) -- (49,19) -- (46,20) -- (45,18);
\fill[black!15](45,18) -- (46,20) -- (46,21) -- (45,18);
\fill[black!15](45,23) -- (45,18) -- (46,21) -- (45,23);
\fill[black!15](47,23) -- (45,23) -- (46,21) -- (47,22) -- (47,23);
\fill[black!15](47,22) -- (48,22) -- (47,23) -- (47,22);
% Draw the boundary of piece5
\draw[black](45,18) -- (48,18) -- (49,19) -- (46,20) -- (46,21) -- (47,22) -- (48,22) -- (47,23) -- (45,23) -- (45,18);
% Draw the convex parts of piece8
\fill[black!15](5,12) -- (4,12) -- (6,11) -- (5,12);
\fill[black!15](6,13) -- (5,12) -- (6,11) -- (6,13);
\fill[black!15](6,11) -- (7,14) -- (6,14) -- (6,13) -- (6,11);
\fill[black!15](4,9) -- (6,9) -- (5,10) -- (4,9);
\fill[black!15](4,9) -- (5,10) -- (5,11) -- (4,9);
\fill[black!15](5,11) -- (6,11) -- (4,12) -- (5,11);
\fill[black!15](4,9) -- (5,11) -- (4,12) -- (4,9);
\fill[black!15](4,9) -- (4,12) -- (4,13) -- (3,13) -- (4,9);
% Draw the boundary of piece8
\draw[black](4,9) -- (6,9) -- (5,10) -- (5,11) -- (6,11) -- (7,14) -- (6,14) -- (6,13) -- (5,12) -- (4,12) -- (4,13) -- (3,13) -- (4,9);
% Draw the convex parts of piece8
\fill[black!15](10,15) -- (9,15) -- (11,14) -- (10,15);
\fill[black!15](11,16) -- (10,15) -- (11,14) -- (11,16);
\fill[black!15](11,14) -- (12,17) -- (11,17) -- (11,16) -- (11,14);
\fill[black!15](9,12) -- (11,12) -- (10,13) -- (9,12);
\fill[black!15](9,12) -- (10,13) -- (10,14) -- (9,12);
\fill[black!15](10,14) -- (11,14) -- (9,15) -- (10,14);
\fill[black!15](9,12) -- (10,14) -- (9,15) -- (9,12);
\fill[black!15](9,12) -- (9,15) -- (9,16) -- (8,16) -- (9,12);
% Draw the boundary of piece8
\draw[black](9,12) -- (11,12) -- (10,13) -- (10,14) -- (11,14) -- (12,17) -- (11,17) -- (11,16) -- (10,15) -- (9,15) -- (9,16) -- (8,16) -- (9,12);
% Draw the convex parts of piece8
\fill[black!15](26,11) -- (25,11) -- (27,10) -- (26,11);
\fill[black!15](27,12) -- (26,11) -- (27,10) -- (27,12);
\fill[black!15](27,10) -- (28,13) -- (27,13) -- (27,12) -- (27,10);
\fill[black!15](25,8) -- (27,8) -- (26,9) -- (25,8);
\fill[black!15](25,8) -- (26,9) -- (26,10) -- (25,8);
\fill[black!15](26,10) -- (27,10) -- (25,11) -- (26,10);
\fill[black!15](25,8) -- (26,10) -- (25,11) -- (25,8);
\fill[black!15](25,8) -- (25,11) -- (25,12) -- (24,12) -- (25,8);
% Draw the boundary of piece8
\draw[black](25,8) -- (27,8) -- (26,9) -- (26,10) -- (27,10) -- (28,13) -- (27,13) -- (27,12) -- (26,11) -- (25,11) -- (25,12) -- (24,12) -- (25,8);
% Draw the convex parts of piece8
\fill[black!15](31,6) -- (30,6) -- (32,5) -- (31,6);
\fill[black!15](32,7) -- (31,6) -- (32,5) -- (32,7);
\fill[black!15](32,5) -- (33,8) -- (32,8) -- (32,7) -- (32,5);
\fill[black!15](30,3) -- (32,3) -- (31,4) -- (30,3);
\fill[black!15](30,3) -- (31,4) -- (31,5) -- (30,3);
\fill[black!15](31,5) -- (32,5) -- (30,6) -- (31,5);
\fill[black!15](30,3) -- (31,5) -- (30,6) -- (30,3);
\fill[black!15](30,3) -- (30,6) -- (30,7) -- (29,7) -- (30,3);
% Draw the boundary of piece8
\draw[black](30,3) -- (32,3) -- (31,4) -- (31,5) -- (32,5) -- (33,8) -- (32,8) -- (32,7) -- (31,6) -- (30,6) -- (30,7) -- (29,7) -- (30,3);
% Draw the convex parts of piece8
\fill[black!15](33,19) -- (32,19) -- (34,18) -- (33,19);
\fill[black!15](34,20) -- (33,19) -- (34,18) -- (34,20);
\fill[black!15](34,18) -- (35,21) -- (34,21) -- (34,20) -- (34,18);
\fill[black!15](32,16) -- (34,16) -- (33,17) -- (32,16);
\fill[black!15](32,16) -- (33,17) -- (33,18) -- (32,16);
\fill[black!15](33,18) -- (34,18) -- (32,19) -- (33,18);
\fill[black!15](32,16) -- (33,18) -- (32,19) -- (32,16);
\fill[black!15](32,16) -- (32,19) -- (32,20) -- (31,20) -- (32,16);
% Draw the boundary of piece8
\draw[black](32,16) -- (34,16) -- (33,17) -- (33,18) -- (34,18) -- (35,21) -- (34,21) -- (34,20) -- (33,19) -- (32,19) -- (32,20) -- (31,20) -- (32,16);
% Draw the convex parts of piece8
\fill[black!15](37,17) -- (36,17) -- (38,16) -- (37,17);
\fill[black!15](38,18) -- (37,17) -- (38,16) -- (38,18);
\fill[black!15](38,16) -- (39,19) -- (38,19) -- (38,18) -- (38,16);
\fill[black!15](36,14) -- (38,14) -- (37,15) -- (36,14);
\fill[black!15](36,14) -- (37,15) -- (37,16) -- (36,14);
\fill[black!15](37,16) -- (38,16) -- (36,17) -- (37,16);
\fill[black!15](36,14) -- (37,16) -- (36,17) -- (36,14);
\fill[black!15](36,14) -- (36,17) -- (36,18) -- (35,18) -- (36,14);
% Draw the boundary of piece8
\draw[black](36,14) -- (38,14) -- (37,15) -- (37,16) -- (38,16) -- (39,19) -- (38,19) -- (38,18) -- (37,17) -- (36,17) -- (36,18) -- (35,18) -- (36,14);
% Draw the convex parts of piece8
\fill[black!15](37,21) -- (36,21) -- (38,20) -- (37,21);
\fill[black!15](38,22) -- (37,21) -- (38,20) -- (38,22);
\fill[black!15](38,20) -- (39,23) -- (38,23) -- (38,22) -- (38,20);
\fill[black!15](36,18) -- (38,18) -- (37,19) -- (36,18);
\fill[black!15](36,18) -- (37,19) -- (37,20) -- (36,18);
\fill[black!15](37,20) -- (38,20) -- (36,21) -- (37,20);
\fill[black!15](36,18) -- (37,20) -- (36,21) -- (36,18);
\fill[black!15](36,18) -- (36,21) -- (36,22) -- (35,22) -- (36,18);
% Draw the boundary of piece8
\draw[black](36,18) -- (38,18) -- (37,19) -- (37,20) -- (38,20) -- (39,23) -- (38,23) -- (38,22) -- (37,21) -- (36,21) -- (36,22) -- (35,22) -- (36,18);
% Draw the convex parts of piece8
\fill[black!15](43,8) -- (42,8) -- (44,7) -- (43,8);
\fill[black!15](44,9) -- (43,8) -- (44,7) -- (44,9);
\fill[black!15](44,7) -- (45,10) -- (44,10) -- (44,9) -- (44,7);
\fill[black!15](42,5) -- (44,5) -- (43,6) -- (42,5);
\fill[black!15](42,5) -- (43,6) -- (43,7) -- (42,5);
\fill[black!15](43,7) -- (44,7) -- (42,8) -- (43,7);
\fill[black!15](42,5) -- (43,7) -- (42,8) -- (42,5);
\fill[black!15](42,5) -- (42,8) -- (42,9) -- (41,9) -- (42,5);
% Draw the boundary of piece8
\draw[black](42,5) -- (44,5) -- (43,6) -- (43,7) -- (44,7) -- (45,10) -- (44,10) -- (44,9) -- (43,8) -- (42,8) -- (42,9) -- (41,9) -- (42,5);
% Draw the convex parts of piece4
\fill[black!15](2,1) -- (1,1) -- (3,0) -- (2,1);
\fill[black!15](3,0) -- (3,2) -- (2,2) -- (2,1) -- (3,0);
\fill[black!15](0,0) -- (3,0) -- (1,1) -- (0,0);
\fill[black!15](0,0) -- (1,1) -- (1,2) -- (0,2) -- (0,0);
% Draw the boundary of piece4
\draw[black](0,0) -- (3,0) -- (3,2) -- (2,2) -- (2,1) -- (1,1) -- (1,2) -- (0,2) -- (0,0);
% Draw the convex parts of piece4
\fill[black!15](2,3) -- (1,3) -- (3,2) -- (2,3);
\fill[black!15](3,2) -- (3,4) -- (2,4) -- (2,3) -- (3,2);
\fill[black!15](0,2) -- (3,2) -- (1,3) -- (0,2);
\fill[black!15](0,2) -- (1,3) -- (1,4) -- (0,4) -- (0,2);
% Draw the boundary of piece4
\draw[black](0,2) -- (3,2) -- (3,4) -- (2,4) -- (2,3) -- (1,3) -- (1,4) -- (0,4) -- (0,2);
% Draw the convex parts of piece4
\fill[black!15](2,5) -- (1,5) -- (3,4) -- (2,5);
\fill[black!15](3,4) -- (3,6) -- (2,6) -- (2,5) -- (3,4);
\fill[black!15](0,4) -- (3,4) -- (1,5) -- (0,4);
\fill[black!15](0,4) -- (1,5) -- (1,6) -- (0,6) -- (0,4);
% Draw the boundary of piece4
\draw[black](0,4) -- (3,4) -- (3,6) -- (2,6) -- (2,5) -- (1,5) -- (1,6) -- (0,6) -- (0,4);
% Draw the convex parts of piece4
\fill[black!15](2,7) -- (1,7) -- (3,6) -- (2,7);
\fill[black!15](3,6) -- (3,8) -- (2,8) -- (2,7) -- (3,6);
\fill[black!15](0,6) -- (3,6) -- (1,7) -- (0,6);
\fill[black!15](0,6) -- (1,7) -- (1,8) -- (0,8) -- (0,6);
% Draw the boundary of piece4
\draw[black](0,6) -- (3,6) -- (3,8) -- (2,8) -- (2,7) -- (1,7) -- (1,8) -- (0,8) -- (0,6);
% Draw the convex parts of piece4
\fill[black!15](2,24) -- (1,24) -- (3,23) -- (2,24);
\fill[black!15](3,23) -- (3,25) -- (2,25) -- (2,24) -- (3,23);
\fill[black!15](0,23) -- (3,23) -- (1,24) -- (0,23);
\fill[black!15](0,23) -- (1,24) -- (1,25) -- (0,25) -- (0,23);
% Draw the boundary of piece4
\draw[black](0,23) -- (3,23) -- (3,25) -- (2,25) -- (2,24) -- (1,24) -- (1,25) -- (0,25) -- (0,23);
% Draw the convex parts of piece4
\fill[black!15](2,26) -- (1,26) -- (3,25) -- (2,26);
\fill[black!15](3,25) -- (3,27) -- (2,27) -- (2,26) -- (3,25);
\fill[black!15](0,25) -- (3,25) -- (1,26) -- (0,25);
\fill[black!15](0,25) -- (1,26) -- (1,27) -- (0,27) -- (0,25);
% Draw the boundary of piece4
\draw[black](0,25) -- (3,25) -- (3,27) -- (2,27) -- (2,26) -- (1,26) -- (1,27) -- (0,27) -- (0,25);
% Draw the convex parts of piece4
\fill[black!15](5,1) -- (4,1) -- (6,0) -- (5,1);
\fill[black!15](6,0) -- (6,2) -- (5,2) -- (5,1) -- (6,0);
\fill[black!15](3,0) -- (6,0) -- (4,1) -- (3,0);
\fill[black!15](3,0) -- (4,1) -- (4,2) -- (3,2) -- (3,0);
% Draw the boundary of piece4
\draw[black](3,0) -- (6,0) -- (6,2) -- (5,2) -- (5,1) -- (4,1) -- (4,2) -- (3,2) -- (3,0);
% Draw the convex parts of piece4
\fill[black!15](10,1) -- (9,1) -- (11,0) -- (10,1);
\fill[black!15](11,0) -- (11,2) -- (10,2) -- (10,1) -- (11,0);
\fill[black!15](8,0) -- (11,0) -- (9,1) -- (8,0);
\fill[black!15](8,0) -- (9,1) -- (9,2) -- (8,2) -- (8,0);
% Draw the boundary of piece4
\draw[black](8,0) -- (11,0) -- (11,2) -- (10,2) -- (10,1) -- (9,1) -- (9,2) -- (8,2) -- (8,0);
% Draw the convex parts of piece4
\fill[black!15](16,26) -- (15,26) -- (17,25) -- (16,26);
\fill[black!15](17,25) -- (17,27) -- (16,27) -- (16,26) -- (17,25);
\fill[black!15](14,25) -- (17,25) -- (15,26) -- (14,25);
\fill[black!15](14,25) -- (15,26) -- (15,27) -- (14,27) -- (14,25);
% Draw the boundary of piece4
\draw[black](14,25) -- (17,25) -- (17,27) -- (16,27) -- (16,26) -- (15,26) -- (15,27) -- (14,27) -- (14,25);
% Draw the convex parts of piece4
\fill[black!15](21,26) -- (20,26) -- (22,25) -- (21,26);
\fill[black!15](22,25) -- (22,27) -- (21,27) -- (21,26) -- (22,25);
\fill[black!15](19,25) -- (22,25) -- (20,26) -- (19,25);
\fill[black!15](19,25) -- (20,26) -- (20,27) -- (19,27) -- (19,25);
% Draw the boundary of piece4
\draw[black](19,25) -- (22,25) -- (22,27) -- (21,27) -- (21,26) -- (20,26) -- (20,27) -- (19,27) -- (19,25);
% Draw the convex parts of piece4
\fill[black!15](25,21) -- (24,21) -- (26,20) -- (25,21);
\fill[black!15](26,20) -- (26,22) -- (25,22) -- (25,21) -- (26,20);
\fill[black!15](23,20) -- (26,20) -- (24,21) -- (23,20);
\fill[black!15](23,20) -- (24,21) -- (24,22) -- (23,22) -- (23,20);
% Draw the boundary of piece4
\draw[black](23,20) -- (26,20) -- (26,22) -- (25,22) -- (25,21) -- (24,21) -- (24,22) -- (23,22) -- (23,20);
% Draw the convex parts of piece4
\fill[black!15](28,21) -- (27,21) -- (29,20) -- (28,21);
\fill[black!15](29,20) -- (29,22) -- (28,22) -- (28,21) -- (29,20);
\fill[black!15](26,20) -- (29,20) -- (27,21) -- (26,20);
\fill[black!15](26,20) -- (27,21) -- (27,22) -- (26,22) -- (26,20);
% Draw the boundary of piece4
\draw[black](26,20) -- (29,20) -- (29,22) -- (28,22) -- (28,21) -- (27,21) -- (27,22) -- (26,22) -- (26,20);
% Draw the convex parts of piece4
\fill[black!15](33,21) -- (32,21) -- (34,20) -- (33,21);
\fill[black!15](34,20) -- (34,22) -- (33,22) -- (33,21) -- (34,20);
\fill[black!15](31,20) -- (34,20) -- (32,21) -- (31,20);
\fill[black!15](31,20) -- (32,21) -- (32,22) -- (31,22) -- (31,20);
% Draw the boundary of piece4
\draw[black](31,20) -- (34,20) -- (34,22) -- (33,22) -- (33,21) -- (32,21) -- (32,22) -- (31,22) -- (31,20);
% Draw the convex parts of piece4
\fill[black!15](41,15) -- (40,15) -- (42,14) -- (41,15);
\fill[black!15](42,14) -- (42,16) -- (41,16) -- (41,15) -- (42,14);
\fill[black!15](39,14) -- (42,14) -- (40,15) -- (39,14);
\fill[black!15](39,14) -- (40,15) -- (40,16) -- (39,16) -- (39,14);
% Draw the boundary of piece4
\draw[black](39,14) -- (42,14) -- (42,16) -- (41,16) -- (41,15) -- (40,15) -- (40,16) -- (39,16) -- (39,14);
% Draw the convex parts of piece4
\fill[black!15](41,17) -- (40,17) -- (42,16) -- (41,17);
\fill[black!15](42,16) -- (42,18) -- (41,18) -- (41,17) -- (42,16);
\fill[black!15](39,16) -- (42,16) -- (40,17) -- (39,16);
\fill[black!15](39,16) -- (40,17) -- (40,18) -- (39,18) -- (39,16);
% Draw the boundary of piece4
\draw[black](39,16) -- (42,16) -- (42,18) -- (41,18) -- (41,17) -- (40,17) -- (40,18) -- (39,18) -- (39,16);
% Draw the convex parts of piece1
\fill[black!15](2,20) -- (4,20) -- (4,22) -- (2,22) -- (2,20);
% Draw the boundary of piece1
\draw[black](2,20) -- (4,20) -- (4,22) -- (2,22) -- (2,20);
% Draw the convex parts of piece1
\fill[black!15](4,7) -- (6,7) -- (6,9) -- (4,9) -- (4,7);
% Draw the boundary of piece1
\draw[black](4,7) -- (6,7) -- (6,9) -- (4,9) -- (4,7);
% Draw the convex parts of piece1
\fill[black!15](5,24) -- (7,24) -- (7,26) -- (5,26) -- (5,24);
% Draw the boundary of piece1
\draw[black](5,24) -- (7,24) -- (7,26) -- (5,26) -- (5,24);
% Draw the convex parts of piece1
\fill[black!15](6,0) -- (8,0) -- (8,2) -- (6,2) -- (6,0);
% Draw the boundary of piece1
\draw[black](6,0) -- (8,0) -- (8,2) -- (6,2) -- (6,0);
% Draw the convex parts of piece1
\fill[black!15](8,25) -- (10,25) -- (10,27) -- (8,27) -- (8,25);
% Draw the boundary of piece1
\draw[black](8,25) -- (10,25) -- (10,27) -- (8,27) -- (8,25);
% Draw the convex parts of piece1
\fill[black!15](13,15) -- (15,15) -- (15,17) -- (13,17) -- (13,15);
% Draw the boundary of piece1
\draw[black](13,15) -- (15,15) -- (15,17) -- (13,17) -- (13,15);
% Draw the convex parts of piece1
\fill[black!15](25,22) -- (27,22) -- (27,24) -- (25,24) -- (25,22);
% Draw the boundary of piece1
\draw[black](25,22) -- (27,22) -- (27,24) -- (25,24) -- (25,22);
% Draw the convex parts of piece1
\fill[black!15](28,11) -- (30,11) -- (30,13) -- (28,13) -- (28,11);
% Draw the boundary of piece1
\draw[black](28,11) -- (30,11) -- (30,13) -- (28,13) -- (28,11);
% Draw the convex parts of piece1
\fill[black!15](30,0) -- (32,0) -- (32,2) -- (30,2) -- (30,0);
% Draw the boundary of piece1
\draw[black](30,0) -- (32,0) -- (32,2) -- (30,2) -- (30,0);
% Draw the convex parts of piece1
\fill[black!15](33,5) -- (35,5) -- (35,7) -- (33,7) -- (33,5);
% Draw the boundary of piece1
\draw[black](33,5) -- (35,5) -- (35,7) -- (33,7) -- (33,5);
% Draw the convex parts of piece1
\fill[black!15](34,2) -- (36,2) -- (36,4) -- (34,4) -- (34,2);
% Draw the boundary of piece1
\draw[black](34,2) -- (36,2) -- (36,4) -- (34,4) -- (34,2);
% Draw the convex parts of piece1
\fill[black!15](38,2) -- (40,2) -- (40,4) -- (38,4) -- (38,2);
% Draw the boundary of piece1
\draw[black](38,2) -- (40,2) -- (40,4) -- (38,4) -- (38,2);
% Draw the convex parts of piece1
\fill[black!15](47,7) -- (49,7) -- (49,9) -- (47,9) -- (47,7);
% Draw the boundary of piece1
\draw[black](47,7) -- (49,7) -- (49,9) -- (47,9) -- (47,7);
% Draw the convex parts of piece1
\fill[black!15](47,12) -- (49,12) -- (49,14) -- (47,14) -- (47,12);
% Draw the boundary of piece1
\draw[black](47,12) -- (49,12) -- (49,14) -- (47,14) -- (47,12);
% Draw the convex parts of piece1
\fill[black!15](47,20) -- (49,20) -- (49,22) -- (47,22) -- (47,20);
% Draw the boundary of piece1
\draw[black](47,20) -- (49,20) -- (49,22) -- (47,22) -- (47,20);
% Draw the convex parts of piece6
\fill[black!15](12,0) -- (12,4) -- (11,3) -- (11,1) -- (12,0);
% Draw the boundary of piece6
\draw[black](12,0) -- (12,4) -- (11,3) -- (11,1) -- (12,0);
% Draw the convex parts of piece6
\fill[black!15](12,5) -- (12,9) -- (11,8) -- (11,6) -- (12,5);
% Draw the boundary of piece6
\draw[black](12,5) -- (12,9) -- (11,8) -- (11,6) -- (12,5);
% Draw the convex parts of piece6
\fill[black!15](12,10) -- (12,14) -- (11,13) -- (11,11) -- (12,10);
% Draw the boundary of piece6
\draw[black](12,10) -- (12,14) -- (11,13) -- (11,11) -- (12,10);
% Draw the convex parts of piece6
\fill[black!15](13,15) -- (13,19) -- (12,18) -- (12,16) -- (13,15);
% Draw the boundary of piece6
\draw[black](13,15) -- (13,19) -- (12,18) -- (12,16) -- (13,15);
% Draw the convex parts of piece6
\fill[black!15](13,20) -- (13,24) -- (12,23) -- (12,21) -- (13,20);
% Draw the boundary of piece6
\draw[black](13,20) -- (13,24) -- (12,23) -- (12,21) -- (13,20);
% Draw the convex parts of piece6
\fill[black!15](19,15) -- (19,19) -- (18,18) -- (18,16) -- (19,15);
% Draw the boundary of piece6
\draw[black](19,15) -- (19,19) -- (18,18) -- (18,16) -- (19,15);
% Draw the convex parts of piece6
\fill[black!15](25,16) -- (25,20) -- (24,19) -- (24,17) -- (25,16);
% Draw the boundary of piece6
\draw[black](25,16) -- (25,20) -- (24,19) -- (24,17) -- (25,16);
% Draw the convex parts of piece6
\fill[black!15](28,7) -- (28,11) -- (27,10) -- (27,8) -- (28,7);
% Draw the boundary of piece6
\draw[black](28,7) -- (28,11) -- (27,10) -- (27,8) -- (28,7);
% Draw the convex parts of piece6
\fill[black!15](35,7) -- (35,11) -- (34,10) -- (34,8) -- (35,7);
% Draw the boundary of piece6
\draw[black](35,7) -- (35,11) -- (34,10) -- (34,8) -- (35,7);
% Draw the convex parts of piece6
\fill[black!15](35,15) -- (35,19) -- (34,18) -- (34,16) -- (35,15);
% Draw the boundary of piece6
\draw[black](35,15) -- (35,19) -- (34,18) -- (34,16) -- (35,15);
% Draw the convex parts of piece6
\fill[black!15](36,8) -- (36,12) -- (35,11) -- (35,9) -- (36,8);
% Draw the boundary of piece6
\draw[black](36,8) -- (36,12) -- (35,11) -- (35,9) -- (36,8);
% Draw the convex parts of piece6
\fill[black!15](39,13) -- (39,17) -- (38,16) -- (38,14) -- (39,13);
% Draw the boundary of piece6
\draw[black](39,13) -- (39,17) -- (38,16) -- (38,14) -- (39,13);
% Draw the convex parts of piece6
\fill[black!15](48,14) -- (48,18) -- (47,17) -- (47,15) -- (48,14);
% Draw the boundary of piece6
\draw[black](48,14) -- (48,18) -- (47,17) -- (47,15) -- (48,14);
% Draw the convex parts of piece6
\fill[black!15](49,14) -- (49,18) -- (48,17) -- (48,15) -- (49,14);
% Draw the boundary of piece6
\draw[black](49,14) -- (49,18) -- (48,17) -- (48,15) -- (49,14);
% Draw the convex parts of piece6
\fill[black!15](49,23) -- (49,27) -- (48,26) -- (48,24) -- (49,23);
% Draw the boundary of piece6
\draw[black](49,23) -- (49,27) -- (48,26) -- (48,24) -- (49,23);
% Draw the bounding box of the solution
\draw[black](0,0) -- (49,0) -- (49,27)-- (0,27) -- (0,0);
\draw (24.5,-2) node {artif};
\draw (24.5,-4.5) node {L = 49};
% Draw the bounding box of the solution
\end{tikzpicture}
% SOLUTION for problem instance = shirts4_8
% Scaling factor for coordinates = 1
\begin{tikzpicture}[scale = 0.08]
% Draw the convex parts of piece3
\fill[black!15](3,1) -- (3,0) -- (5,1) -- (3,1);
\fill[black!15](8,1) -- (11,0) -- (12,5) -- (8,1);
\fill[black!15](7,8) -- (0,8) -- (0,1) -- (3,1) -- (7,8);
\fill[black!15](3,1) -- (5,1) -- (8,1) -- (12,5) -- (11,9) -- (7,8) -- (3,1);
% Draw the boundary of piece3
\draw[black](0,1) -- (3,1) -- (3,0) -- (5,1) -- (8,1) -- (11,0) -- (12,5) -- (11,9) -- (7,8) -- (0,8) -- (0,1);
% Draw the convex parts of piece3
\fill[black!15](3,10) -- (3,9) -- (5,10) -- (3,10);
\fill[black!15](8,10) -- (11,9) -- (12,14) -- (8,10);
\fill[black!15](7,17) -- (0,17) -- (0,10) -- (3,10) -- (7,17);
\fill[black!15](3,10) -- (5,10) -- (8,10) -- (12,14) -- (11,18) -- (7,17) -- (3,10);
% Draw the boundary of piece3
\draw[black](0,10) -- (3,10) -- (3,9) -- (5,10) -- (8,10) -- (11,9) -- (12,14) -- (11,18) -- (7,17) -- (0,17) -- (0,10);
% Draw the convex parts of piece3
\fill[black!15](3,19) -- (3,18) -- (5,19) -- (3,19);
\fill[black!15](8,19) -- (11,18) -- (12,23) -- (8,19);
\fill[black!15](7,26) -- (0,26) -- (0,19) -- (3,19) -- (7,26);
\fill[black!15](3,19) -- (5,19) -- (8,19) -- (12,23) -- (11,27) -- (7,26) -- (3,19);
% Draw the boundary of piece3
\draw[black](0,19) -- (3,19) -- (3,18) -- (5,19) -- (8,19) -- (11,18) -- (12,23) -- (11,27) -- (7,26) -- (0,26) -- (0,19);
% Draw the convex parts of piece3
\fill[black!15](3,28) -- (3,27) -- (5,28) -- (3,28);
\fill[black!15](8,28) -- (11,27) -- (12,32) -- (8,28);
\fill[black!15](7,35) -- (0,35) -- (0,28) -- (3,28) -- (7,35);
\fill[black!15](3,28) -- (5,28) -- (8,28) -- (12,32) -- (11,36) -- (7,35) -- (3,28);
% Draw the boundary of piece3
\draw[black](0,28) -- (3,28) -- (3,27) -- (5,28) -- (8,28) -- (11,27) -- (12,32) -- (11,36) -- (7,35) -- (0,35) -- (0,28);
% Draw the convex parts of piece2
\fill[black!15](18,11) -- (15,12) -- (12,11) -- (18,11);
\fill[black!15](12,6) -- (23,6) -- (24,8) -- (23,9) -- (12,6);
\fill[black!15](18,11) -- (12,11) -- (11,9) -- (12,6) -- (23,9) -- (18,11);
\fill[black!15](23,9) -- (23,11) -- (22,12) -- (18,11) -- (23,9);
% Draw the boundary of piece2
\draw[black](12,6) -- (23,6) -- (24,8) -- (23,9) -- (23,11) -- (22,12) -- (18,11) -- (15,12) -- (12,11) -- (11,9) -- (12,6);
% Draw the convex parts of piece2
\fill[black!15](18,20) -- (15,21) -- (12,20) -- (18,20);
\fill[black!15](12,15) -- (23,15) -- (24,17) -- (23,18) -- (12,15);
\fill[black!15](18,20) -- (12,20) -- (11,18) -- (12,15) -- (23,18) -- (18,20);
\fill[black!15](23,18) -- (23,20) -- (22,21) -- (18,20) -- (23,18);
% Draw the boundary of piece2
\draw[black](12,15) -- (23,15) -- (24,17) -- (23,18) -- (23,20) -- (22,21) -- (18,20) -- (15,21) -- (12,20) -- (11,18) -- (12,15);
% Draw the convex parts of piece2
\fill[black!15](18,29) -- (15,30) -- (12,29) -- (18,29);
\fill[black!15](12,24) -- (23,24) -- (24,26) -- (23,27) -- (12,24);
\fill[black!15](18,29) -- (12,29) -- (11,27) -- (12,24) -- (23,27) -- (18,29);
\fill[black!15](23,27) -- (23,29) -- (22,30) -- (18,29) -- (23,27);
% Draw the boundary of piece2
\draw[black](12,24) -- (23,24) -- (24,26) -- (23,27) -- (23,29) -- (22,30) -- (18,29) -- (15,30) -- (12,29) -- (11,27) -- (12,24);
% Draw the convex parts of piece2
\fill[black!15](18,38) -- (15,39) -- (12,38) -- (18,38);
\fill[black!15](12,33) -- (23,33) -- (24,35) -- (23,36) -- (12,33);
\fill[black!15](18,38) -- (12,38) -- (11,36) -- (12,33) -- (23,36) -- (18,38);
\fill[black!15](23,36) -- (23,38) -- (22,39) -- (18,38) -- (23,36);
% Draw the boundary of piece2
\draw[black](12,33) -- (23,33) -- (24,35) -- (23,36) -- (23,38) -- (22,39) -- (18,38) -- (15,39) -- (12,38) -- (11,36) -- (12,33);
% Draw the convex parts of piece1
\fill[black!15](23,11) -- (24,9) -- (31,10) -- (23,11);
\fill[black!15](23,13) -- (22,12) -- (23,11) -- (31,10) -- (23,13);
\fill[black!15](31,10) -- (31,14) -- (24,16) -- (23,14) -- (23,13) -- (31,10);
% Draw the boundary of piece1
\draw[black](24,9) -- (31,10) -- (31,14) -- (24,16) -- (23,14) -- (23,13) -- (22,12) -- (23,11) -- (24,9);
% Draw the convex parts of piece1
\fill[black!15](24,4) -- (25,2) -- (32,3) -- (24,4);
\fill[black!15](24,6) -- (23,5) -- (24,4) -- (32,3) -- (24,6);
\fill[black!15](32,3) -- (32,7) -- (25,9) -- (24,7) -- (24,6) -- (32,3);
% Draw the boundary of piece1
\draw[black](25,2) -- (32,3) -- (32,7) -- (25,9) -- (24,7) -- (24,6) -- (23,5) -- (24,4) -- (25,2);
% Draw the convex parts of piece1
\fill[black!15](24,18) -- (25,16) -- (32,17) -- (24,18);
\fill[black!15](24,20) -- (23,19) -- (24,18) -- (32,17) -- (24,20);
\fill[black!15](32,17) -- (32,21) -- (25,23) -- (24,21) -- (24,20) -- (32,17);
% Draw the boundary of piece1
\draw[black](25,16) -- (32,17) -- (32,21) -- (25,23) -- (24,21) -- (24,20) -- (23,19) -- (24,18) -- (25,16);
% Draw the convex parts of piece1
\fill[black!15](24,32) -- (25,30) -- (32,31) -- (24,32);
\fill[black!15](24,34) -- (23,33) -- (24,32) -- (32,31) -- (24,34);
\fill[black!15](32,31) -- (32,35) -- (25,37) -- (24,35) -- (24,34) -- (32,31);
% Draw the boundary of piece1
\draw[black](25,30) -- (32,31) -- (32,35) -- (25,37) -- (24,35) -- (24,34) -- (23,33) -- (24,32) -- (25,30);
% Draw the convex parts of piece8
\fill[black!15](5,37) -- (4,37) -- (9,36) -- (5,37);
\fill[black!15](9,36) -- (9,38) -- (7,38) -- (5,37) -- (9,36);
\fill[black!15](1,36) -- (9,36) -- (4,37) -- (1,36);
\fill[black!15](1,36) -- (4,37) -- (3,38) -- (0,38) -- (1,36);
% Draw the boundary of piece8
\draw[black](1,36) -- (9,36) -- (9,38) -- (7,38) -- (5,37) -- (4,37) -- (3,38) -- (0,38) -- (1,36);
% Draw the convex parts of piece8
\fill[black!15](17,1) -- (16,1) -- (21,0) -- (17,1);
\fill[black!15](21,0) -- (21,2) -- (19,2) -- (17,1) -- (21,0);
\fill[black!15](13,0) -- (21,0) -- (16,1) -- (13,0);
\fill[black!15](13,0) -- (16,1) -- (15,2) -- (12,2) -- (13,0);
% Draw the boundary of piece8
\draw[black](13,0) -- (21,0) -- (21,2) -- (19,2) -- (17,1) -- (16,1) -- (15,2) -- (12,2) -- (13,0);
% Draw the convex parts of piece8
\fill[black!15](17,5) -- (16,5) -- (21,4) -- (17,5);
\fill[black!15](21,4) -- (21,6) -- (19,6) -- (17,5) -- (21,4);
\fill[black!15](13,4) -- (21,4) -- (16,5) -- (13,4);
\fill[black!15](13,4) -- (16,5) -- (15,6) -- (12,6) -- (13,4);
% Draw the boundary of piece8
\draw[black](13,4) -- (21,4) -- (21,6) -- (19,6) -- (17,5) -- (16,5) -- (15,6) -- (12,6) -- (13,4);
% Draw the convex parts of piece8
\fill[black!15](17,14) -- (16,14) -- (21,13) -- (17,14);
\fill[black!15](21,13) -- (21,15) -- (19,15) -- (17,14) -- (21,13);
\fill[black!15](13,13) -- (21,13) -- (16,14) -- (13,13);
\fill[black!15](13,13) -- (16,14) -- (15,15) -- (12,15) -- (13,13);
% Draw the boundary of piece8
\draw[black](13,13) -- (21,13) -- (21,15) -- (19,15) -- (17,14) -- (16,14) -- (15,15) -- (12,15) -- (13,13);
% Draw the convex parts of piece8
\fill[black!15](26,1) -- (25,1) -- (30,0) -- (26,1);
\fill[black!15](30,0) -- (30,2) -- (28,2) -- (26,1) -- (30,0);
\fill[black!15](22,0) -- (30,0) -- (25,1) -- (22,0);
\fill[black!15](22,0) -- (25,1) -- (24,2) -- (21,2) -- (22,0);
% Draw the boundary of piece8
\draw[black](22,0) -- (30,0) -- (30,2) -- (28,2) -- (26,1) -- (25,1) -- (24,2) -- (21,2) -- (22,0);
% Draw the convex parts of piece8
\fill[black!15](28,27) -- (27,27) -- (32,26) -- (28,27);
\fill[black!15](32,26) -- (32,28) -- (30,28) -- (28,27) -- (32,26);
\fill[black!15](24,26) -- (32,26) -- (27,27) -- (24,26);
\fill[black!15](24,26) -- (27,27) -- (26,28) -- (23,28) -- (24,26);
% Draw the boundary of piece8
\draw[black](24,26) -- (32,26) -- (32,28) -- (30,28) -- (28,27) -- (27,27) -- (26,28) -- (23,28) -- (24,26);
% Draw the convex parts of piece8
\fill[black!15](28,29) -- (27,29) -- (32,28) -- (28,29);
\fill[black!15](32,28) -- (32,30) -- (30,30) -- (28,29) -- (32,28);
\fill[black!15](24,28) -- (32,28) -- (27,29) -- (24,28);
\fill[black!15](24,28) -- (27,29) -- (26,30) -- (23,30) -- (24,28);
% Draw the boundary of piece8
\draw[black](24,28) -- (32,28) -- (32,30) -- (30,30) -- (28,29) -- (27,29) -- (26,30) -- (23,30) -- (24,28);
% Draw the convex parts of piece8
\fill[black!15](28,38) -- (27,38) -- (32,37) -- (28,38);
\fill[black!15](32,37) -- (32,39) -- (30,39) -- (28,38) -- (32,37);
\fill[black!15](24,37) -- (32,37) -- (27,38) -- (24,37);
\fill[black!15](24,37) -- (27,38) -- (26,39) -- (23,39) -- (24,37);
% Draw the boundary of piece8
\draw[black](24,37) -- (32,37) -- (32,39) -- (30,39) -- (28,38) -- (27,38) -- (26,39) -- (23,39) -- (24,37);
% Draw the convex parts of piece4
\fill[black!15](12,21) -- (16,21) -- (16,23) -- (15,24) -- (12,24) -- (12,21);
% Draw the boundary of piece4
\draw[black](12,21) -- (16,21) -- (16,23) -- (15,24) -- (12,24) -- (12,21);
% Draw the convex parts of piece4
\fill[black!15](12,30) -- (16,30) -- (16,32) -- (15,33) -- (12,33) -- (12,30);
% Draw the boundary of piece4
\draw[black](12,30) -- (16,30) -- (16,32) -- (15,33) -- (12,33) -- (12,30);
% Draw the convex parts of piece4
\fill[black!15](16,21) -- (20,21) -- (20,23) -- (19,24) -- (16,24) -- (16,21);
% Draw the boundary of piece4
\draw[black](16,21) -- (20,21) -- (20,23) -- (19,24) -- (16,24) -- (16,21);
% Draw the convex parts of piece4
\fill[black!15](16,30) -- (20,30) -- (20,32) -- (19,33) -- (16,33) -- (16,30);
% Draw the boundary of piece4
\draw[black](16,30) -- (20,30) -- (20,32) -- (19,33) -- (16,33) -- (16,30);
% Draw the convex parts of piece4
\fill[black!15](20,21) -- (24,21) -- (24,23) -- (23,24) -- (20,24) -- (20,21);
% Draw the boundary of piece4
\draw[black](20,21) -- (24,21) -- (24,23) -- (23,24) -- (20,24) -- (20,21);
% Draw the convex parts of piece4
\fill[black!15](20,30) -- (24,30) -- (24,32) -- (23,33) -- (20,33) -- (20,30);
% Draw the boundary of piece4
\draw[black](20,30) -- (24,30) -- (24,32) -- (23,33) -- (20,33) -- (20,30);
% Draw the convex parts of piece4
\fill[black!15](24,23) -- (28,23) -- (28,25) -- (27,26) -- (24,26) -- (24,23);
% Draw the boundary of piece4
\draw[black](24,23) -- (28,23) -- (28,25) -- (27,26) -- (24,26) -- (24,23);
% Draw the convex parts of piece4
\fill[black!15](28,23) -- (32,23) -- (32,25) -- (31,26) -- (28,26) -- (28,23);
% Draw the boundary of piece4
\draw[black](28,23) -- (32,23) -- (32,25) -- (31,26) -- (28,26) -- (28,23);
% Draw the convex parts of piece5
\fill[black!15](0,38) -- (8,38) -- (7,39) -- (1,39) -- (0,38);
% Draw the boundary of piece5
\draw[black](0,38) -- (8,38) -- (7,39) -- (1,39) -- (0,38);
% Draw the convex parts of piece5
\fill[black!15](0,39) -- (8,39) -- (7,40) -- (1,40) -- (0,39);
% Draw the boundary of piece5
\draw[black](0,39) -- (8,39) -- (7,40) -- (1,40) -- (0,39);
% Draw the convex parts of piece5
\fill[black!15](8,39) -- (16,39) -- (15,40) -- (9,40) -- (8,39);
% Draw the boundary of piece5
\draw[black](8,39) -- (16,39) -- (15,40) -- (9,40) -- (8,39);
% Draw the convex parts of piece5
\fill[black!15](12,12) -- (20,12) -- (19,13) -- (13,13) -- (12,12);
% Draw the boundary of piece5
\draw[black](12,12) -- (20,12) -- (19,13) -- (13,13) -- (12,12);
% Draw the convex parts of piece5
\fill[black!15](16,3) -- (24,3) -- (23,4) -- (17,4) -- (16,3);
% Draw the boundary of piece5
\draw[black](16,3) -- (24,3) -- (23,4) -- (17,4) -- (16,3);
% Draw the convex parts of piece5
\fill[black!15](16,39) -- (24,39) -- (23,40) -- (17,40) -- (16,39);
% Draw the boundary of piece5
\draw[black](16,39) -- (24,39) -- (23,40) -- (17,40) -- (16,39);
% Draw the convex parts of piece5
\fill[black!15](17,2) -- (25,2) -- (24,3) -- (18,3) -- (17,2);
% Draw the boundary of piece5
\draw[black](17,2) -- (25,2) -- (24,3) -- (18,3) -- (17,2);
% Draw the convex parts of piece5
\fill[black!15](24,39) -- (32,39) -- (31,40) -- (25,40) -- (24,39);
% Draw the boundary of piece5
\draw[black](24,39) -- (32,39) -- (31,40) -- (25,40) -- (24,39);
% Draw the convex parts of piece6
\fill[black!15](0,8) -- (4,8) -- (4,9) -- (0,9) -- (0,8);
% Draw the boundary of piece6
\draw[black](0,8) -- (4,8) -- (4,9) -- (0,9) -- (0,8);
% Draw the convex parts of piece6
\fill[black!15](0,26) -- (4,26) -- (4,27) -- (0,27) -- (0,26);
% Draw the boundary of piece6
\draw[black](0,26) -- (4,26) -- (4,27) -- (0,27) -- (0,26);
% Draw the convex parts of piece6
\fill[black!15](0,35) -- (4,35) -- (4,36) -- (0,36) -- (0,35);
% Draw the boundary of piece6
\draw[black](0,35) -- (4,35) -- (4,36) -- (0,36) -- (0,35);
% Draw the convex parts of piece6
\fill[black!15](3,17) -- (7,17) -- (7,18) -- (3,18) -- (3,17);
% Draw the boundary of piece6
\draw[black](3,17) -- (7,17) -- (7,18) -- (3,18) -- (3,17);
% Draw the convex parts of piece6
\fill[black!15](8,38) -- (12,38) -- (12,39) -- (8,39) -- (8,38);
% Draw the boundary of piece6
\draw[black](8,38) -- (12,38) -- (12,39) -- (8,39) -- (8,38);
% Draw the convex parts of piece6
\fill[black!15](12,2) -- (16,2) -- (16,3) -- (12,3) -- (12,2);
% Draw the boundary of piece6
\draw[black](12,2) -- (16,2) -- (16,3) -- (12,3) -- (12,2);
% Draw the convex parts of piece6
\fill[black!15](12,3) -- (16,3) -- (16,4) -- (12,4) -- (12,3);
% Draw the boundary of piece6
\draw[black](12,3) -- (16,3) -- (16,4) -- (12,4) -- (12,3);
% Draw the convex parts of piece6
\fill[black!15](28,15) -- (32,15) -- (32,16) -- (28,16) -- (28,15);
% Draw the boundary of piece6
\draw[black](28,15) -- (32,15) -- (32,16) -- (28,16) -- (28,15);
% Draw the convex parts of piece7
\fill[black!15](0,0) -- (3,0) -- (3,1) -- (0,1) -- (0,0);
% Draw the boundary of piece7
\draw[black](0,0) -- (3,0) -- (3,1) -- (0,1) -- (0,0);
% Draw the convex parts of piece7
\fill[black!15](0,9) -- (3,9) -- (3,10) -- (0,10) -- (0,9);
% Draw the boundary of piece7
\draw[black](0,9) -- (3,9) -- (3,10) -- (0,10) -- (0,9);
% Draw the convex parts of piece7
\fill[black!15](0,18) -- (3,18) -- (3,19) -- (0,19) -- (0,18);
% Draw the boundary of piece7
\draw[black](0,18) -- (3,18) -- (3,19) -- (0,19) -- (0,18);
% Draw the convex parts of piece7
\fill[black!15](0,27) -- (3,27) -- (3,28) -- (0,28) -- (0,27);
% Draw the boundary of piece7
\draw[black](0,27) -- (3,27) -- (3,28) -- (0,28) -- (0,27);
% Draw the convex parts of piece7
\fill[black!15](5,0) -- (8,0) -- (8,1) -- (5,1) -- (5,0);
% Draw the boundary of piece7
\draw[black](5,0) -- (8,0) -- (8,1) -- (5,1) -- (5,0);
% Draw the convex parts of piece7
\fill[black!15](5,9) -- (8,9) -- (8,10) -- (5,10) -- (5,9);
% Draw the boundary of piece7
\draw[black](5,9) -- (8,9) -- (8,10) -- (5,10) -- (5,9);
% Draw the convex parts of piece7
\fill[black!15](5,18) -- (8,18) -- (8,19) -- (5,19) -- (5,18);
% Draw the boundary of piece7
\draw[black](5,18) -- (8,18) -- (8,19) -- (5,19) -- (5,18);
% Draw the convex parts of piece7
\fill[black!15](5,27) -- (8,27) -- (8,28) -- (5,28) -- (5,27);
% Draw the boundary of piece7
\draw[black](5,27) -- (8,27) -- (8,28) -- (5,28) -- (5,27);
% Draw the bounding box of the solution
\draw[black](0,0) -- (32,0) -- (32,40)-- (0,40) -- (0,0);
\draw (16,-2) node {shirts4\_8};
\draw (16,-5) node {L = 32};
% Draw the bounding box of the solution
\end{tikzpicture}
% SOLUTION for problem instance = shirts5_10
% Scaling factor for coordinates = 1
\begin{tikzpicture}[scale = 0.08]
% Draw the convex parts of piece3
\fill[black!15](3,1) -- (3,0) -- (5,1) -- (3,1);
\fill[black!15](8,1) -- (11,0) -- (12,5) -- (8,1);
\fill[black!15](7,8) -- (0,8) -- (0,1) -- (3,1) -- (7,8);
\fill[black!15](3,1) -- (5,1) -- (8,1) -- (12,5) -- (11,9) -- (7,8) -- (3,1);
% Draw the boundary of piece3
\draw[black](0,1) -- (3,1) -- (3,0) -- (5,1) -- (8,1) -- (11,0) -- (12,5) -- (11,9) -- (7,8) -- (0,8) -- (0,1);
% Draw the convex parts of piece3
\fill[black!15](3,10) -- (3,9) -- (5,10) -- (3,10);
\fill[black!15](8,10) -- (11,9) -- (12,14) -- (8,10);
\fill[black!15](7,17) -- (0,17) -- (0,10) -- (3,10) -- (7,17);
\fill[black!15](3,10) -- (5,10) -- (8,10) -- (12,14) -- (11,18) -- (7,17) -- (3,10);
% Draw the boundary of piece3
\draw[black](0,10) -- (3,10) -- (3,9) -- (5,10) -- (8,10) -- (11,9) -- (12,14) -- (11,18) -- (7,17) -- (0,17) -- (0,10);
% Draw the convex parts of piece3
\fill[black!15](3,19) -- (3,18) -- (5,19) -- (3,19);
\fill[black!15](8,19) -- (11,18) -- (12,23) -- (8,19);
\fill[black!15](7,26) -- (0,26) -- (0,19) -- (3,19) -- (7,26);
\fill[black!15](3,19) -- (5,19) -- (8,19) -- (12,23) -- (11,27) -- (7,26) -- (3,19);
% Draw the boundary of piece3
\draw[black](0,19) -- (3,19) -- (3,18) -- (5,19) -- (8,19) -- (11,18) -- (12,23) -- (11,27) -- (7,26) -- (0,26) -- (0,19);
% Draw the convex parts of piece3
\fill[black!15](3,28) -- (3,27) -- (5,28) -- (3,28);
\fill[black!15](8,28) -- (11,27) -- (12,32) -- (8,28);
\fill[black!15](7,35) -- (0,35) -- (0,28) -- (3,28) -- (7,35);
\fill[black!15](3,28) -- (5,28) -- (8,28) -- (12,32) -- (11,36) -- (7,35) -- (3,28);
% Draw the boundary of piece3
\draw[black](0,28) -- (3,28) -- (3,27) -- (5,28) -- (8,28) -- (11,27) -- (12,32) -- (11,36) -- (7,35) -- (0,35) -- (0,28);
% Draw the convex parts of piece3
\fill[black!15](15,1) -- (15,0) -- (17,1) -- (15,1);
\fill[black!15](20,1) -- (23,0) -- (24,5) -- (20,1);
\fill[black!15](19,8) -- (12,8) -- (12,1) -- (15,1) -- (19,8);
\fill[black!15](15,1) -- (17,1) -- (20,1) -- (24,5) -- (23,9) -- (19,8) -- (15,1);
% Draw the boundary of piece3
\draw[black](12,1) -- (15,1) -- (15,0) -- (17,1) -- (20,1) -- (23,0) -- (24,5) -- (23,9) -- (19,8) -- (12,8) -- (12,1);
% Draw the convex parts of piece2
\fill[black!15](18,20) -- (15,21) -- (12,20) -- (18,20);
\fill[black!15](12,15) -- (23,15) -- (24,17) -- (23,18) -- (12,15);
\fill[black!15](18,20) -- (12,20) -- (11,18) -- (12,15) -- (23,18) -- (18,20);
\fill[black!15](23,18) -- (23,20) -- (22,21) -- (18,20) -- (23,18);
% Draw the boundary of piece2
\draw[black](12,15) -- (23,15) -- (24,17) -- (23,18) -- (23,20) -- (22,21) -- (18,20) -- (15,21) -- (12,20) -- (11,18) -- (12,15);
% Draw the convex parts of piece2
\fill[black!15](18,29) -- (15,30) -- (12,29) -- (18,29);
\fill[black!15](12,24) -- (23,24) -- (24,26) -- (23,27) -- (12,24);
\fill[black!15](18,29) -- (12,29) -- (11,27) -- (12,24) -- (23,27) -- (18,29);
\fill[black!15](23,27) -- (23,29) -- (22,30) -- (18,29) -- (23,27);
% Draw the boundary of piece2
\draw[black](12,24) -- (23,24) -- (24,26) -- (23,27) -- (23,29) -- (22,30) -- (18,29) -- (15,30) -- (12,29) -- (11,27) -- (12,24);
% Draw the convex parts of piece2
\fill[black!15](18,38) -- (15,39) -- (12,38) -- (18,38);
\fill[black!15](12,33) -- (23,33) -- (24,35) -- (23,36) -- (12,33);
\fill[black!15](18,38) -- (12,38) -- (11,36) -- (12,33) -- (23,36) -- (18,38);
\fill[black!15](23,36) -- (23,38) -- (22,39) -- (18,38) -- (23,36);
% Draw the boundary of piece2
\draw[black](12,33) -- (23,33) -- (24,35) -- (23,36) -- (23,38) -- (22,39) -- (18,38) -- (15,39) -- (12,38) -- (11,36) -- (12,33);
% Draw the convex parts of piece2
\fill[black!15](19,14) -- (16,15) -- (13,14) -- (19,14);
\fill[black!15](13,9) -- (24,9) -- (25,11) -- (24,12) -- (13,9);
\fill[black!15](19,14) -- (13,14) -- (12,12) -- (13,9) -- (24,12) -- (19,14);
\fill[black!15](24,12) -- (24,14) -- (23,15) -- (19,14) -- (24,12);
% Draw the boundary of piece2
\draw[black](13,9) -- (24,9) -- (25,11) -- (24,12) -- (24,14) -- (23,15) -- (19,14) -- (16,15) -- (13,14) -- (12,12) -- (13,9);
% Draw the convex parts of piece2
\fill[black!15](30,22) -- (27,23) -- (24,22) -- (30,22);
\fill[black!15](24,17) -- (35,17) -- (36,19) -- (35,20) -- (24,17);
\fill[black!15](30,22) -- (24,22) -- (23,20) -- (24,17) -- (35,20) -- (30,22);
\fill[black!15](35,20) -- (35,22) -- (34,23) -- (30,22) -- (35,20);
% Draw the boundary of piece2
\draw[black](24,17) -- (35,17) -- (36,19) -- (35,20) -- (35,22) -- (34,23) -- (30,22) -- (27,23) -- (24,22) -- (23,20) -- (24,17);
% Draw the convex parts of piece1
\fill[black!15](23,29) -- (24,27) -- (31,28) -- (23,29);
\fill[black!15](23,31) -- (22,30) -- (23,29) -- (31,28) -- (23,31);
\fill[black!15](31,28) -- (31,32) -- (24,34) -- (23,32) -- (23,31) -- (31,28);
% Draw the boundary of piece1
\draw[black](24,27) -- (31,28) -- (31,32) -- (24,34) -- (23,32) -- (23,31) -- (22,30) -- (23,29) -- (24,27);
% Draw the convex parts of piece1
\fill[black!15](25,2) -- (26,0) -- (33,1) -- (25,2);
\fill[black!15](25,4) -- (24,3) -- (25,2) -- (33,1) -- (25,4);
\fill[black!15](33,1) -- (33,5) -- (26,7) -- (25,5) -- (25,4) -- (33,1);
% Draw the boundary of piece1
\draw[black](26,0) -- (33,1) -- (33,5) -- (26,7) -- (25,5) -- (25,4) -- (24,3) -- (25,2) -- (26,0);
% Draw the convex parts of piece1
\fill[black!15](25,11) -- (26,9) -- (33,10) -- (25,11);
\fill[black!15](25,13) -- (24,12) -- (25,11) -- (33,10) -- (25,13);
\fill[black!15](33,10) -- (33,14) -- (26,16) -- (25,14) -- (25,13) -- (33,10);
% Draw the boundary of piece1
\draw[black](26,9) -- (33,10) -- (33,14) -- (26,16) -- (25,14) -- (25,13) -- (24,12) -- (25,11) -- (26,9);
% Draw the convex parts of piece1
\fill[black!15](31,25) -- (32,23) -- (39,24) -- (31,25);
\fill[black!15](31,27) -- (30,26) -- (31,25) -- (39,24) -- (31,27);
\fill[black!15](39,24) -- (39,28) -- (32,30) -- (31,28) -- (31,27) -- (39,24);
% Draw the boundary of piece1
\draw[black](32,23) -- (39,24) -- (39,28) -- (32,30) -- (31,28) -- (31,27) -- (30,26) -- (31,25) -- (32,23);
% Draw the convex parts of piece1
\fill[black!15](33,8) -- (34,6) -- (41,7) -- (33,8);
\fill[black!15](33,10) -- (32,9) -- (33,8) -- (41,7) -- (33,10);
\fill[black!15](41,7) -- (41,11) -- (34,13) -- (33,11) -- (33,10) -- (41,7);
% Draw the boundary of piece1
\draw[black](34,6) -- (41,7) -- (41,11) -- (34,13) -- (33,11) -- (33,10) -- (32,9) -- (33,8) -- (34,6);
% Draw the convex parts of piece8
\fill[black!15](5,37) -- (4,37) -- (9,36) -- (5,37);
\fill[black!15](9,36) -- (9,38) -- (7,38) -- (5,37) -- (9,36);
\fill[black!15](1,36) -- (9,36) -- (4,37) -- (1,36);
\fill[black!15](1,36) -- (4,37) -- (3,38) -- (0,38) -- (1,36);
% Draw the boundary of piece8
\draw[black](1,36) -- (9,36) -- (9,38) -- (7,38) -- (5,37) -- (4,37) -- (3,38) -- (0,38) -- (1,36);
% Draw the convex parts of piece8
\fill[black!15](5,39) -- (4,39) -- (9,38) -- (5,39);
\fill[black!15](9,38) -- (9,40) -- (7,40) -- (5,39) -- (9,38);
\fill[black!15](1,38) -- (9,38) -- (4,39) -- (1,38);
\fill[black!15](1,38) -- (4,39) -- (3,40) -- (0,40) -- (1,38);
% Draw the boundary of piece8
\draw[black](1,38) -- (9,38) -- (9,40) -- (7,40) -- (5,39) -- (4,39) -- (3,40) -- (0,40) -- (1,38);
% Draw the convex parts of piece8
\fill[black!15](28,8) -- (27,8) -- (32,7) -- (28,8);
\fill[black!15](32,7) -- (32,9) -- (30,9) -- (28,8) -- (32,7);
\fill[black!15](24,7) -- (32,7) -- (27,8) -- (24,7);
\fill[black!15](24,7) -- (27,8) -- (26,9) -- (23,9) -- (24,7);
% Draw the boundary of piece8
\draw[black](24,7) -- (32,7) -- (32,9) -- (30,9) -- (28,8) -- (27,8) -- (26,9) -- (23,9) -- (24,7);
% Draw the convex parts of piece8
\fill[black!15](28,36) -- (27,36) -- (32,35) -- (28,36);
\fill[black!15](32,35) -- (32,37) -- (30,37) -- (28,36) -- (32,35);
\fill[black!15](24,35) -- (32,35) -- (27,36) -- (24,35);
\fill[black!15](24,35) -- (27,36) -- (26,37) -- (23,37) -- (24,35);
% Draw the boundary of piece8
\draw[black](24,35) -- (32,35) -- (32,37) -- (30,37) -- (28,36) -- (27,36) -- (26,37) -- (23,37) -- (24,35);
% Draw the convex parts of piece8
\fill[black!15](28,38) -- (27,38) -- (32,37) -- (28,38);
\fill[black!15](32,37) -- (32,39) -- (30,39) -- (28,38) -- (32,37);
\fill[black!15](24,37) -- (32,37) -- (27,38) -- (24,37);
\fill[black!15](24,37) -- (27,38) -- (26,39) -- (23,39) -- (24,37);
% Draw the boundary of piece8
\draw[black](24,37) -- (32,37) -- (32,39) -- (30,39) -- (28,38) -- (27,38) -- (26,39) -- (23,39) -- (24,37);
% Draw the convex parts of piece8
\fill[black!15](37,15) -- (36,15) -- (41,14) -- (37,15);
\fill[black!15](41,14) -- (41,16) -- (39,16) -- (37,15) -- (41,14);
\fill[black!15](33,14) -- (41,14) -- (36,15) -- (33,14);
\fill[black!15](33,14) -- (36,15) -- (35,16) -- (32,16) -- (33,14);
% Draw the boundary of piece8
\draw[black](33,14) -- (41,14) -- (41,16) -- (39,16) -- (37,15) -- (36,15) -- (35,16) -- (32,16) -- (33,14);
% Draw the convex parts of piece8
\fill[black!15](37,33) -- (36,33) -- (41,32) -- (37,33);
\fill[black!15](41,32) -- (41,34) -- (39,34) -- (37,33) -- (41,32);
\fill[black!15](33,32) -- (41,32) -- (36,33) -- (33,32);
\fill[black!15](33,32) -- (36,33) -- (35,34) -- (32,34) -- (33,32);
% Draw the boundary of piece8
\draw[black](33,32) -- (41,32) -- (41,34) -- (39,34) -- (37,33) -- (36,33) -- (35,34) -- (32,34) -- (33,32);
% Draw the convex parts of piece8
\fill[black!15](37,35) -- (36,35) -- (41,34) -- (37,35);
\fill[black!15](41,34) -- (41,36) -- (39,36) -- (37,35) -- (41,34);
\fill[black!15](33,34) -- (41,34) -- (36,35) -- (33,34);
\fill[black!15](33,34) -- (36,35) -- (35,36) -- (32,36) -- (33,34);
% Draw the boundary of piece8
\draw[black](33,34) -- (41,34) -- (41,36) -- (39,36) -- (37,35) -- (36,35) -- (35,36) -- (32,36) -- (33,34);
% Draw the convex parts of piece8
\fill[black!15](37,37) -- (36,37) -- (41,36) -- (37,37);
\fill[black!15](41,36) -- (41,38) -- (39,38) -- (37,37) -- (41,36);
\fill[black!15](33,36) -- (41,36) -- (36,37) -- (33,36);
\fill[black!15](33,36) -- (36,37) -- (35,38) -- (32,38) -- (33,36);
% Draw the boundary of piece8
\draw[black](33,36) -- (41,36) -- (41,38) -- (39,38) -- (37,37) -- (36,37) -- (35,38) -- (32,38) -- (33,36);
% Draw the convex parts of piece8
\fill[black!15](37,39) -- (36,39) -- (41,38) -- (37,39);
\fill[black!15](41,38) -- (41,40) -- (39,40) -- (37,39) -- (41,38);
\fill[black!15](33,38) -- (41,38) -- (36,39) -- (33,38);
\fill[black!15](33,38) -- (36,39) -- (35,40) -- (32,40) -- (33,38);
% Draw the boundary of piece8
\draw[black](33,38) -- (41,38) -- (41,40) -- (39,40) -- (37,39) -- (36,39) -- (35,40) -- (32,40) -- (33,38);
% Draw the convex parts of piece4
\fill[black!15](12,21) -- (16,21) -- (16,23) -- (15,24) -- (12,24) -- (12,21);
% Draw the boundary of piece4
\draw[black](12,21) -- (16,21) -- (16,23) -- (15,24) -- (12,24) -- (12,21);
% Draw the convex parts of piece4
\fill[black!15](12,30) -- (16,30) -- (16,32) -- (15,33) -- (12,33) -- (12,30);
% Draw the boundary of piece4
\draw[black](12,30) -- (16,30) -- (16,32) -- (15,33) -- (12,33) -- (12,30);
% Draw the convex parts of piece4
\fill[black!15](16,21) -- (20,21) -- (20,23) -- (19,24) -- (16,24) -- (16,21);
% Draw the boundary of piece4
\draw[black](16,21) -- (20,21) -- (20,23) -- (19,24) -- (16,24) -- (16,21);
% Draw the convex parts of piece4
\fill[black!15](17,30) -- (21,30) -- (21,32) -- (20,33) -- (17,33) -- (17,30);
% Draw the boundary of piece4
\draw[black](17,30) -- (21,30) -- (21,32) -- (20,33) -- (17,33) -- (17,30);
% Draw the convex parts of piece4
\fill[black!15](24,24) -- (28,24) -- (28,26) -- (27,27) -- (24,27) -- (24,24);
% Draw the boundary of piece4
\draw[black](24,24) -- (28,24) -- (28,26) -- (27,27) -- (24,27) -- (24,24);
% Draw the convex parts of piece4
\fill[black!15](33,3) -- (37,3) -- (37,5) -- (36,6) -- (33,6) -- (33,3);
% Draw the boundary of piece4
\draw[black](33,3) -- (37,3) -- (37,5) -- (36,6) -- (33,6) -- (33,3);
% Draw the convex parts of piece4
\fill[black!15](35,20) -- (39,20) -- (39,22) -- (38,23) -- (35,23) -- (35,20);
% Draw the boundary of piece4
\draw[black](35,20) -- (39,20) -- (39,22) -- (38,23) -- (35,23) -- (35,20);
% Draw the convex parts of piece4
\fill[black!15](37,3) -- (41,3) -- (41,5) -- (40,6) -- (37,6) -- (37,3);
% Draw the boundary of piece4
\draw[black](37,3) -- (41,3) -- (41,5) -- (40,6) -- (37,6) -- (37,3);
% Draw the convex parts of piece4
\fill[black!15](37,17) -- (41,17) -- (41,19) -- (40,20) -- (37,20) -- (37,17);
% Draw the boundary of piece4
\draw[black](37,17) -- (41,17) -- (41,19) -- (40,20) -- (37,20) -- (37,17);
% Draw the convex parts of piece4
\fill[black!15](37,29) -- (41,29) -- (41,31) -- (40,32) -- (37,32) -- (37,29);
% Draw the boundary of piece4
\draw[black](37,29) -- (41,29) -- (41,31) -- (40,32) -- (37,32) -- (37,29);
% Draw the convex parts of piece5
\fill[black!15](9,39) -- (17,39) -- (16,40) -- (10,40) -- (9,39);
% Draw the boundary of piece5
\draw[black](9,39) -- (17,39) -- (16,40) -- (10,40) -- (9,39);
% Draw the convex parts of piece5
\fill[black!15](20,39) -- (28,39) -- (27,40) -- (21,40) -- (20,39);
% Draw the boundary of piece5
\draw[black](20,39) -- (28,39) -- (27,40) -- (21,40) -- (20,39);
% Draw the convex parts of piece5
\fill[black!15](24,16) -- (32,16) -- (31,17) -- (25,17) -- (24,16);
% Draw the boundary of piece5
\draw[black](24,16) -- (32,16) -- (31,17) -- (25,17) -- (24,16);
% Draw the convex parts of piece5
\fill[black!15](24,23) -- (32,23) -- (31,24) -- (25,24) -- (24,23);
% Draw the boundary of piece5
\draw[black](24,23) -- (32,23) -- (31,24) -- (25,24) -- (24,23);
% Draw the convex parts of piece5
\fill[black!15](24,34) -- (32,34) -- (31,35) -- (25,35) -- (24,34);
% Draw the boundary of piece5
\draw[black](24,34) -- (32,34) -- (31,35) -- (25,35) -- (24,34);
% Draw the convex parts of piece5
\fill[black!15](32,0) -- (40,0) -- (39,1) -- (33,1) -- (32,0);
% Draw the boundary of piece5
\draw[black](32,0) -- (40,0) -- (39,1) -- (33,1) -- (32,0);
% Draw the convex parts of piece5
\fill[black!15](32,16) -- (40,16) -- (39,17) -- (33,17) -- (32,16);
% Draw the boundary of piece5
\draw[black](32,16) -- (40,16) -- (39,17) -- (33,17) -- (32,16);
% Draw the convex parts of piece5
\fill[black!15](33,1) -- (41,1) -- (40,2) -- (34,2) -- (33,1);
% Draw the boundary of piece5
\draw[black](33,1) -- (41,1) -- (40,2) -- (34,2) -- (33,1);
% Draw the convex parts of piece5
\fill[black!15](33,2) -- (41,2) -- (40,3) -- (34,3) -- (33,2);
% Draw the boundary of piece5
\draw[black](33,2) -- (41,2) -- (40,3) -- (34,3) -- (33,2);
% Draw the convex parts of piece5
\fill[black!15](33,13) -- (41,13) -- (40,14) -- (34,14) -- (33,13);
% Draw the boundary of piece5
\draw[black](33,13) -- (41,13) -- (40,14) -- (34,14) -- (33,13);
% Draw the convex parts of piece6
\fill[black!15](0,17) -- (4,17) -- (4,18) -- (0,18) -- (0,17);
% Draw the boundary of piece6
\draw[black](0,17) -- (4,17) -- (4,18) -- (0,18) -- (0,17);
% Draw the convex parts of piece6
\fill[black!15](1,8) -- (5,8) -- (5,9) -- (1,9) -- (1,8);
% Draw the boundary of piece6
\draw[black](1,8) -- (5,8) -- (5,9) -- (1,9) -- (1,8);
% Draw the convex parts of piece6
\fill[black!15](2,35) -- (6,35) -- (6,36) -- (2,36) -- (2,35);
% Draw the boundary of piece6
\draw[black](2,35) -- (6,35) -- (6,36) -- (2,36) -- (2,35);
% Draw the convex parts of piece6
\fill[black!15](3,26) -- (7,26) -- (7,27) -- (3,27) -- (3,26);
% Draw the boundary of piece6
\draw[black](3,26) -- (7,26) -- (7,27) -- (3,27) -- (3,26);
% Draw the convex parts of piece6
\fill[black!15](12,8) -- (16,8) -- (16,9) -- (12,9) -- (12,8);
% Draw the boundary of piece6
\draw[black](12,8) -- (16,8) -- (16,9) -- (12,9) -- (12,8);
% Draw the convex parts of piece6
\fill[black!15](20,22) -- (24,22) -- (24,23) -- (20,23) -- (20,22);
% Draw the boundary of piece6
\draw[black](20,22) -- (24,22) -- (24,23) -- (20,23) -- (20,22);
% Draw the convex parts of piece6
\fill[black!15](20,23) -- (24,23) -- (24,24) -- (20,24) -- (20,23);
% Draw the boundary of piece6
\draw[black](20,23) -- (24,23) -- (24,24) -- (20,24) -- (20,23);
% Draw the convex parts of piece6
\fill[black!15](28,33) -- (32,33) -- (32,34) -- (28,34) -- (28,33);
% Draw the boundary of piece6
\draw[black](28,33) -- (32,33) -- (32,34) -- (28,34) -- (28,33);
% Draw the convex parts of piece6
\fill[black!15](28,39) -- (32,39) -- (32,40) -- (28,40) -- (28,39);
% Draw the boundary of piece6
\draw[black](28,39) -- (32,39) -- (32,40) -- (28,40) -- (28,39);
% Draw the convex parts of piece6
\fill[black!15](32,30) -- (36,30) -- (36,31) -- (32,31) -- (32,30);
% Draw the boundary of piece6
\draw[black](32,30) -- (36,30) -- (36,31) -- (32,31) -- (32,30);
% Draw the convex parts of piece7
\fill[black!15](0,9) -- (3,9) -- (3,10) -- (0,10) -- (0,9);
% Draw the boundary of piece7
\draw[black](0,9) -- (3,9) -- (3,10) -- (0,10) -- (0,9);
% Draw the convex parts of piece7
\fill[black!15](0,18) -- (3,18) -- (3,19) -- (0,19) -- (0,18);
% Draw the boundary of piece7
\draw[black](0,18) -- (3,18) -- (3,19) -- (0,19) -- (0,18);
% Draw the convex parts of piece7
\fill[black!15](0,27) -- (3,27) -- (3,28) -- (0,28) -- (0,27);
% Draw the boundary of piece7
\draw[black](0,27) -- (3,27) -- (3,28) -- (0,28) -- (0,27);
% Draw the convex parts of piece7
\fill[black!15](5,0) -- (8,0) -- (8,1) -- (5,1) -- (5,0);
% Draw the boundary of piece7
\draw[black](5,0) -- (8,0) -- (8,1) -- (5,1) -- (5,0);
% Draw the convex parts of piece7
\fill[black!15](5,9) -- (8,9) -- (8,10) -- (5,10) -- (5,9);
% Draw the boundary of piece7
\draw[black](5,9) -- (8,9) -- (8,10) -- (5,10) -- (5,9);
% Draw the convex parts of piece7
\fill[black!15](5,18) -- (8,18) -- (8,19) -- (5,19) -- (5,18);
% Draw the boundary of piece7
\draw[black](5,18) -- (8,18) -- (8,19) -- (5,19) -- (5,18);
% Draw the convex parts of piece7
\fill[black!15](5,27) -- (8,27) -- (8,28) -- (5,28) -- (5,27);
% Draw the boundary of piece7
\draw[black](5,27) -- (8,27) -- (8,28) -- (5,28) -- (5,27);
% Draw the convex parts of piece7
\fill[black!15](12,0) -- (15,0) -- (15,1) -- (12,1) -- (12,0);
% Draw the boundary of piece7
\draw[black](12,0) -- (15,0) -- (15,1) -- (12,1) -- (12,0);
% Draw the convex parts of piece7
\fill[black!15](17,0) -- (20,0) -- (20,1) -- (17,1) -- (17,0);
% Draw the boundary of piece7
\draw[black](17,0) -- (20,0) -- (20,1) -- (17,1) -- (17,0);
% Draw the convex parts of piece7
\fill[black!15](28,24) -- (31,24) -- (31,25) -- (28,25) -- (28,24);
% Draw the boundary of piece7
\draw[black](28,24) -- (31,24) -- (31,25) -- (28,25) -- (28,24);
% Draw the bounding box of the solution
\draw[black](0,0) -- (41,0) -- (41,40)-- (0,40) -- (0,0);
\draw (20.5,-2) node {shirts5\_10};
\draw (20.5,-5) node {L = 41};
% Draw the bounding box of the solution
\end{tikzpicture}
\caption{Ten open problem instances whose upper bounds are improved in this work. Eight of these improvements are obtained by the DB-CC \cite{Rodrigues2017-zl} model solved with Gurobi 11.0.3, whilst two others are obtained by our 0-1 DB-CC model solved with Google CPSat 9.11.4210. Most of the solutions above are likely to be the optimal despite it could not be proven.}
\label{fig:improved_UB_open_instances}
\end{figure}
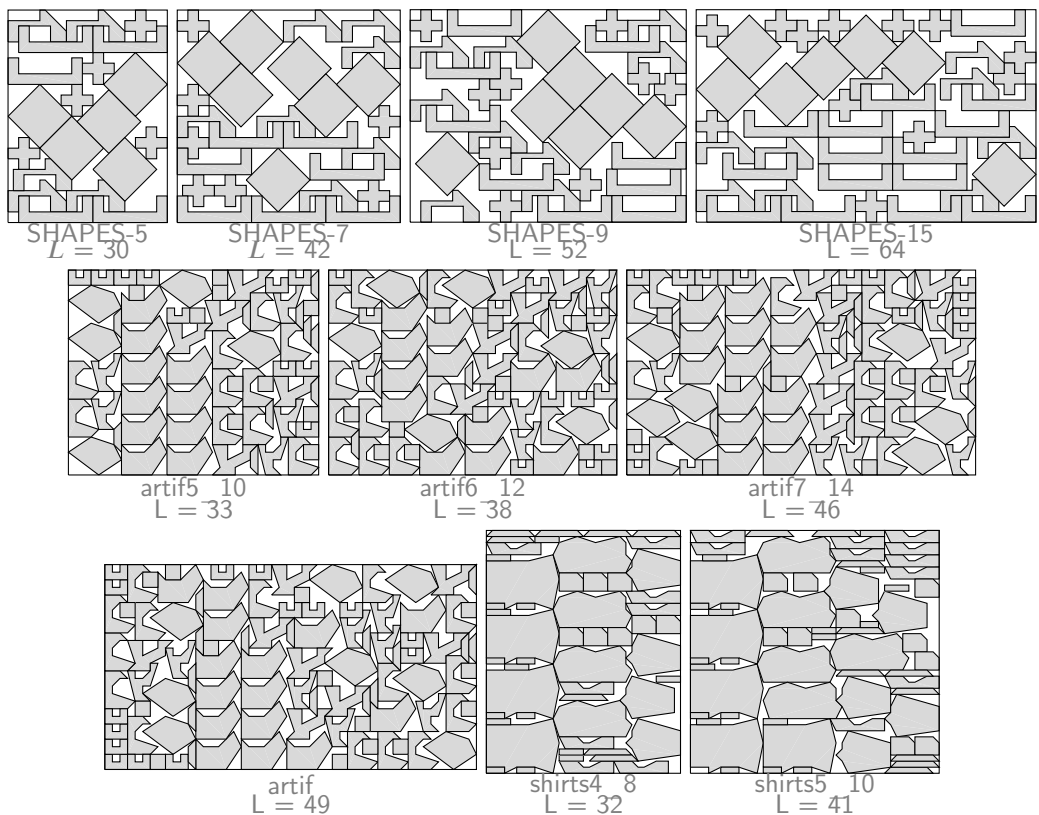

Our new Dotted-Board Parallel Backtracking (DB-PB) algorithm sets the new state of the art among the family of exact methods for the discrete nesting problem, as shown in Tables \ref{tab:first_results} and \ref{tab:second_results}, and Figures \ref{fig:performance_profile_total_time} and \ref{fig:performance_profile_total_time}. However, our DB-PB algorithm lacks of a feasibility pump and primal heuristics that help it to get better feasible solutions on those cases in which it cannot solve the problem instances up to optimality, as shown in Tables \ref{tab:first_results} and \ref{tab:second_results}.

\begin{table}[h]
\centering
\caption{Comparison of DB-CC model \cite{Rodrigues2017-zl} solved with Gurobi and our DB-PB algorithm with TL = 3600 seconds.}
\label{tab:first_results}
\tiny
\begin{tabular}{l|cccccc|cccccc}
& \multicolumn{6}{c}{\textbf{Dotted-Board Clique Covering (DB-CC) model} \cite{Rodrigues2017-zl}} & \multicolumn{6}{c}{\textbf{Dotted-Board Parallel-Backtracking (DB-PB) method (this work)}} \\
\textbf{Instance} & \textbf{UB} & \textbf{LB} & \textbf{Gap (\%)} & \textbf{Prep (secs)} & \textbf{Solve (secs)} & \textbf{Total (secs)} & \textbf{UB} & \textbf{LB} & \textbf{Gap} & \textbf{Prep (secs)} & \textbf{Solve (secs)} & \textbf{Total (secs)} \\
\hline
three & 6 & 6 & 0 & 0.598 & 0.019 & 0.617 & 6 & 6 & 0 & 0.027 & 0.006 & \textcolor{red}{\textbf{0.033}} \\ 
threep2 & 10 & 10 & 0 & 0.77 & 0.021 & 0.791 & 10 & 10 & 0 & 0.003 & 0 & \textcolor{red}{\textbf{0.003}} \\ 
threep2w9 & 8 & 8 & 0 & 0.767 & 0.019 & 0.786 & 8 & 8 & 0 & 0.004 & 0 & \textcolor{red}{\textbf{0.004}} \\ 
threep3 & 14 & 14 & 0 & 0.829 & 0.036 & 0.865 & 14 & 14 & 0 & 0.005 & 0.004 & \textcolor{red}{\textbf{0.009}} \\ 
threep3w9 & 12 & 12 & 0 & 0.919 & 0.097 & 1.016 & 12 & 12 & 0 & 0.006 & 0.108 & \textcolor{red}{\textbf{0.114}} \\ 
Shapes4 & 24 & 24 & 0 & 1.889 & 0.067 & 1.956 & 24 & 24 & 0 & 0.031 & 0 & \textcolor{red}{\textbf{0.031}} \\ 
Shapes8 & 26 & 26 & 0 & 4.812 & 7.188 & 12 & 26 & 26 & 0 & 0.085 & 1.31 & \textcolor{red}{\textbf{1.395}} \\ 
SHAPES-2 & 14 & 14 & 0 & 3.161 & 0.913 & 4.074 & 14 & 14 & 0 & 0.101 & 0.001 & \textcolor{red}{\textbf{0.102}} \\ 
SHAPES-4 & \textcolor{red}{\textbf{25}} & 23 & 8 & 17.47 & TL & TL & 28 & 21 & 25 & 0.282 & TL & TL \\ 
SHAPES-5 & \textcolor{red}{\textbf{30}} & 28 & 6.67 & 28.71 & TL & TL & 35 & 25 & 28.57 & 0.617 & TL & TL \\ 
SHAPES-7 & \textcolor{red}{\textbf{42}} & 35 & 16.67 & 63.56 & TL & TL & X & 34 & X & 0.819 & TL & TL \\ 
SHAPES-9 & X & 38 & X & 81.24 & TL & TL & X & 38 & X & 0.914 & TL & TL \\ 
SHAPES-15 & X & 47 & X & 126.3 & TL & TL & X & 47 & X & 1.08 & TL & TL \\ 
Blasz2 & 26 & 26 & 0 & 3.335 & 35.07 & \textcolor{red}{\textbf{38.4}} & 27 & 24 & 11.11 & 0.085 & TL & TL \\ 
BLAZEWICZ1 & 8 & 8 & 0 & 1.315 & 0.081 & 1.396 & 8 & 8 & 0 & 0.027 & 0 & \textcolor{red}{\textbf{0.027}} \\ 
BLAZEWICZ2 & 14 & 14 & 0 & 3.68 & 14.06 & \textcolor{red}{\textbf{17.74}} & 14 & 14 & 0 & 0.045 & 18.75 & 18.8 \\ 
BLAZEWICZ3 & 20 & 20 & 0 & 5.164 & 170.8 & \textcolor{red}{\textbf{176}} & 21 & 19 & 9.52 & 0.079 & TL & TL \\ 
BLAZEWICZ4 & \textcolor{red}{\textbf{28}} & 26 & 7.14 & 7.317 & TL & TL & 29 & 25 & 13.79 & 0.136 & TL & TL \\ 
BLAZEWICZ5 & \textcolor{red}{\textbf{34}} & 31 & 8.82 & 9.993 & TL & TL & X & 31 & X & 0.181 & TL & TL \\ 
RCO1 & 8 & 8 & 0 & 1.477 & 0.075 & 1.552 & 8 & 8 & 0 & 0.047 & 0 & \textcolor{red}{\textbf{0.047}} \\ 
RCO2 & 15 & 15 & 0 & 3.726 & 7.525 & \textcolor{red}{\textbf{11.25}} & 15 & 15 & 0 & 0.059 & 37.79 & 37.85 \\ 
RCO3 & 22 & 22 & 0 & 6.01 & 133 & \textcolor{red}{\textbf{139}} & 22 & 21 & 4.55 & 0.062 & TL & TL \\ 
RCO4 & \textcolor{red}{\textbf{29}} & 28 & 3.45 & 7.77 & TL & TL & X & 27 & X & 0.144 & TL & TL \\ 
RC05 & \textcolor{red}{\textbf{36}} & 35 & 2.78 & 11.6 & TL & TL & 37 & 34 & 8.11 & 0.372 & TL & TL \\ 
artif1\_2 & 7 & 7 & 0 & 2.298 & 0.208 & 2.506 & 7 & 7 & 0 & 0.068 & 0.003 & \textcolor{red}{\textbf{0.071}} \\ 
artif2\_4 & \textcolor{red}{\textbf{13}} & 12 & 7.69 & 7.28 & TL & TL & \textcolor{red}{\textbf{13}} & 12 & 7.69 & 0.17 & TL & TL \\ 
artif3\_6 & \textcolor{red}{\textbf{19}} & 17 & 10.53 & 15.43 & TL & TL & 20 & 17 & 15 & 0.214 & TL & TL \\ 
artif4\_8 & \textcolor{red}{\textbf{26}} & 22 & 15.38 & 26.64 & TL & TL & \textcolor{red}{\textbf{26}} & 22 & 15.38 & 0.279 & TL & TL \\ 
artif5\_10 & \textcolor{red}{\textbf{33}} & 27 & 18.18 & 38.69 & TL & TL & 34 & 27 & 20.59 & 0.333 & TL & TL \\ 
artif6\_12 & \textcolor{red}{\textbf{38}} & 33 & 13.16 & 54.8 & TL & TL & 41 & 33 & 19.51 & 0.427 & TL & TL \\ 
artif7\_14 & \textcolor{red}{\textbf{46}} & 38 & 17.39 & 76.9 & TL & TL & 47 & 38 & 19.15 & 0.486 & TL & TL \\ 
artif & \textcolor{red}{\textbf{49}} & 41 & 16.33 & 94.06 & TL & TL & 52 & 41 & 21.15 & 0.579 & TL & TL \\ 
shirts1\_2 & 13 & 13 & 0 & 4.638 & 0.021 & 4.659 & 13 & 13 & 0 & 0.095 & 0 & \textcolor{red}{\textbf{0.095}} \\ 
shirts2\_4 & 17 & 17 & 0 & 24.15 & 24.81 & 48.96 & 17 & 17 & 0 & 0.246 & 0.707 & \textcolor{red}{\textbf{0.953}} \\ 
shirts3\_6 & \textcolor{red}{\textbf{24}} & 24 & 0 & 55.67 & 1022 & 1078 & 25 & 24 & 4 & 0.496 & TL & TL \\ 
shirts4\_8 & \textcolor{red}{\textbf{32}} & 31 & 3.12 & 128.6 & TL & TL & 33 & 31 & 6.06 & 1.104 & TL & TL \\ 
shirts5\_10 & \textcolor{red}{\textbf{41}} & 38 & 7.32 & 205 & TL & TL & 42 & 37 & 11.9 & 1.244 & TL & TL \\ 
Dagli1 & 23 & 23 & 0 & 296.3 & 3.647 & 299.9 & 23 & 23 & 0 & 1.963 & 0.003 & \textcolor{red}{\textbf{1.966}} \\ 
fu5 & 18 & 18 & 0 & 3.146 & 0.139 & 3.285 & 18 & 18 & 0 & 0.057 & 0 & \textcolor{red}{\textbf{0.057}} \\ 
fu6 & 23 & 23 & 0 & 14.98 & 0.642 & 15.62 & 23 & 23 & 0 & 0.158 & 0.007 & \textcolor{red}{\textbf{0.165}} \\ 
fu7 & 24 & 24 & 0 & 25.74 & 0.231 & 25.98 & 24 & 24 & 0 & 0.266 & 0.016 & \textcolor{red}{\textbf{0.282}} \\ 
fu8 & 24 & 24 & 0 & 40.38 & 1.809 & 42.19 & 24 & 24 & 0 & 0.44 & 0.035 & \textcolor{red}{\textbf{0.475}} \\ 
fu9 & 25 & 25 & 0 & 120 & 105.5 & 225.5 & 25 & 25 & 0 & 1.297 & 4.714 & \textcolor{red}{\textbf{6.011}} \\ 
fu10 & 29 & 29 & 0 & 267.1 & 3182 & 3449 & 29 & 29 & 0 & 2.013 & 256 & \textcolor{red}{\textbf{258}} \\ 
fu & \textcolor{red}{\textbf{34}} & 31 & 8.82 & 636 & TL & TL & \textcolor{red}{\textbf{34}} & 31 & 8.82 & 3.407 & TL & TL \\ 
J1-10-10-0 & 18 & 18 & 0 & 3.408 & 0.681 & 4.089 & 18 & 18 & 0 & 0.047 & 0.001 & \textcolor{red}{\textbf{0.048}} \\ 
J1-10-10-1 & 17 & 17 & 0 & 3.793 & 2.4 & 6.193 & 17 & 17 & 0 & 0.048 & 0.001 & \textcolor{red}{\textbf{0.049}} \\ 
J1-10-10-2 & 20 & 20 & 0 & 3.593 & 2.217 & 5.81 & 20 & 20 & 0 & 0.045 & 0.001 & \textcolor{red}{\textbf{0.046}} \\ 
J1-10-10-3 & 21 & 21 & 0 & 5.717 & 15.52 & 21.24 & 21 & 21 & 0 & 0.066 & 0.07 & \textcolor{red}{\textbf{0.136}} \\ 
J1-10-10-4 & 13 & 13 & 0 & 3.332 & 1.46 & 4.792 & 13 & 13 & 0 & 0.042 & 0.004 & \textcolor{red}{\textbf{0.046}} \\ 
J1-12-20-0 & 12 & 12 & 0 & 6.365 & 3.874 & 10.24 & 12 & 12 & 0 & 0.104 & 0.003 & \textcolor{red}{\textbf{0.107}} \\ 
J1-12-20-1 & 10 & 10 & 0 & 4.653 & 3.43 & 8.083 & 10 & 10 & 0 & 0.074 & 0.008 & \textcolor{red}{\textbf{0.082}} \\ 
J1-12-20-2 & 12 & 12 & 0 & 10.57 & 7.575 & 18.15 & 12 & 12 & 0 & 0.176 & 0.01 & \textcolor{red}{\textbf{0.186}} \\ 
J1-12-20-3 & 8 & 8 & 0 & 4.191 & 1.131 & 5.322 & 8 & 8 & 0 & 0.072 & 0.005 & \textcolor{red}{\textbf{0.077}} \\ 
J1-12-20-4 & 13 & 13 & 0 & 14.41 & 23.86 & 38.27 & 13 & 13 & 0 & 0.197 & 0.049 & \textcolor{red}{\textbf{0.246}} \\ 
J1-14-20-0 & 12 & 12 & 0 & 12.42 & 24.83 & 37.25 & 12 & 12 & 0 & 0.187 & 0.03 & \textcolor{red}{\textbf{0.217}} \\ 
J1-14-20-1 & 12 & 12 & 0 & 10.72 & 42.47 & 53.18 & 12 & 12 & 0 & 0.154 & 0.936 & \textcolor{red}{\textbf{1.09}} \\ 
J1-14-20-2 & 14 & 14 & 0 & 18.75 & 105.4 & 124.2 & 14 & 14 & 0 & 0.265 & 2.271 & \textcolor{red}{\textbf{2.536}} \\ 
J1-14-20-3 & 10 & 10 & 0 & 8.146 & 6.53 & 14.68 & 10 & 10 & 0 & 0.125 & 0.013 & \textcolor{red}{\textbf{0.138}} \\ 
J1-14-20-4 & 14 & 14 & 0 & 15.97 & 61.71 & 77.69 & 14 & 14 & 0 & 0.215 & 0.319 & \textcolor{red}{\textbf{0.534}} \\ 
J2-10-35-0* & 24 & 24 & 0 & 121.6 & 1847 & 1968 & 24 & 24 & 0 & 1.499 & 5.572 & \textcolor{red}{\textbf{7.071}} \\ 
J2-10-35-1* & 24 & 22 & 8.33 & 126.2 & TL & TL & 24 & 24 & 0 & 1.123 & 18.22 & \textcolor{red}{\textbf{19.34}} \\ 
J2-10-35-2* & 23 & 23 & 0 & 139.8 & 3598 & 3738 & 23 & 23 & 0 & 1.351 & 3.609 & \textcolor{red}{\textbf{4.96}} \\ 
J2-10-35-3* & 20 & 20 & 0 & 97.43 & 192.6 & 290 & 20 & 20 & 0 & 1.187 & 0.636 & \textcolor{red}{\textbf{1.823}} \\ 
J2-10-35-4* & 19 & 19 & 0 & 87.12 & 1603 & 1691 & 19 & 19 & 0 & 1.193 & 1.493 & \textcolor{red}{\textbf{2.686}} \\ 
J2-12-35-0* & 30 & 26 & 13.33 & 285.7 & TL & TL & 29 & 29 & 0 & 2.304 & 428.3 & \textcolor{red}{\textbf{430.6}} \\ 
J2-12-35-1* & 26 & 23 & 11.54 & 228 & TL & TL & 25 & 25 & 0 & 1.822 & 196.7 & \textcolor{red}{\textbf{198.6}} \\ 
J2-12-35-2* & 24 & 22 & 8.33 & 202.6 & TL & TL & 24 & 24 & 0 & 1.606 & 205.1 & \textcolor{red}{\textbf{206.7}} \\ 
J2-12-35-3* & 22 & 20 & 9.09 & 210.3 & TL & TL & 21 & 21 & 0 & 1.217 & 143.5 & \textcolor{red}{\textbf{144.8}} \\ 
J2-12-35-4* & 26 & 24 & 7.69 & 239.4 & TL & TL & 26 & 26 & 0 & 1.805 & 279.3 & \textcolor{red}{\textbf{281.1}} \\ 
J2-14-35-0* & 30 & 26 & 13.33 & 592.1 & TL & TL & 30 & 30 & 0 & 3.242 & 2192 & \textcolor{red}{\textbf{2195}} \\ 
J2-14-35-1 & 30 & 25 & 16.67 & 537.9 & TL & TL & \textcolor{red}{\textbf{29}} & 25 & 13.79 & 2.911 & TL & TL \\ 
J2-14-35-2* & 26 & 23 & 11.54 & 353.6 & TL & TL & 26 & 26 & 0 & 2.209 & 2354 & \textcolor{red}{\textbf{2356}} \\ 
J2-14-35-3 & 27 & 22 & 18.52 & 389.1 & TL & TL & \textcolor{red}{\textbf{26}} & 22 & 15.38 & 2.811 & TL & TL \\ 
J2-14-35-4* & 27 & 24 & 11.11 & 578.4 & TL & TL & 26 & 26 & 0 & 4.244 & 1815 & \textcolor{red}{\textbf{1820}} \\ 
Poly1a* & 17 & 14 & 17.65 & 136.5 & TL & TL & 15 & 15 & 0 & 1.035 & 2422 & \textcolor{red}{\textbf{2423}} \\ 
Poly1b* & 19 & 16 & 15.79 & 240.8 & TL & TL & 18 & 18 & 0 & 1.985 & 53.9 & \textcolor{red}{\textbf{55.88}} \\ 
Poly1c & 13 & 13 & 0 & 50.21 & 3.793 & 54 & 13 & 13 & 0 & 0.405 & 0 & \textcolor{red}{\textbf{0.405}} \\ 
Poly1d* & 13 & 11 & 15.38 & 67.92 & TL & TL & 12 & 12 & 0 & 0.608 & 36.08 & \textcolor{red}{\textbf{36.69}} \\ 
Poly1e* & 12 & 11 & 8.33 & 39.58 & TL & TL & 12 & 12 & 0 & 0.351 & 527.9 & \textcolor{red}{\textbf{528.2}} \\ 
Jakobs1 & 11 & 11 & 0 & 137.7 & 164.4 & \textcolor{red}{\textbf{302}} & 11 & 11 & 0 & 0.969 & 791.1 & 792.1 \\ 
\hline
\# Feasible & \textcolor{red}{\textbf{79}} & & & & & &76 & & & &  \\
\# Optimal & & & 47 &  &  & & & & \textcolor{red}{\textbf{56}} & & \\
\# Best total &  & & & & & 6 &  & & & & & \textcolor{red}{\textbf{53}} \\
\# Best UB/ties & \textcolor{red}{\textbf{18}} & & & & &  & 5 & & & &  \\
Preprocessing  &  & & & 7518.224 & &  &  & & & \textcolor{red}{\textbf{58.247}} &
\end{tabular}
\end{table}
\clearpage

\begin{table}[h]
\centering
\caption{Comparison of our 0-1 DB-CC model solved with Google CPSat and our DB-PB algorithm with TL = 3600 secs.}
\label{tab:second_results}
\tiny
\begin{tabular}{l|cccccc|cccccc}
& \multicolumn{6}{c}{\textbf{Binary DB Clique Covering (0-1 DB-CC) model (this work)}} & \multicolumn{6}{c}{\textbf{DB Parallel-Backtracking (DB-PB) method (this work)}} \\
\textbf{Instance} & \textbf{UB} & \textbf{LB} & \textbf{Gap (\%)} & \textbf{Prep (secs)} & \textbf{Solve (secs)} & \textbf{Total (secs)} & \textbf{UB} & \textbf{LB} & \textbf{Gap} & \textbf{Prep (secs)} & \textbf{Solve (secs)} & \textbf{Total (secs)} \\
\hline
three & 6 & 6 & 0 & 0.578 & 0.02 & 0.598 & 6 & 6 & 0 & 0.027 & 0.006 & \textcolor{red}{\textbf{0.033}} \\ 
threep2 & 10 & 10 & 0 & 0.642 & 0.043 & 0.685 & 10 & 10 & 0 & 0.003 & 0 & \textcolor{red}{\textbf{0.003}} \\ 
threep2w9 & 8 & 8 & 0 & 0.62 & 0.04 & 0.66 & 8 & 8 & 0 & 0.004 & 0 & \textcolor{red}{\textbf{0.004}} \\ 
threep3 & 14 & 14 & 0 & 0.778 & 0.121 & 0.899 & 14 & 14 & 0 & 0.005 & 0.004 & \textcolor{red}{\textbf{0.009}} \\ 
threep3w9 & 12 & 12 & 0 & 0.682 & 0.637 & 1.319 & 12 & 12 & 0 & 0.006 & 0.108 & \textcolor{red}{\textbf{0.114}} \\ 
Shapes4 & 24 & 24 & 0 & 1.373 & 0.035 & 1.408 & 24 & 24 & 0 & 0.031 & 0 & \textcolor{red}{\textbf{0.031}} \\ 
Shapes8 & 26 & 26 & 0 & 4.111 & 2.928 & 7.039 & 26 & 26 & 0 & 0.085 & 1.31 & \textcolor{red}{\textbf{1.395}} \\ 
SHAPES-2 & 14 & 14 & 0 & 2.367 & 0.479 & 2.846 & 14 & 14 & 0 & 0.101 & 0.001 & \textcolor{red}{\textbf{0.102}} \\ 
SHAPES-4 & \textcolor{red}{\textbf{25}} & 23 & 8 & 14.84 & TL & TL & 28 & 21 & 25 & 0.282 & TL & TL \\ 
SHAPES-5 & \textcolor{red}{\textbf{31}} & 25 & 19.35 & 29.06 & TL & TL & 35 & 25 & 28.57 & 0.617 & TL & TL \\ 
SHAPES-7 & \textcolor{red}{\textbf{43}} & 34 & 20.93 & 63.29 & TL & TL & X & 34 & X & 0.819 & TL & TL \\ 
SHAPES-9 & \textcolor{red}{\textbf{52}} & 38 & 26.92 & 81.16 & TL & TL & X & 38 & X & 0.914 & TL & TL \\ 
SHAPES-15 & \textcolor{red}{\textbf{64}} & 47 & 26.56 & 126.5 & TL & TL & X & 47 & X & 1.08 & TL & TL \\ 
Blasz2 & \textcolor{red}{\textbf{26}} & 24 & 7.69 & 3.475 & TL & TL & 27 & 24 & 11.11 & 0.085 & TL & TL \\ 
BLAZEWICZ1 & 8 & 8 & 0 & 1.772 & 0.125 & 1.897 & 8 & 8 & 0 & 0.027 & 0 & \textcolor{red}{\textbf{0.027}} \\ 
BLAZEWICZ2 & 14 & 14 & 0 & 4.136 & 1099 & 1103 & 14 & 14 & 0 & 0.045 & 18.75 & \textcolor{red}{\textbf{18.8}} \\ 
BLAZEWICZ3 & \textcolor{red}{\textbf{21}} & 19 & 9.52 & 6.492 & TL & TL & \textcolor{red}{\textbf{21}} & 19 & 9.52 & 0.079 & TL & TL \\ 
BLAZEWICZ4 & \textcolor{red}{\textbf{28}} & 25 & 10.71 & 8.568 & TL & TL & 29 & 25 & 13.79 & 0.136 & TL & TL \\ 
BLAZEWICZ5 & \textcolor{red}{\textbf{35}} & 31 & 11.43 & 10.56 & TL & TL & X & 31 & X & 0.181 & TL & TL \\ 
RCO1 & 8 & 8 & 0 & 1.827 & 0.087 & 1.914 & 8 & 8 & 0 & 0.047 & 0 & \textcolor{red}{\textbf{0.047}} \\ 
RCO2 & 15 & 15 & 0 & 4.473 & 237.4 & 241.9 & 15 & 15 & 0 & 0.059 & 37.79 & \textcolor{red}{\textbf{37.85}} \\ 
RCO3 & \textcolor{red}{\textbf{22}} & 21 & 4.55 & 7.08 & TL & TL & \textcolor{red}{\textbf{22}} & 21 & 4.55 & 0.062 & TL & TL \\ 
RCO4 & \textcolor{red}{\textbf{29}} & 27 & 6.9 & 7.23 & TL & TL & X & 27 & X & 0.144 & TL & TL \\ 
RC05 & \textcolor{red}{\textbf{37}} & 34 & 8.11 & 12.81 & TL & TL & \textcolor{red}{\textbf{37}} & 34 & 8.11 & 0.372 & TL & TL \\ 
artif1\_2 & 7 & 7 & 0 & 2.371 & 0.251 & 2.622 & 7 & 7 & 0 & 0.068 & 0.003 & \textcolor{red}{\textbf{0.071}} \\ 
artif2\_4 & \textcolor{red}{\textbf{13}} & 12 & 7.69 & 7.059 & TL & TL & \textcolor{red}{\textbf{13}} & 12 & 7.69 & 0.17 & TL & TL \\ 
artif3\_6 & \textcolor{red}{\textbf{20}} & 17 & 15 & 15.18 & TL & TL & \textcolor{red}{\textbf{20}} & 17 & 15 & 0.214 & TL & TL \\ 
artif4\_8 & \textcolor{red}{\textbf{26}} & 22 & 15.38 & 25.75 & TL & TL & \textcolor{red}{\textbf{26}} & 22 & 15.38 & 0.279 & TL & TL \\ 
artif5\_10 & \textcolor{red}{\textbf{33}} & 27 & 18.18 & 50.58 & TL & TL & 34 & 27 & 20.59 & 0.333 & TL & TL \\ 
artif6\_12 & \textcolor{red}{\textbf{39}} & 33 & 15.38 & 71.35 & TL & TL & 41 & 33 & 19.51 & 0.427 & TL & TL \\ 
artif7\_14 & \textcolor{red}{\textbf{46}} & 38 & 17.39 & 98.39 & TL & TL & 47 & 38 & 19.15 & 0.486 & TL & TL \\ 
artif & \textcolor{red}{\textbf{50}} & 41 & 18 & 111.2 & TL & TL & 52 & 41 & 21.15 & 0.579 & TL & TL \\ 
shirts1\_2 & 13 & 13 & 0 & 6.587 & 0.473 & 7.06 & 13 & 13 & 0 & 0.095 & 0 & \textcolor{red}{\textbf{0.095}} \\ 
shirts2\_4 & 17 & 17 & 0 & 33.8 & 20.16 & 53.96 & 17 & 17 & 0 & 0.246 & 0.707 & \textcolor{red}{\textbf{0.953}} \\ 
shirts3\_6 & \textcolor{red}{\textbf{25}} & 24 & 4 & 74.01 & TL & TL & \textcolor{red}{\textbf{25}} & 24 & 4 & 0.496 & TL & TL \\ 
shirts4\_8 & \textcolor{red}{\textbf{33}} & 31 & 6.06 & 181.1 & TL & TL & \textcolor{red}{\textbf{33}} & 31 & 6.06 & 1.104 & TL & TL \\ 
shirts5\_10 & \textcolor{red}{\textbf{41}} & 37 & 9.76 & 268 & TL & TL & 42 & 37 & 11.9 & 1.244 & TL & TL \\ 
Dagli1 & 23 & 23 & 0 & 335.2 & 10.67 & 345.8 & 23 & 23 & 0 & 1.963 & 0.003 & \textcolor{red}{\textbf{1.966}} \\ 
fu5 & 18 & 18 & 0 & 3.355 & 0.089 & 3.444 & 18 & 18 & 0 & 0.057 & 0 & \textcolor{red}{\textbf{0.057}} \\ 
fu6 & 23 & 23 & 0 & 14.43 & 2.245 & 16.67 & 23 & 23 & 0 & 0.158 & 0.007 & \textcolor{red}{\textbf{0.165}} \\ 
fu7 & 24 & 24 & 0 & 23.78 & 1.284 & 25.06 & 24 & 24 & 0 & 0.266 & 0.016 & \textcolor{red}{\textbf{0.282}} \\ 
fu8 & 24 & 24 & 0 & 40.96 & 3.221 & 44.18 & 24 & 24 & 0 & 0.44 & 0.035 & \textcolor{red}{\textbf{0.475}} \\ 
fu9 & 25 & 25 & 0 & 109.5 & 80.78 & 190.3 & 25 & 25 & 0 & 1.297 & 4.714 & \textcolor{red}{\textbf{6.011}} \\ 
fu10 & 29 & 28 & 3.45 & 261.1 & TL & TL & 29 & 29 & 0 & 2.013 & 256 & \textcolor{red}{\textbf{258}} \\ 
fu & \textcolor{red}{\textbf{34}} & 31 & 8.82 & 642.9 & TL & TL & \textcolor{red}{\textbf{34}} & 31 & 8.82 & 3.407 & TL & TL \\ 
J1-10-10-0 & 18 & 18 & 0 & 4.151 & 0.288 & 4.439 & 18 & 18 & 0 & 0.047 & 0.001 & \textcolor{red}{\textbf{0.048}} \\ 
J1-10-10-1 & 17 & 17 & 0 & 4.009 & 1.072 & 5.081 & 17 & 17 & 0 & 0.048 & 0.001 & \textcolor{red}{\textbf{0.049}} \\ 
J1-10-10-2 & 20 & 20 & 0 & 4.039 & 0.708 & 4.747 & 20 & 20 & 0 & 0.045 & 0.001 & \textcolor{red}{\textbf{0.046}} \\ 
J1-10-10-3 & 21 & 21 & 0 & 5.403 & 10.73 & 16.13 & 21 & 21 & 0 & 0.066 & 0.07 & \textcolor{red}{\textbf{0.136}} \\ 
J1-10-10-4 & 13 & 13 & 0 & 2.842 & 1.244 & 4.086 & 13 & 13 & 0 & 0.042 & 0.004 & \textcolor{red}{\textbf{0.046}} \\ 
J1-12-20-0 & 12 & 12 & 0 & 6.789 & 1.642 & 8.431 & 12 & 12 & 0 & 0.104 & 0.003 & \textcolor{red}{\textbf{0.107}} \\ 
J1-12-20-1 & 10 & 10 & 0 & 6.333 & 0.756 & 7.089 & 10 & 10 & 0 & 0.074 & 0.008 & \textcolor{red}{\textbf{0.082}} \\ 
J1-12-20-2 & 12 & 12 & 0 & 11.5 & 0.668 & 12.17 & 12 & 12 & 0 & 0.176 & 0.01 & \textcolor{red}{\textbf{0.186}} \\ 
J1-12-20-3 & 8 & 8 & 0 & 4.352 & 0.656 & 5.008 & 8 & 8 & 0 & 0.072 & 0.005 & \textcolor{red}{\textbf{0.077}} \\ 
J1-12-20-4 & 13 & 13 & 0 & 14.8 & 13.57 & 28.36 & 13 & 13 & 0 & 0.197 & 0.049 & \textcolor{red}{\textbf{0.246}} \\ 
J1-14-20-0 & 12 & 12 & 0 & 12.29 & 31.92 & 44.21 & 12 & 12 & 0 & 0.187 & 0.03 & \textcolor{red}{\textbf{0.217}} \\ 
J1-14-20-1 & 12 & 12 & 0 & 11.5 & 182.9 & 194.4 & 12 & 12 & 0 & 0.154 & 0.936 & \textcolor{red}{\textbf{1.09}} \\ 
J1-14-20-2 & 14 & 14 & 0 & 21.02 & 675.8 & 696.8 & 14 & 14 & 0 & 0.265 & 2.271 & \textcolor{red}{\textbf{2.536}} \\ 
J1-14-20-3 & 10 & 10 & 0 & 9.891 & 2.554 & 12.44 & 10 & 10 & 0 & 0.125 & 0.013 & \textcolor{red}{\textbf{0.138}} \\ 
J1-14-20-4 & 14 & 14 & 0 & 17.58 & 72.76 & 90.34 & 14 & 14 & 0 & 0.215 & 0.319 & \textcolor{red}{\textbf{0.534}} \\ 
J2-10-35-0* & 24 & 24 & 0 & 121.7 & 640.8 & 762.5 & 24 & 24 & 0 & 1.499 & 5.572 & \textcolor{red}{\textbf{7.071}} \\ 
J2-10-35-1* & 24 & 23 & 4.17 & 126.4 & TL & TL & 24 & 24 & 0 & 1.123 & 18.22 & \textcolor{red}{\textbf{19.34}} \\ 
J2-10-35-2* & 23 & 23 & 0 & 148.6 & 1491 & 1640 & 23 & 23 & 0 & 1.351 & 3.609 & \textcolor{red}{\textbf{4.96}} \\ 
J2-10-35-3* & 20 & 20 & 0 & 100.7 & 148.7 & 249.4 & 20 & 20 & 0 & 1.187 & 0.636 & \textcolor{red}{\textbf{1.823}} \\ 
J2-10-35-4* & 19 & 19 & 0 & 84.72 & 180.8 & 265.5 & 19 & 19 & 0 & 1.193 & 1.493 & \textcolor{red}{\textbf{2.686}} \\ 
J2-12-35-0* & 29 & 25 & 13.79 & 285.4 & TL & TL & 29 & 29 & 0 & 2.304 & 428.3 & \textcolor{red}{\textbf{430.6}} \\ 
J2-12-35-1* & 25 & 23 & 8 & 226.3 & TL & TL & 25 & 25 & 0 & 1.822 & 196.7 & \textcolor{red}{\textbf{198.6}} \\ 
J2-12-35-2* & 24 & 22 & 8.33 & 201.6 & TL & TL & 24 & 24 & 0 & 1.606 & 205.1 & \textcolor{red}{\textbf{206.7}} \\ 
J2-12-35-3* & 21 & 20 & 4.76 & 162.3 & TL & TL & 21 & 21 & 0 & 1.217 & 143.5 & \textcolor{red}{\textbf{144.8}} \\ 
J2-12-35-4* & 26 & 24 & 7.69 & 244.4 & TL & TL & 26 & 26 & 0 & 1.805 & 279.3 & \textcolor{red}{\textbf{281.1}} \\ 
J2-14-35-0* & 30 & 26 & 13.33 & 592.3 & TL & TL & 30 & 30 & 0 & 3.242 & 2192 & \textcolor{red}{\textbf{2195}} \\ 
J2-14-35-1 & 30 & 25 & 16.67 & 532 & TL & TL & \textcolor{red}{\textbf{29}} & 25 & 13.79 & 2.911 & TL & TL \\ 
J2-14-35-2* & 26 & 23 & 11.54 & 338.9 & TL & TL & 26 & 26 & 0 & 2.209 & 2354 & \textcolor{red}{\textbf{2356}} \\ 
J2-14-35-3 & \textcolor{red}{\textbf{26}} & 22 & 15.38 & 382.5 & TL & TL & \textcolor{red}{\textbf{26}} & 22 & 15.38 & 2.811 & TL & TL \\ 
J2-14-35-4* & 27 & 24 & 11.11 & 570.8 & TL & TL & 26 & 26 & 0 & 4.244 & 1815 & \textcolor{red}{\textbf{1820}} \\ 
Poly1a* & 16 & 14 & 12.5 & 192.5 & TL & TL & 15 & 15 & 0 & 1.035 & 2422 & \textcolor{red}{\textbf{2423}} \\ 
Poly1b* & 18 & 17 & 5.56 & 233.8 & TL & TL & 18 & 18 & 0 & 1.985 & 53.9 & \textcolor{red}{\textbf{55.88}} \\ 
Poly1c & 13 & 13 & 0 & 52.14 & 3.63 & 55.77 & 13 & 13 & 0 & 0.405 & 0 & \textcolor{red}{\textbf{0.405}} \\ 
Poly1d* & 12 & 11 & 8.33 & 66.74 & TL & TL & 12 & 12 & 0 & 0.608 & 36.08 & \textcolor{red}{\textbf{36.69}} \\ 
Poly1e* & 12 & 11 & 8.33 & 44.5 & TL & TL & 12 & 12 & 0 & 0.351 & 527.9 & \textcolor{red}{\textbf{528.2}} \\ 
Jakobs1 & 11 & 11 & 0 & 140.2 & 134.6 & \textcolor{red}{\textbf{274.9}} & 11 & 11 & 0 & 0.969 & 791.1 & 792.1 \\ 
\hline
\# Feasible & \textcolor{red}{\textbf{81}} & & & & & & 76 & & & &  \\
\# Optimal & & & 42 &  &  & & & & \textcolor{red}{\textbf{56}} & & \\
\# Best time &  & & & & & 1 &  & & & & & \textcolor{red}{\textbf{55}} \\
\# Best UB/ties & \textcolor{red}{\textbf{24}} & & & & &  & 11 & & & &  \\
Preprocessing  &  & & & 7755.892 & &  &  & & & \textcolor{red}{\textbf{58.247}} &
\end{tabular}
\end{table}
\clearpage

\begin{table}[h]
\centering
\caption{Number of binary decision variables and constraints defined by each exact model and algorithm.}
\label{tab:binary_variables}    
\tiny
\begin{tabular}{l|ccc|ccc}
 & \multicolumn{3}{c}{\textbf{Total number of binary variables per instance}} & \multicolumn{3}{|c}{\textbf{Total number of constraints per instance}} \\
\textbf{Instance} & \textbf{DB-CC} \cite{Rodrigues2017-zl} & \makecell[c]{\textbf{0-1 DB-CC} \\ \textbf{(this work)}} & \makecell[c]{\textbf{DB-PB} \\ \textbf{(this work)}} & \textbf{DB-CC} \cite{Rodrigues2017-zl}& \makecell[c]{\textbf{0-1 DB-CC} \\ \textbf{(this work)}} & \makecell[c]{\textbf{DB-PB} \\ \textbf{(this work)}} \\
\hline
three &   61 &   62 &   61 &   42 &   34 & 1266 \\ 
threep2 &  117 &  118 &  117 &   98 &   75 & 3006 \\ 
threep2w9 &  127 &  128 &  127 &   96 &   84 & 3394 \\ 
threep3 &  187 &  189 &  187 &  151 &  139 & 5181 \\ 
threep3w9 &  207 &  209 &  207 &  178 &  166 & 6150 \\ 
Shapes4 &  389 &  390 &  389 &  326 &  323 & 46060 \\ 
Shapes8 & 1008 & 1013 & 1008 & 1104 & 1062 & 201702 \\ 
SHAPES-2 &  848 &  850 &  848 &  960 &  906 & 92766 \\ 
SHAPES-4 & 2468 & 2475 & 2468 & 3137 & 3006 & 604259 \\ 
SHAPES-5 & 3413 & 3423 & 3413 & 4455 & 4253 & 939916 \\ 
SHAPES-7 & 5168 & 5182 & 5168 & 6740 & 6456 & 1563279 \\ 
SHAPES-9 & 5978 & 5994 & 5978 & 7719 & 7634 & 1850985 \\ 
SHAPES-15 & 7733 & 7753 & 7733 & 9944 & 9821 & 2474348 \\ 
Blasz2 & 1080 & 1083 & 1080 &  876 &  854 & 91670 \\ 
BLAZEWICZ1 &  432 &  432 &  432 &  280 &  276 & 35462 \\ 
BLAZEWICZ2 & 1096 & 1099 & 1096 &  854 &  761 & 120814 \\ 
BLAZEWICZ3 & 1594 & 1597 & 1594 & 1199 & 1137 & 184828 \\ 
BLAZEWICZ4 & 2175 & 2179 & 2175 & 1673 & 1589 & 259511 \\ 
BLAZEWICZ5 & 2756 & 2761 & 2756 & 2079 & 2040 & 334194 \\ 
RCO1 &  432 &  432 &  432 &  199 &  202 & 36049 \\ 
RCO2 & 1179 & 1182 & 1179 &  679 &  609 & 133924 \\ 
RCO3 & 1843 & 1847 & 1843 & 1054 &  973 & 220924 \\ 
RCO4 & 2175 & 2177 & 2175 & 1185 & 1124 & 264424 \\ 
RC05 & 3171 & 3178 & 3171 & 1784 & 1714 & 394924 \\ 
artif1\_2 &  947 &  948 &  947 &  604 &  584 & 68540 \\ 
artif2\_4 & 2099 & 2101 & 2099 & 1687 & 1585 & 250220 \\ 
artif3\_6 & 3251 & 3254 & 3251 & 2767 & 2675 & 437228 \\ 
artif4\_8 & 4787 & 4793 & 4787 & 4205 & 4129 & 686572 \\ 
artif5\_10 & 5939 & 5946 & 5939 & 5289 & 5274 & 873580 \\ 
artif6\_12 & 7283 & 7291 & 7283 & 6538 & 6444 & 1091756 \\ 
artif7\_14 & 8627 & 8637 & 8627 & 7767 & 7609 & 1309932 \\ 
artif & 9587 & 9599 & 9587 & 8676 & 8552 & 1465772 \\ 
shirts1\_2 & 1924 & 1924 & 1924 &  603 &  609 & 147101 \\ 
shirts2\_4 & 4010 & 4013 & 4010 & 1910 & 1822 & 736227 \\ 
shirts3\_6 & 5798 & 5800 & 5798 & 3121 & 3041 & 1364984 \\ 
shirts4\_8 & 8480 & 8484 & 8480 & 5085 & 4957 & 2319578 \\ 
shirts5\_10 & 10566 & 10571 & 10566 & 6563 & 6445 & 3062040 \\ 
Dagli1 & 6911 & 6913 & 6911 & 5452 & 5211 & 6170439 \\ 
fu5 &  721 &  721 &  721 &  226 &  230 & 120663 \\ 
fu6 & 1722 & 1723 & 1722 &  762 &  736 & 734345 \\ 
fu7 & 2172 & 2172 & 2172 &  849 &  837 & 1124614 \\ 
fu8 & 2772 & 2772 & 2772 & 1374 & 1317 & 1616281 \\ 
fu9 & 4202 & 4207 & 4202 & 2555 & 2210 & 3446231 \\ 
fu10 & 5982 & 5988 & 5982 & 3440 & 2942 & 6365641 \\ 
fu & 8642 & 8649 & 8642 & 6739 & 6086 & 11355356 \\ 
J1-10-10-0 & 1063 & 1065 & 1063 &  581 &  525 & 118778 \\ 
J1-10-10-1 & 1065 & 1068 & 1065 &  766 &  680 & 154905 \\ 
J1-10-10-2 & 1073 & 1075 & 1073 &  576 &  500 & 140430 \\ 
J1-10-10-3 & 1284 & 1287 & 1284 &  702 &  624 & 232436 \\ 
J1-10-10-4 &  896 &  899 &  896 &  403 &  340 & 96289 \\ 
J1-12-20-0 & 1597 & 1598 & 1597 &  727 &  685 & 260436 \\ 
J1-12-20-1 & 1377 & 1378 & 1377 &  752 &  728 & 207624 \\ 
J1-12-20-2 & 2039 & 2041 & 2039 &  938 &  838 & 445512 \\ 
J1-12-20-3 & 1363 & 1365 & 1363 &  621 &  554 & 172873 \\ 
J1-12-20-4 & 2408 & 2412 & 2408 & 1186 & 1034 & 580191 \\ 
J1-14-20-0 & 2344 & 2346 & 2344 &  934 &  863 & 469901 \\ 
J1-14-20-1 & 2319 & 2322 & 2319 & 1227 & 1056 & 459468 \\ 
J1-14-20-2 & 2781 & 2784 & 2781 & 1133 &  969 & 776783 \\ 
J1-14-20-3 & 1979 & 1981 & 1979 &  651 &  537 & 348312 \\ 
J1-14-20-4 & 2631 & 2634 & 2631 & 1481 & 1326 & 634924 \\ 
J2-10-35-0 & 4486 & 4492 & 4486 & 4787 & 4353 & 3505202 \\ 
J2-10-35-1 & 4572 & 4578 & 4572 & 4807 & 4416 & 3566129 \\ 
J2-10-35-2 & 4892 & 4898 & 4892 & 4744 & 4225 & 3811214 \\ 
J2-10-35-3 & 4328 & 4335 & 4328 & 3335 & 2930 & 2872285 \\ 
J2-10-35-4 & 3992 & 3997 & 3992 & 3717 & 3353 & 2678511 \\ 
J2-12-35-0 & 6240 & 6246 & 6240 & 5947 & 5416 & 6493492 \\ 
J2-12-35-1 & 6100 & 6106 & 6100 & 6103 & 5653 & 5234224 \\ 
J2-12-35-2 & 5852 & 5861 & 5852 & 5878 & 5366 & 4829349 \\ 
J2-12-35-3 & 5420 & 5426 & 5420 & 4729 & 4291 & 4143829 \\ 
J2-12-35-4 & 6256 & 6261 & 6256 & 4804 & 4376 & 5439275 \\ 
J2-14-35-0 & 8718 & 8726 & 8718 & 8024 & 7353 & 10155708 \\ 
J2-14-35-1 & 8620 & 8628 & 8620 & 7768 & 7140 & 9308609 \\ 
J2-14-35-2 & 7484 & 7494 & 7484 & 6888 & 6247 & 6813466 \\ 
J2-14-35-3 & 7552 & 7559 & 7552 & 6362 & 5687 & 7557702 \\ 
J2-14-35-4 & 8056 & 8063 & 8056 & 6826 & 6202 & 8157073 \\ 
Poly1a & 5620 & 5624 & 5620 & 5104 & 4815 & 3283976 \\ 
Poly1b & 6858 & 6863 & 6858 & 6507 & 6261 & 4922042 \\ 
Poly1c & 4437 & 4438 & 4437 & 3356 & 3224 & 1513472 \\ 
Poly1d & 4527 & 4530 & 4527 & 3639 & 3412 & 1946264 \\ 
Poly1e & 3638 & 3640 & 3638 & 2956 & 2873 & 1301254 \\ 
Jakobs1 & 7460 & 7462 & 7460 & 2472 & 2354 & 2846369
\end{tabular}
\end{table}
\clearpage

\section{Conclusions and future work}
\label{sec:conclusions}

We introduce a parallel exact branch-and-bound-and-prune algorithm requiring no pre-processing to solve efficiently the discrete nesting problem with rotations, together with two binary linear programming models for the same problem based on two reformulations of the family of Dotted-Board models \cite{Toledo2013-oi, Rodrigues2017-zl, Cherri2018-rq} and a new lower bound algorithm as by-products. Finally, we introduce a detailed reproducibility protocol and dataset to allow the exact replication of all our models, experiments, and results.

We reproduce the same benchmark introduced by \citet{Rodrigues2017-zl} to compare the current state-of-the-art DB-CC \cite{Rodrigues2017-zl} model with our 0-1 DB-CC model and our ad-hoc exact DB-PB algorithm by replicating and implementing all methods evaluated herein into the same software and hardware platform. Our experiments shows that DB-PB significantly outperforms the DB-CC and 0-1 DB-CC models both in number of problem instances solved up to optimality and solving time. Our ad-hoc exact DB-PB algorithm outperforms the current state-of-the-art DB-CC model and our new 0-1 DB-CC reformulation by two orders of magnitude in most of problem instances without requiring their costly preprocessing time to compute the edge and vertex clique coverings. Our experiments positively confirm our main hypothesis claiming that the preprocessing time required by the DB-CC is unnecessary and counterproductive, and that an ad-hoc exact algorithm might be a better option to solve the discrete nesting problem in comparison with costly state-of-the-art Branch-and-Cut, Constraint Programming, or hybrid CP-SAT approaches. Finally, our experiments show that our DP-PB-LB lower bound algorithm outperforms the iterative lower algorithm of \citet{Rodrigues2017-zl} by solving 51 problem instances up to optimality within 10 minutes, that means 4 problem instances more than DB-CC within 1 hour. On the other hand, the DB-CC \cite{Rodrigues2017-zl} model obtains tighter upper bounds that our DB-PB algorithm in most of unsolved problem instances and two feasible solutions more, which we attribute to the set of feasibility pump and primal heuristics implemented by Gurobi and Google CPSat, as well as the lack of these in our new DB-PB algorithm. Finally, our experiments solve seventeen open problem instances within 1 hour for the first time in the literature and improve the best known upper bound for other ten instances.

As forthcoming activities. we plan to devise ad-hoc feasibility pump and primal heuristics to improve the capability of our current DB-PB algorithm of finding high-quality solutions for all problem instances, specially when they cannot be solved up to optimality within the time limits. Likewise, we plan to devise new nesting heuristics based on our DB-PB framework, which can compete with current state-of-the-art heuristics for the continuous nesting problem.

\section*{Acknowledgments}

We thank Hervé Gateuil for the access to a high-performance Lectra computer to validate our algorithms.

\appendix

\section{Appendix A: Best solutions obtained for all problem instances}
\label{sec:appendix_A}

This appendix introduces the best solutions for all problem instances in our experiments.

\section{Appendix B: The reproducible experiments on irregular strip-packing}
\label{sec:appendix_B}

This appendix introduces a detailed reproducibility protocol and dataset \cite{Lastra-Diaz2025-gp} providing our raw output data and a collection of software and data resources to allow the replication of all our experiments and results.

\bibliographystyle{cas-model2-names}

%\bibliography{biblio.bib}

\end{document}